\documentclass[11pt,oneside]{article}

\usepackage{amsmath,amsthm,amsfonts,amssymb}
\usepackage{mathrsfs}
\usepackage{mathtools}
\usepackage{enumitem}
\usepackage{tabularx}
\usepackage{physics}
\usepackage{caption}
\usepackage{subcaption}
\usepackage{contour}
\usepackage{xcolor}
\usepackage{tikz}
\usepackage[linesnumbered,ruled,vlined]{algorithm2e}
\usepackage[%dvips,
CJKbookmarks=true,
bookmarksnumbered=true,
bookmarksopen=true,
colorlinks=true,
citecolor=red,
linkcolor=blue,
anchorcolor=red,
urlcolor=blue,
linktocpage=true
]{hyperref}

\usepackage{fullpage}

\usepackage{tocloft}
\RequirePackage{amsthm,amsmath,amsfonts,amssymb}
\RequirePackage[numbers,sort]{natbib}
\RequirePackage{graphicx}

\usetikzlibrary{matrix,positioning,calc}
\NewDocumentCommand{\tens}{t_}
{%
	\IfBooleanTF{#1}
	{\tensop}
	{\otimes}%
}
\def\linebelowrow#1#2{%
	\pgfmathsetmacro{\nextrow}{int(#2+1)}
	\coordinate (aux) at ($(#1-#2-1.center)!.5!(#1-\nextrow-1.center)$);
	\draw[dashed] (#1.west|-aux) -- (#1.east|-aux);
}

\def\lineaftercolumn#1#2{%
	\pgfmathsetmacro{\nextcol}{int(#2+1)}
	\coordinate (aux) at ($(#1-1-#2.center)!.5!(#1-1-\nextcol.center)$);
	\draw[dashed] (#1.north-|aux) -- (#1.south-|aux);
}

\providecommand{\keywords}[1]{\textbf{\textit{Keywords:}} #1}
\newcommand{\bs}{\boldsymbol}
\newcommand{\prob}{\mathbb{P}}
\newcommand{\err}{\textrm{err}}

\newcommand{\expect}{\mathbb{E}}
\newcommand{\stoptime}{\mathcal{T}}
\DeclareMathAlphabet{\mathpzc}{OT1}{pzc}{m}{it}
\newcommand{\Ffunc}{\mathpzc{F}}

\newcommand{\Ffuncbold}{\contour[2]{black}{$\mathpzc{F}$}}
\newcommand{\myquad}[1][1]{\hspace*{#1em}\ignorespaces}
\DeclareMathOperator*{\argmin}{arg\,min}

\DeclarePairedDelimiter\floor{\lfloor}{\rfloor}

\makeatletter
\newcommand{\Biggg}{\bBigg@{3}}

\makeatother
\newcommand\Tstrut{\rule{0pt}{2.6ex}}         % = `top' strut
\newcommand\DTstrut{\rule{0pt}{3.6ex}}
\newcommand\Bstrut{\rule[-1ex]{0pt}{0pt}}   % = `bottom' strut
\newcommand\DBstrut{\rule[-1.5ex]{0pt}{0pt}}   % = `bottom' strut
\newcommand\DDBstrut{\rule[-11.5ex]{0pt}{0pt}}
\SetKwComment{Comment}{//}{ }

\theoremstyle{plain}

\newtheorem{theorem}{Theorem}[section]
\newtheorem{lemma}[theorem]{Lemma}
\newtheorem{corollary}[theorem]{Corollary}

\theoremstyle{remark}
\newtheorem{definition}[theorem]{Definition}
\newtheorem{remark}[theorem]{Remark}

\newtheorem{point}{Point}

\usepackage[noblocks,affil-sl]{authblk}

\author[1]{Mehrdad Moharrami}
\author[2]{Vijay Subramanian}
\author[3]{Mingyan Liu}
\author[4]{Marc Lelarge}
\affil[1]{Coordinated Science Lab, University of Illinois at Urbana Champaign\vspace{1.5ex}}
\affil[2,3]{Electrical and Computer Engineering, University of Michigan\vspace{1.5ex}}
\affil[4]{ENS/INRIA\vspace{1.5ex}}

\begin{document}
		\title{Impact of Community Structure on Cascades}
		\date{}
		\maketitle

	\begin{abstract}
		We study cascades under the threshold model on sparse random graphs with community structure. In this model, individuals adopt the new behavior based on how many neighbors have already chosen it. Specifically, we consider the permanent adoption model wherein individuals that have adopted the new behavior (or opinion) cannot change their state. We present a differential-equation-based tight approximation to the stochastic process of adoption and prove the validity of the mean-field equations. In addition, we characterize both necessary and sufficient conditions for contagion to happen no matter how small the set of initial adopters is. Finally, we study the problem of optimum seeding given budget constraints and propose a gradient-based heuristic seeding strategy. Our algorithm, numerically, dispels commonly held beliefs in the literature that suggest the best seeding strategy is to seed over the vertices with the highest number of neighbors.
	\end{abstract}


	\tableofcontents

	\section{Introduction}\label{sec:intro}
	In this paper, we investigate a type of cascade problem on graphs that has been used to study the spread of new technology or opinions in social networks, see e.g.,~\cite{Granovetter1978,Schelling1978,Naylor1990,Watts2007,Vega2007,Easley2010}.  This spread is also referred to as a contagion in networks.  The underlying model typically consists of a few (selected) initial adopters (vertices in the network) or ``seeds'' and a particular adoption model that determines the condition under which a vertex will choose to adopt given the states of its neighbors.  A commonly studied model here is the threshold model~\cite{Morris2000,Watts2002}, whereby individuals adopt the new technology (or opinion) based on how many neighbors have already chosen it.

Prior work in this area has generally focused on analyzing what happens when the underlying network consists of a single community modeled as a sparse random graph, either heuristically, see e.g., \cite{Watts2002,Lopez2008}, or rigorously, see e.g., \cite{Balogh2007,Amini2010, Lelarge2012,Como2019}. In this work, we instead consider graphs with a type of community structure (also known as modular networks), whereby multiple sparse random graphs are weakly interconnected. This could model, for instance, segments of the population (e.g., different age or ethnic groups), where members of a single segment are more strongly connected (with a relatively high vertex degree) and cross-segment connections are weak, i.e., fewer members are connected to those from a different segment. This would be a more realistic and interesting model for many practical scenarios and serves as a natural next step beyond the studies with a single community. We are particularly interested in whether the existence of communities affects the number of individuals who eventually adopt the new technology. Also of interest is whether seeding in all communities is a better strategy in terms of maximizing the number of eventual adopters than exclusively in one community or, in particular, the optimum seeding strategy given budget constraints. While earlier works have looked at this problem using heuristic methods, see e.g., \cite{Galstyan2007,Gleeson2007,Gleeson2008,Lopez2008,Galstyan2009}, we present a mathematically rigorous analysis of this problem.

Specifically, we consider the permanent adoption model where vertices that have adopted the new technology/behavior/opinion---called active vertices---cannot change their state. Our analysis in Sections \ref{sec:markovproc}-\ref{sec:odeanalysisinf} presents a differential-equation-based tight approximation to the stochastic process of adoption under the threshold model. While the approach is similar to the analysis of contagions in a single community in the case of $d$-regular random graphs \cite{Balogh2007} and random graphs \cite{Amini2010}, the additional community structure requires significant technical development to establish the validity of this approach in the new setting. We also present a probabilistic approach to solve (in an intuitive manner) the associated system of ordinary differential equations (ODEs) in Section \ref{sec:odeprobsol}, which provides an intuitive explanation to the ``surprising'' dimension reduction observed in \cite{Balogh2007,Amini2010}. This dimension reduction is crucial to developing a comprehensive understanding of the contagion process.

Analyzing the trajectory of the ODEs in Sections \ref{sec:odeanalysisfinite}-\ref{sec:odeanalysisinf}, we propose a fixed point equation whose solution can be used to exactly determine the final fraction of the population that are eventual adopters, i.e., the size of the cascade---Theorem \ref{thm:odesmaininf} and Corollary \ref{cor:odesmaininf}. In particular, we prove the validity of the mean-field analysis of the contagion process over infinite trees, presented in Section \ref{sec:meanfield}. Furthermore, when the fixed point equation has multiple solutions, we identify the correct solution and provide an algorithmic means to calculate it. For general thresholds, we also provide a sharp characterization of the contagion threshold---the condition on the thresholds for which a contagion occurs with a finite set of seed vertices---in terms of the Perron-Frobenius eigenvalue of an associated matrix---Theorem \ref{thm:contagion}. Specializing to Poisson degree distributions with symmetric community structure and linear thresholds of \cite{Morris2000,Watts2002}, we prove that the existence of communities does not matter for global properties like the contagion threshold---Corollary \ref{cor:poissym}. This last set of results are presented as a rigorous counterpart to the many heuristic and empirical results in the literature \cite{Gleeson2008,Galstyan2007} for Poisson degree distributions and linear thresholds.

Using the fixed point characterization of the size of the cascade, we then study the impact of the community structure on the viral seeding of vertices in Section \ref{sec:numres}. We develop a gradient-based heuristic seeding strategy to maximize the size of the cascade given budget constraints. Empirically, our algorithm suggests that commonly held beliefs in the literature, which point to the best strategy being to seed over vertices with the highest number of neighbors, may be misguided. Most notably, we can demonstrate many cases wherein our seeding algorithm achieves a global cascade reaching almost all nodes while seeding over the vertices with the highest degree fails to spread much further from the seeds.

\medskip
\noindent{\bf Proof Technique:}
The basic idea behind the proof is to couple the evolution of the cascade with the realization of the random graph \cite{Balogh2007,Amini2010,Lelarge2012}. This is done either by exploring all neighbors of a uniformly selected active vertex \cite{Lelarge2012} or by realizing their connections one by one \cite{Balogh2007,Amini2010}. We adopt the latter approach but note that the resulting process evolves slower as we explore edges instead of vertices. This results in a less correlated structure that makes the analysis possible in the presence of community structure. We present the details of this coupling in Section \ref{sec:markovproc}.

Given the above coupling, the next step is to approximate the evolution of the process. Following the ideas in \cite{Wormald1995,Kurtz1970} and similar to \cite{Balogh2007,Amini2010}, one may attempt to trace the cascade using a system of ODEs. However, in the presence of community structure, the resulted ODEs are intractable due to the interconnected nature of the problem. In particular, all variables of the associated ODEs depend on each other, and the dimension of ODEs increases unboundedly as the number of vertices goes to infinity. This is in contrast with the analysis of single community \cite{Balogh2007,Amini2010} in which the trajectory of the evolution of each variable can be studied separately using a natural ordering.

To resolve this issue, we study the following truncated versions of the problem: $(1)$ all vertices with large degrees are initial adopters, and $(2)$ vertices with large degrees that are not seeded initially will never adopt the new technology. As it will become clear in our exposition, the dimension of the associated ODEs will remain bounded for truncated processes. Using a natural coupling, we then show that the final fraction of adopters in the original process is sandwiched between the same quantities given for these two truncated versions. Hence, we only need to study truncated processes to characterize the asymptotic behavior of the cascade. See Section \ref{sec:odeapprox} for details.

The next hurdle we address is the solution of the ODEs associated with a truncated process in the presence of community structure. In the case of one community, the corresponding ODEs are surprisingly simple: their solution is characterized by the solution of a one-dimensional ODE \cite{Balogh2007,Amini2010,Lelarge2012}. In Section \ref{sec:odeprobsol}, we present an intuitive probabilistic approach to solve the ODEs, which also explains the dimension reduction observed in the case of one community. Specifically, we show that the solution of the ODEs can be obtained by solving a much simpler $k^2$-dimensional set of ODEs, where $k$ is the number of communities. This dimension reduction is crucial to developing a comprehensive understanding of the contagion process.

The final and most critical part of our work is to establish the connection between the final fraction of adopters and the equilibrium point of the ODEs rigorously. For any population of size $n < \infty$, one can use the ODEs to approximate the evolution of the cascade in the corresponding truncated process. This approximation is valid before getting too close to the boundary of the region where Lipschitzness holds and is applicable only for a constant number of updates. Hence, using the ODEs to characterize the final fraction of adopters and its asymptotic behavior as $n$ increases without bound needs extra care; this type of analysis is prone to an unjustified interchange of limits. We would also like to point out that the analysis of \cite{Amini2010} appears to suffer from this issue. In particular, the author did not properly address the interchange of limits: while the analysis shows that for any finite $n$, the fraction of adopters gets close to a specific fixed point of an associated ODE, characterizing the final proportion of adopters and its asymptotic behavior need more work.

In Sections \ref{sec:odeanalysisfinite}-\ref{sec:odeanalysisinf}, we rigorously establish the connection between the final proportion of adopters in a truncated process and the fixed point of the associated ODEs. The core idea is to augment the graph by adding two active vertices with high degrees after running the process for some time; the degrees are proportional to the Perron-Frobenius eigenvector of the Jacobian matrix of the associated set of ODEs at its equilibrium point. We then pair some of these newly added active half-edges with other half-edges and approximate the state of the augmented process using a new set of ODEs. Coupling the truncated process with the augmented process, we characterize the state of the coupled truncated process after removing these newly added half-edges. In particular, we show that all active half-edges of the coupled truncated process have already been explored with high probability, given the fixed point of the ODEs associated with it is stable. This argument results in a probabilistic bound for the stopping time of the process for any finite value of $n$. Analyzing the asymptotic behavior of these bounds, we show a concentration of the stopping time of the truncated process and hence, the connection between asymptotics of the truncated process and the fixed point of the associated ODEs. This is the main result of the paper, and it is presented in Theorem \ref{thm:odesmaininf}.

\medskip
\noindent{\bf Related Works:} The threshold model~\cite{Granovetter1978,Schelling1978,Naylor1990,Watts2007,Moore2014} is a well accepted model for explaining the adoption of a new technology, opinion or behavior in a population that interacts via a social network. The linear threshold model, where the threshold is a function of the degree, was analyzed for the contagion threshold for specific graphs in~\cite{Morris2000}, and using heuristically derived formulae for single community random graphs in~\cite{Watts2002,Lopez2008}. The results on the single community random graphs were rigorously proved using branching processes in~\cite{Lelarge2012}, where the importance of pivotal players (those whose degree is low enough that one neighbor will make them adopt the new behavior) was identified and studied. Similar results were derived using the differential equation method in~\cite{Balogh2007,Amini2010}, and in~\cite{Como2019} for the non-permanent adoption model.

The threshold model has been studied for networks with communities, but using heuristically derived mean-field approximations and approximate differential equations~\cite{Galstyan2007,Gleeson2007,Gleeson2008,Lopez2008,Galstyan2009}. In these studies, it was numerically shown in~\cite{Gleeson2008,Galstyan2007} for the linear threshold model that the community structure leads to a different dynamic in terms of the evolution of the cascade itself. It is important to note that the authors in these works postulate both the mean-field equation and the differential equations in an \textit{ad hoc} manner without a formal proof. This is particularly the case for the multi-community work in~\cite{Gleeson2008,Galstyan2007} where the authors combine the adoption processes in the different communities without proper mathematical justification.

The problem of maximizing influence propagation in networks, by targeting certain influential vertices that have the potential to influence many others, has been an important follow-up problem~\cite{Granovetter1978,Schelling1978,Naylor1990,Watts2007} once the impact of a social network on behavior adoption was discovered.
%From the algorithmic standpoint, this selection problem can be stated as follows~\cite{Domingos:2001:MNV:502512.502525}: given a social network, an influence model (say the threshold model), and a set of vertices $S$, let $\sigma(S)$ be the expected number of vertices that will be activated by the end of the influence propagation process. Then, for a given budget $b$, the influence maximization problem is concerned with finding the set $S$ of size $b$ that maximizes the return $\sigma(S)$.
%The objective is to find the optimal subset of vertices with a pre-specified size to target initially so as to maximize the eventual influence propagation in the entire population.
While this problem is known to be NP hard for many influence models, several approximate methods have been designed, see e.g.,~\cite{Kempe2003,Mossel2010}.
% An important result established in~\cite{Kempe:2003:MSI:956750.956769,doi:10.1137/080714452} states that for a class of models, a simple greedy algorithm, which selects the next best candidate vertex, yields a solution that is guaranteed to be within $(1-1/e)$ of the optimal, under a few crucial assumptions.
%Although this result is general, it is very sensitive to some crucial assumptions. To be more specific, let us focus on the so called linear threshold models where a vertex $v$ is activated whenever the fraction of its active neighbors exceeds some predefined threshold $\theta_v$. In order to apply the results in~\cite{Kempe:2003:MSI:956750.956769,doi:10.1137/080714452}, we need to assume that $\theta_v$ are i.i.d. uniform in the interval $[0,1]$. In the Appendix~\ref{app:greed} we show by an example how these results break if $\theta_v$ is assumed to be fixed.
%
A contrasting strategy to identifying and targeting influential vertices is to use viral marketing~\cite{Schelling1978,Naylor1990,Vega2007}. % by seeding a new behavior in a population. %, i.e., targeting a certain fraction of the population.
%Often this can be network-structure agnostic, i.e., target a certain fraction of the population, but it can also use properties such as the degrees of the vertices, without the knowledge of the fine or local network structure, i.e., the location of the vertices.
A randomized version of viral marketing, also referred to as seeding or advertising in the paper, was studied in~\cite{Lelarge2012,Amini2010} where the resulting cascade was precisely identified. The results in~\cite{Lelarge2012} also suggested that targeting higher degree vertices is a better seeding strategy over degree-unaware random seeding. With community structure,~\cite{Gleeson2007,Gleeson2008,Galstyan2009} showed using heuristic analysis methods that the seeding strategies could be dramatically different from the one-community optimal strategies. Typically asymmetric seeding strategies, wherein the seeding is principally carried out in one community over another, were shown to perform better than more uniform (over the communities) seeding strategies. %For e.g.,~\cite{gleeson2008cascades} shows a four community example where seeding in community 1 produces a global cascade but seeding in community 4 doesn't. Again it should be emphasized that since the analysis of~\cite{gleeson2008cascades,galstyan2007cascading}  is not mathematically rigorous, the results can be questioned when not accompanied by large-scale simulations.

\noindent{\bf Organization:} The remainder of this paper is organized as follows. We present our model in Section \ref{sec:model}.  In Section \ref{sec:meanfield}, we present a mean-field approximation of the adoption process, whose validity is then established in Sections \ref{sec:markovproc} through \ref{sec:odeanalysisinf}: in Section \ref{sec:markovproc}, we construct a Markov process coupling the evolution of the adoption process with the process generating the random graph; in Section \ref{sec:odeapprox}, we present two truncated versions of this process which are then approximated using a set of ODEs; we then provide a probabilistic approach to solve this set of ODEs in Section \ref{sec:odeprobsol}; the analysis of the trajectory of the ODEs is presented in Sections \ref{sec:odeanalysisfinite}; we establish the connection between the asymptotics of the cascade process and the ODEs in Section \ref{sec:odeanalysisinf}. We discuss the results on the contagion threshold for general thresholds in Section \ref{sec:contagion}. Many results are then specialized to the case of Poisson degree distributions in Section \ref{sec:poissdeg}. We present numerical results and discuss the optimal seeding strategy using a heuristic policy in Section \ref{sec:numres}. %Finally, we present some open problems in Section \ref{sec:openprob}.% and conclude in Section \ref{sec:conclusion}.

\medskip
\noindent{\bf Notation:} Random variables are denoted by capital letters (sometimes using a bold typeset, too); realizations or deterministic quantities are in small letters. Vectors are denoted by using a bold typeset and individual components without it. Adhering to game-theoretic notation, a vertex's community is denoted by $j\in\{1,2\}$ and the other community by $-j=\{1,2\}\setminus\{j\}$. Multigraphs are denoted by an asterisk and simple graphs without one. The words ``community'' and ``side'' are used interchangeably. $\mathbb{R}_+$ denote the set of non-negative real numbers, $\mathbb{Z}_+$ denote the set of non-negative integers, and $\mathbb{Z}_{++} = \{1,2,\cdots,\}$. We say that the set of events $\{A_n\}_{n\in\mathbb{Z}_{++}}$ holds with high probability if $\lim_{n\to\infty}\prob(A_n)=1$. The superscript/subscript $(j'\gets j)$ on a parameter denotes that the parameter is associated with half-edges/vertices in community $j$ that can be paired with half-edges/vertices in community $j'$. The superscript $(j)$ is used to denote that the parameter is associated with vertices in community $j$.
	\section{Mathematical Model}\label{sec:model}
	Consider a set $[n]=\{1,\dotsc,n\}$ of agents that are organized into two communities, community 1 $\{1,2,\dotsc, n_1\}$ and community 2 $\{n_1+1,\dotsc, n\}$ with $n_2:=n-n_1$ individuals. Assume that we are given three sequences of non-negative integers: $\mathbf{d}^n_1=(d_{1,i}^{n})_1^{n_1}$, $\mathbf{d}^n_2=(d_{2,i}^{n})_{n_1+1}^n$, and $\mathbf{d}^n_m=(d_{m,i}^{n})_1^n$, which satisfy the following conditions:
1) $\sum_{i=1}^{n_1} d_{1,i}^n$ is even;
2) $\sum_{i=n_1+1}^{n} d_{2,i}^n$ is even; and
3) $\sum_{i=1}^{n_1} d_{m,i}^n = \sum_{i=n_1+1}^{n} d_{m,i}^n$. The sequence $\mathbf{d}^n_j$ is the degree sequence of the sub-graph for community $j$ for $j\in\{1,2\}$ and $\mathbf{d}^n_m$ is the degree sequence of the bipartite graph connecting the two communities.

Construct a two-community random multigraph (allowing for self-loops and multiple links) with given degree sequences $\mathbf{d}^n_1$, $\mathbf{d}^n_2$ and $\mathbf{d}^n_m$ generated by the configuration model~\cite{Bollobas2001} as the concatenation of $G^*(n_1,\mathbf{d}^n_1)$, $G^*(n_2,\mathbf{d}^n_2)$ (both generated via the configuration model) and a random bipartite multigraph $G^*(n_1,n_2,\mathbf{d}^n_m)$: generate half-edges for each vertex corresponding to the different degree sequences and combine the half-edges into edges by a uniform random matching of the set of half-edges of each sequence. Conditioned on the random multigraphs and the random bipartite graph being simple graphs, we obtain uniformly distributed random graphs $G(n_1,\mathbf{d}^n_1)$, $G(n_2,\mathbf{d}^n_2)$, and $G(n_1,n_2,\mathbf{d}^n_m)$ with the given degree sequences.
The concatenation of these produces a simple two-community graph $G(n,\mathbf{d}^n_1,\mathbf{d}^n_2,\mathbf{d}^n_m)$ with the desired distributions. In Definition~\ref{def:regcon_graph} we impose standard regularity assumptions~\cite{Janson2014} on the degree sequences so that the resulted multigraphs are simple with positive probability. We assume that $\beta_1(n)\coloneqq n_1/n \xrightarrow{n\to\infty} \beta$ (equivalently $\beta_2(n) \coloneqq n_2/n  \xrightarrow{n\to\infty} 1-\beta$). The stochastic block model~\cite{Mossel2012} is a prototypical example of a two-community graph.

Following Lelarge~\cite{Lelarge2012}, we analyze the threshold model of Morris~\cite{Morris2000} and Watts~\cite{Watts2002} on the two-community random graph model described above. In this model, vertices have the choice between two types of opinions/technologies, A and B; we often also use ``inactive'' to denote type A and ``active'' to denote type B. All vertices initially start in type A, i.e., are inactive. Each vertex has a threshold that is a function of its community and degrees (in the same community and across to the other community); the value of the threshold is fixed and allowed to be any non-negative real number. If a vertex finds that the number of its neighbors (across both communities) who have chosen type B is greater than its threshold, then it will permanently choose to switch to type B. Again following \cite{Lelarge2012} we initially seed vertices with type B using a Bernoulli random variable (1 implying that a vertex gets seeded with type B) that is independently chosen with the mean depending on the vertex's parameters, namely, community and degrees. Note that a degree and/or community-unaware seeding strategy would imply an appropriate uniformity in the means of the seeding random variables. After the seeding process is completed, the remaining vertices then react to the seed vertices and decide whether to adopt type B.  This process continues until a final state of the vertices is reached. A cascade is said to happen if the number of vertices adopting type B is substantially greater than the seed set.

	%\section{Literature Review}\label{sec:review}
	%\input{Sections/review}

	\section{Mean-Field Approximation}\label{sec:meanfield}
	We start by presenting a mean-field approximation of the process of adoption of type B, i.e., becoming active, in a typical simple graph generated through the configuration model described in Section~\ref{sec:model}. We comment that the goal of this section is not to pin down the final proportion of the adopters; instead, we aim to provide an approximation based on a heuristic argument.

The graphs that we consider are locally tree-like \cite{Remco2016}  so that the structure up to any finite depth when viewed from a uniformly selected vertex of the graph is a tree with high probability. Therefore, the local structure of a simple graph produced by the configuration model converges~\cite{Aldous2007,Aldous2004,Bordenave2016} to a rooted unimodular Galton-Watson Multi-type Tree ($\text{GWMT}_*$).
In a rooted unimodular Galton-Watson Tree, the degree distribution of any non-root vertex is the size-biased/sampling-biased distribution of the root: for a random variable $D\in\mathbb{Z}_+$ with distribution $\prob(\cdot)$ and finite mean $\expect[D]$, the size-biased/sampling-biased distribution $\prob_*(\cdot)$ is given by $\prob_*(d)\coloneqq d \prob(d)/\expect[D]$ for all $d\in \mathbb{Z}_+$. In the case of $\text{GWMT}_*$, the degree distribution of each child depends on the community of its parent. The joint degree distribution of the root vertex is $\prob_{j,m}$ if the community of the root is $j\in \{1,2\}$. The degree distributions of each child is then given by the size-biased/sampling-biased distribution for the community of the parent and the regular distribution for the other community. In particular, if the parent is in community $j\in \{1,2\}$ and the child vertex is in community $j$ too, then the joint degree distribution is the size-biased distribution $\prob_{j*,m}$ given by $\prob_{j*,m}(d_j,d_{-j}) \coloneqq d_j\prob_{j,m}(d_j,d_{-j})/\sum_{k,k'} k \prob_{j,m}(k,k')$ for all $d_j,d_{-j}\in\mathbb{Z}_+$; on the other hand, if the parent is in community $j\in \{1,2\}$ and the child vertex is in community $-j\in\{1,2\}\setminus\{j\}$, then the joint degree distribution of the child is the size-biased distribution $\prob_{-j,m*}$ given by $\prob_{-j,m*}(d_{-j},d_{j}) \coloneqq d_j\prob_{-j,m}(d_{-j},d_{j})/\sum_{k,k'} k'\prob_{-j,m}(k,k')$ for all $d_j,d_{-j}\in\mathbb{Z}_+$.
We denote a random variable with the size-biased distribution by $D^*+1$ where $D^*$ takes values in $\mathbb{Z}_+$. For a Poisson random variable with parameter $\lambda > 0$, i.e., $D\sim \mathrm{Poi}(\lambda)$, we have $D^* \sim \mathrm{Poi}(\lambda)$, so that the size-biased/sampling-biased distribution is a shifted Poisson distribution. This is the only distribution with this property.

\begin{figure}[htbp]
\begin{center}
 \includegraphics[width=0.9\textwidth]{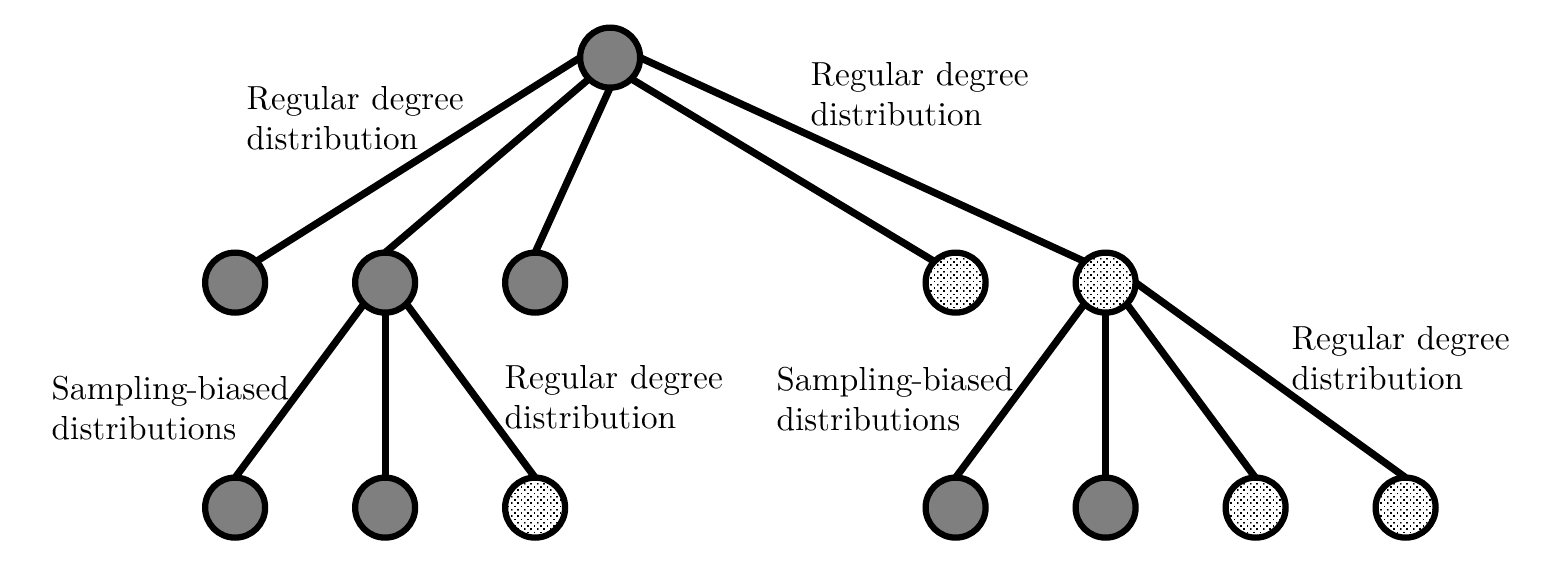}
\caption{Illustration of the limiting rooted unimodular Galton-Watson Multi-type Tree. Solid circles denote vertices in community 1, and dotted circles denote vertices in community 2.}
\label{fig:GWMIT}
\end{center}
\end{figure}

An example of the limiting rooted $\text{GWMT}_*$ is shown in Figure~\ref{fig:GWMIT} where the root vertex is in community 1.

Assume that we have a rooted $\text{GWMT}_*$ (with root vertex $\psi$) denoted by $T_\psi$. For a vertex $l\neq\psi$ let $l_p$ be its parent, indicated by $(l_p\gets l)$, and $T_{(l_p\gets l)}$ be the sub-tree rooted at $l$ when the link $(l_p,l)$ is excised. Then assuming that $l_p$ is inactive, state of vertex $l$ only depends on the state of her children in sub-tree $T_{(l_p \gets l )}$. Next, we define a few random variables that will aid in describing the mean-field approximation. \\
$X^{(j)}_\psi$: Bernoulli r.v; $=1$ if root vertex $\psi$ of the rooted $\text{GWMT}_*$ is on side $j$ and inactive. \\
$Y^{(j\gets j)}_l$: Bernoulli r.v; $=1$ if vertex $l$($\neq \psi$) and its parent ${l_p}$ are both on side $j$ and vertex $l$ is inactive on $T_{({l_p} \gets l)}$.\\
$Y^{(j\gets-j)}_l$: Bernoulli r.v; $=1$ if vertex $l$($\neq \psi$) is on side $-j$ and its parent $l_p$ is on side $j$ and vertex $l$ is inactive on $T_{({l_p}\gets l)}$.\\
$\bs{\alpha}^{(j)}_l$: Bernoulli r.v; $=1$ if vertex $l$ on side $j$ is a seed vertex.\\
$K_l^{(j)}$: Threshold of vertex $l$ on side $j$ that is determined by number of its neighbors in either community, i.e., by the degrees of the vertex $l$ on side $j$ and $-j$.\\
Then we can write down the following equations:\\
\noindent (i) A non-root vertex $l\neq \psi$ remains inactive on $T_{({l_p}\gets l)}$, if it is not seeded initially and the number of her children who are active does not exceed her threshold, i.e.,
\begin{align}
Y^{(j\gets j)}_l &= \left(1 - \bs{\alpha}^{(j)}_l\right) \mathbf{1}\Big\{\sum_{i\longrightarrow l} \left(1-Y^{(j\gets j)}_i\right) + \sum_{i\longrightarrow l} \left(1-Y^{(j\gets-j)}_i\right) \leq K^{(j)}_l\Big\}, \label{eq:meanfiled_Yjj}\displaybreak[0]\\
Y^{(j\gets-j)}_l &= \left(1 - \bs{\alpha}^{(-j)}_l\right) \mathbf{1}\Big\{\sum_{i\longrightarrow l} \left(1-Y^{(-j\gets -j)}_i\right) + \sum_{i\longrightarrow l} \left(1-Y^{(-j\gets j)}_i\right) \leq K^{(-j)}_l\Big\} \label{eq:meanfiled_Yj-j},
\end{align}
where $\mathbf{1}\{O\}$ is the indicator function of set $O$. \\
\noindent (ii) Root vertex $\psi$ (on side $j$) remains inactive if it is not seeded initially and the number of her active children falls below her threshold, i.e.,
\begin{align}
X^{(j)}_\psi &= \left(1 - \bs{\alpha}^{(j)}_\psi\right)\mathbf{1}\Big\{\sum_{i\longrightarrow \psi} \left(1-Y^{(j\gets j)}_i\right)  + \sum_{i\longrightarrow \psi} \left(1-Y^{(j\gets-j)}_i\right) \leq K^{(j)}_\psi\Big\} \label{eq:meanfield_Xj}.
\end{align}

For the mean-field approximation it is assumed that the random variables $Y^{(1\gets 1)}_l$, $Y^{(1\gets 2)}_l$, $Y^{(2\gets 1)}_l$, and $Y^{(2\gets 2)}_l$ for $l\neq\psi$ are, respectively, identically distributed when considering $l$ as the variable and keeping $(j\gets j)$ or $(j\gets -j)$ fixed. Moreover, it is assumed that all these random variables are mutually independent. These random variables are then related via the following Recursive Distributional Equations (RDEs), where equality below should be interpreted in terms of distribution.
\begin{align}
\begin{split}
& \widetilde{Y}^{(j\gets j)} \overset{d}{=} \left(1 - \bs{\bar{\alpha}}^{(j)}(D_j^*+1,D_m)\right) \mathbf{1}\Big\{\sum_{i=1}^{D^*_j} \left(1-\widetilde{Y}^{(j\gets j)}_i\right)
\\
& \qquad \qquad
+ \sum_{i=1}^{D_m} \left(1-\widetilde{Y}^{(j\gets-j)}_i\right) \leq K^{(j)}(D_j^*+1,D_m)\Big\},
\end{split} \label{eq:meanfiled_Yjjdist}\displaybreak[0]\\
\begin{split}
& \widetilde{Y}^{(j\gets-j)} \overset{d}{=} \left(1 - \bs{\bar{\alpha}}^{(-j)}(D_{-j},D_m^*+1)\right) \mathbf{1}\Big\{ \sum_{i=1}^{D_{-j}} \left(1-\widetilde{Y}^{(-j\gets -j)}_i\right)
\\
& \qquad \qquad
+ \sum_{i=1}^{D^*_m} \left(1-\widetilde{Y}^{(-j\gets j)}_i\right) \leq K^{(-j)}(D_{-j},D_m^*+1)\Big\},
\end{split}\label{eq:meanfiled_Yj-jdist}\displaybreak[0]
\end{align}
where for every $j\in \{1,2\}$, $\widetilde{Y}^{(j\gets j)}$ and $\widetilde{Y}^{(j\gets j)}_i$ as well as $\widetilde{Y}^{(j\gets-j)}$ and $\widetilde{Y}^{(j\gets-j)}_i$ are {\it i.i.d.} copies (Bernoulli random variables with unknown parameters). We also have a set of random variables: $D_j$ is a random variable with the community $j$ degree distribution, $D_j^*+1$ is a random variable with the size-biased distribution of $D_j$, $D_m$ has inter-community degree distribution, and $D_m^*+1$ is a random variable with the size-biased distribution of $D_m$; the joint distribution of $(D^*_j+1,D_m)$ is given by $\prob_{j*,m}$, and the joint distribution of $(D_j,D^*_m+1)$ is given by $\prob_{j,m*}$ (for all $d_j,d_m \in\mathbb{Z}_+$, we have $\prob_{j*,m}(d_j,d_m) = d_j\prob_{j,m}(d_j,d_m)/\expect[D_j]$ and $\prob_{j,m*}(d_j,d_m) = d_m\prob_{j,m}(d_j,d_m)/\expect[D_m]$).
We have also assumed, without loss of generality, that the seeding Bernoulli random variables have means that depend on the community and the degrees of the vertex, namely, $\alpha_j(d_j,d_{-j})$ for $j\in \{1,2\}$ and $d_j, d_{-j} \in \mathbb{Z}_+$. We also assume that threshold random variables are deterministic functions of the community and degrees of the vertex, namely, $K_j(d_j,d_{-j})$ for $j\in \{1,2\}$ and $d_j, d_{-j} \in \mathbb{Z}_+$. These are then used to construct the random variables $\bs{\bar{\alpha}}^{(j)}(D_j^*+1,D_m)$, $\bs{\bar{\alpha}}^{(-j)}(D_{-j},D_m^*+1)$, $K^{(j)}(D_j^*+1,D_m)$ and $K^{(-j)}(D_{-j},D_m^*+1)$.

Since we have RDEs with Bernoulli random variables, we can equivalently obtain the solutions by taking expectations and solving for the means of the underlying random variables.  We set $\expect[X^{(j)}_\psi]= \phi_j$, $\expect[\widetilde{Y}^{(j\gets j)}]= \mu^{(j\gets j)}$ and $\expect[\widetilde{Y}^{(j\gets-j)}] = \mu^{(j\gets-j)}$. Taking expectation in \eqref{eq:meanfiled_Yjjdist}-\eqref{eq:meanfiled_Yj-jdist} and then \eqref{eq:meanfield_Xj} yields
\begin{align}
\begin{split}
\mu^{(j\gets j)} &= \sum_{u_j+u_{-j} \leq
K_j(d_j,d_{-j}) } \prob_{j*,m}(d_j,d_{-j})(1-\alpha_j(d_j,d_{-j}))\\
& \qquad \qquad \qquad \qquad \times Bi(u_j;d_j-1,1 -\mu^{(j\gets j)})Bi(u_{-j};d_{-j},1 - \mu^{(j\gets-j)}) ,
\end{split}
 \label{eq:meanfield_mujj}\displaybreak[0]\\
\begin{split}
\mu^{(j\gets-j)}& = \sum_{u_j+u_{-j} \leq { K_{-j}(d_{-j},d_j) } } \prob_{-j,m*}(d_{-j},d_j)(1-\alpha_{-j}(d_{-j},d_j))\\
& \qquad \qquad \qquad \qquad \times Bi(u_j;d_j-1,1-\mu^{(-j\gets j)})Bi(u_{-j};d_{-j},1-\mu^{(-j\gets -j)}),
\end{split}
 \label{eq:meanfield_muj-j}\displaybreak[0]\\
\begin{split}
\phi_j &= \sum_{u_j+u_{-j} \leq K_j(d_j,d_{-j}) } \prob_{j,m}(d_j,d_{-j})(1-\alpha_j(d_j,d_{-j}))\\
& \qquad \qquad \qquad \qquad \times Bi(u_j;d_j,1 -\mu^{(j\gets j)})Bi(u_{-j};d_{-j},1 - \mu^{(j\gets-j)}),
\end{split}
  \label{eq:meanfield_phij}
\end{align}
where $Bi(k;n,p):={n \choose k}p^k(1-p)^{n-k}$ is the probability mass function of the binomial distribution.

To find the probability of a vertex in community $j\in \{1,2\}$ remaining inactive, i.e. $\phi_j=1$, one needs to first solve the fixed point equations \eqref{eq:meanfield_mujj}-\eqref{eq:meanfield_muj-j}, and then substitute the result into \eqref{eq:meanfield_phij}. For ease of understanding we write equations \eqref{eq:meanfield_mujj}-\eqref{eq:meanfield_phij} as follows:
\begin{align}
\bs{\mu} = \bs{F}(\bs{\mu}), \text{ and } \bs{\phi}=\bs{\Phi}(\bs{\mu}), \label{eq:meanfield_mu}
\end{align}
for functions $\bs{F}(\cdot)$ and $\bs{\Phi}(\cdot)$ defined component-wise via the right-hand sides of \eqref{eq:meanfield_mujj}-\eqref{eq:meanfield_muj-j}, and \eqref{eq:meanfield_phij}, respectively.

A basic question at this point is whether one can rigorously justify \eqref{eq:meanfield_mu}, particularly given the various independence and uniformity assumptions for the derivation. A few other questions also arise: i) Does a solution to \eqref{eq:meanfield_mu} exist? ii) Are there multiple solutions to \eqref{eq:meanfield_mu}? Numerically, we observed that there are many cases where \eqref{eq:meanfield_mu} has multiple solutions; and iii) Which solution should one pick if there are multiple solutions? Note that for every $\bs{\mu} \in [0,1]^4$ and $j\in \{1,2\}$, we have
\begin{align}
\begin{split}
\phi_j &= \sum_{d_j,d_{-j}} \prob_{j,m}(d_j,d_{-j})(1-\alpha_j(d_j,d_{-j})) \times \\
&~~~~ \sum_{u_j+u_{-j} \leq K_j(d_j,d_{-j})} Bi(u_j;d_j-1,1 -\mu^{(j\gets j)})Bi(u_{-j};d_{-j},1 - \mu^{(j\gets-j)})
\\
&\leq \sum_{d_j,d_{-j}} \prob_{j,m}(d_j,d_{-j})(1-\alpha_j(d_j,d_{-j})) =\prob(\bs{\alpha}^{(j)}_{\psi} = 0),
\end{split}
\end{align}
so that the seeding distribution gets automatically accounted in any solution of \eqref{eq:meanfield_mu}, and the final population of active vertices includes at least the seed vertices.

Before proceeding, we should again point out that equations of a similar form were heuristically postulated in the literature~\cite{Galstyan2007,Gleeson2007,Gleeson2008,Lopez2008,Galstyan2009}. An important contribution of our work is thus to rigorously prove the validity of \eqref{eq:meanfield_mu}, and to identify the correct solution to choose. As discussed in \cite{Lelarge2012}, the existence of multiple solutions and a lack of ``monotonicity" makes it extremely challenging to use the techniques developed in \cite{Aldous2004,Aldous2007} to prove the needed results.

	\section{Markov Process of Adoption}\label{sec:markovproc}
	As the first step to proving the validity of the mean-field equations, we construct a Markov process that couples the evolution of the adoption process with the process of generating the random graph using the configuration model.

The typical way to create a random graph with a given degree sequence $(d^{n}_i)_{1}^n$ using the configuration model is as follows: first label vertices of the graph $1,2,\dots,n$ such that vertex $i$ has $d_i^{n}$ half-edges sticking out of it; next iterate through all the unpaired half-edges so that at each step, two half-edges are paired randomly; and declare the final graph as the desired random graph. In our setting, when there are two communities, the basic idea of generating the random graph using the configuration model is exactly the same.
However, to analyze the adoption process, we work a little differently. We start by realizing the early adopter vertices using the seeding random variables. We set the early adopters to be active and make all their half-edges active. Any other vertex and its half-edges will initially be counted as inactive. We then run the adoption process and draw the random graph simultaneously by iterating through the active half-edges (if any).  At each iteration, we pick an active half-edge, i.e., a half-edge connected to an active vertex, and connect it to some other half-edge that belongs to the appropriate community. Then we remove both half-edges from the graph. Moreover, if the second half-edge belongs to an inactive vertex, we reduce its threshold by one. If the threshold of the inactive vertex becomes minus one after this change, we activate this vertex and also all the half-edges that are still connected to this vertex. Note that this process stops when all active half-edges have been omitted. In particular, the remainder of the graph (containing only inactive half-edges) is not realized (or can be realized but will not influence the contagion process).
This process is described in Algorithm~\ref{algorithm}.
%\normalem
\begin{algorithm}
	\SetAlgoLined
	\KwData{degree sequences and early adopters}
	\KwResult{sub-graph of the final random graph that contains all active vertices}
	initialization\;
	\While{there is an active half-edge}{
		randomly choose an active half-edge\;
		randomly choose another half-edge belongs to proper community \;
		omit two selected half-edges from the set of half-edges\;
		update the state of the inactive vertex (if any)\;
	}
	\caption{process of jointly generating the random graph and running the adoption process.}
	\label{algorithm}
\end{algorithm}
%\ULforem

We keep track of active half-edges, inactive vertices, and number of times that the process described in Algorithm~\ref{algorithm} picks half-edges from each community. The random variables associated with these quantities are given as follows:\\
$A_{j}(k)$: Number of active half-edges belonging entirely to community $j$ at time k.\\
$A^{(j)}_{m}(k)$: Number of active half-edges between the two communities belonging to vertices in community $j$ at time $k$.\\
$T_{j}(k)$: Number of times the algorithm visits community $j$ up to time k where a visit means removing two half-edges within the same community.\\
$I^{(j)}_{d_j,d_{-j},u_j,u_{-j}}(k)$: Number of inactive vertices in community $j$ with $d_j$ initially assigned half-edges corresponding to community $j$ where $u_j$ of them have been removed by time $k$, and similarly, $d_{-j}$ initially assigned half-edges corresponding to community $-j$ where $u_{-j}$ of them have been removed by time $k$. Note that $u_j +u_{-j}\leq K_j(d_j,d_{-j})$, $0\leq u_j\leq d_j$, and $0\leq u_{-j}\leq d_{-j}$.\\
It is easily verified that $\{ X^n(k) \}_{k\in \mathbb{Z}_+}$ is a discrete-time Markov chain, where
\begin{align*}
X^n(k):=(A_j(k), A^{(j)}_m(k), T_j(k), I^{(j)}_{d_j,d_{-j},u_j,u_{-j}}(k), \dotsc),
\end{align*}
and $j\in\{1,2\}$. For ease of exposition we denote the number of edges entirely in community $j$ by $m_j(n)$ and the number of edges between the two communities by $m_m(n)$; these can be determined once the degrees have been realized.
\begin{remark}\label{rem:notationsimp}
	We use the following notational convention throughout the paper: whenever $I^{(j)}_{d_j,d_{-j},u_j,u_{-j}}$ appears as a member of some sequence, it represents all $I^{(j)}_{d_j,d_{-j},u_j,u_{-j}}$ for $j\in\{1,2\}$, $(d_j,u_j)\in\mathbb{Z}_+^2$, $(d_{-j},u_{-j})\in\mathbb{Z}_+^2$, $u_j\leq d_j$, $u_{-j}\leq d_{-j}$, and $u_j+u_{-j} \leq K_j(d_j,d_{-j})$. We use the same convention for $A_j(k)$, $A^{(j)}_m(k)$, and $T_j(k)$.
\end{remark}

The mean-field analysis~\cite{Wormald1995,Molloy1998,Wormald1999,Amini2010} proceeds by scaling both space and time by $n$ and considering the one-step drift of the scaled process. We will now present the one-step drift analysis of our Markov chain (for the unscaled variables). At each iteration, one of the following events will happen:
\begin{enumerate}[itemindent=0pt,label*=\arabic*.]
	\item Two active half-edges will be omitted. This event results in the half-edges being ``wasted", in a manner of speaking. Here two sub-cases are possible:
	\begin{enumerate}[label*=\roman*.]
		\item Both half-edges belongs to community $j$. This event happens with probability $$\frac{A_j(k)\left(A_j(k)-1\right)}{\left(A_1(k) + A_2(k) + A^{(1)}_m(k) + A^{(2)}_m(k)\right)\left(2m_j(n) - 2T_j(k) - 1\right)}.$$
		In this case, we should update the corresponding variables as follows:
		\begin{align*}
			A_j(k+1) = A_j(k) - 2, \;T_j(k+1) = T_j(k) + 1.
		\end{align*}
		\item Half-edges belongs to different sides. This event happens with probability $$\frac{2A^{(1)}_m(k) A^{(2)}_m(k)}{\left(A_1(k) + A_2(k) + A^{(1)}_m(k) + A^{(2)}_m(k)\right)\left(m_m(n)-(k-T_1(k)-T_2(k))\right)}.$$
		In this case, we should update the variables as follows:
		\begin{align*}
			A^{(2)}_m(k+1) = A^{(2)}_m(k) - 1, \;A^{(1)}_m(k+1) = A^{(1)}_m(k) - 1.
		\end{align*}
	\end{enumerate}
	\item One active half-edge and one inactive half-edge will be omitted, while the inactive half-edge belongs to an inactive vertex in community $j$. Four sub-cases arise here:
	\begin{enumerate}[label*=\roman*.]
		\item The inactive vertex belongs to $I^{(j)}_{d_j,d_{-j},u_j,u_{-j}}(k)$ and the active half-edges belongs to community $j$, while $K_j(d_j,d_{-j}) \geq u_j + u_{-j} + 1$. This event results in the threshold of the inactive vertex in community $j$ being lowered by $1$ owing to a vertex within its own community.
		This occurs with probability $$\frac{A_j(k)}{A_1(k) + A_2(k) + A^{(1)}_m(k) + A^{(2)}_m(k)}\times\frac{\left(d_j-u_j\right)I^{(j)}_{d_j,d_{-j},u_j,u_{-j}}(k)}{2m_j(n) - 2T_j(k) - 1}.$$
		In this case, we should update variables as follows:
		\begin{align*}
			& A_j(k+1) = A_j(k) - 1, \; I^{(j)}_{d_j,d_{-j},u_j,u_{-j}}(k+1) =I^{(j)}_{d_j,d_{-j},u_j,u_{-j}}(k) - 1, \\
			& T_j(k+1) = T_j(k) + 1, \; I^{(j)}_{d_j,d_{-j},u_j+1,u_{-j}}(k+1) =I^{(j)}_{d_j,d_{-j},u_j+1,u_{-j}}(k) + 1.
		\end{align*}
		\item The inactive vertex belongs to $I^{(j)}_{d_j,d_{-j},u_j,u_{-j}}(k)$ and the active half-edge belongs to community $j$, while $K_j(d_j,d_{-j}) = u_j + u_{-j}$. During this event, the inactive vertex and all its remaining half-edges become active. This is an important growth event for our process.
		This occurs with probability $$\frac{A_j(k)}{A_1(k) + A_2(k) + A^{(1)}_m(k) + A^{(2)}_m(k)}\times\frac{\left(d_j-u_j\right)I^{(j)}_{d_j,d_{-j},u_j,u_{-j}}(k)}{2m_j(n) - 2T_j(k) - 1}.$$
		Here we update the variables as follows:
		\begin{align*}
			& A_j(k+1) = A_j(k) - 1 + d_j - u_j - 1, \; T_j(k+1) = T_j(k) + 1,\\
			& A^{(j)}_m(k+1) = A^{(j)}_m(k) + d_{-j} - u_{-j},\;
			%\\ &
			I^{(j)}_{d_j,d_{-j},u_j,u_{-j}}(k+1) =I^{(j)}_{d_j,d_{-j},u_j,u_{-j}}(k) - 1.
		\end{align*}
		\item The inactive vertex belongs to $I^{(j)}_{d_j,d_{-j},u_j,u_{-j}}(k)$ and the active half-edge comes from the other community, while $K_j(d_j,d_{-j}) \geq u_j + u_{-j} + 1$. During this event, the threshold of an inactive vertex is reduced by a vertex from the other community. This occurs with probability $$\frac{A^{(-j)}_m(k)}{A_1(k) + A_2(k) + A^{(1)}_m(k) + A^{(2)}_m(k)} \times\frac{\left(d_{-j}-u_{-j}\right)I^{(j)}_{d_j,d_{-j},u_j,u_{-j}}(k)}{m_m(n) - (k - T_1(k) - T_2(k)) }.$$
		Here we update the variables as follows:
		\begin{align*}
			& A^{(-j)}_m(k+1) = A^{(-j)}_m(k) - 1,\;
			%\\ &
			I^{(j)}_{d_j,d_{-j},u_j,u_{-j}}(k+1) =I^{(j)}_{d_j,d_{-j},u_j,u_{-j}}(k) - 1,\\
			& I^{(j)}_{d_j,d_{-j},u_j,u_{-j}+1}(k+1) =I^{(j)}_{d_j,d_{-j},u_j,u_{-j}+1}(k) + 1.
		\end{align*}
		\item The inactive vertex belongs to $I^{(j)}_{d_j,d_{-j},u_j,u_{-j}}(k)$ and the active half-edge comes from the other community, while $K_j(d_j,d_{-j}) = u_j + u_{-j}$. This is another important growth event for our process wherein an inactive vertex becomes active owing to a vertex from the other community.
		This occurs with probability $$\frac{A^{(-j)}_m(k)}{A_1(k) + A_2(k) + A^{(1)}_m(k) + A^{(2)}_m(k)}\times\frac{\left(d_{-j}-u_{-j}\right)I^{(j)}_{d_j,d_{-j},u_j,u_{-j}}(k)}{m_m(n) - (k - T_1(k) - T_2(k)) }.$$
		Here we update the variables as follows:
		\begin{align*}
			& A^{(-j)}_m(k+1) = A^{(-j)}_m(k) - 1,\; A_j(k+1) = A_j(k) + d_j - u_j, \\
			& A^{(j)}_m(k+1) = A^{(j)}_m(k) + d_{-j} - u_{-j} - 1, \;
			%\\ &
			I^{(j)}_{d_j,d_{-j},u_j,u_{-j}}(k+1) =I^{(j)}_{d_j,d_{-j},u_j,u_{-j}}(k) - 1.
		\end{align*}
	\end{enumerate}
\end{enumerate}

Finally, note that these random variables satisfy balance equations given by the realization of degrees. For $j\in \{1, 2\}$ we have
\begin{align}
\begin{split}
&~A_j(k) +\!\! \sum_{u_j+u_{-j} \leq K_j(d_j,d_{-j})}\!\! (d_j-u_j)I^{(j)}_{d_j,d_{-j},u_j,u_{-j}}(k) = 2m_j(n)- 2T_j(k), \\
&~A^{(j)}_m(k) +\!\! \sum_{u_j+u_{-j} \leq K_j(d_j,d_{-j})} \!\!(d_{-j}-u_{-j})I^{(j)}_{d_j,d_{-j},u_j,u_{-j}}(k) = m_m(n) - (k - T_j(k) - T_{-j}(k)),
\end{split}
\label{eq:balanceeq}
\end{align}
where the summations above are understood to be over both the degrees $(d_{j}, d_{-j})\in\mathbb{Z}_+^2$ and the used half-edges $(u_{j},u_{-j})\in\mathbb{Z}_+^2$ meeting the constraint listed underneath.
\begin{remark}\label{rem:trackvar}
	The the coordinates of the discrete-time Markov chain $\{X^n(k)\}_{k\in\mathbb{Z}_+}$ are dependent because of the balance equations. Hence, we only need to keep track of $T_j$ and $I^{(j)}_{d_j,d_{-j},u_j,u_{-j}}$ for $j\in\{1,2\}$ and different values of $d_j$, $d_{-j}$, $u_j$, and $u_{-j}$.
\end{remark}
The one-step drifts of the unscaled random variables are obtained by summing over all possible events, given the current state of the Markov chain. The details can be found in Appendix \ref{app:onestepdrift}.

We conclude this section by stating the regularity conditions on the degree sequences and some consequences of these conditions. Recall that $\mathbf{d}^n_1=(d_{1,i}^{n})_1^{n_1}$ and $\mathbf{d}^n_2=(d_{2,i}^{n})_{n_1+1}^n$ are the corresponding degree sequence of the sub-graphs for community $1$ and $2$, respectively, and $\mathbf{d}^n_m=(d_{m,i}^{n})_1^n$ is the degree sequence of the bipartite graph connecting the two communities; these sequences satisfy the following basic conditions:
1) $\sum_{i=1}^{n_1} d_{1,i}^n$ is even;
2) $\sum_{i=n_1+1}^{n} d_{2,i}^n$ is even; and
3) $\sum_{i=1}^{n_1} d_{m,i}^n = \sum_{i=n_1+1}^{n} d_{m,i}^n$.
The size of the community $1$ is $n_1$ and the size of the community $2$ is $n_2= n-n_1$. We assume that $\beta_1(n)\coloneqq  n_1/n \xrightarrow{n\to\infty} \beta$, and $\beta_2(n) =  n_2/n \xrightarrow{n\to\infty} 1-\beta$. We also assume that $m_j(n)$ for $j\in\{1,2\}$ and $m_m(n)$ grows unboundedly as $n\to\infty$, where $m_j(n)$ denotes the number of edges in community $j$ and $m_m(n)$ denotes the number of edges between the two communities.

\begin{definition}\label{def:regcon_deg}
	We say the degree sequences $\mathbf{d}^n_1$, $\mathbf{d}^n_2$, and $\mathbf{d}^n_{m}$ satisfy the {\bf degree regularity conditions} if the following conditions hold:
	\begin{enumerate}[wide]
		\setlength{\itemindent}{\parindent}
		\item The empirical degree distribution of vertices in community $j\in\{1,2\}$  converges to some joint probability distribution:
		\begin{align*}
		|\{i:d^{n}_{j,i} = r, d^{n}_{m,i} = s\,\text{and } i\in\text{community }j\}|/n_j \to \prob_{j,m}(r,s)\text{ for every }r,s\geq 0.
		\end{align*}
		\item The average degree of vertices, and the ratio of vertices in each community converges:
		\begin{gather*}
			\lambda_{j,j}(n) \coloneqq 2m_j(n)/n_j\to\lambda_{j,j},\qquad \lambda_{j,m}(n)\coloneqq m_m(n)/n_j \to \lambda_{j,m} ,\\
			 \beta_j(n)\coloneqq n_j/n \to \beta_j,
		\end{gather*}
		where $\beta_1 = \beta$, $\beta_2 = 1-\beta$, and following the first condition, $\lambda_{j,j}$ and $\lambda_{j,m}$ for $j\in\{1,2\}$ are given as follows:
		\begin{align*}
		\lambda_{j,j} = \sum_{r\geq 0}r\,\prob_{j,m}(r,\mathbb{Z}_+) \in (0,\infty), \qquad
		\lambda_{j,m}= \sum_{s\geq 0}s\,\prob_{j,m}(\mathbb{Z}_+,s) \in (0,\infty),
		\end{align*}
		where $\prob_{j,m}(\mathbb{Z}_+,s) \coloneqq \sum_r \prob_{j,m}(r,s)$, and $\prob_{j,m}(r,\mathbb{Z}_+) \coloneqq \sum_s \prob_{j,m}(r,s)$ for $j\in\{1,2\}$.
		\item The second moment of the empirical degree distribution grows smaller than $n$:
		\begin{align*}
		&\frac{\sum_{i}(d^{n}_{j,i})^2}{n_j} = o(n_j) \text{ for }j\in\{1,2\},\allowdisplaybreaks\\
		&\frac{\sum_{i=1}^{n_1}(d^{n}_{m,i})^2}{n_1} = o(n_1),\text{ and }\frac{\sum_{i=n_1+1}^{n}(d^{n}_{m,i})^2}{n_2} = o(n_2)
		\end{align*}
	\end{enumerate}
\end{definition}
\begin{definition}\label{def:regcon_graph}
	We say the degree sequences $\mathbf{d}^n_1$, $\mathbf{d}^n_2$, and $\mathbf{d}^n_{m}$ satisfy the {\bf graph regularity conditions} if they satisfy the degree regularity conditions and the following conditions also hold:
	\begin{enumerate}[wide]
		\setlength{\itemindent}{\parindent}
		\item The random multigraphs $G^*(n_1,\mathbf{d}^n_1)$ and $G^*(n_2,\mathbf{d}^n_2)$ are asymptotically simple random graphs with positive probability. Mathematically, the condition is
		\begin{align*}
		\sum_{i}(d^{n}_{j,i})^2 = O(m_j(n)) \text{ for }j\in\{1,2\},
		\end{align*}
		which implies $\liminf_{n\to\infty} \prob(G^*(n_j,\mathbf{d}^n_j)\text{ is simple}) > 0$ for $j\in\{1,2\}$~\cite[Theorem 1.1]{Janson2014}.
		\item The random bipartite multigraph $G^*(n_1,n_2,\mathbf{d}^n_m)$ is asymptotically a simple random bipartite graph with positive probability. Mathematically, the conditions are
		\begin{align*}
		&(\text{i})~\sum_{i=1}^{n_1}\sum_{i'=n_1+1}^n d^{n}_{m,i}(d^{n}_{m,i}-1) d^{n}_{m,i'}(d^{n}_{m,i'}-1) = O((m_m(n))^2),\\
		&(\text{ii})~\text{for any $M\geq 1$, }\\
		&\begin{aligned}
		\myquad[1]\sum_{i=\min(d_{2}^{\max},M)}^{n_1} d^{n}_{m,(i)} = \Omega(m_m(n)) \text{ and } \sum_{i=n_1+\min(d_{1}^{\max},M)}^{n} d^{n}_{m,(i)} = \Omega(m_m(n)),
		\end{aligned}
		\end{align*}
		where $(d^{n}_{m,(i)})_{i=1}^{n_1}$ is the descending-sorted version of $(d_{m,i}^{n})_{i=1}^{n_1}$ , $(d^{n}_{m,(i)})_{i=n_1+1}^{n}$ is the descending-sorted version of $(d_{m,i}^{n})_{i=n_1+1}^{n}$, $d_{1}^{\max} = d^{n}_{m,(1)}$, and $d_{2}^{\max} = d^{n}_{m,(n_1+1)}$.
		This implies that~\cite[Theorem 6.1]{Janson2014}$$\liminf_{n\to\infty} \prob(G^*(n_1,n_2,\mathbf{d}^n_m)\text{ is simple}) > 0.$$Recall that $x = \Omega(N(n))$ means $\liminf_{n\to\infty} x/N(n) > 0$.
	\end{enumerate}
\end{definition}
Janson in~\cite{Janson2009} proved that the probability of the event ``the random multigraph generated by the configuration model is simple'' is strictly positive if and only if the first assumption in the graph regularity condition holds. Various sufficient conditions for this property were given by the authors of~\cite{Bender1978,Bollobas1980,Mckay1991,Bollobas2001}. The final result in~\cite{Janson2009} provides a necessary and sufficient condition. Blanchet and Stauffer in~\cite{Blanchet2013} proved a similar result for bipartite random graphs: ``the random bipartite multigraph generated by the configuration model is simple'' with strictly positive probability if and only if the second graph regularity condition holds. Both results are presented in~\cite{Janson2014}.
\begin{remark}
	In Definition \ref{def:regcon_graph}, the part (i) of the second condition is similar to the first condition . Also, note that part (ii) of the second condition follows from the first condition if $d_1^{\max}=o(m_m(n))$ and $d_2^{\max}=o(m_m(n))$~\cite[Remark 6.1]{Janson2014}.
\end{remark}
\begin{remark}
	Almost all the results presented in this paper only need the degree regularity conditions. The additional assumptions in the graph regularity conditions extend these results to uniformly sampled simple graphs.
	Note that the second and third assumptions in Definition \ref{def:regcon_deg} are required to study the asymptotic behavior of the Markov process of adoption and are not invoked until Section \ref{sec:odeanalysisinf}.
\end{remark}
\begin{remark}\label{rem:newlambda}
	Suppose the degree regularity conditions hold. Since $\sum_{i=1}^{n_1}d^{n}_{m,i} = \sum_{i=n_1+1}^{n}d^{n}_{m,i}$, we have $\lambda_{1,m}(n) \beta_1(n) = \lambda_{2,m}(n) \beta_2(n)$. We denote this quantity by $\lambda_m(n)\coloneqq \lambda_{m,1}(n) \beta_1(n) = \lambda_{m,2}(n) \beta_2(n)$. Similarly, we define $\lambda_1(n) \coloneqq \lambda_{1,1}(n) \beta_1(n)$ and $\lambda_2(n) \coloneqq \lambda_{2,2}(n) \beta_2(n)$. Note that $\lambda_1(n) = 2m_1(n)/n$, $\lambda_m(n) = m_m(n)/n$, and $\lambda_2 = 2m_2(n)/n.$
	Following the same notation, we set $\lambda_j\coloneqq \lim_{n\to\infty} \lambda_j(n) = \lambda_{j,j}\beta_j$, and $\lambda_m\coloneqq \lim_{n\to\infty}\lambda_m(n) = \lambda_{1,m} \beta_1 = \lambda_{2,m}\beta_2$.
\end{remark}

	\section{Convergence to ODEs}\label{sec:odeapprox}
	In this section, we use techniques developed for the mean-field analysis~\cite{Wormald1995,Molloy1998,Wormald1999,Amini2010} of the resulting population density-dependent Markov processes to approximate the process by a system of ODEs. Approximating the vanilla form of the Markov process of adoption from Section~\ref{sec:markovproc} can be problematic as the dimension of the resulted ODEs may grow unboundedly as $n\to\infty$. Hence, we introduce two truncated versions of the Markov process of adoption, which sandwich the vanilla version. Focusing on the truncated Markov processes, we then approximate a scaled-version of this process by continuous functions obtained from the solution of a set of ODEs. We start by highlighting why the analysis is non-trivial and why truncation is necessary:

\begin{point}\label{point:1}
 The first point concerns some of the terms that appear in the one-step drift. Notice that we have many terms like
\begin{align*}
\frac{A_j(k)}{A_1(k) + A_2(k) + A^{(1)}_m(k) + A^{(2)}_m(k)}.
\end{align*}
In terms of the scaled variables, these terms are not Lipschitz unless there is a lower bound on the value of the (scaled) denominator. Owing to this, in our ODEs approximation, we will have to stop the Markov process of adoption just before the sum of these scaled variables hits zero (corresponding to the denominator above), i.e., before all the active half-edges have been omitted; it is important that this be the sum and not the individual components. For the same reason, we have to stop the process before we run out of half-edges in any of the two communities or between the communities.
\end{point}
\begin{point}\label{point:2}
The second point is regarding the one-step drift of variables like $A_j(k)$, i.e., the number of active half-edges based on the community structure. The one-step drift can be unbounded as the increase can equal the number of vertices (in the appropriate community) minus one. However, owing to the balance equations, as we pointed out in Remark \ref{rem:trackvar}, there is no need to keep track of the random variables associated with the number of active half-edges.

There is, however, another technical issue with the one-step drift of other quantities as they depend on all terms $I^{(j)}_{d_j,d_{-j},u_j,u_{-j}}(k)$ through a sum associated with $A_1(k)+A_2(k)+A_m^{(1)}(k)+A_m^{(2)}(k)$.
For any finite $n$, we only need to account for a finite number of terms, but in the limit, we have a countable number of terms leading to a similar property for the functions associated with these variables. Since the coefficients of these variables are increasing without bound, the associated functions are not Lipschitz continuous. This precludes the direct application of the results of~\cite{Wormald1995,Wormald1999}. On the other hand, given degree regularity conditions, this should be a superficial problem as the total number of half-edges associated with vertices with high degrees is small.

To address the specific scenario outlined above, we bound the original Markov process of adoption, from above and below, using two truncated versions of the process. We denote these Markov chains with $X_{U,\delta}^n$ and $X_{L,\delta}^n$ respectively, where $\delta> 0$ is the tuning parameter. Fix $\delta > 0$ and set $d_\delta>0$ large enough so that the following inequalities hold for all $n$:
\begin{align*}
&\sum_{i = 1}^{n_1} (d_{1,i}^n + d_{m,i}^n) \bs{1}\{d_{1,i}^n + d_{m,i}^n > d_\delta\}  \leq \delta n_1,  \\
&\sum_{i = n_1+1}^{n} (d_{2,i}^n + d_{m,i}^n) \bs{1}\{d_{2,i}^n + d_{m,i}^n > d_\delta\}  \leq \delta n_2.
\end{align*}
$X_{U,\delta}^n$ is defined by activating all vertices with total degree larger than $d_{\delta}$, i.e., by setting $\alpha_j(d_j,d_{-j}) = 1$ for all $d_j+d_{-j} > d_{\delta}$ and $j\in\{1,2\}$. $X_{L,\delta}^n$ is defined by assuming $K_j(d_j,d_{-j}) = d_j + d_{-j}$ for all $j\in\{1,2\}$ and all $d_j + d_{-j} > d_{\delta}$ so that these nodes can never be activated during the contagion. Note that by degree regularity conditions (Definition \ref{def:regcon_deg}) such a $d_\delta > 0$ exists.

Since none of the inactive vertices with degree higher than $d_{\delta}$ can be activated, instead of tracking the random variables $I^{(j)}_{d_j,d_{-j},u_j,u_{-j}}$ for these vertices, we track the total number of half-edges associated with these random variables.

Intuitively speaking, the number of inactive vertices of any degree at the natural stopping time of $X^n$ is bounded between the same quantities for $X_{U,\delta}^n$ and $X_{L,\delta}^n$. In particular, there is a natural coupling between $X_{U,\delta}^n$, $X^n$ and $X_{L,\delta}^n$ such that for all $j\in\{1,2\}$:
\begin{align} \label{eq:naturalstop}
&I^{(j)}_{d_j,d_{-j},u_j,u_{-j}}(\stoptime^n_{U,\delta}) \leq I^{(j)}_{d_j,d_{-j},u_j,u_{-j}}(\stoptime^n) \leq I^{(j)}_{d_j,d_{-j},u_j,u_{-j}}(\stoptime^n_{L,\delta})\qquad \forall d_j + d_{-j} \leq d_\delta,
\end{align}
where $\stoptime^n_{U,\delta}$, $\stoptime^n$, and $\stoptime^n_{L,\delta}$ are the natural stopping times of $X_{U,\delta}^n$, $X^n$, and $X_{L,\delta}^n$ respectively (we are abusing notation here, i.e., $I^{(j)}_{d_j,d_{-j},u_j,u_{-j}}(\stoptime^n_{U,\delta})$ is the number of inactive vertices with certain parameters at the natural stopping time of $X_{U,\delta}^n$ etc.). Next, we present this natural coupling.

Let us consider a realization of the Markov process of adoption $X_{U,\delta}^n$ up to its natural stopping time. We couple this realization with a realization of $X^n$ as follows: at each time of the random process $X^n$, we pick an active half-edge uniformly at random and pair it with the same half-edge as in $X_{U,\delta}^n$ (remove both half-edges) so that the same vertices appear as neighbors in $X^n$. Similarly, we realize $X_{L,\delta}^n$. Note that the order in which we pick an active half-edge to be paired with a random half-edge does not affect the state of the Markov processes of adoption at its natural stopping time.

It is important to note that this bound only works for the stopping time of the processes and not the whole trajectory. In the following subsection, we use the techniques developed by Wormald~\cite{Wormald1995,Wormald1999} to approximate the trajectories of $X_{L,\delta}^n$ and $X_{U,\delta}^n$. In Section \ref{sec:odeanalysisinf}, we first approximate the state of these random processes at their natural stopping time. Then we show that the difference between these two approximations can be made arbitrarily small by tuning the parameter $\delta > 0$.
\end{point}

\subsection{Convergence to ODEs for the Truncated Processes}
In this subsection, we focus on the case where inactive vertices with degree higher than some constant $d_{\max}$ cannot be activated. Let us define a new set of random variables to keep track of half-edges associated with these inactive vertices:\\
$W_{j}(k)$: Number of remaining half-edges belonging to inactive vertices with degree higher than $d_{\max}$ in community $j$ at time k.\\
$W^{(j)}_{m}(k)$: Number of remaining half-edges between the two communities belonging to inactive vertices with degree higher than $d_{\max}$ in community $j$ at time $k$.\\
At each step of the Markov process of adoption, the value of $W_{j}(k)$ or $W^{(j)}_{m}(k)$ can reduce at most by one.
\begin{remark}\label{rem:newrand}
	The one-step drifts of these random variables are given as follows:
	\begin{align*}
	&\expect[W_j(k+1) - W_j(k)|X^n(k)] =  \\
	&\myquad[4]\frac{A_j(k)}{A_1(k) + A_2(k) + A^{(1)}_m(k) + A^{(2)}_m(k)} \times \frac{ -W_j(k)}{2m_j(n) - 2T_j(n) - 1},\\
	&\expect[W^{(j)}_{m}(k+1) - W^{(j)}_{m}(k)|X^n(k)] =  \\
	&\myquad[4]\frac{A^{(-j)}_m(k)}{A_1(k) + A_2(k) + A^{(1)}_m(k) + A^{(2)}_m(k)} \times \frac{ -W^{(j)}_{m}(k)}{m_m(n) - (k-T_1(k)-T_2(k))}.
	\end{align*}
	Note that these random variables do not change the dynamic of the Markov process of adoption; hence, the one-step drifts of all other random variables remain the same. Also, note that the sum in the balance equations \eqref{eq:balanceeq} now has a finite number of summands (less than $(d_{\max}+1)^4$ many summands) as the other terms are replaced with either $W_{j}(k)$ or $W^{(j)}_{m}(k)$:
	\begin{align}
	\begin{split}
	&~A_j(k) + \sum_{\substack{u_j+u_{-j} \leq K_j(d_j,d_{-j})\\d_j + d_{-j} \leq d_{\max}}}  (d_j-u_j)I^{(j)}_{d_j,d_{-j},u_j,u_{-j}}(k) + W_j(k)= 2m_j(n)- 2T_j(k), \\
	&~A^{(j)}_m(k) + \sum_{\substack{u_j+u_{-j} \leq K_j(d_j,d_{-j})\\d_j + d_{-j} \leq d_{\max}}} (d_{-j}-u_{-j})I^{(j)}_{d_j,d_{-j},u_j,u_{-j}}(k) + W^{(j)}_{m}(k) = m_m(n) - T_m(k),
	\end{split}
	\label{eq:balanceeq_W}
	\end{align}
	where $T_m(k) \coloneqq k - T_j(k) - T_{-j}(k)$.
\end{remark}
Recall that $m_j(n)$ denotes the total number of edges on side $j\in\{1,2\}$, and $m_m(n)$ denotes the total number of edges between the two communities. Recall also that $\lambda_1(n) = 2m_1(n)/n$, $\lambda_m = m_m(n)/n$, and $\lambda_2(n)=2m_2(n)/n$. Then the ODEs follow by defining the real functions $\tau_j(t)$, $i^{(j)}_{d_j,d_{-j},u_j,u_{-j}}(t)$, $w_{j}(t)$, and $w^{(j)}_{m}(t)$ to model the behavior of their discrete counterpart, i.e., intuitively speaking:
\begin{equation}\label{eq:scalevar}
\begin{aligned}
&\tau_j(t)=\lim_{n\rightarrow\infty} \frac{1}{n}T_j(tn), \;
&&i^{(j)}_{d_j,d_{-j},u_j,u_{-j}}(t)=\lim_{n\rightarrow\infty} \frac{1}{n}I^{(j)}_{d_j,d_{-j},u_j,u_{-j}}(tn),\\
&w_{j}(t)=\lim_{n\rightarrow\infty} \frac{1}{n}W_{j}(tn), \;
&&w^{(j)}_{m}(t)=\lim_{n\rightarrow\infty} \frac{1}{n}W^{(j)}_{m}(tn), \;
\end{aligned}
\end{equation}
where all the limits are in probability and sample-path-wise.
We can then use the one-step drifts from Appendix~\ref{app:onestepdrift} and Remark \ref{rem:newrand} to derive the ODEs. The details are in Appendix~\ref{app:odederive}.

As we mentioned in Point \ref{point:1}, we have to stop the process just before we run out of half-edges within each community or between the two communities, and before we run out of active half-edges. Invoking \cite[Theorem 5.1]{Wormald1999}, we can approximate the trajectory of the truncated process using the solution of the system of ODEs \eqref{eq:diffeq_i}-\eqref{eq:ode_ic}, as long as the functions associated with the ODEs are Lipschitz continuous.

\begin{theorem}\label{thm:odesfinite}
Fix $\varepsilon> 0$ small enough, and consider the Markov process of adoption. Assume there is a constant $d_{\max} > 0$, independent of $n$, such that $K_j(d_j,d_{-j}) = d_j+d_{-j}$ for all $d_j+d_{-j} > d_{\max}$ and $j\in\{1,2\}$. Consider a realization of the initial condition \eqref{eq:ode_ic} given in Appendix \ref{app:onestepdrift}, for which the total number of active half-edges at the beginning of the process is greater than $2\varepsilon n$. Let $\theta = O(n^{-\gamma})$ for some $\gamma < 1/3$.
Then, with probability $1 -  O\left(\theta^{-1}\exp(-n\theta^3)  \right)$, we have
\begin{align*}
&\left|I^{(j)}_{d_j,d_{-j},u_j,u_{-j}}(t) - n i^{(j)}_{d_j,d_{-j},u_j,u_{-j}}(t/n)\right| = O(n\theta),\\
&\left| T_j(t) - n\tau_j(t/n)\right| = O(n\theta),\\
&\left| W_j(t) - nw_j(t/n)\right| = O(n\theta),\\
&\left| W^{(j)}_m(t) - nw^{(j)}_m(t/n)\right| = O(n\theta),
\end{align*}
uniformly for $0 \leq t\leq \sigma_\varepsilon n$, where $w_j$, $w^{(j)}_m$, $i^{(j)}_{d_j,d_{-j},u_j,u_{-j}}$ and $\tau_j$ are the solution of the ODEs given in Appendix \ref{app:odederive} with the realized initial condition, and $\sigma_\varepsilon = \sigma_\varepsilon(n)$ is the supremum of those $x$ to which the solution of the ODEs \eqref{eq:diffeq_i}-\eqref{eq:ode_ic} can be extended before reaching within $l^{\infty}$-distance $C\theta$ of the boundary of $\widehat{\mathcal{D}}_{\varepsilon,n}$, for a sufficiently large constant $C$. The open connected set $\widehat{\mathcal{D}}_{\varepsilon,n}$ is defined as follows:
\begin{align*}
\widehat{\mathcal{D}}_{\varepsilon,n} \coloneqq &\bigg\{ (t,\tau_1,\tau_2,w_1,w_2,w^{(1)}_m,w^{(2)}_m,i^{(1)}_{d_1,d_2,u_1,u_2},i^{(2)}_{d_2,d_1,u_2,u_1})\in R^{K}:\\
&\myquad[2]  -\varepsilon < t - \tau_1 - \tau_2 <\lambda_m(n) - \varepsilon,\\
&\myquad[2]\text{for $j\in\{1,2\}$}:
\begin{aligned}
	&-\frac{\varepsilon}{2} < \tau_j <\frac{ \lambda_j(n) - \varepsilon }{2},\\
	& -\varepsilon < w_j < 4\lambda_j(n),~ -\varepsilon < w^{(j)}_m < 2\lambda_m(n),
\end{aligned}\\
&\myquad[2] \text{for $ j\in\{1,2\}$, $u_j\leq d_j$, $u_{-j}\leq d_{-j}$, $d_j + d_{-j} \leq d_{\max}$}:-\varepsilon < i^{(j)}_{d_j,d_{-j},u_j,u_{-j}} < 2,\\
&\myquad[2]  -\varepsilon< a_1,a_2, a^{(1)}_m,a^{(2)}_m , \text{ and } \\
&\myquad[2]\varepsilon < a_1 + a_2 + a^{(1)}_m + a^{(2)}_m < 4(\lambda_1(n)+\lambda_2(n) + \lambda_m(n))\bigg\}
\end{align*}
where $K \leq 6 + (d_{\max}+1)^4$ is a constant, and
\begin{align*}
&a_j \coloneqq -\sum_{\substack{u_j+u_{-j} \leq K_j(d_j,d_{-j})\\d_j + d_{-j} \leq d_{\max}}} (d_j-u_j)i^{(j)}_{d_j,d_{-j},u_j,u_{-j}} + \lambda_j(n) - 2\tau_j - w_j,\\
&a^{(j)}_m \coloneqq -\sum_{\substack{u_j+u_{-j} \leq K_j(d_j,d_{-j})\\d_j + d_{-j} \leq d_{\max}}} (d_{-j}-u_{-j})i^{(j)}_{d_j,d_{-j},u_j,u_{-j}} + \lambda_m(n) -\tau_m- w^{(j)}_m,
\end{align*}
and $\tau_m \coloneqq t - \tau_1 - \tau_2$.
\end{theorem}
\begin{proof}
See Appendix \ref{proof:odesfinite}.
\end{proof}

Using  the balance equations \eqref{eq:balanceeq_W}, we get the following corollary.
\begin{corollary}\label{cor:odesfinite}
	In the setting of Theorem \ref{thm:odesfinite}, with probability $1 -  O\left(\theta^{-1}\exp(-n\theta^3)  \right)$, we also have
	\begin{align*}
	&\left|A_j(t) - n a_j(t/n)\right| = O(n\theta),\\
	&\left|A^{(j)}_m(t) - n a^{(j)}_m(t/n)\right| = O(n\theta),
	\end{align*}
	uniformly for $0 \leq t\leq \sigma_\varepsilon n$, where $a_j$ and $a^{(j)}_m$ are given by equations \eqref{eq:diffeq_aj} and \eqref{eq:diffeq_amj} respectively, in Appendix \ref{app:odederive}.
\end{corollary}
We comment that the initial condition \eqref{eq:ode_ic} given in Appendix \ref{app:onestepdrift} is random and that the only source of randomness is the state of the vertices at time $0$ (active or inactive). The above statements hold for any realization of the initial condition in which the total number of active half-edges is greater than $2\varepsilon n$ (so that the state of the Markov processes of adoption at time $0$ is an interior point of $\widehat{\mathcal{D}}_{\varepsilon,n}$).

	\section{A Probabilistic Argument to Solve the ODEs}\label{sec:odeprobsol}
	In this section, we present a probabilistic heuristic to derive the form of the solution of the ODEs given in Appendix \ref{app:odederive}. Our sketch also provides an intuitive answer to the ``surprising simplification'' that has been observed in the solution of the ODEs in~\cite{Amini2010,Balogh2007}. The formal proof is in the following section.

Each iteration of the Markov process of adoption given by Algorithm \ref{algorithm} has two important phases: first, we pick an active half-edge, and then, we pair it with a random half-edge in the proper community. Consider a fixed half-edge $e$ in community $j$ that belongs to an inactive vertex $v$ at the beginning of the Markov process of adoption. We want to estimate the probability that $e$ has not been paired with any other half-edge up to time $k$.

Condition on the event that $v$ is still inactive at time $k$. Then $e$ has not been paired with any other half-edge if it has not been picked in the second phase of any iteration of Algorithm \ref{algorithm} up to time $k$. Any additional dependencies introduced by the conditioning should fade away as $n\to\infty$, and we will proceed by ignoring them. Since half-edges are chosen uniformly at random in the second phase of each iteration, we have
\begin{equation*}
\begin{aligned}
&\prob(\{e\text{ has not been paired with any other half-edge}\}|\{v \text{ is inactive}\})\approx\\
&\myquad[2]\left(1 - \frac{1}{2m_j(n) - 1}\right)\times \left(1 - \frac{1}{2m_j(n) - 2 - 1}\right) \times \cdots\times \left(1 - \frac{1}{2m_j(n) - 2T_j(k) -1}\right).
\end{aligned}
\end{equation*}
Using the simple approximation $1 - x \approx \mathrm{e}^{-x}$ for small values of $x$, we get
\begin{align*}
&\prob(\{e\text{ has not been paired with any other half-edge}\}|\{v \text{ is inactive}\})\\
&\myquad[1]\approx\exp\left( -\sum_{i=1}^{T_j(k)} \frac{1}{2m_j(n) - 2i-1}\right)\approx\exp\left( -\frac{1}{2}\int_{m_j(n)-T_j(k)}^{m_j(n)} \frac{1}{z}\,dz\right) = \left(1 - \frac{T_j(k)}{m_j(n)}\right)^\frac{1}{2}.
\end{align*}
Next, we use the same argument for the half-edges that are supposed to connect to the vertices in the other community. However, there is an important distinction here as two phases happen in different communities: if in the first phase we pick an active half-edge from community $j$, in the second phase we pick a random half-edge from community $-j$. This makes the direct use of the above argument almost impossible. To fix it, we track two half-edges dangling from two vertices in communities $1$ and $2$. Consider two fixed half-edges $e_1$ and $e_2$ between the communities such that that $e_j$ belongs to an inactive vertex $v_j$ in community $j$ for $j\in\{1\gets 2\}$. Now, we can use the same argument as before:
\begin{equation*}
\begin{aligned}
&\prob(\{e_1 \text{ and }e_2\text{ have not been paired with any other half-edges}\}|\{v_1 \text{ and }v_2 \text{ are inactive}\})\\
&\myquad[2]\approx\left(1 - \frac{1}{m_m(n) }\right)\times \left(1 - \frac{1}{m_m(n) - 1 - 1}\right) \times \cdots\times \left(1 - \frac{1}{m_m(n) - T_m(k) -1}\right)\\
&\myquad[2]\approx \left(1 - \frac{T_m(k)}{m_m(n)}\right).
\end{aligned}
\end{equation*}
where $T_m(k)$ denotes the number of times the algorithm removes one half-edge from each community. Note that $T_m(k) = k - T_1(k)-T_2(k)$. Intuitively speaking, for large values of $n$, the events $\{e_j \text{ is not paired with }\allowbreak \text{any other half-edges}\}$ for $j\in\{1\gets 2\}$ are independent, and we can write:
\begin{equation*}
\begin{aligned}
&\prob(\{e_1 \text{ and }e_2\text{ have not been paired with any other half-edges}\}|\{v_1 \text{ and }v_2 \text{ are inactive}\})\\
&\myquad[4]\approx\prob(\{e_1\text{ has not been paired with any other half-edge}\}|\{v_1 \text{ is inactive}\})\\
&\myquad[6]\times\prob(\{e_2\text{ has not been paired with any other half-edge}\}|\{v_2 \text{ is inactive}\})
\end{aligned}
\end{equation*}
If we denote the two quantities on the RHS as $Z_1(k)$ and $Z_2(k)$ respectively, we expect to have
\begin{align*}
Z_1(k)Z_2(k) \approx \left(1 - \frac{T_m(k)}{m_m(n)}\right).
\end{align*}
Note that $1-Z_1(k)$ is the probability that the half-edge $e_1$ has been paired with an active half-edge in community $2$. Considering $v_1$ as the root vertex, $1-Z_1(k)$ is the probability that the descendant of $v_1$ in the second community through the half-edge $e_1$ is active, even if this link is excised. In particular, following the discussion of Section \ref{sec:meanfield}, it is natural to expect that $1-Z_1(tn) \approx 1-\mu^{(1\gets 2)}(t)$.

Next, consider a vertex $v$ in community $j$ with $d_j$ half-edges in community $j$ and $d_{-j}$ half-edges between the two communities. The above sketch, together with some independence assumptions (which can be justified as $n\to\infty$), suggests the following approximate equality at time $k$:
\begin{align*}
&\prob\Big(\Big\{
\begin{minipage}{0.6\textwidth}
$u_j$ out of $d_j$ half-edges and $u_{-j}$ out of $d_{-j}$ half-edges of $v$ have been removed where $u_j+u_{-j} \leq K_j(d_j,d_{-j})$
\end{minipage}
\Big\}\Big) \approx \\
&\myquad[6]Bi\left(u_j,d_j;1-\left(1 - \frac{T_j(k)}{m_j(n)}\right)^{\frac{1}{2}}\right) \times
Bi\left(u_{-j},d_{-j};1-Z_j(k)\right),
\end{align*}
which further implies
\begin{equation}\label{eq:probsol_i}
\begin{aligned}
	&\expect\big[I^{(j)}_{d_j,d_{-j},u_j,u_{-j}}(k) \big| I^{(j)}_{d_j,d_{-j},0,0}(0) \big]\approx \\
	&\myquad[2]I^{(j)}_{d_j,d_{-j},0,0}(0) \times
	Bi\left(u_j,d_j;1-\left(1 - \frac{T_j(k)}{m_j(n)}\right)^{\frac{1}{2}}\right) \times
	Bi\left(u_{-j},d_{-j};1-Z_j(k)\right).
\end{aligned}
\end{equation}

Now, we are ready to simplify the differential equations given in Appendix \ref{app:odederive}. Let $\mu^{(j\gets j)}$ and $\mu^{(j\gets -j)}$ model the limiting value of the probabilities defined above; intuitively speaking, we have
\begin{align}\label{eq:scalevar_mu}
\mu^{(j\gets j)}(t)=\lim_{n\to\infty}\left(1 - \frac{T_j(tn)}{m_j(n)}\right)^\frac{1}{2},\text{ and }
\mu^{(j\gets -j)}(t)=\lim_{n\to\infty}Z_j(tn).
\end{align}
Then the heuristic equality \eqref{eq:probsol_i} suggests that
\begin{align}
&i^{(j)}_{d_j,d_{-j},u_j,u_{-j}}(t) = i^{(j)}_{d_j,d_{-j},0,0}(0)\, Bi(u_j;d_j,1-\mu^{(j\gets j)}(t))\,Bi(u_{-j};d_{-j},1-\mu^{(j\gets -j)}(t)). \label{eq:sol_i}
\end{align}
Also, by the definition of the scaled variables
\begin{align}
&\tau_j(t) =\frac{\lambda_j(n)}{2}\left( 1 - {\mu^{(j\gets j)}(t)}^2  \right),\label{eq:sol_tau} \\
&\tau_m(t) =\lambda_m(n)\left( 1 - \mu^{(1\gets 2)}(t)\mu^{(2\gets 1)}(t)  \right),\label{eq:sol_taum}
\end{align}
where $\tau_m$ is the continuous counterpart of $T_m$. Moreover, we expect the following equality to hold:
\begin{align}
&\frac{\lambda_1(n)}{2} \left(\mu^{(1\gets 1)}(t)\right)^2 + \frac{\lambda_2(n)}{2} \left(\mu^{(2\gets 2)}(t)\right)^2 +\lambda_m(n)\mu^{(1\gets 2)}(t)\mu^{(2\gets 1)}(t)   \\
&\myquad[20]= \lambda_m(n) + \frac{\lambda_1(n)}{2} + \frac{\lambda_2(n)}{2} - t, \label{eq:sol_t}
\end{align}
as it is equivalent to the equality $T_1(k) + T_2(k)+T_m(k) = k$. The importance of the above equality is that given the vector $\left(\mu^{(1\gets 1)}(t),\mu^{(1\gets 2)}(t),\mu^{(2\gets 1)}(t),\mu^{(2\gets 2)}(t)\right)$, the value of $t$ is uniquely determined. Furthermore, all the other relevant quantities, like $i^{(j)}_{d_j, d_{-j}, u_j, u_{-j}}(t)$, are also determined.

\begin{remark}
	We emphasize that the above argument is not meant to be rigorous, but the result is surprisingly correct. For example, there might be no inactive vertex in community $j\in\{1,2\}$ to begin with; however, equations \eqref{eq:sol_i}-\eqref{eq:sol_t} are still valid.
\end{remark}

	\section{Analysis of the ODEs}\label{sec:odeanalysisfinite}
	Following the heuristic argument provided in Section \ref{sec:odeprobsol}, we start by presenting the solution of the ODEs in Appendix \ref{app:odederive}.
The following lemma characterizes the solution of the differential equations that (with high probability) approximate the adoption process.
\begin{lemma} \label{lem:diffeq_sol}
The solution of differential equations \eqref{eq:diffeq_i}--\eqref{eq:diffeq_wmj} with initial condition \eqref{eq:ode_ic} in $\widehat{\mathcal{D}}_{\varepsilon,n}$, is given by equations \eqref{eq:sol_i}, \eqref{eq:sol_tau}, \eqref{eq:sol_taum}, and the following:
\begin{align}\label{eq:sol_w}
w_j(t) = w_j(0) \mu^{(j\gets j)}(t), \myquad[4]w^{(j)}_m(t) = w^{(j)}_m(0) \mu^{(j\gets - j)}(t),
\end{align}
for $j\in\{1,2\}$, where $\left(\mu^{(1\gets 1)}(0),\mu^{(1\gets 2)}(0),\mu^{(2\gets 1)}(0),\mu^{(2\gets 2)}(0)\right)$ is the unique solution of the following four-dimensional differential equations
\begin{align}
&\frac{-a_j(t)}{a_1(t) + a_2(t) + a^{(1)}_m(t) + a^{(2)}_m(t)}
= \lambda_j(n) \frac{d\mu^{(j\gets j)}}{dt}\left(\mu^{(j\gets j)}(t)\right) \label{eq:sol_mujj},\\
&\frac{-a^{(-j)}_m(t)}{a_1(t) + a_2(t) + a^{(1)}_m(t) + a^{(2)}_m(t)}
=\lambda_m(n) \frac{d\mu^{(j\gets - j)}}{dt}\left(\mu^{(-j\gets j)}(t)\right)\label{eq:sol_muj-j},
\end{align}
with the initial condition given by
\begin{align}
&\left(\mu^{(1\gets 1)}(0),\mu^{(1\gets 2)}(0),\mu^{(2\gets 1)}(0),\mu^{(2\gets 2)}(0)\right) = (1,1,1,1),\label{eq:sol_mu0}
\end{align}
and $\left(\mu^{(1\gets 1)}(t),\mu^{(1\gets 2)}(t),\mu^{(2\gets 1)}(t),\mu^{(2\gets 2)}(t)\right) \in \mathcal{D}_{\varepsilon,n}$. The set $\mathcal{D}_{\varepsilon,n}$ is defined as follows:
\begin{align*}
\mathcal{D}_{\varepsilon,n} \coloneqq &\bigg\{\bs{\mu} =
(\mu^{(1\gets 1)},\mu^{(1\gets 2)},\mu^{(2\gets 1)},\mu^{(2\gets 2)}) \in [0,1]^4:\\
&\myquad[4]\text{for all $j\in\{1,2\}$:}~ \sqrt{\frac{\varepsilon}{\lambda_j(n)}}< \mu^{(j\gets j)},~ \sqrt{\frac{\varepsilon}{\lambda_m(n)}} < \mu^{(j\gets - j)} \\
&\myquad[4]  -\varepsilon< a_1(\bs{\mu}),a_2(\bs{\mu}), a^{(1)}_m(\bs{\mu}),a^{(2)}_m(\bs{\mu}), \\
&\myquad[4] \varepsilon < a_1(\bs{\mu}) + a_2(\bs{\mu}) + a^{(1)}_m(\bs{\mu}) + a^{(2)}_m(\bs{\mu}) \bigg\}.
\end{align*}
The functions $a_j(t)$ and $a^{(j)}_m(t)$ are given as follows:
\begin{align*}
&a_j(t) = -\sum_{\substack{u_j+u_{-j} \leq K_j(d_j,d_{-j})\\d_j+d_{-j} \leq d_{\max}}} (d_j-u_j)i^{(j)}_{d_j,d_{-j},u_j,u_{-j}}(t) + \lambda_j(n) - 2\tau_j(t) - w_j(t),\\[5pt]
&a^{(j)}_m(t) = -\sum_{\substack{u_j+u_{-j} \leq K_j(d_j,d_{-j})\\d_j+d_{-j} \leq d_{\max}}} (d_{-j}-u_{-j})i^{(j)}_{d_j,d_{-j},u_j,u_{-j}}(t) + \lambda_m(n) -\tau_m(t) - w^{(j)}_m(t),
\end{align*}
and $\tau_m(t)$ is given by \eqref{eq:sol_taum}. Abusing notation, the functions $a_j(\bs{\mu})$ and $a^{(j)}_m(\bs{\mu})$ are defined similar to $a_j(t)$ and $a^{(j)}_m(t)$ using \eqref{eq:sol_i}, \eqref{eq:sol_tau}, \eqref{eq:sol_taum} and \eqref{eq:sol_w}. Also, the solution of the four dimensional differential equation satisfies the equality \eqref{eq:sol_t}.
\end{lemma}
\begin{proof}
See Appendix \ref{proof:diffeq_sol}.
\end{proof}
\begin{remark}
	The set $\widehat{\mathcal{D}}_{\varepsilon,n}$ is a fixed set, however, $\mathcal{D}_{\varepsilon,n}$ is a random set since the terms in \eqref{eq:sol_i} and \eqref{eq:sol_w} depend on  the initial condition \eqref{eq:ode_ic}. Note that both of them depend on $n$.
\end{remark}

The significance of this result is in demonstrating that the set of ODEs from Section~\ref{sec:odeapprox} can be reduced to a set of four-dimensional ODEs (which has a unique solution). Note that this dimension reduction applies to the sample path of the adoption process and not just the final population of active vertices as suggested by the mean-field approximation of Section~\ref{sec:meanfield}.

The denominator of all four equations given by \eqref{eq:sol_mujj} and \eqref{eq:sol_muj-j} are the same. Since in $\mathcal{D}_{\varepsilon,n}$ this quantity is bounded away from zero by $\varepsilon>0$, it is safe to remove this term from the denominator of the differential equations for equilibrium analysis. More specifically, if we consider a particle at $(1,1,1,1)$ whose movement is governed by \eqref{eq:sol_mujj} and \eqref{eq:sol_muj-j}, removing the denominator will not change the trajectory of the particle but will affect its speed. Hence, after some simple algebra, we find that the trajectory of \eqref{eq:sol_mujj}-\eqref{eq:sol_muj-j} is the same as the trajectory of the following system of differential equations:
\begin{align}
\frac{d\mu^{(j\gets j)}}{dt} &= \Ffunc_{(j\gets j)}(\mu^{(j\gets j)},\mu^{(j\gets - j)}) - \mu^{(j\gets j)}, \label{eq:altsol_mujj}\\
\frac{d\mu^{(j\gets - j)}}{dt} &= \Ffunc_{(j\gets -j)}(\mu^{(-j\gets - j)},\mu^{(-j\gets j)}) - \mu^{(j\gets - j)}, \label{eq:altsol_muj-j}
\end{align}
for $j\in\{1,2\}$ with the same initial condition \eqref{eq:sol_mu0}, where the functions $\Ffunc_{(j\gets j)}$ and $\Ffunc_{(j\gets -j)}$ are given as follows:
\begin{align}
&\begin{aligned}\label{eq:trunc_Fjj}
&\Ffunc_{(j\gets j)}(\mu^{(j\gets j)},\mu^{(j\gets - j)}) \coloneqq \sum_{\substack{u_j+u_{-j} \leq K_j(d_j,d_{-j})\\d_j+d_{-j} \leq d_{\max}}}  \,\frac{d_j}{\lambda_j(n)} \,i^{(j)}_{d_j,d_{-j},0,0}(0)\,\\
&\myquad[8] Bi(u_j;d_j-1,1-\mu^{(j\gets j)})Bi(u_{-j};d_{-j},1-\mu^{(j\gets - j)})  +  \frac{w_j(0)}{\lambda_j(n)},
\end{aligned}
\hspace{10mm}\\
&\begin{aligned}\label{eq:trunc_Fj-j}
&\Ffunc_{(j\gets -j)}(\mu^{(-j\gets - j)},\mu^{(-j\gets j)}) \coloneqq \sum_{\substack{u_j+u_{-j} \leq K_{-j}(d_{-j},d_{j})\\d_j+d_{-j} \leq d_{\max}}} \,\frac{d_j}{\lambda_m(n)}\,i^{(-j)}_{d_{-j},d_{j},0,0}(0)\,\\
&\myquad[8] Bi(u_j;d_j-1,1-\mu^{(-j\gets j)})Bi(u_{-j};d_{-j},1-\mu^{(-j\gets - j)}) + \frac{w^{(-j)}_m(0)}{\lambda_m(n)}.
\end{aligned}
\end{align}
\begin{remark}\label{rem:atoFcon}
	In derivation of \eqref{eq:altsol_mujj} and \eqref{eq:altsol_muj-j}, we have used the following equalities:
	\begin{align}
	&a_j = \lambda_j(n)\mu^{(j\gets j)}(\mu^{(j\gets j)} - \Ffunc_{(j\gets j)}(\mu^{(j\gets j)},\mu^{(j\gets - j)}) )\label{eq:ajtoFcon}\\
	&a^{(-j)}_m = \lambda_m(n)\mu^{(-j\gets j)}(\mu^{(j\gets - j)} - \Ffunc_{(j\gets -j)}(\mu^{(-j\gets - j)},\mu^{(-j\gets j)}) )\label{eq:a-jtoFcon}
	\end{align}
	These equalities are algebraic and straightforward.
\end{remark}
\begin{remark}\label{rem:FandFconection}
	There is a clear connection between the function $\bs{F}$ given by \eqref{eq:meanfield_mu} and the function $\Ffuncbold \coloneqq (\Ffunc_{(1\gets 1)},\Ffunc_{(1\gets 2)},\Ffunc_{(2\gets 1)},\Ffunc_{(2\gets 2)})$. Note that as $n\to\infty$:
	\begin{align*}
	&i^{(j)}_{d_j,d_{-j},0,0}(0)\,d_j/\lambda_j(n) \xrightarrow{P} \prob_{j*,m}(d_j,d_{-j})(1-\alpha_j(d_j,d_{-j})), \allowdisplaybreaks\\ &i^{(-j)}_{d_{-j},d_{j},0,0}(0)\,d_j/\lambda_m(n)\xrightarrow{P}\prob_{-j,m*}(d_{-j},d_{j})(1-\alpha_{-j}(d_{-j},d_{j})),\allowdisplaybreaks\\
	&w_j(0)/\lambda_j(n) \xrightarrow{P} \sum_{d_j+d_{-j} > d_{\max}}\prob_{j*,m}(d_j,d_{-j})(1-\alpha_j(d_j,d_{-j})),\allowdisplaybreaks\\
	&w^{(-j)}_m(0)/\lambda_m(n)\xrightarrow{P}\sum_{d_j+d_{-j} > d_{\max}}\prob_{-j,m*}(d_{-j},d_{j})(1-\alpha_{-j}(d_{-j},d_{j})).
	\end{align*}
	The same terms appear in \eqref{eq:meanfield_mujj} and \eqref{eq:meanfield_muj-j} if we assume $K_j(d_j,d_{-j}) = d_j+d_{-j}$ for all $d_j+d_{-j}>d_{\max}$ and $j\in\{1,2\}$ (for the proof, see Lemma \ref{lem:initcond_conv}). Almost all properties of $\bs{F}$ and $\Ffuncbold$ are the same, and all the statements in the rest of this section hold for both. Note that the function $\Ffuncbold$, given a realization of the initial condition \eqref{eq:ode_ic}, is not a random function.
\end{remark}
\begin{remark}
	For the sake of notational simplicity, we may write $\Ffunc_{(j\gets j)}(\bs{\mu})$ instead of $\Ffunc_{(j\gets j)}(\mu^{(j\gets j)},\allowbreak\mu^{(j\gets - j)})$, and $\Ffunc_{(j\gets -j)}(\bs{\mu})$ instead of $\Ffunc_{(j\gets -j)}(\mu^{(-j\gets - j)},\mu^{(-j\gets j)})$.
\end{remark}
Note that the function $\Ffuncbold$ is well-defined on $[0,1]^4$, hence, we do not need any restriction on its domain. We continue with some basic properties of the function $\Ffuncbold$. As the first step, we prove that $\Ffuncbold$ is increasing in each of its components.
\begin{lemma} \label{lem:Fprop_increasing}
If $\bs{\mu} \geq \bs{\mu}^\prime$ component-wise with $\bs{\mu}\neq \bs{\mu}^{\prime}$, then $\Ffuncbold(\bs{\mu}) \geq \Ffuncbold(\bs{\mu}^\prime)$ component-wise, and $\Ffuncbold(\bs{\mu})\neq \Ffuncbold(\bs{\mu}^\prime)$.
\end{lemma}
\begin{proof}
See Appendix \ref{proof:Fprop_increasing}.
\end{proof}

To analyze the equilibrium of \eqref{eq:altsol_mujj}-\eqref{eq:altsol_muj-j}, we use the LaSalle Invariance Principle~\cite{Lasalle1960}. The following lemma characterizes the most important properties of $\Ffuncbold$ which then enables us to invoke this principle. All inequalities are interpreted component-wise.
\begin{lemma}\label{lem:Fprop_feasreg}
Let $\mathcal{U}\subseteq [0,1]^4$ be the largest connected set containing $\bs{1}\coloneqq(1,1,1,1)$ such that $\forall \bs{\mu}\in \mathcal{U}$, $\bs{\mu} \geq \Ffuncbold(\bs{\mu})$. Then we have the followings:
\begin{enumerate}[label=(\roman*)]
	\item $\Ffuncbold(\mathcal{U}) \subseteq \mathcal{U}$.
	\item $\mathcal{U}$ is closed and compact.
	\item $\forall \bs{u}\in \mathcal{U}$, $\lim_{k\to\infty}\Ffuncbold^k(\bs{u})$ converges to some point $\bs{u}_*\in \mathcal{U}$, which is a fixed point of $\Ffuncbold$.
	\item If $\bs{u}_*\in [0,1]^4$ is a fixed point of $\Ffuncbold$, then for any $\bs{u} \geq \bs{u}_*$ such that $\bs{u}$ and $\bs{u}_*$ are equal in at least one component, we have $\bs{u} \notin \mathcal{U}$.
\end{enumerate}
\end{lemma}
\begin{proof}
See Appendix \ref{proof:Fprop_feasreg}.
\end{proof}
An immediate and important corollary of the Lemma \ref{lem:Fprop_feasreg} is the following.
\begin{corollary}\label{cor:Fprop_sol}
	Let $\bs{\mu}_*$ denote the closest fixed point of $\Ffuncbold$ to $\bs{1}$ in sup-norm, i.e.,
	\begin{align*}
	\bs{\mu}_* \coloneqq \argmin_{\bs{u}: \Ffuncbold(\bs{u}) = \bs{u}} \norm{\bs{u} - \bs{1}}_\infty,
	\end{align*}
	where $\norm{\bs{x}-\bs{y}}_\infty \coloneqq \max_i\abs{x_i - y_i}$. Then, we have $\bs{\mu}_*\in \mathcal{U}$, and for all $\bs{u}\in \mathcal{U} \cap \{\bs{x}: \bs{1} \geq \bs{x} \geq \bs{\mu}_*\}$:
	\begin{align*}
	\bs{\mu}_* = \lim_{k\to\infty}\Ffuncbold^k(\bs{u}).
	\end{align*}
\end{corollary}
\begin{proof}
	The proof follows by parts (ii), (iii) and (iv) of Lemma \ref{lem:Fprop_feasreg}, and the fact that $\bs{1}\in \mathcal{U}$.
\end{proof}
Finally, we characterize the equilibrium point at which the ODEs \eqref{eq:altsol_mujj}-\eqref{eq:altsol_muj-j} settles starting from $(1,1,1,1)$ and also provide an iterative method to find it.
\begin{theorem}\label{thm:alterode_sol}
Consider the following system of ODEs:
\begin{align}
\frac{d\bs{\mu}}{dt} = \Ffuncbold(\bs{\mu}) - \bs{\mu} \myquad[4]\bs{\mu}(0) = \bs{1}, \myquad[4]\bs{\mu} \in [0,1]^4.
\label{eq:alterode}
\end{align}
The solution to the ODEs settles at $\bs{\mu}_*$, given by Corollary \ref{cor:Fprop_sol}. Moreover, any point $\bs{\mu}_* < \bs{\mu} < \bs{1}$ of the trajectory of the solution is an interior point of $\mathcal{U}$.
\end{theorem}
\begin{proof}
See Appendix~\ref{proof:alterode_sol}.
\end{proof}

Let us, for the moment, pretend that the initial condition of the ODEs in Appendix \ref{app:odederive} is not random. This assumption is not far from being correct since as $n\to\infty$ the initial values given by \eqref{eq:ode_ic} concentrate around their mean (as we also pointed out in Remark \ref{rem:FandFconection}). Now, solving the ODEs \eqref{eq:alterode}, we obtain the trajectory of the ODEs \eqref{eq:sol_mujj}-\eqref{eq:sol_muj-j} as long as $\bs{\mu}(t)\in \mathcal{D}_{\varepsilon,n}$.

\begin{figure}[t!]
	\centering
	\includegraphics[width=0.6\textwidth]{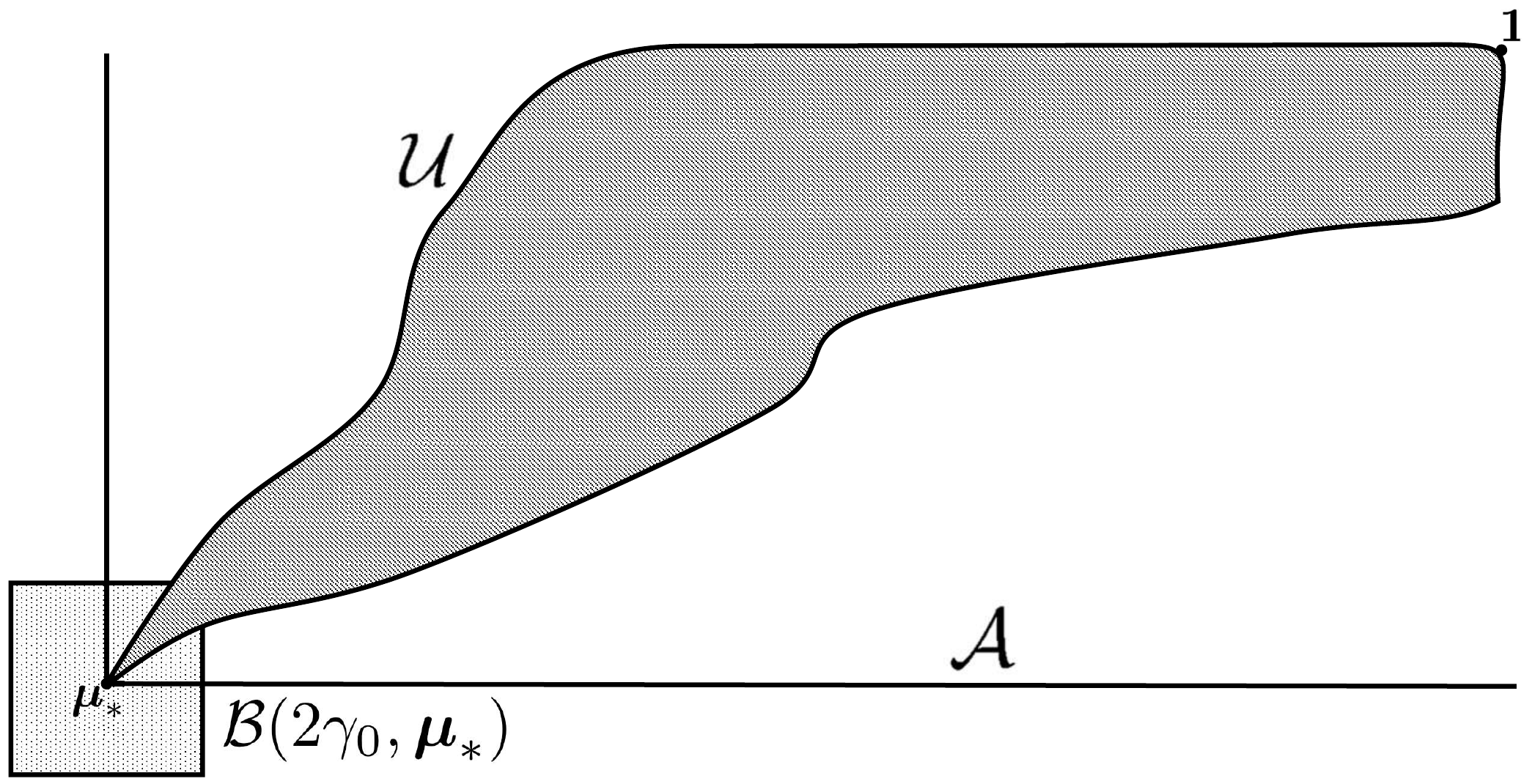}
	\caption{Schematic of $\mathcal{A}$, $\mathcal{U}$, and the choice of $\gamma_0 > 0$ in \eqref{eq:epschoice} .}\label{fig:schematic}
\end{figure}

By part (iv) of Lemma \ref{lem:Fprop_feasreg} and Corollary \ref{cor:Fprop_sol}, the trajectory of ODEs \eqref{eq:alterode} hits the set
\begin{align}
\mathcal{A} \coloneqq \{\bs{x}: \bs{1} \geq \bs{x} \geq \bs{\mu}_* \text{ and }\exists i,j\in\{1,2\}: {\mu}^{(i\gets j)}_* = {x}^{(i\gets j)} \} \label{eq:setA}
\end{align}
at $\bs{\mu}_*$. More specifically, $\mathcal{A}\cap \mathcal{U} = \bs{\mu}_*$. Now, given the fact that both $\mathcal{A}$ and $\mathcal{U}$ are closed and compact, we can pick $\gamma_0>0$ small enough such that
\begin{align}
\{\bs{x}\in\mathcal{U}: \exists y\in\mathcal{A} \text{ such that } \norm{\bs{x} - \bs{y}}_\infty < \gamma_0 \} \subset \mathcal{B}(\bs{\mu}_*,2\gamma_0), \label{eq:epschoice}
\end{align}
where $\mathcal{B}(\bs{\mu}_*,2\gamma_0)$ is a ball of radius $2\gamma_0$ centered at $\bs{\mu}_*$ (in infinity norm). Figure~\ref{fig:schematic} provides a $2$-dimensional schematic for the choice of $\gamma_0>0$.

It is easy to see that the equalities \eqref{eq:ajtoFcon}-\eqref{eq:a-jtoFcon} hold for all $\bs{\mu} \in \mathcal{U}$. Also, note that $\bs{\mu}_*$ is the only point in $\mathcal{U}$ such that $a_1+a_2+a_m^{(1)}+a_m^{(2)} = 0$. This justifies our claim earlier in this section that removing the denominator of the ODEs \eqref{eq:sol_mujj}-\eqref{eq:sol_muj-j} will not affect its trajectory.
Since the function $\Ffuncbold$ is continuous, it is easy to see that there is an $\varepsilon_0 > 0$ such that for all $\varepsilon < \varepsilon_0$,
\begin{align}
\mathcal{U} \cap \{\bs{x}: \bs{1} \geq \bs{x} \geq \bs{\mu}_*\} \setminus \mathcal{B}(\bs{\mu}_*,2\gamma_0) \subset \mathcal{D}_{\varepsilon,n}.% \subset %\mathcal{U} \cap \{\bs{x}: \bs{1} \geq \bs{x} \geq \bs{\mu}_*\}.
\end{align}

Combining the above argument with Lemma \ref{lem:diffeq_sol} and Theorem \ref{thm:odesfinite}, we can track the Markov process of adoption (the truncated version) up to any $\varepsilon$ neighborhood of $\bs{\mu}_*$. Intuitively speaking, if $\bs{\mu}_*$ is a stable equilibrium point of ODEs \eqref{eq:alterode}, then $\bs{\mu}_*$ should correspond to the natural stopping point of the Markov process of adoption. We can also use equation \eqref{eq:sol_t} to estimate the natural stopping time. Define $t_*$ as follows:
\begin{align}
t_* \coloneqq  \lambda_m(n)\left(1-\mu^{(1\gets 2)}_*\mu^{(2\gets 1)}_*\right) + \frac{\lambda_1(n)}{2}\left(1-\left({\mu^{(1\gets 1)}_*}\right)^2 \right)+ \frac{\lambda_2(n)}{2}\left(1-{\left(\mu^{(2\gets 2)}_*\right)}^2\right) . \label{eq:conj_t}
\end{align}
Then, conjecturally, the natural stopping point of the Markov process of adoption is $\sim t_* n$. However, both of these arguments are far from being rigorous. In the following section, we characterize the state of the Markov processes of adoption at its natural stopping time as $n$ goes to infinity.

	\section{Asymptotic Behavior of the Markov Process of Adoption}\label{sec:odeanalysisinf}
	As we mentioned in Point \ref{point:2} of Section \ref{sec:odeapprox}, we have assumed that inactive vertices of degree higher than $d_{\max}$ cannot be activated. In this section, we will study the behavior of the original Markov process of adoption at its natural stopping time for large values of $n$. In what follows, we focus on the truncated version, i.e., vertices with degree higher than $d_{\max}$ cannot be activated; however, this assumption is relaxed in the statement of the main theorem.

Following the discussion at the end of Section \ref{sec:odeanalysisfinite}, for the moment let us pretend that the initial condition of the ODEs in Appendix \ref{app:odederive} is not random. We can track the trajectory of the truncated version of the Markov process of adoption up to $\epsilon$ neighborhood of $\bs{\mu}_*$, after proper scaling. The question is whether the process stops at $\bs{\mu}_*$. If $\bs{\mu}_* = \bs{0}$, then the answer is clear as we have removed almost all the half-edges. However, the same cannot be said if $\bs{\mu}_* \neq \bs{0}$.

Note that there is an essential difference between the case of one community and multiple communities, as the process may move between the communities in the latter case. The first case is much easier and more intuitive to handle, and the latter is more technical. For illustration, we treat these two cases separately, introducing two different but related approaches. Note that the method introduced for the latter case applies to the case of one community and not vice-versa; for this reason, we skip some technical details for the case of one community.

\begin{remark}\label{rem:infnotation}
	We use the index `$\infty$' to refer to terms in the previous sections, when the initial condition \eqref{eq:ode_ic} and the values of $(\lambda_1(n),\lambda_2(n),\lambda_m(n))$ are replaced with their asymptotic counterparts given by Lemma \ref{lem:initcond_conv}. For example, $\bs{\mu}_{*,\infty}$ is the equilibrium point of ODEs \eqref{eq:alterode}, when $\Ffuncbold$ is replaced with $\Ffuncbold_\infty$. Note that $\Ffuncbold_\infty$ is same as the function $\bs{F}$ given by the right-hand side of \eqref{eq:meanfield_mujj}-\eqref{eq:meanfield_muj-j} (mean-field equations) if we assume $K_j(d_j,d_{-j}) = d_j+d_{-j}$ for all $d_j+d_{-j}>d_{\max}$ and $j\in\{1,2\}$ (see Remark \ref{rem:FandFconection}). Also, note that $\Ffuncbold_\infty$ is not random and does not depend on $n$.
\end{remark}

\subsection{Case of One Community}
In this case, the ODEs \eqref{eq:sol_mujj}-\eqref{eq:sol_muj-j} simplify to a one-dimensional ODE as there is only one community from which we can pick an active half-edge uniformly at random. In particular, we have
\begin{align*}
	-1 =  \lambda(n)\frac{d\mu}{dt}\left(\mu(t)\right),\qquad\mu(0) = 1
\end{align*}
which implies that $\mu(t) = \sqrt{1 - 2 t/\lambda(n)}$. In this case, the function $\Ffunc$ is one-dimensional as well and is given as follows:
\begin{align*}
	\Ffunc(\mu) \coloneqq \sum_{\substack{u \leq K(d),\, d\leq d_{\max}}}  \,\frac{d}{\lambda(n)} \,i_{d,0}(0)\, Bi(u;d-1,1-\mu) +  \frac{w(0)}{\lambda(n)},
\end{align*}
where $K(\cdot)$, $\lambda(n)$, $i_{d,0}(0)$ and $w(0)$ are given as before for the case of one community. Now, similar to \eqref{eq:ajtoFcon}-\eqref{eq:a-jtoFcon}, we have
\begin{align*}
	a(t) = \lambda(n) \mu(t) (\mu(t) - \Ffunc(\mu(t))),
\end{align*}
and $\mu_*$ is given by Corollary \ref{cor:Fprop_sol}. Equivalently, as we have the closed-form solution of $\mu(t)$, $\mu_*$ can be written as follows:
\begin{align}
	\mu_* = \mu(t_*), \text{ where } t_*\coloneqq \inf\{t \in [0,\lambda(n)/2]: \mu(t) - \Ffunc(\mu(t)) = 0\} \label{eq:t*_onecom}
\end{align}
Note that the above set is non-empty as $\mu(\lambda(n)/2) - \Ffunc(\mu(\lambda(n)/2)) \leq 0$, $\mu(0) - \Ffunc(\mu(0)) \geq 0$, and $\mu(t) - \Ffunc(t)$ is a continuous function of $t$. Similar results were reported by authors in \cite{Balogh2007} for $d$-regular random graphs and \cite{Amini2010} for random-graphs given by configuration model. Note that we can track the Markov process of adoption up to $\mu_*$ (as the denominator is $1$ in the case of one community); however, to show that the process stops at $\mu_*$, more work needs to be done.

The basic idea is as follows: $(1)$ augment the truncated process by adding an active vertex with high degree, $(2)$ couple the augmented process and the truncated process, and $(3)$ show that the truncated process hits its natural stopping time before the augmented process passes the conjectured stopping time, with high probability. To show the last step, we need to assume that $\mu_*$ is a stable equilibrium point.

Let us denote the truncated process by $X^n_\delta$ (which is either $X^n_{L,\delta}$ or $X^n_{U,\delta}$, see Point \ref{point:2}). We augment $X^n_\delta$ by adding one active vertex $\widetilde{v}$ with $\floor{2\epsilon n}$ half-edges. Let us denote the augmented process by $\widetilde{X}^n_{\delta,\epsilon}$. Since $\widetilde{X}^n_{\delta,\epsilon}$ has more active half-edges than $X^n_\delta$ ($\floor{2\epsilon n}$ more active half-edges at time $0$), we would expect the size of the cascade in the augmented process to be larger. Next, we show that this is indeed the case by constructing a coupling between $X^n_\delta$ and $\widetilde{X}^n_{\delta,\epsilon}$.

As we mentioned earlier in Point \ref{point:2}, the order in which the active-half edges are paired with a random half-edge does not affect the natural stopping time of the process. Also, note that labels of active vertices are irrelevant in the Markov process of adoption, as it only tracks the number of active half-edges. Keeping these two observations in mind, we introduce a refinement to the augmented process $\widetilde{X}^n_{\delta,\epsilon}$ and alter the way half-edges are removed.

\begin{figure}
	\centering
	\includegraphics[width=0.8\textwidth]{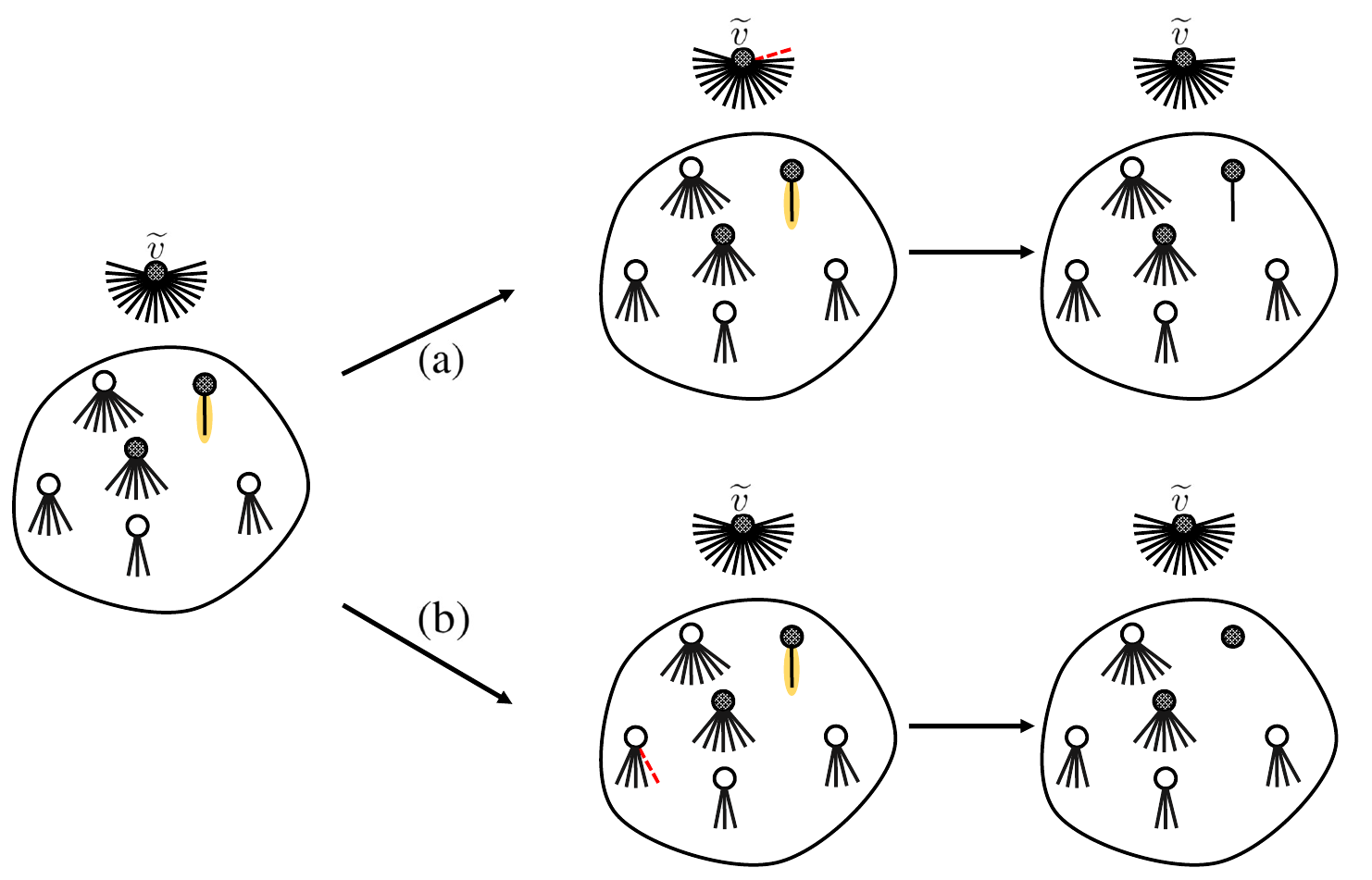}
	\caption{The update rule for the case $M=0$: {\normalfont (a)} the randomly selected half-edge belongs to $\widetilde{v}$; {\normalfont (b)} the randomly selected half-edge does not belong to $\widetilde{v}$. Highlighted active half-edge is denoted by yellow shade, randomly selected half-edge is denoted by dashed-red line, and active vertices are denoted by dotted circles.} \label{fig:coupling_onecomm}
\end{figure}

Recall that during the Markov process of adoption given by Algorithm \ref{algorithm}, we pick one active half-edge and then another half-edge uniformly at random and remove them both. Consider a typical state of the process. Based on the available active half-edges and the choice of random half-edge in the augmented process, we remove different half-edges according to the following rules. To keep track of an important event, we introduce a mark variable $M$, which is initialized to be zero.
\begin{enumerate}[label = (\roman*):]
	\item $M=0$ and there are active half-edges other than the ones that belong to $\widetilde{v}$: highlight one of the active half-edges that does not belong to $\widetilde{v}$ as a potential active half-edge to be removed. Pick another half-edge uniformly at random. If the random half-edge does not belong to $\widetilde{v}$, then remove both the highlighted active half-edge and the random half-edge. Otherwise, remove two active half-edges from $\widetilde{v}$. An example of the update rule for the case $M=0$ is shown in Figure \ref{fig:coupling_onecomm}.
	\item $M=1$ or the only remaining active half-edges (if any) are the ones that belong to $\widetilde{v}$: set $M$ to be one and proceed regularly, i.e., pick one active half-edge and another half-edge uniformly at random and remove them both.
\end{enumerate}

Note that $M = 0$ for some time, and then it becomes $1$ for the rest of the process. It is also possible to have $M=0$ for the entire process, after which we set $M$ to be $1$, as is mentioned above. Also, note that the above refinements do not affect the trajectory of the augmented process $\widetilde{X}^n_{\delta,\epsilon}$.

Next, we couple the truncated and augmented processes as follows: in the truncated process, follow the same order of active half-edges (other than the ones that belong to $\widetilde{v}$) that have been picked in the augmented process and pair each with the same random half-edge. It is easy to see that the resulted coupling has the desired marginal distribution as long as $M=0$. On the other hand, the coupled truncated process hits its natural stopping time whenever the value of $M$ in the augmented processes changes from $0$ to $1$.

Let us denote the equilibrium points of ODE \eqref{eq:alterode} associated with the truncated process and the augmented process with $\mu_*\in\mathbb{R}_+$ and $\widetilde{\mu}_*\in\mathbb{R}_+$ respectively (given by \eqref{eq:t*_onecom}). It is easy to see that $0\leq \widetilde{\mu}_* \leq \mu_*$. Assume that $\mu_* > 0$ is a stable equilibrium point. Hence, $\widetilde{\mu}_*$ and $\mu_*$ can be made arbitrarily close to each other by setting $\epsilon>0$ to be small enough.

\begin{figure}
	\centering
	\includegraphics[width=0.6\textwidth]{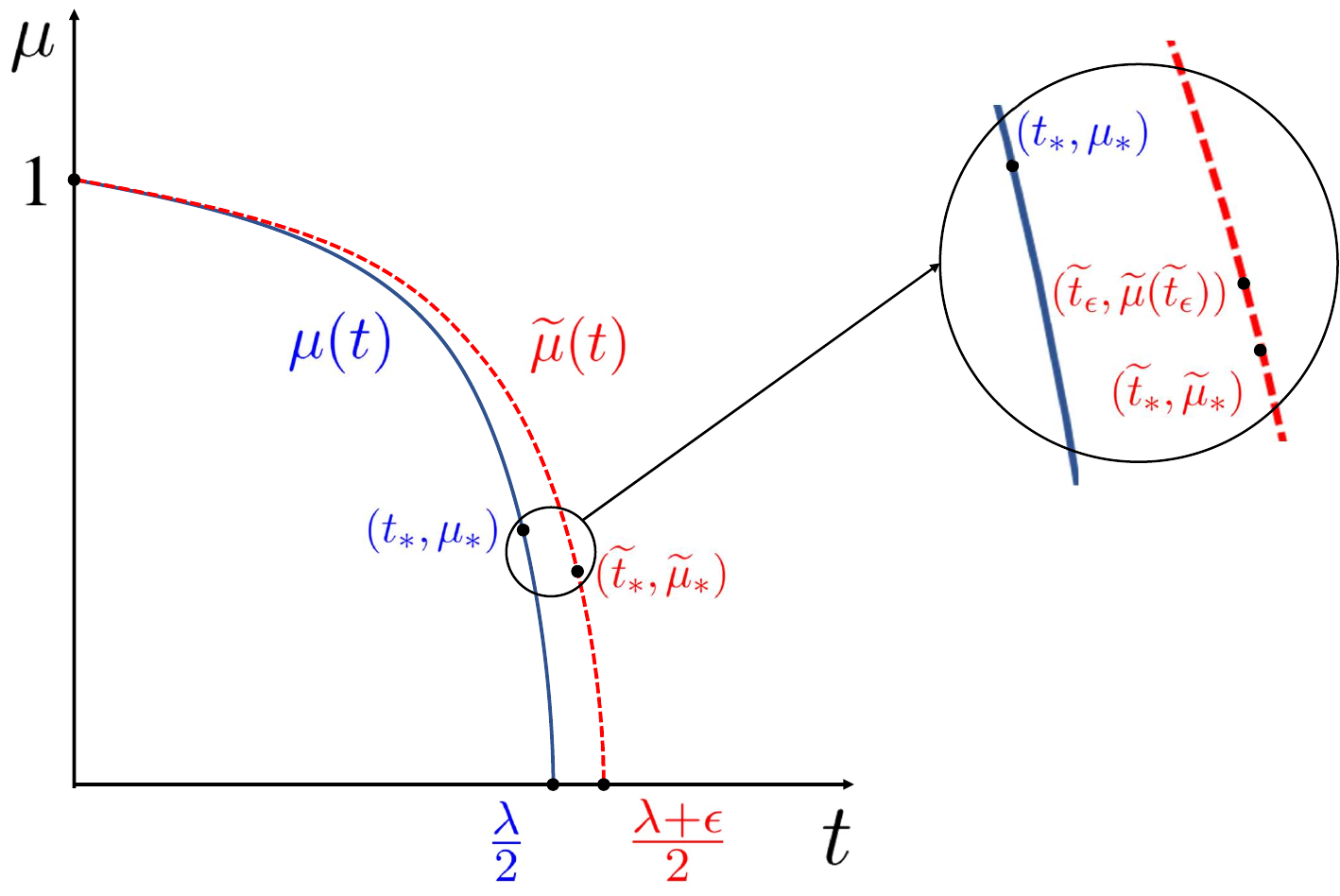}
	\caption{A schematic of the trajectories of the truncated and augmented process, and the value of $\widetilde{t}_\epsilon$, given a fixed realization of the initial seed set; $\widetilde{t}_\epsilon$ is chosen so that $\widetilde{a}(\widetilde{t}_\epsilon) = \lambda \widetilde{\mu}(\widetilde{t}_\epsilon) (\widetilde{\mu}(\widetilde{t}_\epsilon) - \widetilde{\Ffunc}(\widetilde{\mu}(\widetilde{t}_\epsilon))) = \epsilon^2/2$. Solid-blue line is the trajectory of the truncated process, and dashed-red color is the trajectory of the augmented process.
	}\label{fig:teps_onecomm}
\end{figure}

As we mentioned, we can track the augmented process up to points arbitrarily close to $\widetilde{\mu}_*$. In particular, we can track the augmented process up to time  $\widetilde{t}_\epsilon n$, where the total number of active half-edges falls below $\epsilon^2 n$ (incorporating $o(n)$ error), with high probability, before reaching $\widetilde{\mu}_*$. See Figure \ref{fig:teps_onecomm} for a schematic of the choice of $\widetilde{t}_\epsilon$. At this point, we have already used most of the active half-edges associated with $\widetilde{v}$. The claim is that if $\epsilon > 0$ is small enough, then $M=1$ with high probability.

Recall that while $M=0$, we remove two active half-edges from $\widetilde{v}$ if and only if after highlighting an active half-edge, the uniformly selected random half-edge belongs to $\widetilde{v}$. Using an ODE-based approach, similar to what we did in the previous sections, it is easy to see that the number of times such events happen up to time $k$ of the augmented process concentrates around $2\epsilon k/(\lambda(n) + 2\epsilon)$. Moreover, we can pick $\epsilon > 0$ to be small enough such that $\widetilde{\mu}_* > 0$ is in a small neighborhood of $\mu_*$ and the total number of remaining half-edges at time $\widetilde{t}_\epsilon n$ is larger than $(\lambda(n) + 2\epsilon)\epsilon n$. In particular,  $(\lambda(n) + 2\epsilon) n - 2\widetilde{t}_\epsilon n > (\lambda(n) + 2\epsilon)\epsilon n$ which implies that $\widetilde{t}_\epsilon < (\lambda(n) + 2\epsilon) (1-\epsilon)/2$ (note that $\lim_{\epsilon\to 0}(\lambda(n) + 2\epsilon) (1-\epsilon)/2 = \lambda(n)/2$). On the other hand, if at time $\widetilde{t}_\epsilon n$ we have $M=0$, then the number of times that we have removed two active half-edges from $\widetilde{v}$ is smaller than
\begin{align*}
	\frac{2\epsilon}{\lambda(n) + 2\epsilon} \times n(\lambda(n) + 2\epsilon) (1-\epsilon)/2  = \epsilon (1-\epsilon)n
\end{align*}
with high probability. Note that ${2\epsilon}\mathbin{/}{(\lambda(n) + 2\epsilon)}$ is (almost) the ratio of the augmented half-edges to the regular half-edges. Conditioned on sample paths that satisfy the above inequality, the total number of active half-edges that belong to $\widetilde{v}$ (originally) is smaller than $2 \epsilon (1-\epsilon)n + \epsilon^2 n = 2\epsilon n - \epsilon^2 n$, which is a contradiction.

Note that in the above argument, the initial condition is random, $n$ is fixed, and the choice of $\epsilon$ depends on $\mu_*$, which in turn depends on the initial condition; hence, $\epsilon$ is also a random variable. However, as $n\to\infty$, the initial condition converges in probability to its mean (using a similar argument as in Lemma \ref{lem:initcond_conv} for the case of one community) and $\lambda(n)\to\lambda$. Hence, the value of $\epsilon$ also converges in probability to a constant; that is to say, we can pick a non-random small enough $\epsilon > 0$ constant such that $\mu_*$ and $\widetilde{\mu}_*$ are close enough to each other, with high probability.

To summarize, for any small enough $\epsilon > 0$, we have defined an event $\Omega_n(\epsilon)$ with $\lim_{n\to\infty}\! P(\Omega_n(\epsilon))\allowbreak = 0$ such that outside $\Omega_n(\epsilon)$ a scaled-version of the truncated process hits its natural stopping time at time $\left(t_*\pm O(\epsilon)\right)n$. Letting $\epsilon\to 0$, we can characterize the asymptotic behavior of the Markov process of adoption at its natural stopping time (see the statement of Theorem \ref{thm:odesmaininf}). Note that $\Omega_n(\epsilon)$ is defined by considering the union of the complements of finitely many high probability events.

\subsection{Case of Multiple Communities}
As we have mentioned before, we only need to focus on the case of two communities. Note that the same approach as in the previous subsection cannot be used here: the truncated process may run out of active half-edges in one community, while the augmented process keeps picking active half-edges from the same community. The basic idea in the case of multiple communities is to introduce a perturbation to the truncated process when it reaches a small neighborhood of the conjectured stopping time, i.e., $t_* n$ where $t_*$ is given by \eqref{eq:conj_t}.

Let us denote the truncated process by $X^n_\delta$ (which is either $X^n_{L,\delta}$ or $X^n_{U,\delta}$, see Point \ref{point:2}). Recall from Theorem~\ref{thm:odesfinite} that for any fixed $\varepsilon > 0$, we can track the process as long as the ODEs are within $l^{\infty}$-distance $C\theta$ of the boundary of $\widehat{\mathcal{D}}_{\varepsilon,n}$, for a large enough constant $C$ independent of $n$. Equivalently, we can track $X^n_\delta$ using Lemma \ref{lem:diffeq_sol}, and the solution of ODEs \eqref{eq:alterode} as long as the trajectory of the solution is within $C'\theta$ of the boundary of $\mathcal{D}_{\varepsilon,n}$, for some large enough constant $C'$ related to $C$. Also, recall that $\widehat{\mathcal{D}}_{\varepsilon,n}$ is a fixed set and $\mathcal{D}_{\varepsilon,n}$ depends on the initial condition \eqref{eq:ode_ic} which is random.

The rest of this subsection is organized as follows. In Section \ref{subseubsec:constants}, we define a series of constants associated with the asymptotic behavior of the ODEs. In Section \ref{subseubsec:finite}, we present a set of initial conditions \eqref{eq:ode_ic} for which the behavior of the resulting ODEs are not far from the behavior of the asymptotic ODEs (see Remark \ref{rem:infnotation}). We next focus on estimating the stopping time of the truncated process for any such realization of the initial condition. In Section \ref{subseubsec:augment}, we augment the truncated process by adding extra active vertices. We then analyze the sample path of the resulted augmented process. The basic idea is to show that this augmentation cannot initiate a larger cascade. To make this argument rigorous, in Section \ref{subseubsec:twist}, we define a twisted process which is used as a bridge to couple the augmented process and the truncated process. The details of these couplings and their implications are discussed in Section \ref{subseubsec:coupling}, after which we present the main theorem of the paper and discuss the generalization to the case of $k$ communities. Some of the details are postponed to Appendix \ref{app:augment} for readability.

\subsubsection{Constants Associated with Asymptotic ODEs} \label{subseubsec:constants}
We proceed with defining some constants associated with the limiting behavior of the ODEs that will be used later in the proof. We use the notational conventions discussed in Remark \ref{rem:infnotation}. Without loss of generality, we assume that $\bs{\mu}_{*,\infty} > \bs{0}$. Note that the argument is similar for the other cases: $(1)$ we augmented half-edges in community $j\in\{1,2\}$ only if ${\mu}_{*,\infty}^{(j\gets j)} > 0$, $(2)$ we augment half-edges between the two communities only if ${\mu}_{*,\infty}^{(1\gets 2)} {\mu}_{*,\infty}^{(2\gets 1)} > 0$, and $(3)$ we only keep track of changes that the augmentation may cause. In particular, in what follows, we are only interested in rows and columns of the Jacobian matrix of $\Ffuncbold_\infty(\cdot)$ at $\bs{\mu}_{*,\infty}$ in which augmentation has happened.

Let $\bs{J}_{\Ffuncbold_\infty(\cdot)}(\bs{\mu}_{*,\infty})$ denote the Jacobian matrix of $\Ffuncbold_\infty(\cdot)$ at $\bs{\mu}_{*,\infty}$:
\begin{center}
	\begin{tikzpicture}[
		every left delimiter/.style={xshift=.45em},
		every right delimiter/.style={xshift=-.45em},
		]
		\matrix[matrix of math nodes, left delimiter={[},right delimiter={],}, nodes={scale = 1.3, minimum height=5ex, inner sep=4pt},] (M)
		{
			\frac{\partial \Ffunc_{(1\gets 1),\infty}(\bs{\mu}_{*,\infty})}{\partial \mu^{(1\gets 1)}}  &
			\frac{\partial \Ffunc_{(1\gets 1),\infty}(\bs{\mu}_{*,\infty})}{\partial \mu^{(1\gets 2)}}  & 0 & 0 \\

			0 & 0 &  \frac{\partial \Ffunc_{(1\gets 2),\infty}(\bs{\mu}_{*,\infty})}{\partial \mu^{(2\gets 1)}}  &
			\frac{\partial \Ffunc_{(1\gets 2),\infty}(\bs{\mu}_{*,\infty})}{\partial \mu^{(2\gets 2)}} \\

			\frac{\partial \Ffunc_{(2\gets 1),\infty}(\bs{\mu}_{*,\infty})}{\partial \mu^{(1\gets 1)}}  &
			\frac{\partial \Ffunc_{(2\gets 1),\infty}(\bs{\mu}_{*,\infty})}{\partial \mu^{(1\gets 2)}} & 0 & 0\\

			0 & 0 &  \frac{\partial \Ffunc_{(2\gets 2),\infty}(\bs{\mu}_{*,\infty})}{\partial \mu^{(2\gets 1)}}  &
			\frac{\partial \Ffunc_{(2\gets 2),\infty}(\bs{\mu}_{*,\infty})}{\partial \mu^{(2\gets 2)}}\\
		};
		%\node[left=0.5em of M] (equal) {$\bs{J}_{\Ffuncbold_\infty(\cdot)}(\bs{\mu}_{*,\infty})=$} ;
		\lineaftercolumn{M}{1}\lineaftercolumn{M}{2}\lineaftercolumn{M}{3}
		\linebelowrow{M}{1}\linebelowrow{M}{2}\linebelowrow{M}{3}
	\end{tikzpicture}\\
\end{center}
where the elements of $\bs{J}_{\Ffuncbold_\infty(\cdot)}(\bs{\mu}_{*,\infty})$ are given in Appendix \ref{proof:jacobianmatrix}. Let $\zeta_{*,\infty}$ denote its largest eigenvalue. Note that $\bs{J}_{\Ffuncbold_\infty(\cdot)}(\bs{\mu}_{*,\infty})$ is a non-negative aperiodic and irreducible matrix, i.e., a primitive matrix. Hence, by the Perron-Frobenius theorem, $\zeta_{*,\infty}>0$ is a simple eigenvalue, and the corresponding eigenvector $\bs{\nu}_{*,\infty}$ is element-wise positive.

Let us pretend for the moment that the interchange of limits is allowed. Note that elements of $\bs{J}_{\Ffuncbold_\infty(\cdot)}(\bs{\mu}_{*,\infty})$ are related to vertices that are one half-edge away of becoming active at the scaled-time $t_{*,\infty}$, a.k.a. pivotal players at the scaled-time $t_{*,\infty}$. In particular, if $\zeta_{*,\infty} > 1$, then the subgraph restricted to these vertices consists of a giant component, and if $\zeta_{*,\infty} < 1$, then this subgraph is a union of many small components; both statements being true with high probability. This is the essential idea behind the proof: if $\zeta_{*,\infty} < 1$, then the cascade cannot grow much further as the pivotal players are the main drivers of the cascade when there is a lack of active half-edges.

Suppose that $\zeta_{*,\infty} < 1$, which implies that $\bs{\mu}_{*,\infty}$ is a stable equilibrium point of ODEs \eqref{eq:alterode}. For any $\kappa\in(0,1)$, let $\bs{\rho}_{*,\infty}(\kappa) \coloneqq \kappa \bs{\nu}_{*,\infty}$, i.e., for $j\in\{1,2\}$:
\begin{align*}
	\rho_{*,\infty}^{(j\gets j)}(\kappa) \coloneqq \kappa \nu_{*,\infty}^{(j\gets j)} \qquad\text{and} \qquad\rho^{(-j\gets j)}_{*,\infty}(\kappa) \coloneqq \kappa \nu_{*,\infty}^{(-j\gets j)}.
\end{align*}

For any $\kappa\in(0,1)$ and $j\in\{0,1\}$, let
\begin{align*}
	\chi_{*,\infty}^{(j\gets j)}(\kappa) \coloneqq \mu_{*,\infty}^{(j\gets j)} \rho_{*,\infty}^{(j\gets j)}(\kappa) \qquad\text{and} \qquad\chi_{*,\infty}^{(j\gets -j)}(\kappa) \coloneqq \mu_{*,\infty}^{(-j\gets j)} \rho_{*,\infty}^{(j\gets -j)}(\kappa).
\end{align*}
Loosely speaking, $\bs{\chi}_{*,\infty}$ determines the perturbation that we are going to introduce to the graph by adding augmented active vertices. In particular, the number of augmented active half-edges that we are going to add to community $j$ for in-community connections is (almost) proportion to $\chi_{*,\infty}^{(j\gets j)}$, and the number of augmented active half-edges that we are going to add to community $-j$ to be paired with half-edges in community $j$ is (almost) proportion to $\chi_{*,\infty}^{(j\gets -j)}$. The fact that $\mu_{*,\infty}^{(-j\gets j)}$ appears in the definition of $\chi_{*,\infty}^{(j\gets -j)}(\kappa)$ is an artifact of our proof.

For any $\kappa \in (0,1)$, let $\bs{\rho}_{*,\infty}(\kappa,\pm \textrm{err}) \coloneqq \bs{\rho}_{*,\infty}(\kappa)\pm \kappa^2\bs{1}$ and $\bs{\mu}_{*,\infty}(\kappa,\pm \textrm{err}) \coloneqq \bs{\mu}_{*,\infty}\pm \kappa\bs{1}$. Similarly, for $j\in\{1,2\}$ and $\kappa\in(0,1)$, let
\begin{align*}
	&\chi_{*,\infty}^{(j\gets j)}(\kappa,\pm\err) \coloneqq \mu_{*,\infty}^{(j\gets j)}(\kappa,\pm\err) \rho_{*,\infty}^{(j\gets j)}(\kappa,\err),\\
	&\chi_{*,\infty}^{(j\gets -j)}(\kappa,\pm\err) \coloneqq \mu_{*,\infty}^{(-j\gets j)}(\kappa,\pm\err) \rho_{*,\infty}^{(j\gets -j)}(\kappa,\pm\err).
\end{align*}
Let $\kappa_0 \in (0,0.5]$ be the largest constant for which $\bs{\rho}_{*,\infty}(\kappa_0,\pm \textrm{err}) \geq \bs{0}$, $\bs{\mu}_{*,\infty}(\kappa_0,\pm \textrm{err}) \geq \bs{0}$, and all non-zero elements of $\bs{J}_{\Ffuncbold_\infty(\cdot)}(\bs{\mu}_{*,\infty})$ are larger than or equal to $\kappa_0$. These error terms are introduced to offset the effect of finite $n$ and the resulted approximation.

Let $\kappa_1 \in (0,\min(1 - \zeta_{*,\infty},\kappa_0)]$ be the largest constant given which for all $j\in\{1,2\}$, we have:
\begin{align}
	&\zeta_{*,\infty} + \kappa_1 \leq \frac{\left({\mu}_{*,\infty}^{(j\gets j)}(\kappa_1,-\err)\right)^2}{2\chi_{*,\infty}^{(j\gets j)}(\kappa_1,+\err) + \left({\mu}_{*,\infty}^{(j\gets j)}(\kappa_1,+\err)\right)^2},\label{eq:propineq_incomm}\\
	&\begin{aligned}
		&\left(1 - \zeta_{*,\infty} - \kappa_1\right) \chi_{*,\infty}^{(j\gets -j)}(\kappa_1,-\err) \geq \\
		&\myquad[4]\frac{\left(\chi_{*,\infty}^{(1\gets 2)}(\kappa_1,+\err) +  \chi_{*,\infty}^{(2\gets 1)}(\kappa_1,+\err)\right)^2}{\chi_{*,\infty}^{(1\gets 2)}(\kappa_1,-\err) +  \chi_{*,\infty}^{(2\gets 1)}(\kappa_1,-\err) + {\mu}_{*,\infty}^{(1\gets 2)}(\kappa_1,-\err){\mu}_{*,\infty}^{(2\gets 1)}(\kappa_1,-\err)}.
	\end{aligned}\label{eq:propineq_outcomm}
\end{align}
Note that \eqref{eq:propineq_incomm} is strict for $\kappa_1 = 0$, and both side of this inequality are continuous and monotone with respect to $\kappa_1$. Also, note that the same statement holds for \eqref{eq:propineq_outcomm}, after dividing both sides of the inequality by $\kappa_1$. Hence, there exists $\kappa_1 > 0$ such that \eqref{eq:propineq_incomm} and \eqref{eq:propineq_outcomm} are valid. Moreover, the above inequalities are strict for any $\kappa \in (0,\kappa_1)$. For a certain set of initial conditions, the right-hand side of \eqref{eq:propineq_incomm} is a lower bound for the proportion of the regular half-edges over the total number of half-edges, for in-community connections, at the point of augmentation, with high probability. Similarly, the right-hand side of \eqref{eq:propineq_outcomm} is an upper bound for the proportion of the augmented half-edges over the total number of half-edges, between the two communities, at the point of augmentation, with high probability. These inequalities are used in Section \ref{subseubsec:coupling} to bound the size of the cascade in the twisted process. These terms will be defined later in Sections \ref{subseubsec:constants}--\ref{subseubsec:coupling}.

Let $\{\bs{\mu}_\infty(t) \text{ for }t>0 \}$ denote the solution of the ODEs \eqref{eq:alterode}, using notational convention introduced in Remark \ref{rem:infnotation}. For any $\kappa \in(0,\kappa_1)$, let $r(\kappa) > 0$ be the largest constant for which the boundary of the ball of radius $r(\kappa)$ in infinity norm centered at $\bs{\mu}_{*,\infty}$ hits the trajectory of $\{\bs{\mu}_\infty(t) \text{ for }t>0 \}$ at a point $t_{r(\kappa),\infty}$ for which
\begin{align*}
		& \bs{0} < \bs{\mu}_\infty(t_{r(\kappa),\infty}) - \bs{\mu}_{*,\infty} \leq \bs{\rho}_{*,\infty}(\kappa,-\err)/2,\text{ and }\\
		&a_{1,\infty}(t_{r(\kappa),\infty}) + a_{2,\infty}(t_{r(\kappa),\infty}) + a_{m,\infty}^{(1)}(t_{r(\kappa),\infty}) + a_{m,\infty}^{(2)}(t_{r(\kappa),\infty}) \leq \kappa^2/2,
\end{align*}
where $a_{j,\infty}(\cdot)$ and $a_{m,\infty}^{(j)}(\cdot)$ are given by Remark \ref{rem:atoFcon} following the notation convention in Remark \ref{rem:infnotation}. These inequalities are used to ensure that the process has been executed until a time close to the conjectured stopping time before adding the augmented vertice, at which point the total number of remaining active half-edges is small.

For any $\kappa \in(0,\kappa_1)$, let  $\varepsilon(\kappa)  \in \left(0,0.5\right]$ be the largest constant for which
\begin{align*}
	\mathcal{U}_\infty\cap \{\bs{\mu}:\bs{1} \geq \bs{\mu} \geq \bs{\mu}_{*,\infty}\} \setminus \mathcal{B}\left(\bs{\mu}_{*,\infty},r(\kappa)/2\right)  \subset \overline{\mathcal{D}_{2\varepsilon(\kappa),\infty}},
\end{align*}
where $\mathcal{D}_{\cdot,\infty}$ is defined similar to $\mathcal{D}_{\cdot,n}$ (given by Lemma \ref{lem:diffeq_sol}) following Remark \ref{rem:infnotation}.
The existence of such $\varepsilon(\kappa)$ follows from Lemma \ref{lem:Fprop_feasreg}, using a similar argument as is presented at the end of Section \ref{sec:odeanalysisfinite}. An important implication is that the trajectory of $\bs{\mu}_{\infty}(t)$ is almost entirely in $\mathcal{D}_{2\varepsilon(\kappa),\infty}$ except for some parts of the trajectory that are in $\mathcal{B}\left(\bs{\mu}_{*,\infty},r(\kappa)/2\right)$. This ensures the validity of Theorem \ref{thm:odesfinite} for a certain set of initial conditions, so we can approximate the trajectory of the process up to a time close to the conjectured stopping time. The above symbols and their definitions are summarized in Table \ref{table:constants}.
\begin{table}[!htbp]
	\centering
	\begin{tabular}{| l | l |}
		\hline
		Symbol & Definition \Tstrut \Bstrut \\ \hline
		$\bs{J}_{\Ffuncbold_\infty(\cdot)}(\bs{\mu}_{*,\infty})$ & Jacobian matrix of $\Ffuncbold_\infty(\cdot)$ at $\bs{\mu}_{*,\infty}$;\Tstrut \Bstrut \\ \hline
		$\zeta_{*,\infty}$ & largest eigenvalue of $\bs{J}_{\Ffuncbold_\infty(\cdot)}(\bs{\mu}_{*,\infty})$, which is simple and positive;\Tstrut \Bstrut \\ \hline
		$\bs{\nu}_{*,\infty}$ & eigenvector of $\bs{J}_{\Ffuncbold_\infty(\cdot)}(\bs{\mu}_{*,\infty})$ corresponding to eigenvalue $\zeta_{*,\infty}$, which is strictly\Tstrut \Bstrut \\
		& positive;\Tstrut \Bstrut \\ \hline
		$\bs{\rho}_{*,\infty}(\kappa)$ & $\kappa \bs{\nu}_{*,\infty}$;\Tstrut \Bstrut \\ \hline
		$\bs{\chi}_{*,\infty}(\kappa)$ & $\big({\rho}_{*,\infty}^{(1\gets 1)}(\kappa) {\mu}^{(1\gets 1)}_{*,\infty}(\kappa)\!,\!
		{\rho}_{*,\infty}^{(1\gets 2)}(\kappa) {\mu}^{(2\gets 1)}_{*,\infty}(\kappa)\!,\!
		{\rho}_{*,\infty}^{(2\gets 1)}(\kappa) {\mu}^{(1\gets 2)}_{*,\infty}(\kappa)\!,\!
		{\rho}_{*,\infty}^{(2\gets 2)}(\kappa) {\mu}^{(2\gets 2)}_{*,\infty}(\kappa)\big)$;\DTstrut \DBstrut \\ \hline
		$\bs{\rho}_{*,\infty}(\kappa,\pm\err)$ & $\bs{\rho}_{*,\infty}(\kappa)\pm \kappa^2\bs{1}$;\Tstrut \Bstrut \\ \hline
		$\bs{\mu}_{*,\infty}(\kappa,\pm\err)$ & $\bs{\mu}_{*,\infty}\pm \kappa\bs{1}$;\Tstrut \Bstrut \\ \hline
		$\bs{\chi}_{*,\infty}(\kappa,\pm\err)$ & defined similar to $\bs{\chi}_{*,\infty}$, using $\bs{\rho}_{*,\infty}(\kappa,\pm\err)$ and $\bs{\mu}_{*,\infty}(\kappa,\pm\err)$;\Tstrut \Bstrut \\ \hline
		$\kappa_0$	& largest constant in $(0,0.5]$ for which $\bs{\rho}_{*,\infty}(\kappa,\pm \textrm{err}) \geq \bs{0}$, $\bs{\mu}_{*,\infty}(\kappa,\pm \textrm{err}) \geq \bs{0}$, and \Tstrut\Bstrut\\
		& all non-zero elements of $\bs{J}_{\Ffuncbold_\infty(\cdot)}(\bs{\mu}_{*,\infty})$ are greater than $\kappa_0$;\Tstrut\DBstrut\\ \hline
		$\kappa_1$	& largest constant in $(0,\min(1 - \zeta_{*,\infty},\kappa_0)]$ for which the following inequalities\Tstrut\DBstrut \\
		& hold for all $j\in\{1,2\}$:\Tstrut\DBstrut \\
		&\multicolumn{1}{c|}{
		$\begin{aligned}
			&\zeta_{*,\infty} + \kappa_1 \leq \frac{\left({\mu}_{*,\infty}^{(j\gets j)}(\kappa_1,-\err)\right)^2}{2\chi_{*,\infty}^{(j\gets j)}(\kappa_1,+\err) + \left({\mu}_{*,\infty}^{(j\gets j)}(\kappa_1,+\err)\right)^2},\text{ and }\\
			&\left(1 - \zeta_{*,\infty} - \kappa_1\right) \chi_{*,\infty}^{(j\gets -j)}(\kappa_1,+\err) \geq \\
			&\myquad[1]\frac{\left(\chi_{*,\infty}^{(j\gets -j)}(\kappa_1,-\err) +  \chi_{*,\infty}^{(-j\gets j)}(\kappa_1,-\err)\right)^2}{\chi_{*,\infty}^{(j\gets -j)}(\kappa_1,+\err) +  \chi_{*,\infty}^{(-j\gets j)}(\kappa_1,+\err) + {\mu}_{*,\infty}^{(j\gets j)}(\kappa_1,+\err){\mu}_{*,\infty}^{(j\gets j)}(\kappa_1,+\err)};
		\end{aligned}$}\DTstrut\DDBstrut\\\hline
		$r(\kappa)$ & largest positive constant for which \Tstrut\DBstrut\\
		& \multicolumn{1}{c|}{$\begin{aligned}
				&\exists \bs{\mu}_\infty(t_{r(\kappa),\infty}) \in \{\bs{\mu}_\infty(t) \text{ for }t>0 \}\cap \partial{\mathcal{B}(\bs{\mu}_{*,\infty},r(\kappa))}\\
				&\myquad[8]\bs{0} < \bs{\mu}_{*,\infty} - \bs{\mu}_\infty(t_{r(\kappa),\infty}) \leq \bs{\rho}_{*,\infty}(\kappa,-\err)/2;
			\end{aligned}$} \Tstrut\DBstrut\\ \hline
		 $\mathcal{U}_\infty$ & largest connected set in $[0,1]^4$ containing $\bs{1}\coloneqq(1,1,1,1)$ such that \Tstrut\Bstrut\\
		 & \multicolumn{1}{c|}{$\bs{\mu} \geq \Ffuncbold_\infty(\bs{\mu}), \qquad \forall \bs{\mu}\in \mathcal{U}_\infty$;} \Tstrut \DBstrut \\ \hline
		 $\epsilon(\kappa)$& largest constant in $(0,1]$ for which \Tstrut \DBstrut \\
		 & \multicolumn{1}{c|}{$\mathcal{U}_\infty\cap \{\bs{\mu}:\bs{1} \geq \bs{\mu} \geq \bs{\mu}_{*,\infty}\} \setminus \mathcal{B}\left(\bs{\mu}_{*,\infty},r(\kappa)/2\right)  \subset \overline{\mathcal{D}_{2\varepsilon(\kappa),\infty}}$;}\Tstrut\DBstrut\\
		 \hline
	\end{tabular}
	\caption{Symbols and their definitions for asymptotics. Note that $\bs{\mu}_{*,\infty}$ is the equilibrium point of the ODEs \eqref{eq:alterode} and is not random (see Remark \ref{rem:infnotation}).}\label{table:constants}
\end{table}

\subsubsection{Set of Desirable Initial Conditions for Finite ODEs}\label{subseubsec:finite}
Next, we define the set of initial conditions \eqref{eq:ode_ic} for which we can estimate the stopping time of the process with high probability. For $n > 0$ and $\kappa < \kappa_1/2$, let $\mathcal{E}(n,\kappa)$ denote the set of realizations of initial condition \eqref{eq:ode_ic} for which the following hold:
\begin{enumerate}[leftmargin=!,itemindent=1em,align=left,label=Condition $\mathcal{E}.\arabic*$:,ref=$\mathcal{E}.\arabic*$]\setlength\itemsep{0.5em}
	\item({\it approximation of fixed point, Jacobian matrix, and initial values}) $\norm{\bs{\mu}_* - \bs{\mu}_{*,\infty}}_\infty < \kappa$, $\norm{\bs{J}_{\Ffuncbold(\cdot)}(\bs{\mu}_{*})- \bs{J}_{\Ffuncbold_\infty(\cdot)}(\bs{\mu}_{*,\infty})}_\infty < \kappa$, and the initial values \eqref{eq:ode_ic} are in a $\kappa$ neighborhood of their asymptotic values given by Lemma \ref{lem:initcond_conv};\label{def:mathcalE_i}
	\item({\it approximation of Perron-Frobenius eigenvalue and eigenvector}) $\norm{\bs{\nu}_* - \bs{\nu}_{*,\infty}}_\infty < \kappa$, and $\abs{\zeta_* - \zeta_{*,\infty}} < \kappa$, where $\zeta_*$ and $\bs{\nu}_*$ are the largest eigenvalue and the corresponding eigenvector of $\bs{J}_{\Ffuncbold(\cdot)}(\bs{\mu}_{*})$ respectively;\label{def:mathcalE_ii}
	\item({\it vanishing of active half-edges near the fixed point}) the boundary of the ball of radius $r(\kappa)$ in infinity norm centered at $\bs{\mu}_{*}$ hits the trajectory of $\{\bs{\mu}(t) \text{ for }t>0 \}$ at a point $t_{r(\kappa)}$ for which
	\begin{align*}
		& \bs{0} < \bs{\mu}(t_{r(\kappa)}) - \bs{\mu}_{*} < \bs{\rho}_{*}(\kappa),\text{ and}\\
		&a_{1}(t_{r(\kappa)}) + a_{2}(t_{r(\kappa)}) + a_{m}^{(1)}(t_{r(\kappa)}) + a_{m}^{(2)}(t_{r(\kappa)}) \leq \kappa^2/2,
	\end{align*}
		where $\{\bs{\mu}(t) \text{ for }t>0 \}$ denote the solution of the ODEs \eqref{eq:alterode}, $a_{j}(\cdot)$ and $a_{m}^{(j)}(\cdot)$ are given by Remark \ref{rem:atoFcon}, and $\bs{\rho}_{*}(\kappa) \coloneqq \kappa\bs{\nu}_*$;\label{def:mathcalE_iii}
	\item({\it tractability of the process using ODEs}) $\{\bs{\mu}(t) \text{ for }t\geq0\} \setminus \mathcal{B}\left(\bs{\mu}_{*},r(\kappa)/2\right)  \subset \mathcal{D}_{\varepsilon(\kappa),n}$;\label{def:mathcalE_iv}
	\item({\it approximation at the point of augmentation}) the following inequalities hold:\label{def:mathcalE_v}
	\begin{align*}
		\norm{\bs{\mu}_* - \bs{\mu}_{\kappa}}_\infty < \kappa_1/2\qquad\text{and}\qquad\norm{\bs{\nu}_* - \bs{\nu}_{\kappa}}_\infty < \kappa_1/2 \qquad\text{and}\qquad\abs{\zeta_* - \zeta_{\kappa}} < \kappa_1/2,
	\end{align*}
	where $\bs{\mu}_{\kappa}\coloneqq\bs{\mu}(t_{r(\kappa)})$, and $\bs{\nu}_{\kappa}$ and $\zeta_{\kappa}$ are the largest eigenvalue and the corresponding eigenvector of $\bs{J}_{\Ffuncbold(\cdot)}(\bs{\mu}_\kappa)$ respectively.
\end{enumerate}
Note that $\exists \kappa_2 < \kappa_1/2$ small enough so that Condition \ref{def:mathcalE_v} hold for all $\kappa < \kappa_2$ as long as Conditions \ref{def:mathcalE_i}--\ref{def:mathcalE_iv} hold. The intuition behind the above conditions is as follows:
\begin{enumerate}[label=--]
	\item Condition \ref{def:mathcalE_i} ensures that the closest fixed point of the ODEs \eqref{eq:alterode} (associated with the truncated process) is in a small neighborhood of $\bs{\mu}_{*,\infty}$, $\bs{J}_{\Ffuncbold(\cdot)}(\bs{\mu}_{*})$ is primitive, and the constants that appear in Theorem \ref{thm:odesfinite} are uniformly bounded for all realizations of the initial condition that belong to $\mathcal{E}(n,\kappa)$.
	\item Conditions \ref{def:mathcalE_i} and \ref{def:mathcalE_ii} guarantee that the Perron-Frobenius eigenvalue of $\bs{J}_{\Ffuncbold(\cdot)}(\bs{\mu}_{*})$ is positive real and smaller than $1$.
	\item Conditions \ref{def:mathcalE_i}, \ref{def:mathcalE_ii} and \ref{def:mathcalE_iii} guarantee that the proportion of active half-edges at the point of augmentation is small, and that the added augmented half-edges cannot initiate a large cascade, both with high probabilities.
	\item Conditions \ref{def:mathcalE_iii} and \ref{def:mathcalE_iv} makes it possible to use the ODEs to approximate the scaled-version of the truncated process up to points sufficiently close to $\bs{\mu}_*$.
	\item Conditions \ref{def:mathcalE_i} and \ref{def:mathcalE_v} ensures that, inequalities \eqref{eq:propineq_incomm} and \eqref{eq:propineq_outcomm} hold, after replacing $\zeta_{*,\infty}$, $\bs{\chi}_{*,\infty}(\kappa,\pm\err)$ and $\bs{\mu}_{*,\infty}(\kappa,\pm\err)$ with $\zeta_{\kappa}$, $\bs{\chi}_{\kappa}$ and $\bs{\mu}_\kappa$ respectively, where
	\begin{align*}
		\bs{\chi}_{\kappa}\coloneqq\left({\rho}_{\kappa}^{(1\gets 1)} {\mu}^{(1\gets 1)}_\kappa,
		{\rho}_{\kappa}^{(1\gets 2)}{\mu}^{(2\gets 1)}_\kappa,
		{\rho}_{\kappa}^{(2\gets 1)}{\mu}^{(1\gets 2)}_\kappa,
		{\rho}_{\kappa}^{(2\gets 2)}{\mu}^{(2\gets 2)}_\kappa\right),
	\end{align*}
	and $\bs{\rho}_\kappa\coloneqq \kappa \bs{\nu}_\kappa$. These inequalities are crucial to ensure that the truncated process stops near $\bs{\mu}_\kappa$, with high probability.
\end{enumerate}
The above symbols and their definitions are summarized in Table \ref{table:constants_random}.

\begin{table}[!htbp]
	\centering
	\begin{tabular}{| l | l |}
		\hline
		Symbol & Definition \Tstrut \Bstrut \\ \hline
		$\mathcal{E}(n,\kappa)$& set of initial conditions that satisfy Conditions \ref{def:mathcalE_i}-\ref{def:mathcalE_v};\Tstrut \Bstrut \\ \hline
		$\bs{J}_{\Ffuncbold(\cdot)}(\bs{\mu})$ & Jacobian matrix of $\Ffuncbold(\cdot)$ at $\bs{\mu}$;\Tstrut \Bstrut \\ \hline
		$\zeta_{*}$ & largest eigenvalue of $\bs{J}_{\Ffuncbold(\cdot)}(\bs{\mu}_{*})$;\Tstrut \Bstrut \\ \hline
		$\bs{\nu}_{*}$ & eigenvector of $\bs{J}_{\Ffuncbold(\cdot)}(\bs{\mu}_{*})$ corresponding to the eigenvalue $\zeta_{*}$;\Tstrut \Bstrut \\ \hline
		$\bs{\rho}_{*}(\kappa)$ & $\kappa \bs{\nu}_{*}$;\Tstrut \Bstrut \\ \hline
		$\bs{\mu}_\kappa$ & $\bs{\mu}(t_{r(\kappa)})$ which is a point that belongs to the set $\{\bs{\mu}(t) \text{ for }t>0 \}\cap \partial{\mathcal{B}(\bs{\mu}_{*},r(\kappa))}$;\Tstrut\Bstrut\\ \hline
		$\zeta_{\kappa}$ & largest eigenvalue of $\bs{J}_{\Ffuncbold(\cdot)}(\bs{\mu}_\kappa)$;\Tstrut \Bstrut \\ \hline
		$\bs{\nu}_{\kappa}$ & eigenvector of $\bs{J}_{\Ffuncbold(\cdot)}(\bs{\mu}_\kappa)$ corresponding to the eigenvalue $\zeta_\kappa$;\Tstrut \Bstrut \\ \hline
		$\bs{\rho}_\kappa$ & $\kappa \bs{\nu}_\kappa$;\Tstrut \Bstrut \\ \hline
		$\bs{\chi}_{\kappa}$ & $\left({\rho}_{\kappa}^{(1\gets 1)} {\mu}^{(1\gets 1)}_{\kappa},
		{\rho}_{\kappa}^{(1\gets 2)} {\mu}^{(2\gets 1)}_{\kappa},
		{\rho}_{\kappa}^{(2\gets 1)} {\mu}^{(1\gets 2)}_{\kappa},
		{\rho}_{\kappa}^{(2\gets 2)} {\mu}^{(2\gets 2)}_{\kappa}\right)$;\Tstrut \Bstrut \\ \hline
		$\kappa_2$	& largest constant in $(0,\kappa_1/2]$ for which for any $\kappa \leq \kappa_2$, given Conditions \ref{def:mathcalE_i}-\ref{def:mathcalE_iv},\Tstrut\Bstrut\\
		& Condition \ref{def:mathcalE_v} holds;\Bstrut\\\hline
	\end{tabular}
	\caption{Symbols and their definitions for finite values of $n$. Most of symbols defined in this table depend on the realization of the initial condition.}\label{table:constants_random}
\end{table}

It is easy to see that $\{\mathcal{E}(n,\kappa)\}_{n\in\mathbb{N}}$ holds with high probability. This is a direct consequence of Lemma~\ref{lem:initcond_conv} and the continuous mapping theorem. For the rest of this subsection, we focus on a realization of an initial condition which belongs to $\mathcal{E}(n,\kappa)$.
\begin{lemma}\label{lem:highprobevent}
	For any fixed $\kappa>0$, $\prob(\mathcal{E}(n,\kappa))\to 1$ as $n\to\infty$.
\end{lemma}

Let $t_{\kappa}$ denote the time at which the trajectory of the ODEs \eqref{eq:sol_mujj}-\eqref{eq:sol_muj-j} reaches $\bs{\mu}_{\kappa}$. Note that $t_{\kappa}\neq t_{r(\kappa)}$, as the ODEs given by \eqref{eq:sol_mujj}-\eqref{eq:sol_muj-j} and the ODEs given by \eqref{eq:alterode} have the same trajectory but different speed-scales. By \eqref{eq:sol_t}, we have
\begin{align*}
	t_\kappa =\lambda_m + \frac{\lambda_1}{2} + \frac{\lambda_2}{2} - \frac{\lambda_1}{2} \left(\mu_{\kappa}^{(1\gets 1)}\right)^2 - \frac{\lambda_2}{2} \left(\mu_{\kappa}^{(2\gets 2)}\right)^2 -\lambda_m\mu_{\kappa}^{(1\gets 2)}\mu_{\kappa}^{(2\gets 1)}.
\end{align*}
By Theorem \ref{thm:odesfinite}, Corollary \ref{cor:odesfinite} and Lemma \ref{lem:diffeq_sol}, the total number of active half-edges at time $\floor{t_{\kappa}n}$ in the truncated process concentrates around $\mathrm{Const} \times n$ for some $\mathrm{Const} \leq \kappa^2$, with high probability.

\subsubsection{The Augmented Process}\label{subseubsec:augment}
At the beginning of time $\floor{t_{\kappa}n} + 1$, we augment $X^n_\delta$ by adding one active vertex to each community $j\in\{1,2\}$ with $2\floor{n \chi_\kappa^{(j\gets j)} \lambda_j(n) /2}$ half-edges for connections in community $j$ and $\floor{n\chi^{(1\gets 2)} \lambda_m(n)} + \floor{n \chi^{(2\gets 1)} \lambda_m(n)}$ half-edges for connections in community $-j$. Recall that $\chi_\kappa^{(j\gets j)} = \rho_\kappa^{(j\gets j)}\mu_{\kappa}^{(j\gets j)}$ and $\chi_\kappa^{(j\gets -j)} = \rho_\kappa^{(j\gets -j)}\mu_{\kappa}^{(-j\gets j)}$, for $j\in\{1,2\}$. We use $\widetilde{v}_j$ for $j\in\{1,2\}$ to denote these two active vertices. We refer to the newly added active half-edges as the augmented half-edges and all the other half-edges as the regular ones (if necessary, to avoid confusion). Notice that the process may have already run out of active half-edge before time $\floor{t_{\kappa}n}$; still, we can augment the process despite the process being halted for a while.

After augmentation, we change the process and proceed in two phases: during the first phase, we pair a subset of augmented half-edges with random regular half-edges, and during the second phase, we proceed normally by pairing randomly selected active half-edges with random half-edges in the proper community. Note that augmented half-edges cannot be paired with each other during the first phase. We then use the ODEs to approximate the state of the process after the first phase. We refer to this process as the augmented process.

\begin{figure}[!htbp]
\centering
\begin{subfigure}[t]{.3\linewidth}
	\centering
	\includegraphics*[width = 0.75\textwidth]{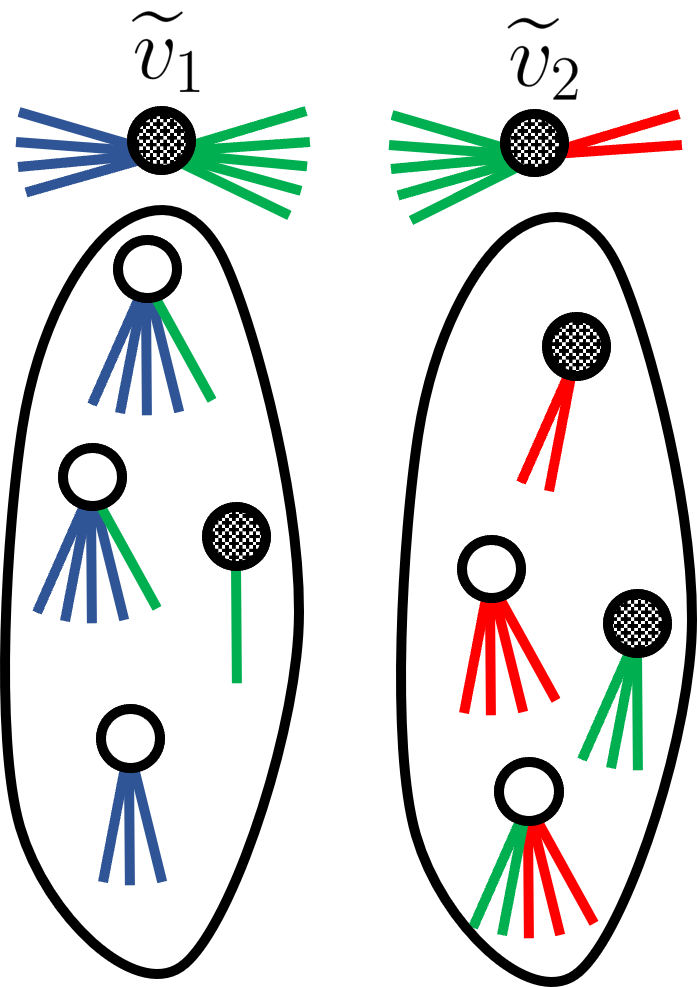}
	\caption*{State of graph at the point of augmentation.}
\end{subfigure}\hfill
\begin{subfigure}[t]{.3\linewidth}
	\centering
	\includegraphics*[width = 0.75\textwidth]{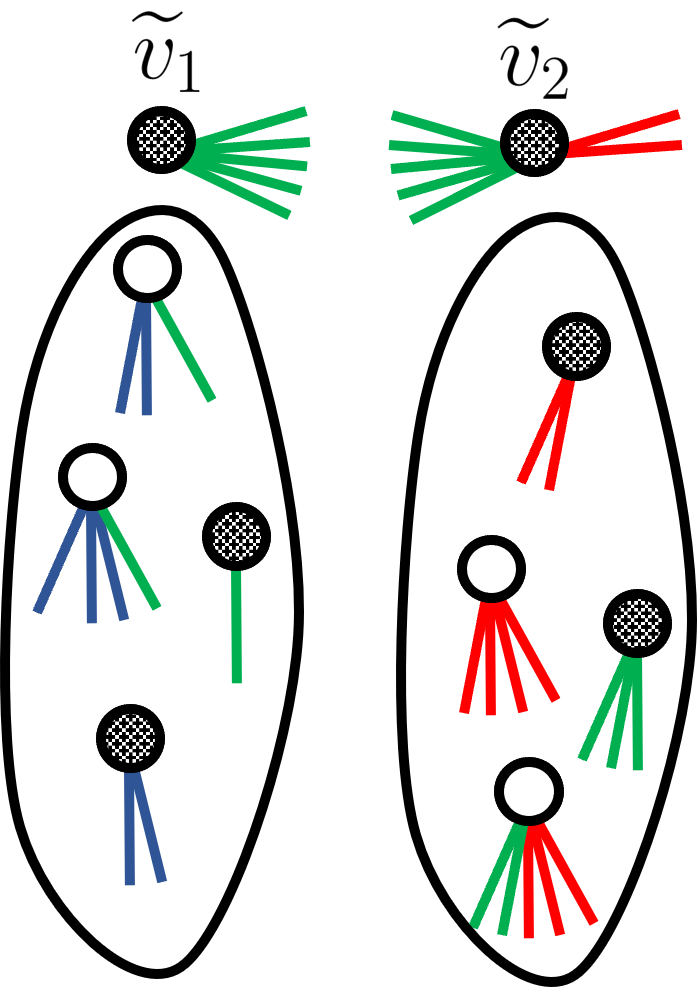}
	\caption*{($i$): $2\floor{n \chi_\kappa^{(1\gets 1)} \lambda_1(n) /2}$ augmented half-edges are paired.}
\end{subfigure}\hfill
\begin{subfigure}[t]{.3\linewidth}
	\centering
	\includegraphics*[width = 0.75\textwidth]{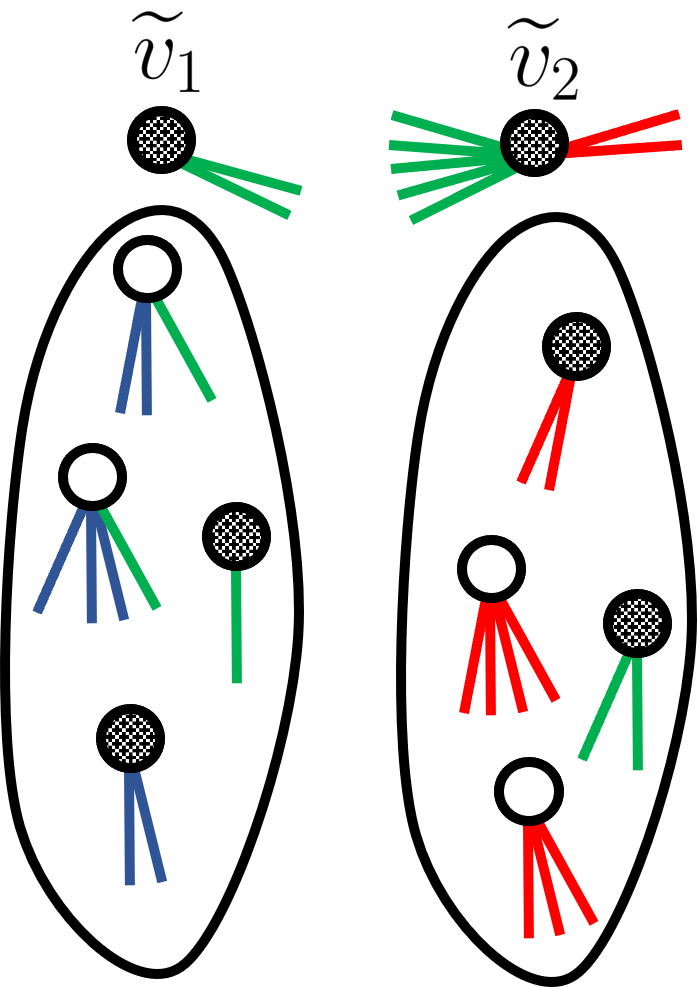}
	\caption*{($ii$): $\floor{n \chi^{(2\gets 1)} \lambda_m(n)}$ augmented half-edges are paired.}
\end{subfigure}\\
\medskip
\begin{subfigure}[t]{.15\linewidth}
\end{subfigure}
\hfill
\begin{subfigure}[t]{.3\linewidth}
	\centering
	\includegraphics*[width = 0.75\textwidth]{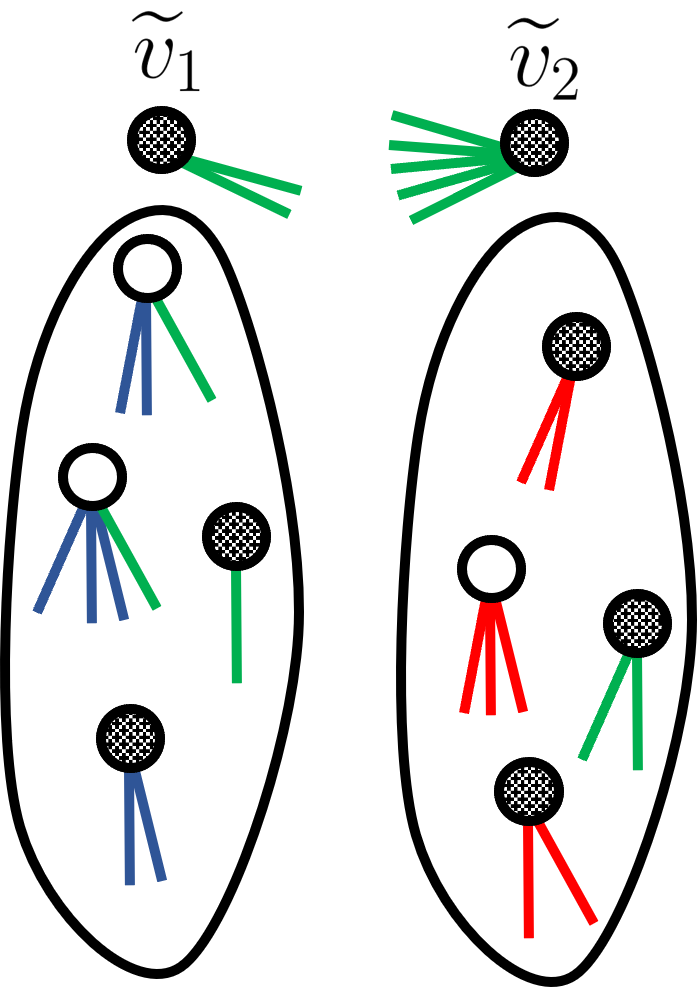}
	\caption*{($iii$): $2\floor{n \chi_\kappa^{(2\gets 2)} \lambda_2(n) /2}$ augmented half-edges are paired.}
\end{subfigure}\hfill
\begin{subfigure}[t]{.3\linewidth}
	\centering
	\includegraphics*[width = 0.75\textwidth]{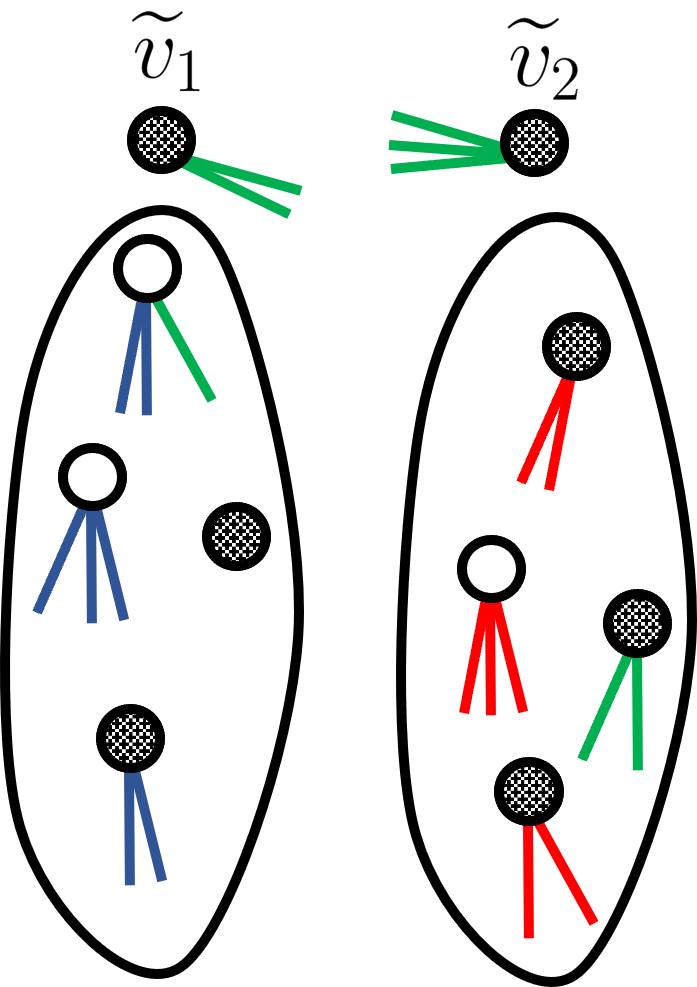}
	\caption*{($iv$): $\floor{n \chi^{(1\gets 2)} \lambda_m(n)}$ augmented half-edges are paired.}
\end{subfigure}\hfill
\begin{subfigure}[t]{.15\linewidth}
\end{subfigure}
\hfill
\caption{A schematic of the timeline of the augmented process during Phase \ref{def:aug_phase1}. Green half-edges are for connections between the communities, blue half-edges are for connections in community $1$, and red half-edges are for connections in community $2$. Active vertices are denoted by dotted circles.}
\label{fig:augmentedp}
\end{figure}

In particular, the augmented process, after adding $\widetilde{v}_j$ for $j\in\{1,2\}$, proceeds as follows:
\begin{enumerate}[leftmargin=!,itemindent=1em,align=left,label=Phase A.$\arabic*$:,ref=A.$\arabic*$]
	\item pairing a subset of augmented half-edges with random half-edge; this phase consists of the following timeline: \label{def:aug_phase1}
	\begin{enumerate}[label = ($\roman*$)]
		\item during the first $2\floor{n \chi_\kappa^{(1\gets 1)} \lambda_1(n) /2}$ times, we pair augmented half-edges in community $1$ with random regular half-edges in community $1$;
		\item during the next $\floor{n \chi^{(2\gets 1)} \lambda_m(n)}$ times, we pair augmented half-edges in community $1$ with random regular half-edges in community $2$;
		\item during the next $2\floor{n \chi_\kappa^{(2\gets 2)} \lambda_2(n) /2}$ times, we pair augmented half-edges in community $2$ with random regular half-edges in community $2$;
		\item during the final $\floor{n \chi^{(1\gets 2)} \lambda_m(n)}$ times, we pair augmented half-edges in community $2$ with random regular half-edges in community $1$;
	\end{enumerate}
	\item normal pairing; this phase proceeds similar to the truncated process, i.e., pairing a randomly selected active half-edge with a random half-edge in the proper community.
\end{enumerate}
Note that at the end of Phase \ref{def:aug_phase1} of the augmented process, $\widetilde{v}_j$ for $j\in\{1,2\}$ has $\floor{n \chi^{(j\gets -j)} \lambda_m(n)}$ augmented half-edges remaining that should be paired with half-edges in community $-j$. Let us denote the augmented process by $\widetilde{X}^n_\delta$. Note that for all $k\leq \floor{t_\kappa n}$, $\widetilde{X}^n_\delta(k) = X^n_\delta(k)$. The symbols associated with the augmented process and their definitions are summarized in \ref{table:constants_augmented}. In Figure \ref{fig:augmentedp}, we illustrate the timeline of the augmented process.

\begin{table}[!htbp]
	\centering
	\begin{tabular}{| l | l |}
		\hline
		Symbol & Definition \Tstrut \Bstrut \\ \hline
		$t_\kappa$& $ \lambda_m(n)\left(1-\mu_{\kappa}^{(1\gets 2)}\mu_{\kappa}^{(2\gets 1)}\right) + \frac{\lambda_1(n)}{2}\left(1-\left({\mu_{\kappa}^{(1\gets 1)}}\right)^2 \right)+ \frac{\lambda_2(n)}{2}\left(1-{\left(\mu_{\kappa}^{(2\gets 2)}\right)}^2\right)$;\DTstrut\DBstrut\\\hline
		$ t_\kappa^{(1\gets 1)}$, $t_\kappa^{(2\gets 1)}$ & $t_\kappa^{(1\gets 1)}= t_\kappa + \chi^{(1\gets 1)}_\kappa \lambda_1(n);\myquad[2] t_\kappa^{(2\gets 1)} = t_\kappa^{(1\gets 1)} + \chi_\kappa^{(2\gets 1)} \lambda_m(n)$;
		\DTstrut\DBstrut\\ \hline
		$ t_\kappa^{(2\gets 2)}$, $t_\kappa^{(1\gets 2)}$ & $t_\kappa^{(2\gets 2)}= t_\kappa^{(2\gets 1)} + \chi_\kappa^{(2\gets 2)} \lambda_2(n);\myquad[2] t_\kappa^{(1\gets 2)} = t_\kappa^{(2\gets 2)} + \chi_\kappa^{(1\gets 2)} \lambda_m(n)$;
		\DTstrut\DBstrut\\ \hline
	\end{tabular}
	\caption{Symbols that appear in the augmented process and their definitions. Note that these values are random and they depend only on the realization of the initial condition.}\label{table:constants_augmented}
\end{table}

Using the same approach as in Sections \ref{sec:markovproc}-\ref{sec:odeanalysisfinite}, we can approximate the augmented process by a system of differential equations. In particular, we are only interested in the state of the augmented process at the end of the Phase \ref{def:aug_phase1}. The details of the one-step drift, the derivation of the corresponding ODEs, and its solution are presented in Appendix \ref{app:augment}.

\begin{remark}
	Note that $\bs{\mu}_{\kappa}$ and $t_\kappa$ are random variables since they depend on the initial condition~\eqref{eq:ode_ic}. Hence, the system of ODEs that we use to approximate the augmented process is random; however, we can still invoke Wormald's theorem as the only source of randomness is the initial condition. Note that both the ``\textit{Trend hypothesis}'' and the ``\textit{Lipschitz hypothesis}'' of Wormald's Theorem~\cite[Theorem 5.1]{Wormald1999} hold for any realization of the initial condition (with a uniform Lipschitz constant, similar to the argument in Appendix \ref{proof:odesfinite}).
\end{remark}

By the analysis of Appendix \ref{app:augment}, at the end of Phase \ref{def:aug_phase1}, the proportion of active half-edges for connections in community $j\in\{1,2\}$ concentrates at
\begin{align*}
	&\widetilde{a}_j(t_\kappa^{(1\gets 2)}) \leq a_j(t_\kappa) +  \chi^{(j\gets j)} \lambda_j(n) \zeta_\kappa + O(\kappa^2),
\end{align*}
where $a_j(t_\kappa)$ is related to the proportion of the same type of active half-edges before the augmentation. Similarly, the proportion of active half-edges in community $-j$ that should be connected to random half-edges in community $j$ concentrates at
\begin{align*}
	&\widetilde{a}^{(-j)}_m(t_\kappa^{(1\gets 2)}) \leq {a}^{(-j)}_m(t_\kappa) +  \chi^{(-j\gets j)}\lambda_m(n) +\chi^{(j\gets -j)}\lambda_m(n) \zeta_\kappa+ O(\kappa^2).
\end{align*}
Note that by Condition \ref{def:mathcalE_iii}, for $j\in\{1,2\}$, we have
\begin{align*}
	&\widetilde{a}_j(t_\kappa^{(1\gets 2)}) \leq \kappa^2/2 +  \chi^{(j\gets j)} \lambda_j(n) \zeta_\kappa + O(\kappa^2),\\
	&\widetilde{a}^{(-j)}_m(t_\kappa^{(1\gets 2)}) \leq \kappa^2/2 +  \chi^{(-j\gets j)}\lambda_m(n) +\chi^{(j\gets -j)}\lambda_m(n) \zeta_\kappa+ O(\kappa^2).
\end{align*}

\subsubsection{The Twisted Process}\label{subseubsec:twist}
\begin{figure}[!htbp]
	\centering
	\begin{subfigure}[t]{0.3\linewidth}
		\centering
		\includegraphics*[width = 0.65\textwidth]{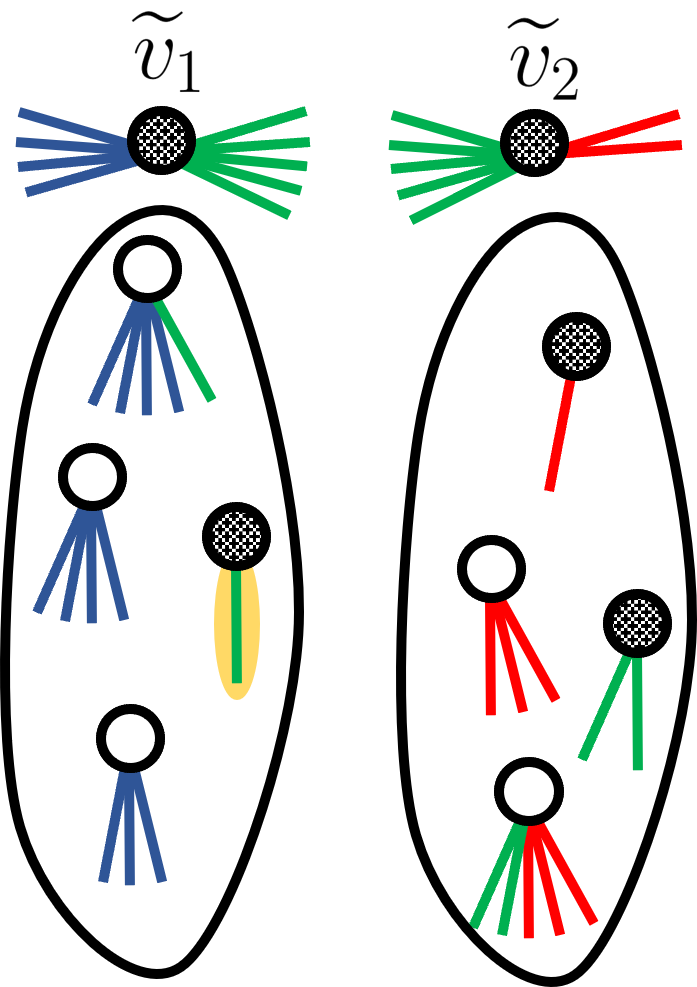}
		\caption*{$1$: Pick a regular active half-edge uniformly at random.}
	\end{subfigure}\hfill
	\begin{subfigure}[t]{0.3\linewidth}
		\centering
		\includegraphics*[width = 0.65\textwidth]{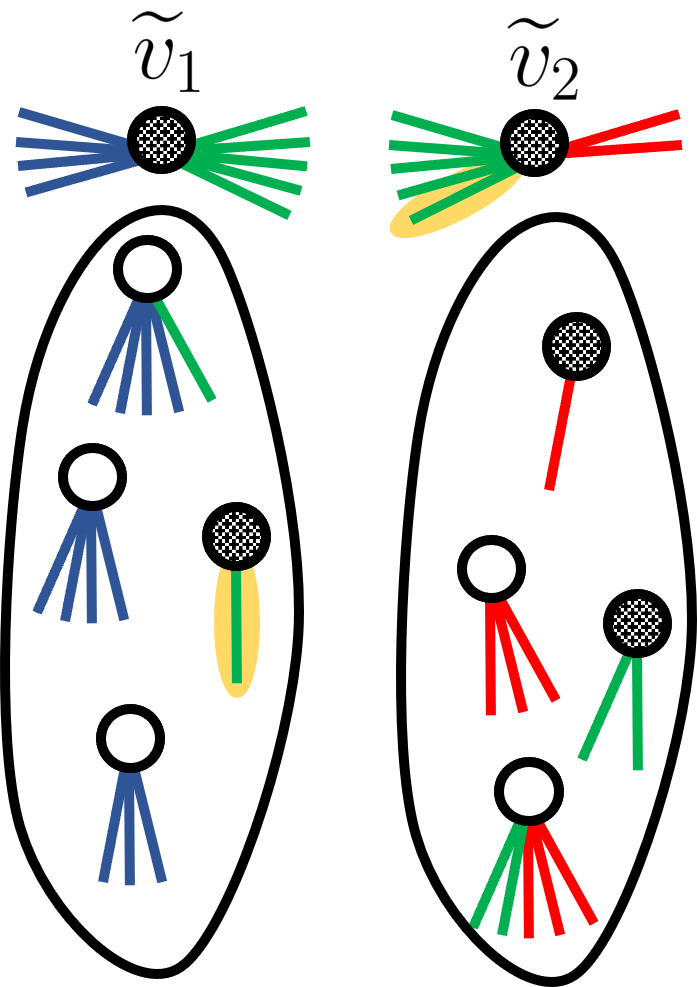}
		\caption*{$2$: Pick another half-edge uniformly at random.}
	\end{subfigure}\hfill
	\begin{subfigure}[t]{0.3\linewidth}
		\centering
		\includegraphics*[width = 0.65\textwidth]{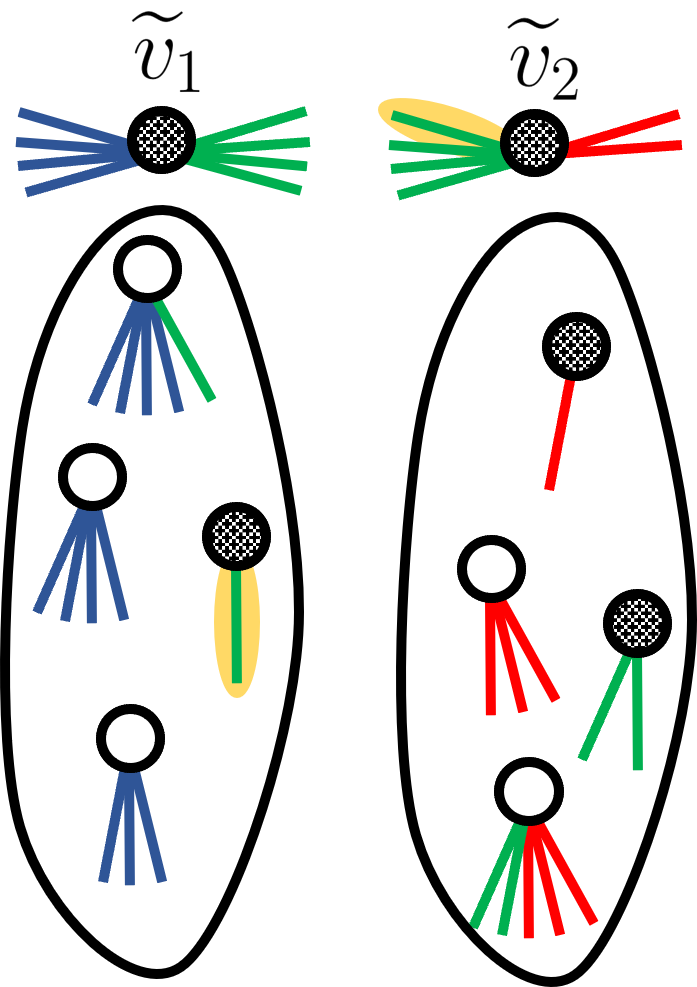}
		\caption*{$3$: Remove two augmented half-edges and pick another random half-edge.}
	\end{subfigure}\\
	\medskip
	\begin{subfigure}[t]{.15\linewidth}
	\end{subfigure}
	\hfill
	\begin{subfigure}[t]{0.3\linewidth}
		\centering
		\includegraphics*[width = 0.65\textwidth]{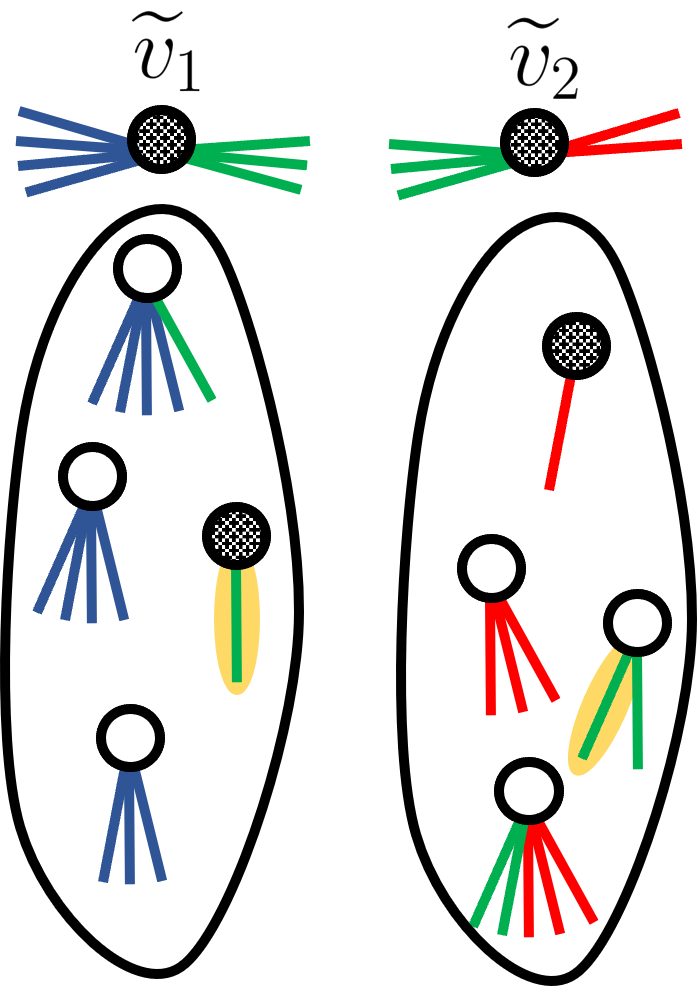}
		\caption*{$4$: Remove two augmented half-edges and pick another random half-edge.}
	\end{subfigure}\hfill
	\begin{subfigure}[t]{0.3\linewidth}
		\centering
		\includegraphics*[width = 0.65\textwidth]{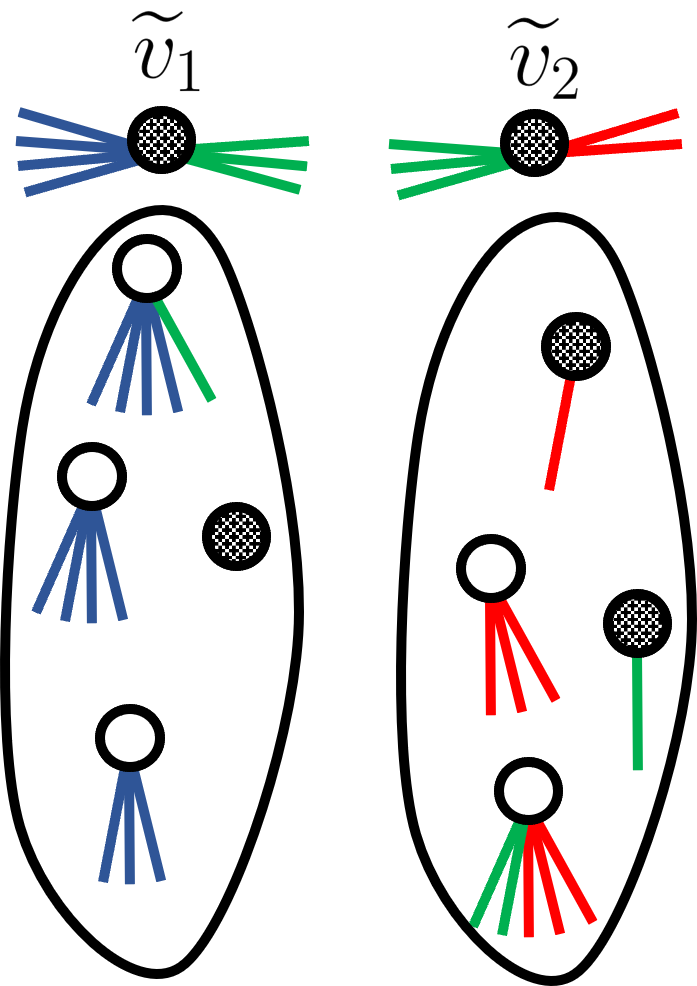}
		\caption*{$5$: Remove the two regular half-edges.}
	\end{subfigure}\hfill
	\begin{subfigure}[t]{.15\linewidth}
	\end{subfigure}
	\hfill
	\caption{A schematic of the timeline of one iteration during Phase \ref{def:twist_phase1} of the twisted process. Green half-edges are for connections between the communities, blue half-edges are for connections in community $1$, and red half-edges are for connections in community $2$. Active vertices are denoted by dotted circles. Selected half-edges during this iteration are highlighted.}\label{fig:twistp}
\end{figure}

To compare the augmented and truncated processes, we define a new process that we call the twisted process. The twisted process and the augmented process have two key differences: $(1)$ augmented half-edges in the twisted process can be paired with each other, and $(2)$ in the twisted process, half-edges are paired in a different order. Note that the sample paths of all these processes (truncated, augmented, and twisted) are the same up to time $\floor{t_\kappa n}$. However, after adding the augmented half-edges, the twisted process proceeds as follows:
\begin{enumerate}[leftmargin=!,itemindent=1em,align=left,label=Phase T.$\arabic*$:,ref=T.$\arabic*$]
	\item pairing regular active half-edges with random regular half-edges; this phase consists of some iterations, and during each iteration, we do the following:\label{def:twist_phase1}
	\begin{enumerate}[label=${(\roman*)}$]
		\item pick a regular active half-edge uniformly at random, and a random half-edge in the proper community;
		\item while the last selected random half-edge is augmented, repeat: $(1)$ remove this half-edge, $(2)$ remove another augmented half-edge from the same community as the initial regular active half-edge, $(3)$ pick another random half-edge from the proper community;
		\item remove the active half-edge, and the last selected regular active half-edge;
	\end{enumerate}
	\item pairing random active half-edges with random regular half-edges; this phase proceeds with pairing all the remaining active half-edges with random half-edges. \label{def:twist_phase2}
\end{enumerate}
Note that Phase \ref{def:twist_phase1} of the twisted process ends when we run out of regular active half-edges. In Figure \ref{fig:twistp}, we illustrate the timeline of an iteration of the Phase \ref{def:twist_phase1} of the twisted process.

%Also, note that the order in which active regular half-edges are removed during Phase \ref{def:twist_phase1} does not matter as we are only concerned with the state of the twisted process at the end of Phase \ref{def:twist_phase1}.

\subsubsection{A Coupling Argument}\label{subseubsec:coupling}
During each iteration of the Phase \ref{def:twist_phase1} of the twisted process, we remove a few half-edges; however, exactly two regular half-edges are removed: the initial active half-edge and the last selected regular half-edge. Hence, there is a natural coupling between the truncated process and the first phase of the twisted process. Note that at the end of Phase \ref{def:twist_phase1}, the coupled truncated process hits its natural stopping time.

Note that augmented half-edges in the twisted process can be paired with each other. Hence, one may expect the twisted process to be more wasteful than the augmented process since more active half-edges can be paired with each other in the twisted process. Next, following this intuition, we propose a coupling between the augmented process and the twisted process to show that Phase \ref{def:twist_phase1} of the twisted process ends after at most $2\floor{n\chi^{(1\gets 1)}\lambda_1(n)/2}+2\floor{n\chi^{(2\gets 2)} \lambda_2(n)/2}+\floor{n\chi^{(1\gets 2)} \lambda_m(n)} + \floor{n\chi^{(2\gets 1)} \lambda_m(n)} = n\,O(\kappa)$ many time (not iteration), with high probability. Note that during this time, at most $n\,O(\kappa)$ vertices can be activated.
\begin{remark}
	During the rest of this section, we use ``$\,\,\widetilde{~}\,\,$'' and ``$\,\,\widehat{~}\,\,$'' to distinguish between quantities that are associated with the augmented and twisted processes, respectively. The only exception is the newly added vertices and label of their half-edges,  which are the same for both processes.
\end{remark}

Both the augmented process and the truncated process follow the same sample path up to time $\floor{t_\kappa n}$ before adding the high-degree active vertices $\widetilde{v}_1$ and $\widetilde{v}_2$. Suppose that vertices of community $j\in\{1,2\}$ are labeled as $\{v_1^{(j)},v_2^{(j)},\cdots,v_{n_j}^{(j)}\}$. Let us label half-edges after adding $\widetilde{v}_1$ and $\widetilde{v}_2$ as follows:
\begin{enumerate}[label = --]
	\item half-edges of the vertex $v_l^{(j)}$ that can be paired with half-edges in community $j$ are labeled as $\{(v_l^{(j)},i)\}_{i=1}^{d_j - u_j}$, where $d_j - u_j$ is the number of remaining such half-edges;
	\item half-edges of the vertex $v_l^{(j)}$ that can be paired with half-edges in community $-j$ are labeled as $\{(v_l^{(j)},-i)\}_{i=1}^{d_{-j}-u_{-j}}$, where $d_{-j} - u_{-j}$ is the number of remaining such half-edges;
	\item half-edges of the augmented vertex $\widetilde{v}_j$ that can be paired with half-edges in community $j$ are labeled as $\{(\widetilde{v}_j,i)\}_{i=1}^{2\floor{n\chi_\kappa^{(j\gets j)} \lambda_j(n)/2}}$;
	\item half-edges of the augmented vertex $\widetilde{v}_j$ that can be paired with half-edges in community $-j$ are labeled as $\{(\widetilde{v}_j,-i)\}_{i=1}^{\floor{n \chi^{(1\gets 2)} \lambda_m(n)} + \floor{n \chi^{(2\gets 1)} \lambda_m(n)}}$.
\end{enumerate}

Consider a realization of the augmented process from time $\floor{t_\kappa n}$ up to time $\floor{{t}^{(1\gets 2)}_\kappa n}$. This realization results in the following sequences:
\begin{enumerate}
	\item The sequence of regular half-edges that are removed within community $j\in\{1,2\}$:
	\begin{align*}
		\widetilde{S}_{(j\gets j)} = \left(\widetilde{e}_{(j\gets j)}(1),\widetilde{e}_{(j\gets j)}(2),\widetilde{e}_{(j\gets j)}(3),\cdots,\widetilde{e}_{(j\gets j)}(2\floor{n\chi_\kappa^{(j\gets j)} \lambda_j(n)/2}) \right),
	\end{align*}
	where $\widetilde{e}_{(j\gets j)}(k)$ is the label of the $k$th regular half-edge in community $j$ that has been paired with an augmented half-edge in community $j$. For $k\leq 2\floor{n\chi_\kappa^{(j\gets j)} \lambda_j(n)/2}$, define $\widetilde{S}_{(j\gets j)}(k) \coloneqq \left(\widetilde{e}_{(j\gets j)}(1),\widetilde{e}_{(j\gets j)}(2),\widetilde{e}_{(j\gets j)}(3),\cdots,\widetilde{e}_{(j\gets j)}(k) \right)$.

	\item The sequence of regular half-edges that are removed between the communities, where the augmented half-edge belongs to community $j\in\{1,2\}$:
	\begin{align*}
		\widetilde{S}_{(j\gets -j)} = \left(\widetilde{e}_{(j\gets   -j)}(1),\widetilde{e}_{(j\gets   -j)}(2),\cdots,\widetilde{e}_{(j\gets   -j)}(\floor{n \chi^{(-j\gets j)} \lambda_m(n)}) \right),
	\end{align*}
	where $\widetilde{e}_{(j\gets   -j)}(k)$ is the label of the $k$th regular half-edge in community $-j$ that has been paired with an augmented half-edge in community $j$.  For $k\leq \floor{n \chi^{(-j\gets j)} \lambda_m(n)}$, define $\widetilde{S}_{(j\gets -j)}(k) \coloneqq \left(\widetilde{e}_{(j\gets   -j)}(1),\widetilde{e}_{(j\gets   -j)}(2),\cdots,\widetilde{e}_{(j\gets   -j)}(k) \right)$.
\end{enumerate}
Let ${E}_{(j\gets j)}$ denote the set of {\em regular half-edges} within community $j\in\{1,2\}$ at the beginning of Phase \ref{def:aug_phase1}. Similarly, let ${E}_{(j\gets   -j)}$ denote the set of {\em regular half-edges} in community $-j\in\{1,2\}$ that can be paired with half-edges in the other community, at the beginning of Phase \ref{def:aug_phase1}. Note that ${E}_{(j\gets j)}$ and ${E}_{(j\gets   -j)}$ are the same in both twisted and augmented processes.

\begin{figure}
	\centering
	\includegraphics[width=0.4\textwidth]{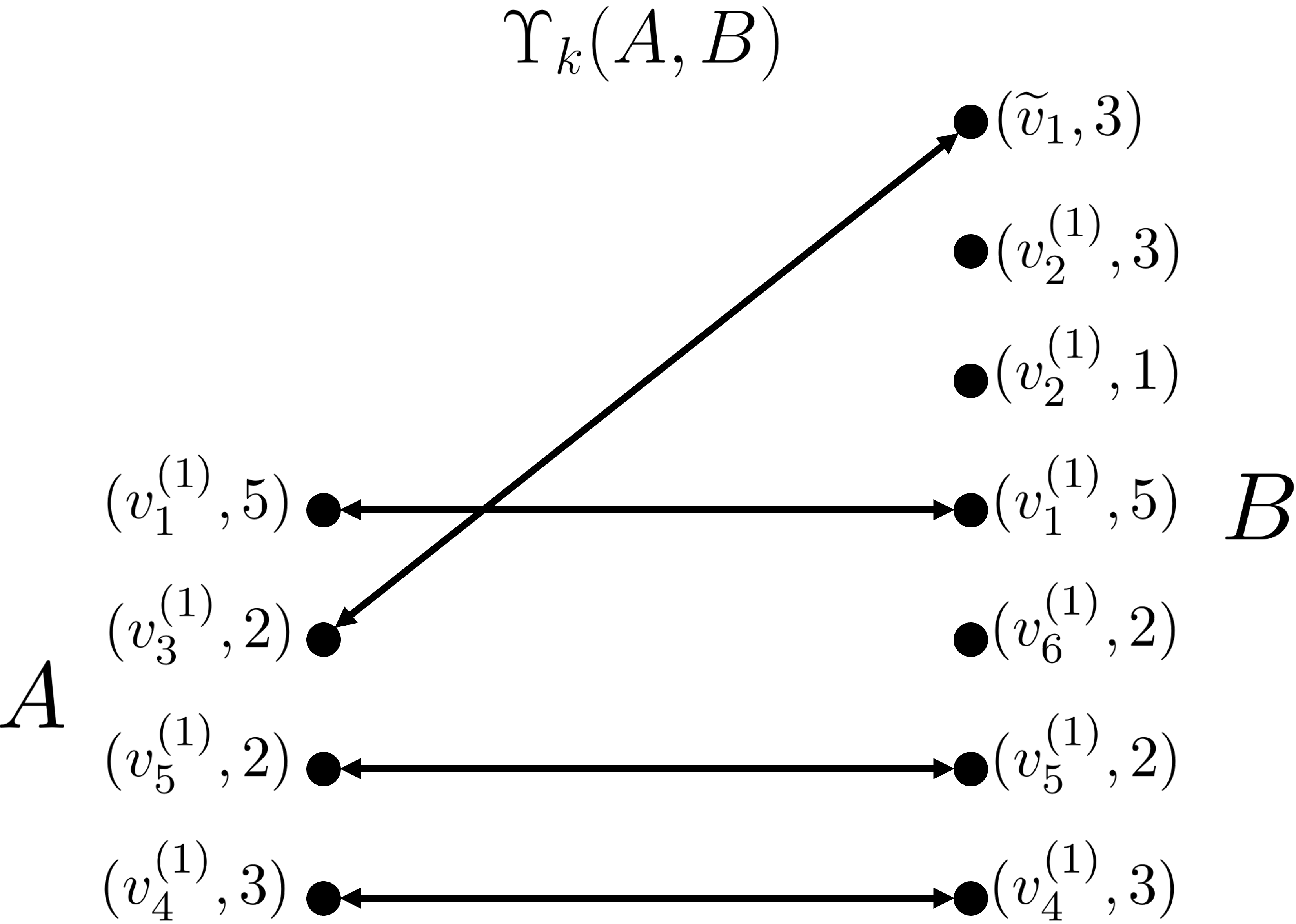}
	\caption{A realization of the relabeling function $\Upsilon_k(A,B)$ at time $k$ that maps $A$ into a subset of $B$. Note that their intersection $A\cap B = \left\{({v}_1^{(1)},5),({v}_5^{(1)},2),({v}_4^{(1)},3)\right\}$ remains intact.} \label{fig:labelreassign}
\end{figure}

Given the above realization of the augmented process up to time $\floor{{t}^{(1\gets 2)}_\kappa n}$, we realize a sequence of pairings for the twisted process. As we mentioned before, we are only concerned with the natural stopping time of Phase \ref{def:twist_phase1}, and hence, the order in which half-edges in different communities are paired does not have any impact. In the following coupling, at each time $k$, we use a random one-to-one relabeling function $\Upsilon_k$, that takes two sets of labels $A$ and $B$ such that $A\subset B$, and maps $A$ to a subset of $B$ uniformly at random while keeping the labels in $A\cap B$ intact. See Figure \ref{fig:labelreassign} for an example. Abusing notation, we use $\Upsilon_k(A,B)(e)$ to denote the label in $B$ that $e\in A$ has been mapped into, and $\Upsilon_k(A,B)$ as the subset of $B$ to which $A$ has been mapped into. 

Given the above sequences, we couple the twisted process and the augmented process. The coupling is based on exploring the sequences $\widetilde{S}_{(1\gets 1)}$, $\widetilde{S}_{(2\gets 2)}$, $\widetilde{S}_{(2\gets 1)}$ and $\widetilde{S}_{(1\gets 2)}$. Hence, we continue the coupling as long as none of these sequences are fully explored and the twisted process is still in Phase \ref{def:twist_phase1}. Specifically, we pair half-edges in the twisted process using the above sequences as long as the following conditions are met:
\begin{enumerate}[leftmargin=!,itemindent=1em,align=left,label = Condition S.$\arabic*$:, ref=S.$\arabic*$]
	\item({\it $\widetilde{S}_{(1\gets 1)}$ and $\widetilde{S}_{(2\gets 2)}$ have members to be explored}) total number of half-edges, after adding the new vertices, that are removed from within community $j\in\{1,2\}$ is smaller than $4\floor{n\chi_\kappa^{(j\gets j)} \lambda_j(n)/2}$, i.e., $\widetilde{S}_{(j\gets j)}$ has not been fully explored yet; \label{def:couplingcond_1}
	\item({\it $\widetilde{S}_{(1\gets 2)}$ and $\widetilde{S}_{(2\gets 1)}$ have members to be explored}) total number of times, after adding the new vertices, that we picked an active half-edge (regular or augmented) from community $j\in\{1,2\}$ and then paired it with a random half-edge in community $-j$ is smaller than $\floor{n \chi^{(-j\gets j)} \lambda_m(n)}$, i.e., $\widetilde{S}_{(j\gets -j)}$ has not been fully explored yet;\label{def:couplingcond_2}
	\item({\it coupled truncated process has not hit its stopping time}) twisted process is in Phase \ref{def:twist_phase1}. \label{def:couplingcond_3}
\end{enumerate}
Whenever one of the above conditions are violated, we decouple the two processes and proceed with the twisted process, independent of the augmented process. We will show that with high probability, Condition \ref{def:couplingcond_3} will be violated first, which in turn implies that the truncated process stops after at most $O(\kappa) n$ many time steps.

\begin{figure}[!htbp]
	\centering
	\begin{subfigure}[t]{.5\linewidth}
		\centering
		\includegraphics*[width = 0.75\textwidth]{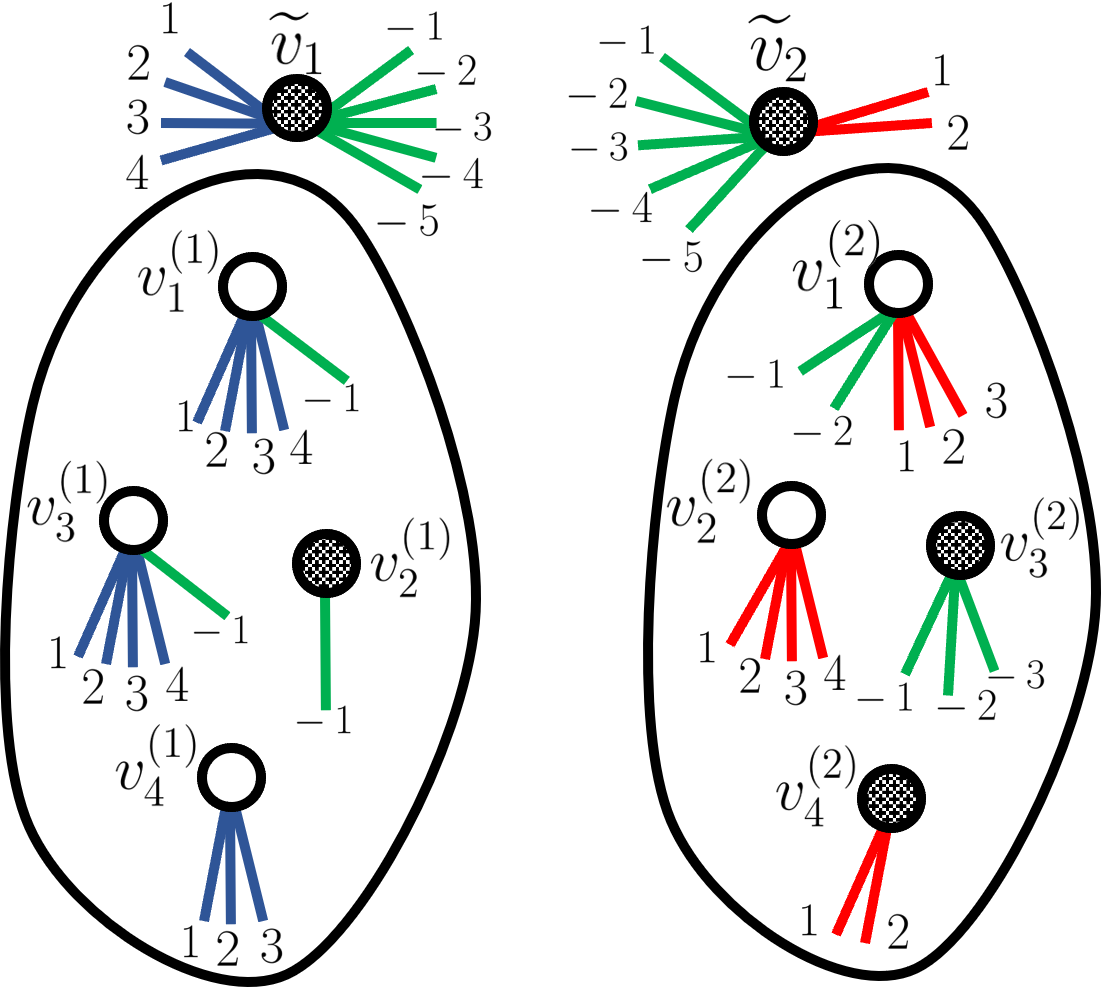}
		\caption*{Labels of vertices and edges.}
	\end{subfigure}\hfill
	\begin{subfigure}[t]{.5\linewidth}
		\centering
		\includegraphics*[width = 0.75\textwidth]{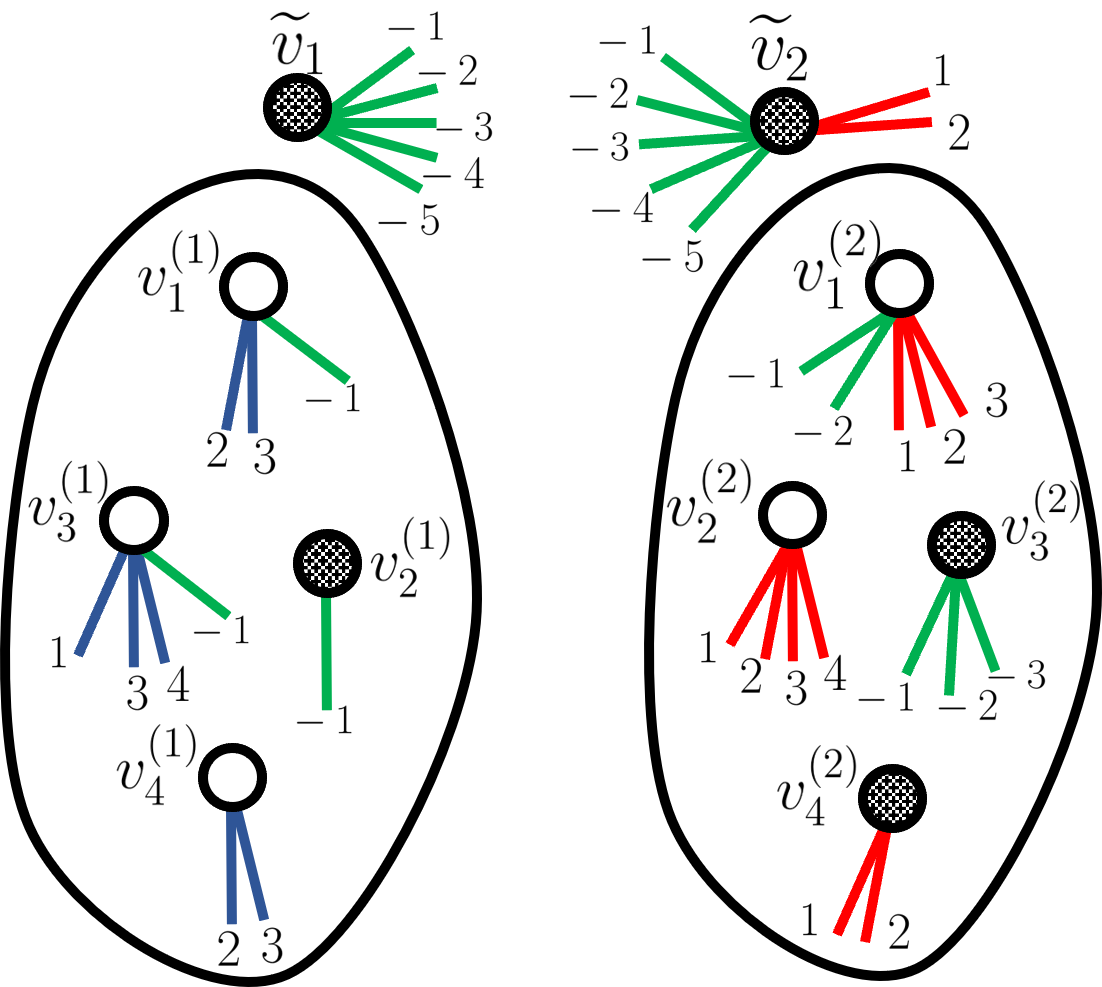}
		\caption*{($i$): $\widetilde{S}_{(1\gets 1)}\!=\!((v_1^{(1)}\!,1),\!(v_3^{(1)}\!,2),\!(v_4^{(1)}\!,1),\! (v_1^{(1)}\!,4))$.}
	\end{subfigure}\\
	\medskip
	\begin{subfigure}[t]{.5\linewidth}
		\centering
		\includegraphics*[width = 0.75\textwidth]{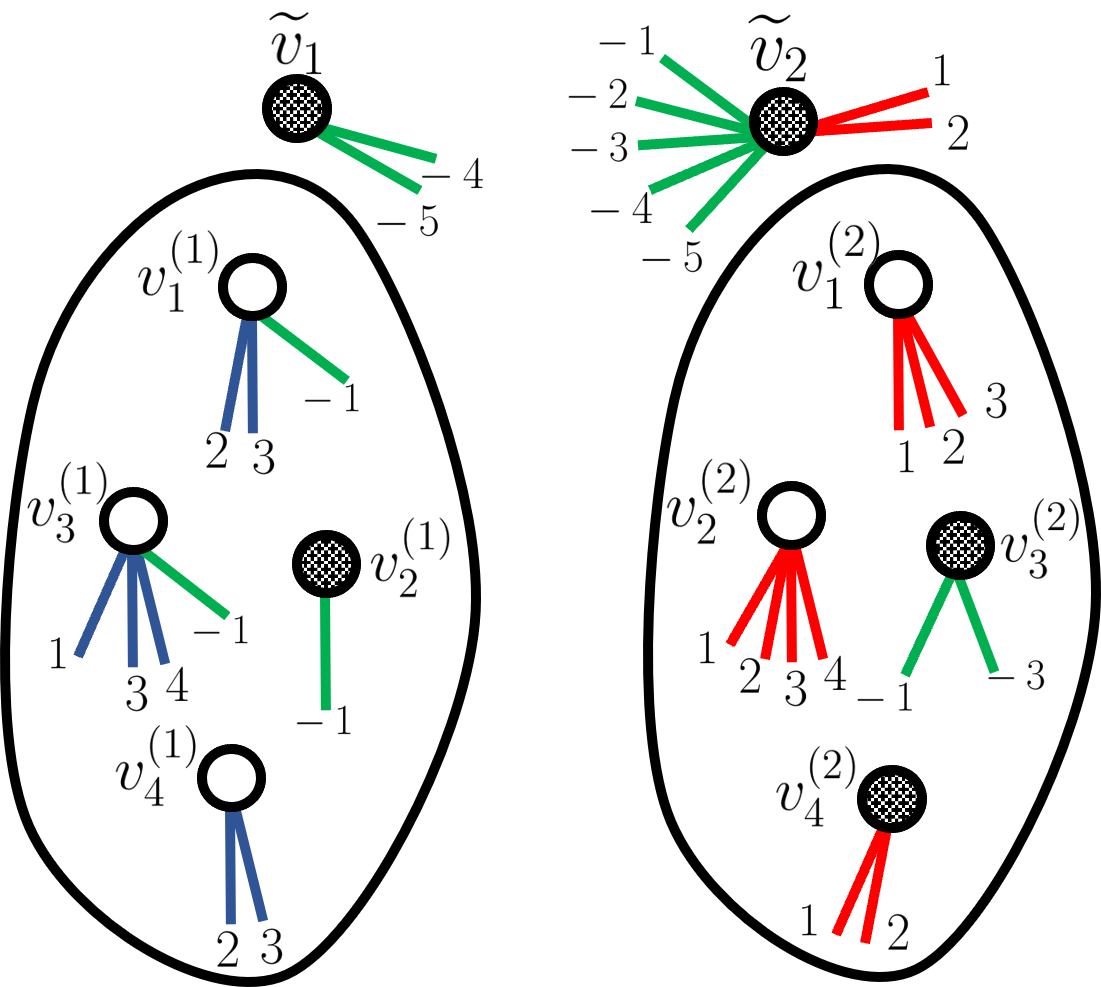}
		\caption*{($ii$): \small$\widetilde{S}_{(1\gets 2)} = ((v_1^{(2)},-2),(v_3^{(2)},-2),(v_1^{(2)},-1))$.}
	\end{subfigure}\hfill
	\begin{subfigure}[t]{.5\linewidth}
		\centering
		\includegraphics*[width = 0.75\textwidth]{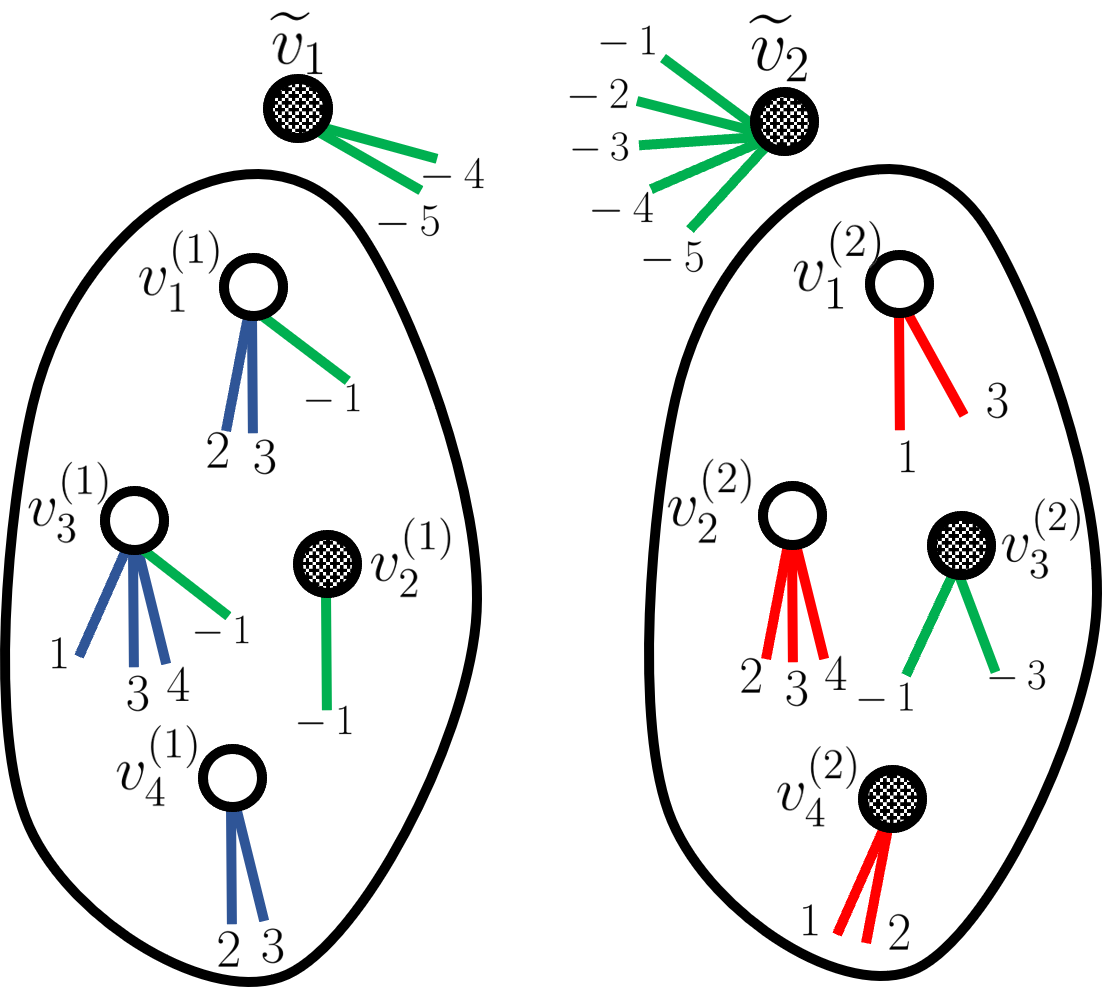}
		\caption*{($iii$): $\widetilde{S}_{(2\gets 2)} = ((v_2^{(2)},1), (v_1^{(2)},2))$.}
	\end{subfigure}\\
	\medskip
	\begin{subfigure}[t]{.5\linewidth}
		\centering
		\includegraphics*[width = 0.75\textwidth]{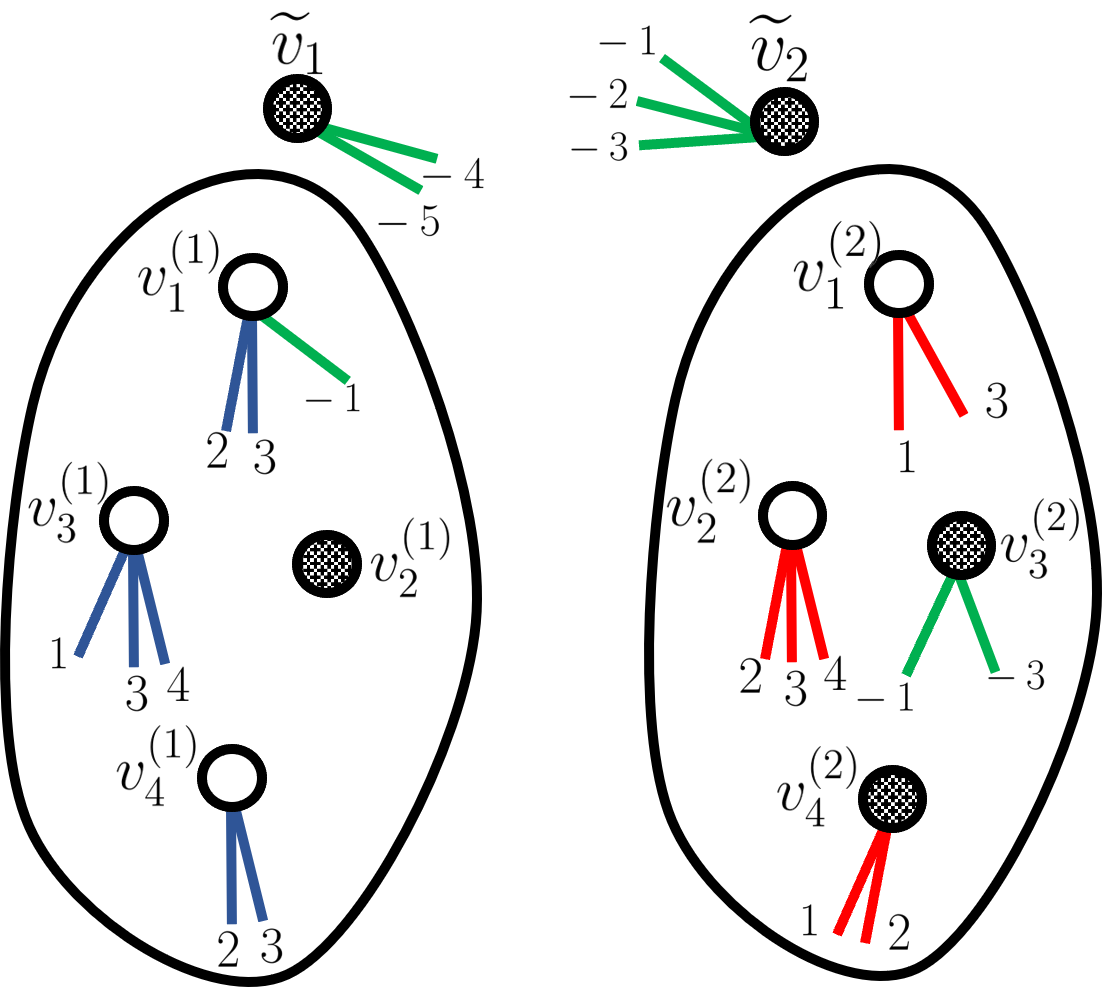}
		\caption*{($iv$): $\widetilde{S}_{(2\gets 1)} = ((v_3^{(1)},-1), (v_2^{(1)},-1))$.}
	\end{subfigure}
	\caption{A realization of sequences $\widetilde{S}_{(1\gets 1)}$, $\widetilde{S}_{(2\gets 2)}$, $\widetilde{S}_{(2\gets 1)}$ and $\widetilde{S}_{(1\gets 2)}$. We use the same configuration and same convention as in Figure \ref{fig:augmentedp}.}
	\label{fig:coupling_multi_aug}
\end{figure}

\begin{figure}[!htbp]
	\centering
	\begin{subfigure}[t]{.45\linewidth}
		\centering
		\includegraphics*[width = 0.75\textwidth]{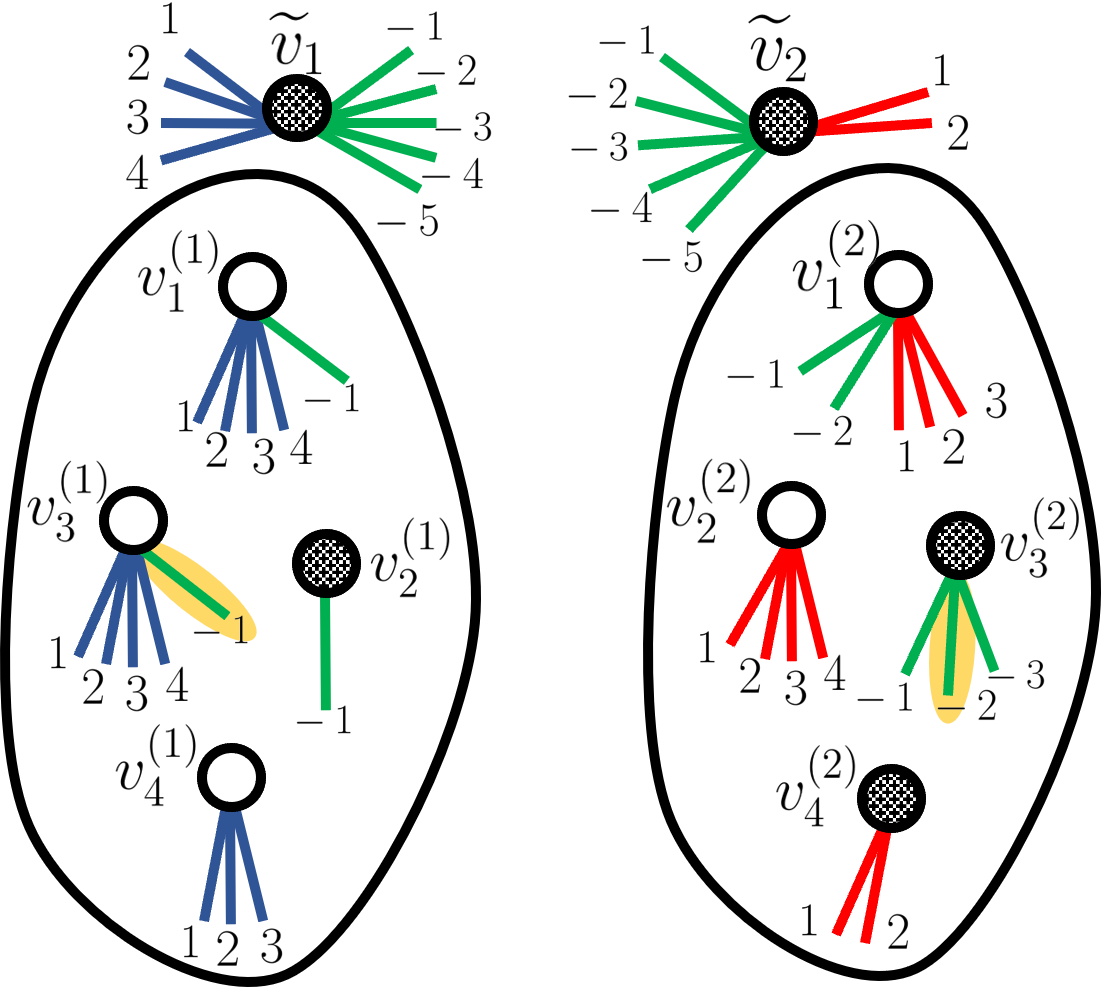}
		\caption{$\widehat{e}({\floor{t_\kappa n}}) = (v_3^{(2)},-2)$, coin toss is head, $\Upsilon_{\floor{t_\kappa n}}$ is identity map, and $e = (v_3^{(1)},-1)$.}
	\end{subfigure}\hfill
	\begin{subfigure}[t]{.45\linewidth}
		\centering
		\includegraphics*[width = 0.75\textwidth]{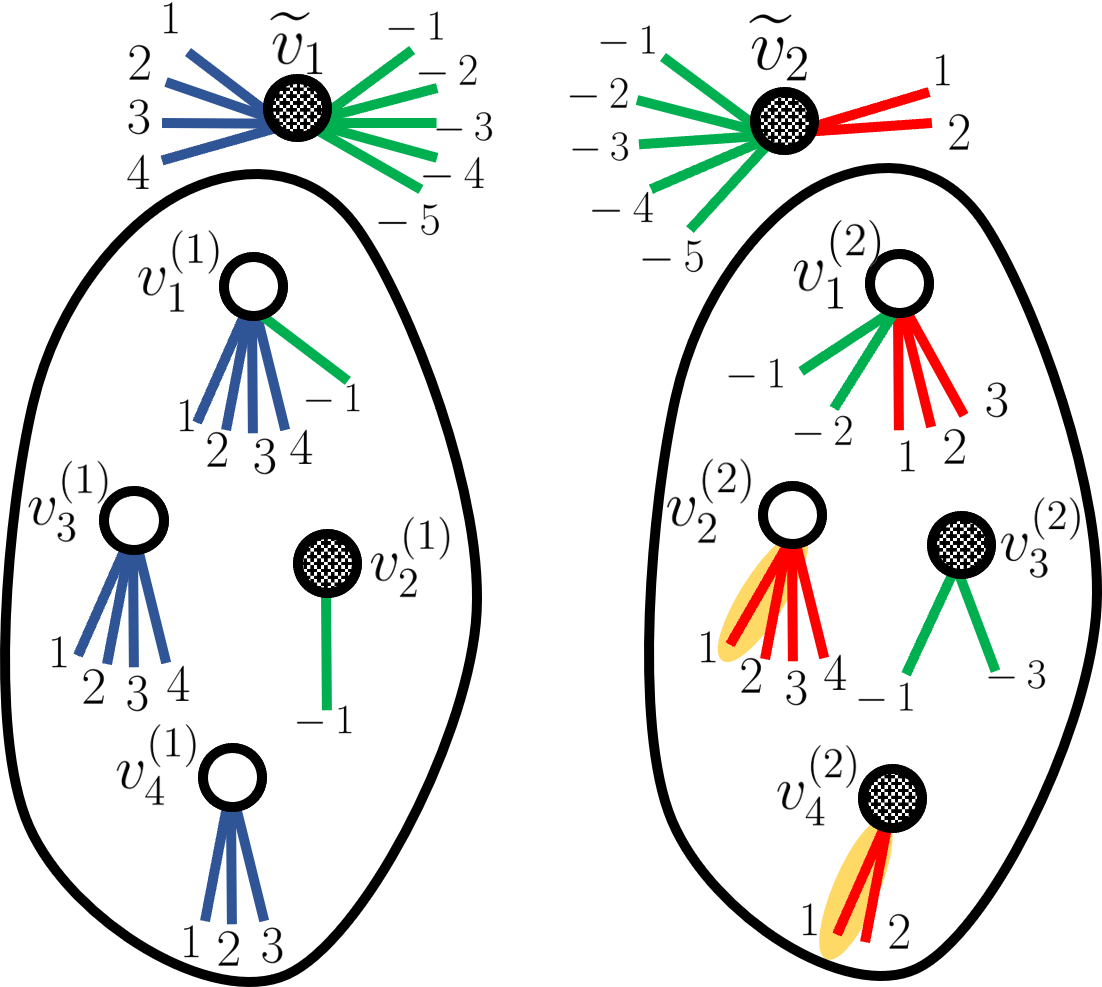}
		\caption{$\widehat{e}({\floor{t_\kappa n}}+1) = (v_4^{(2)},1)$, coin toss is head, $\Upsilon_{\floor{t_\kappa n}+1}$ is not identity map, and $e = (v_2^{(2)},1)$.}
	\end{subfigure}\\
	\medskip
	\begin{subfigure}[t]{.45\linewidth}
		\centering
		\includegraphics*[width = 0.75\textwidth]{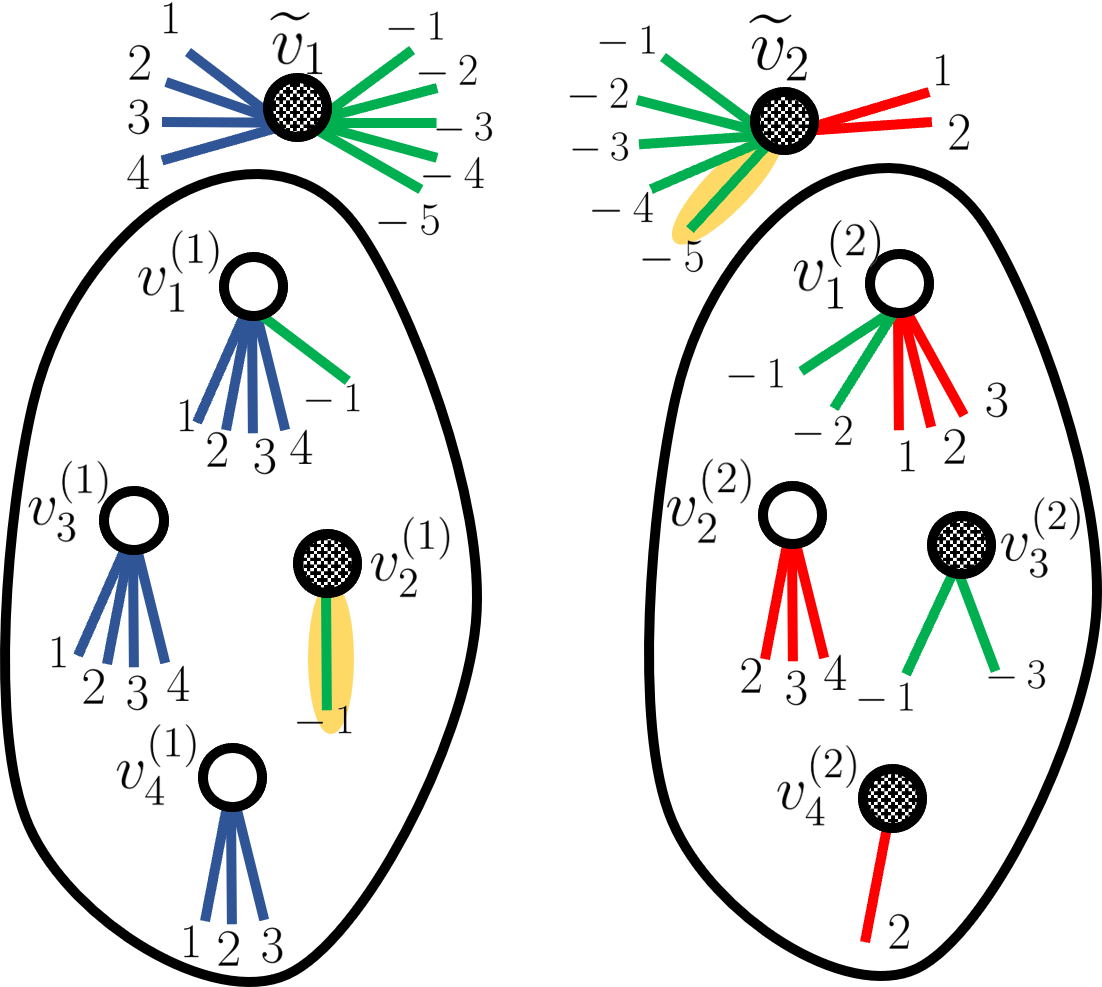}
		\caption{$\widehat{e}({\floor{t_\kappa n}}+2) = (v_2^{(1)},-1)$, coin toss is tail, $\Upsilon_{\floor{t_\kappa n}+2}$ is not identity map, and $e = (\widetilde{v}_2,-5)$.}
	\end{subfigure}\hfill
	\begin{subfigure}[t]{.45\linewidth}
		\centering
		\includegraphics*[width = 0.75\textwidth]{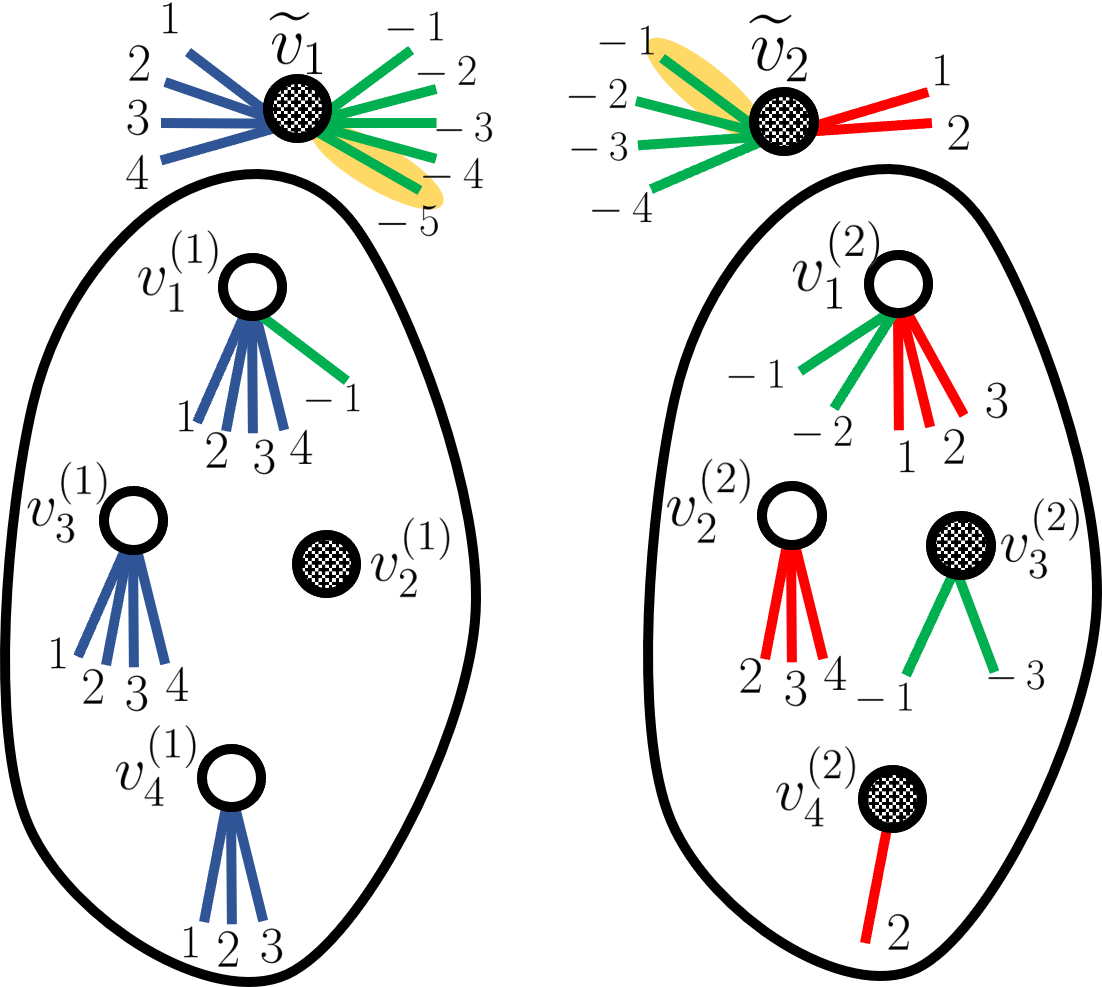}
		\caption{$\widehat{e}({\floor{t_\kappa n}}+3) = (\widetilde{v}_1,-5)$, coin toss is head, $\Upsilon_{\floor{t_\kappa n}+3}$ is not identity map, and $e = (\widetilde{v}_2,-1)$.}
	\end{subfigure}\\
	\medskip
	\begin{subfigure}[t]{.45\linewidth}
		\centering
		\includegraphics*[width = 0.75\textwidth]{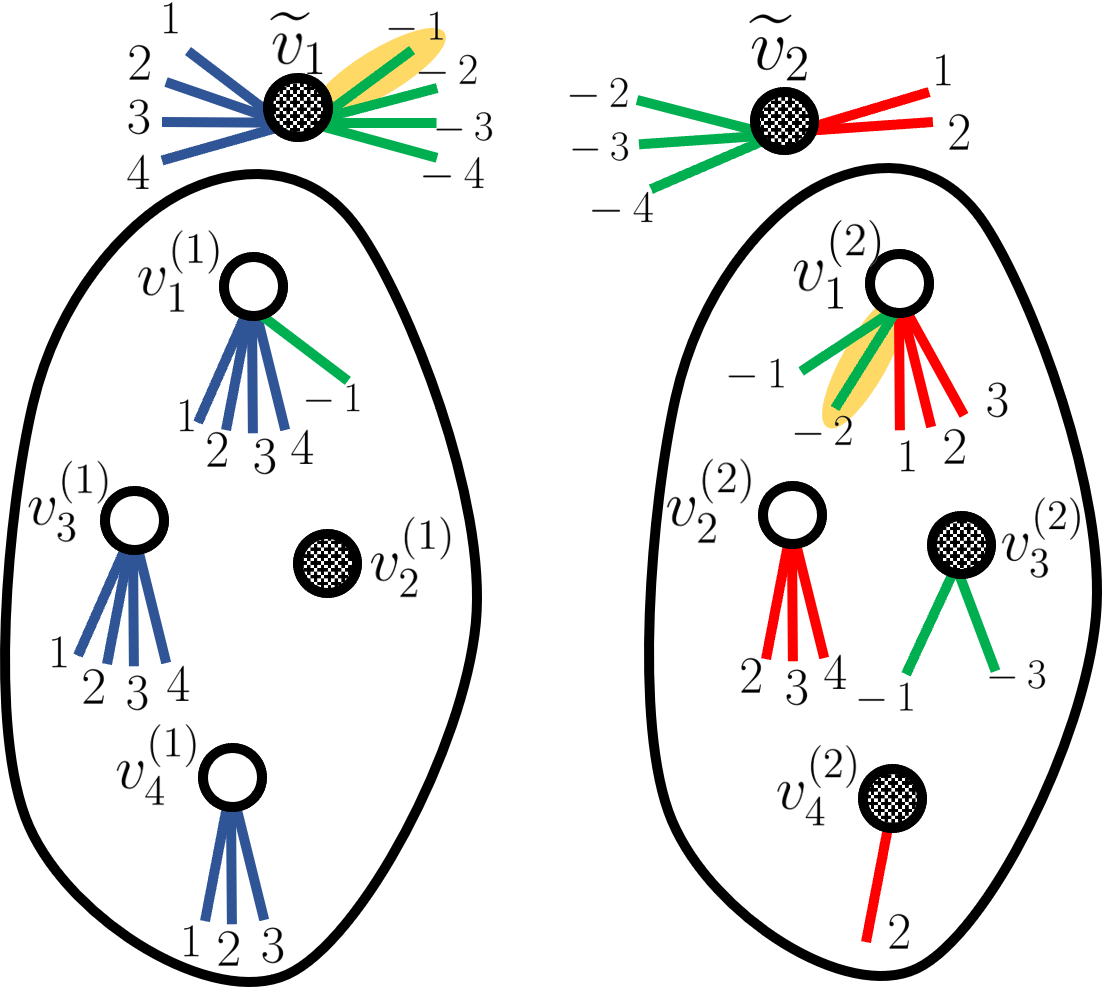}
		\caption{$\widehat{e}({\floor{t_\kappa n}}+4) = (\widetilde{v}_1,-1)$, coin toss is tail, $\Upsilon_{\floor{t_\kappa n}+4}$ is not identity map, and $e \!= \!({v}_1^{(2)}\!,\!-2)$.}
	\end{subfigure}\hfill
	\begin{subfigure}[t]{.45\linewidth}
		\centering
		\includegraphics*[width = 0.75\textwidth]{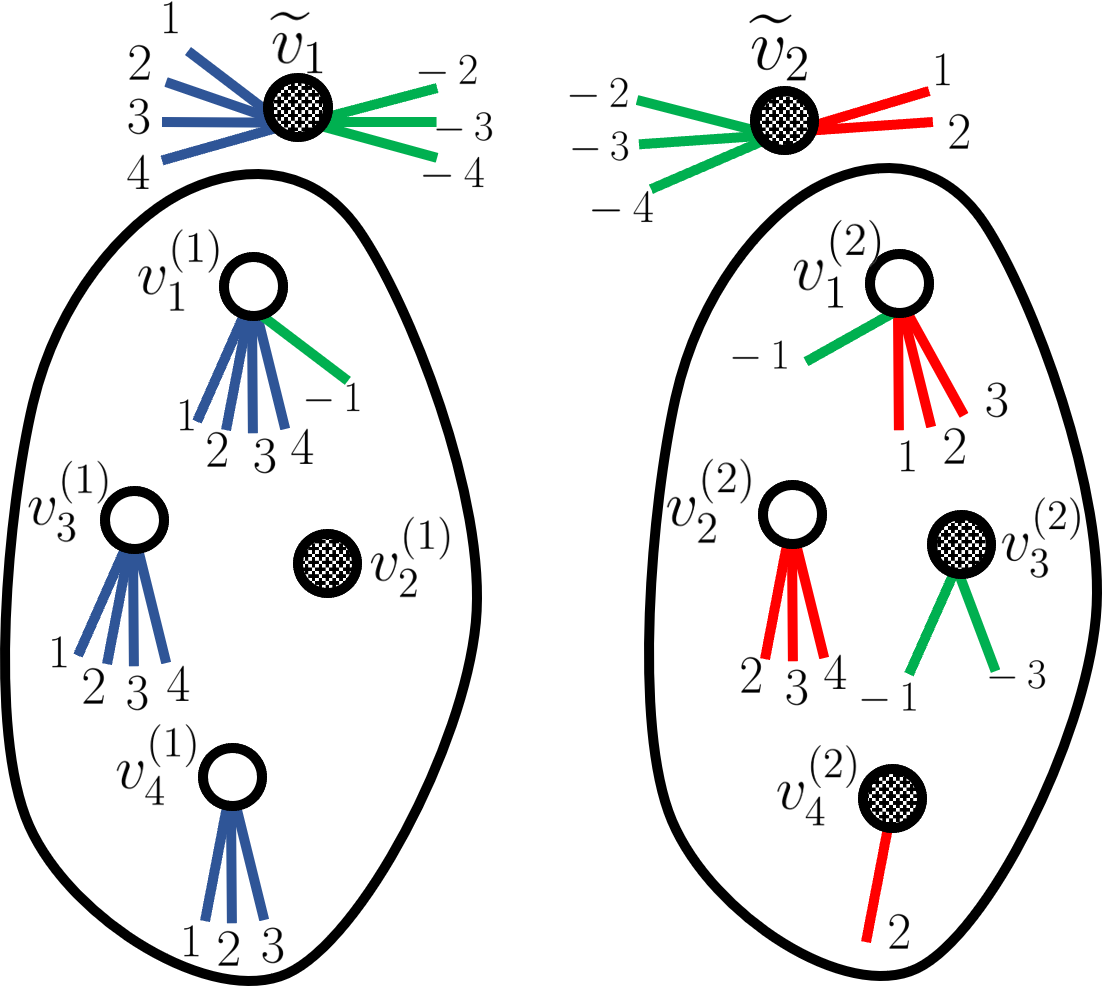}
		\caption{The two processes are decoupled, since $\widehat{T}_{(1\gets   2)}(k) = \floor{n \chi^{(2\gets 1)} \lambda_m(n)} = 3$.}
	\end{subfigure}
	\caption{Coupling of the twisted process and the augmented process given by Figure \ref{fig:coupling_multi_aug}. We use the same convention as in Figure \ref{fig:twistp}. Note that $(c)$ to $(f)$ appears in both figures. Realization of the map $\Upsilon_k$ is denoted in Figure \ref{fig:Upsilonrel}.}
	\label{fig:coupling_multi_twist}
\end{figure}

\begin{figure}[!htbp]
	\centering
	\begin{subfigure}[t]{0.45\linewidth}
		\centering
		\includegraphics*[width = 0.9\textwidth]{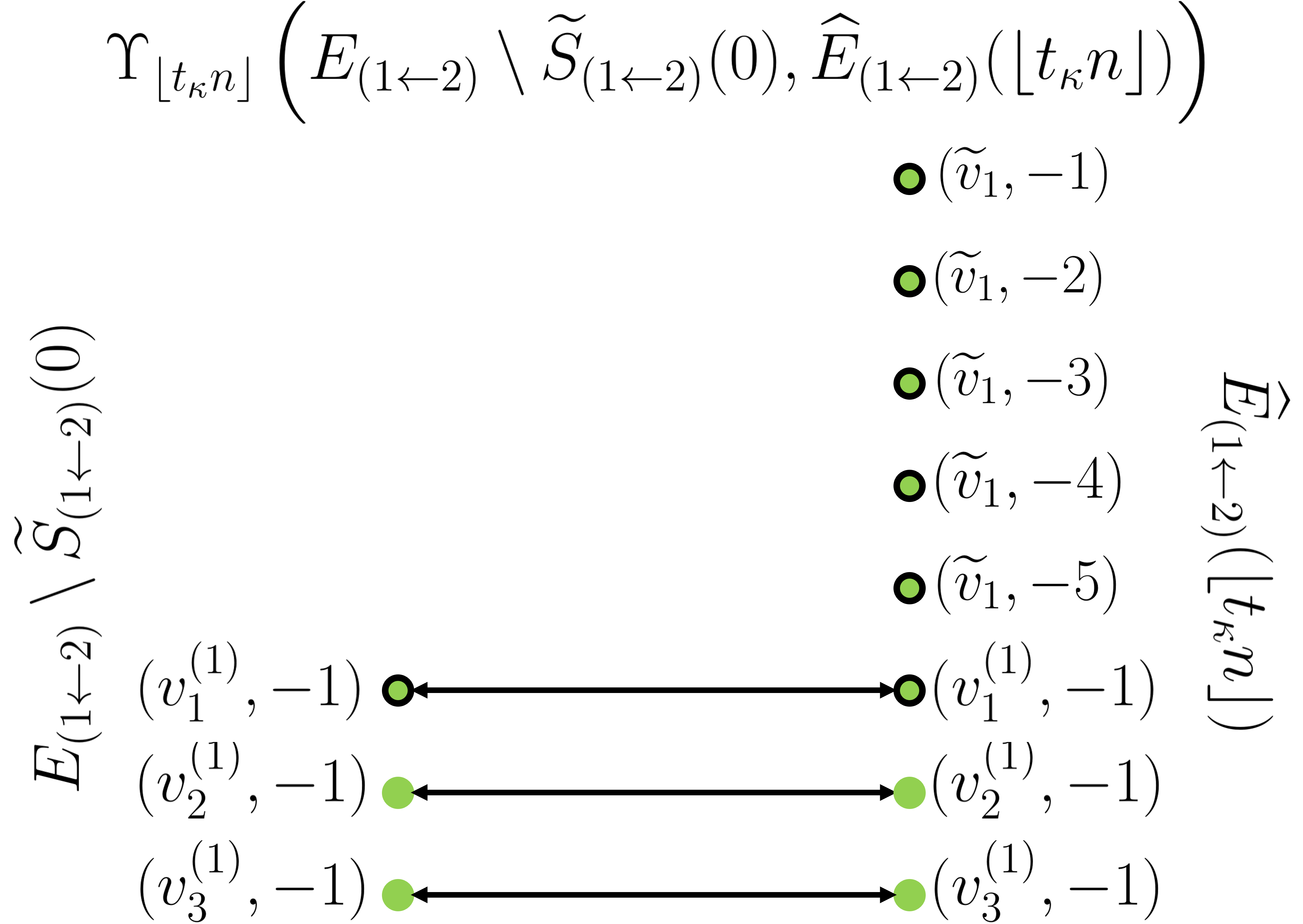}
		\caption{$\widetilde{e}_{(2\gets   1)}(1) = (v_3^{(1)},-1)$}
	\end{subfigure}\hfill
	\begin{subfigure}[t]{0.45\linewidth}
		\centering
		\includegraphics*[width = 0.9\textwidth]{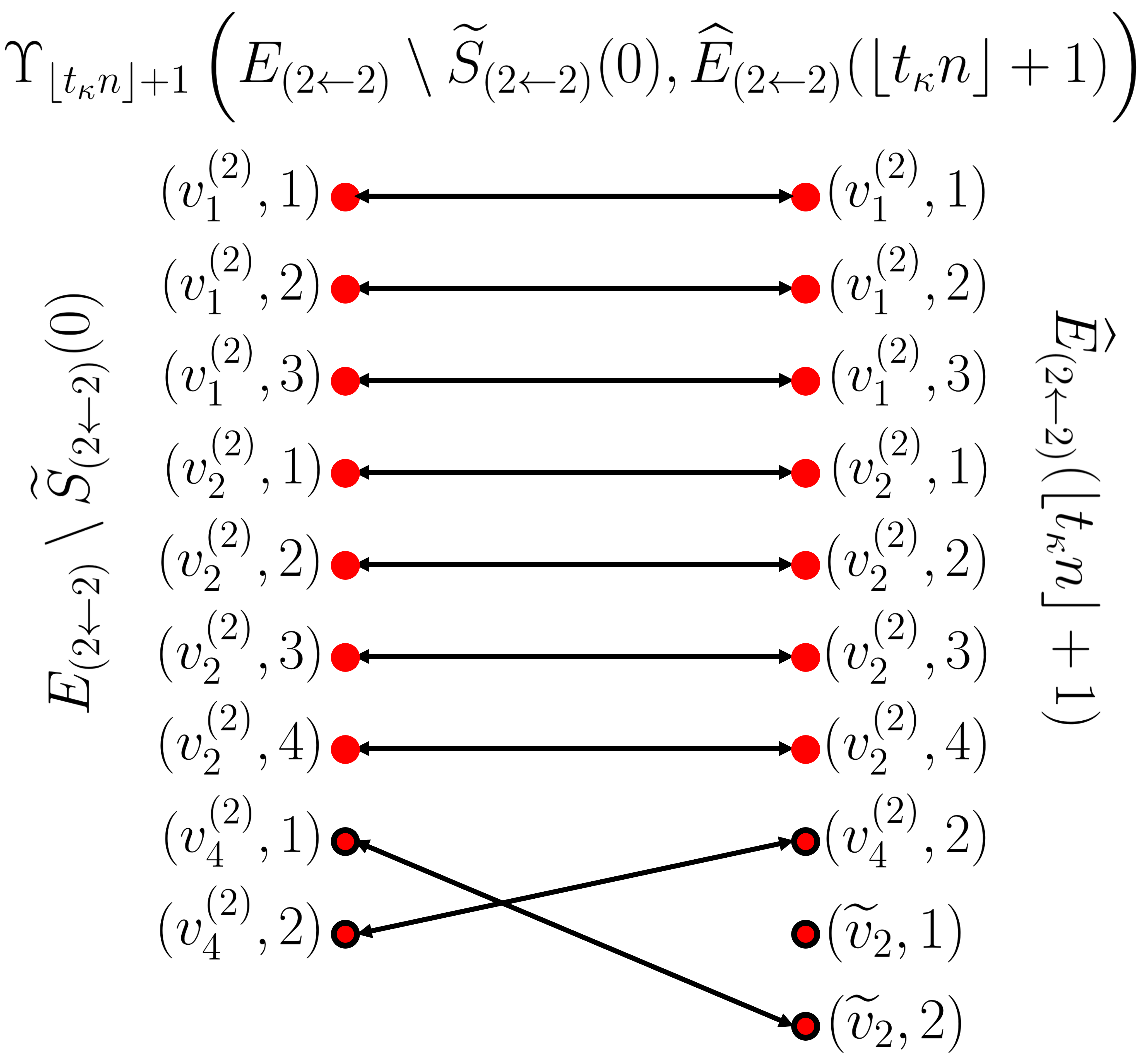}
		\caption{$\widetilde{e}_{2\gets 2}(1) = (v_2^{(2)},1)$.}
	\end{subfigure}\\
	\medskip
	\begin{subfigure}[t]{0.45\linewidth}
		\centering
		\includegraphics*[width = 0.9\textwidth]{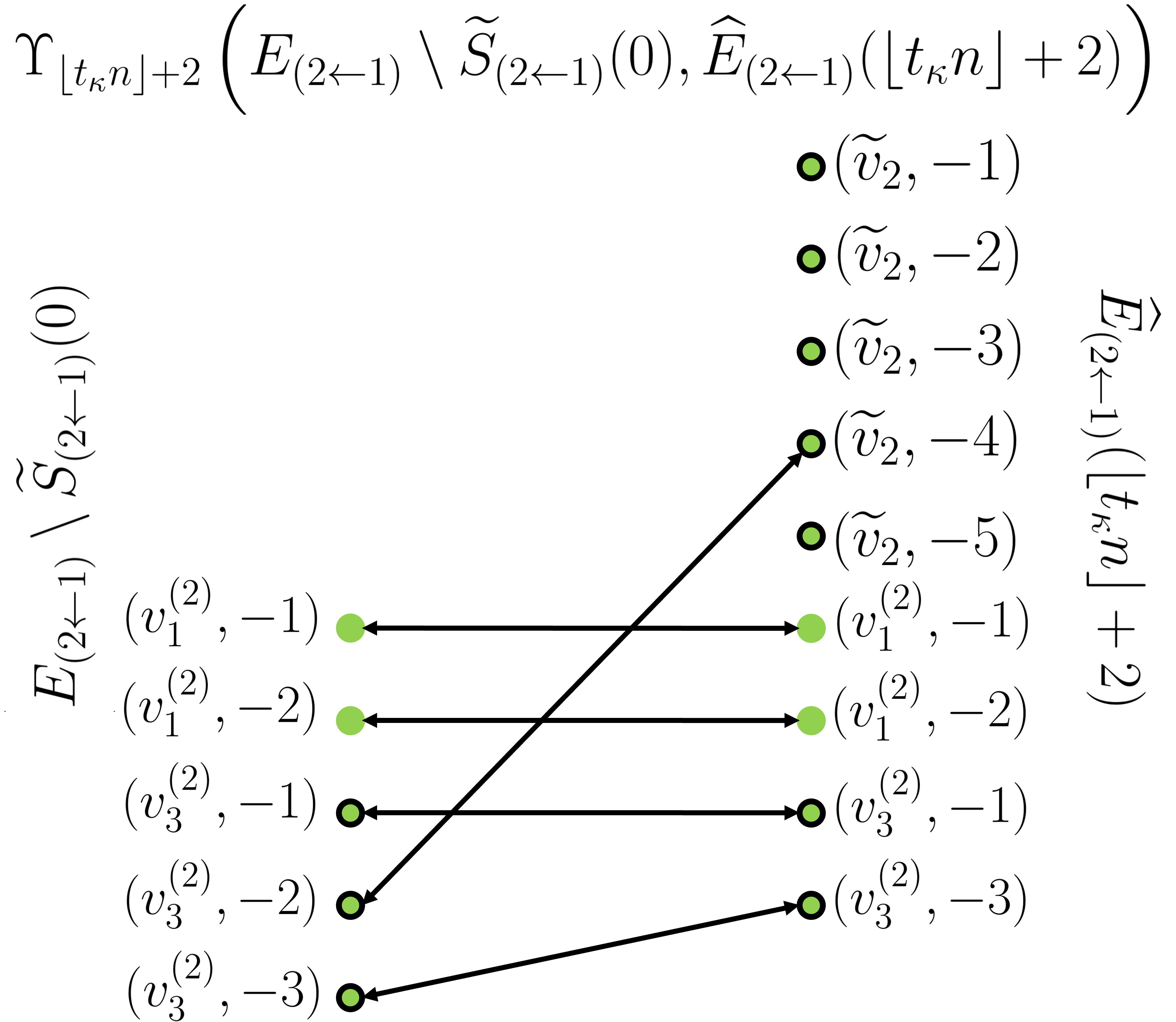}
		\caption{$\widetilde{e}_{(1\gets   2)}(1) = (v_1^{(2)},-2)$.}
	\end{subfigure}\hfill
	\begin{subfigure}[t]{0.45\linewidth}
		\centering
		\includegraphics*[width = 0.9\textwidth]{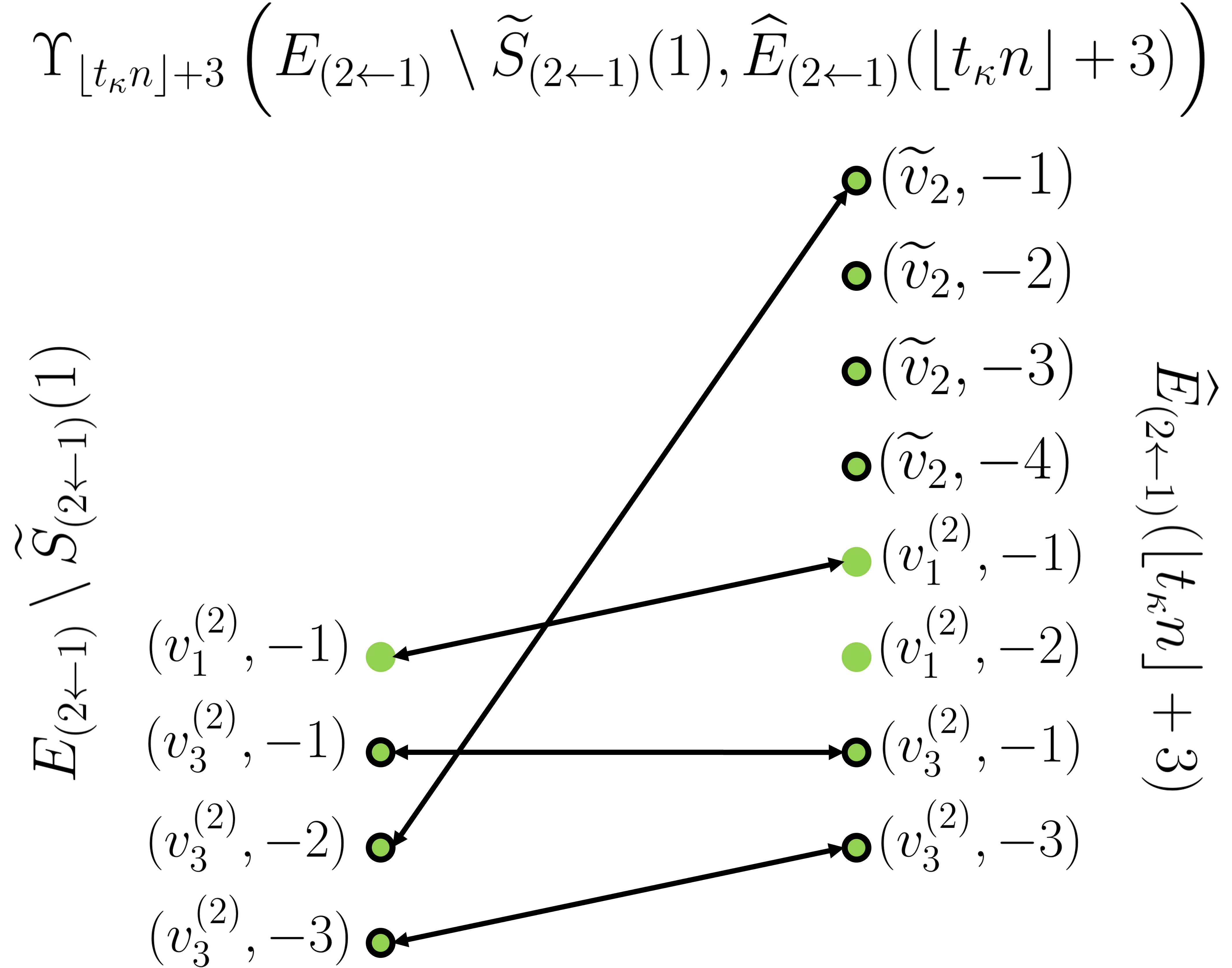}
		\caption{$\widetilde{e}_{(1\gets   2)}(2) = (v_3^{(2)},-2)$.}
	\end{subfigure}\\
	\medskip
	\begin{subfigure}[t]{0.45\linewidth}
		\centering
		\includegraphics*[width = 0.9\textwidth]{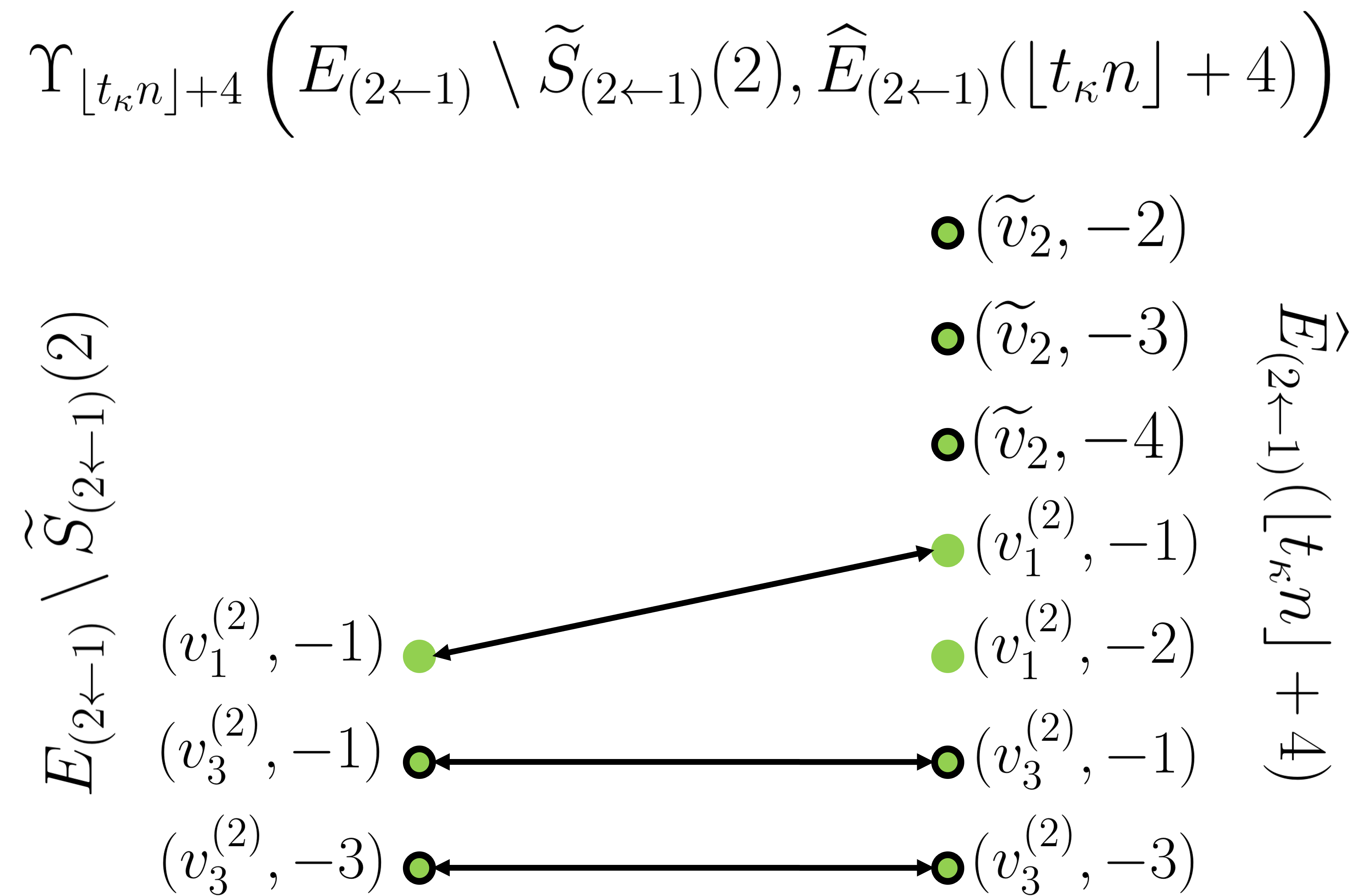}
		\caption{$\widetilde{e}_{(1\gets   2)}(3) = (v_1^{(2)},-1)$.}
	\end{subfigure}
	\caption{Realization of the relabeling function in the coupling given by Figure \ref{fig:coupling_multi_twist}. Color of each node matches the type of associated half-edge. Outline of nodes associated with active half-edges are black.}
	\label{fig:Upsilonrel}
\end{figure}

\begin{figure}[!]
	\centering
	\begin{subfigure}[t]{0.45\linewidth}
		\centering
		\includegraphics*[width = 0.75\textwidth]{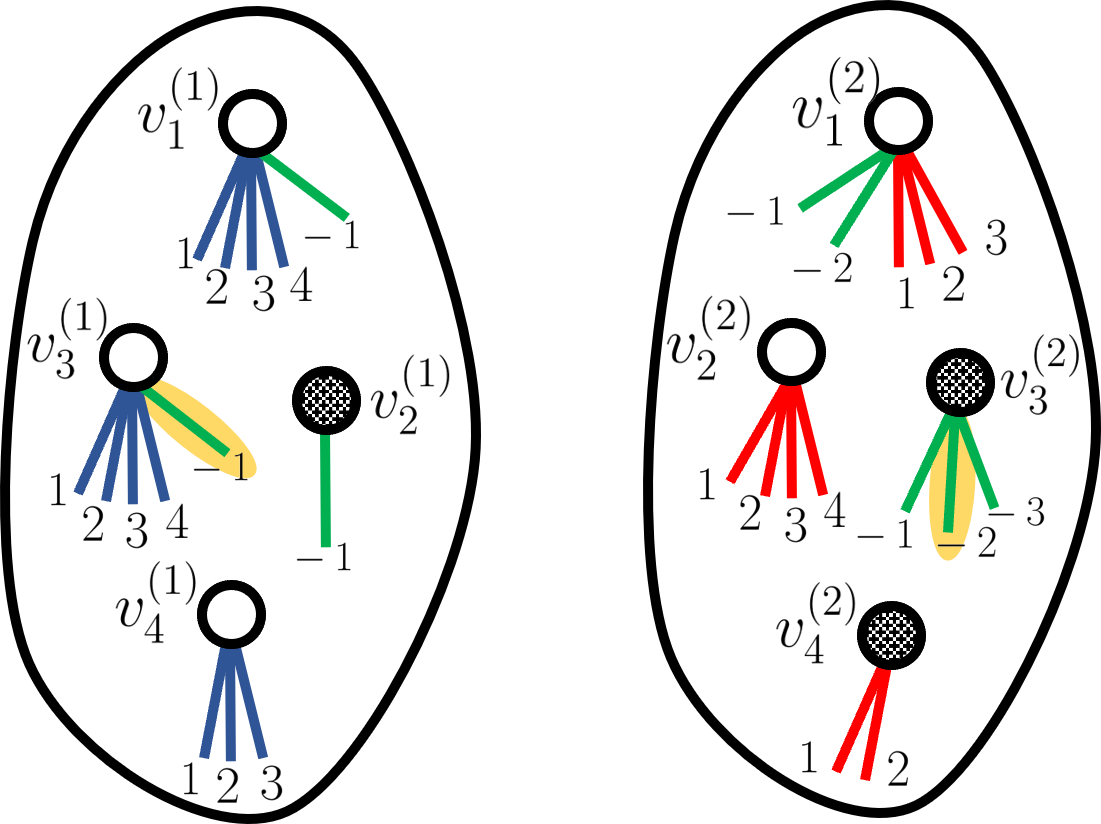}
		\caption{pair $(v_3^{(2)},-2)$ and $(v_3^{(1)},-1)$.}
	\end{subfigure}\hfill
	\begin{subfigure}[t]{0.45\linewidth}
		\centering
		\includegraphics*[width = 0.75\textwidth]{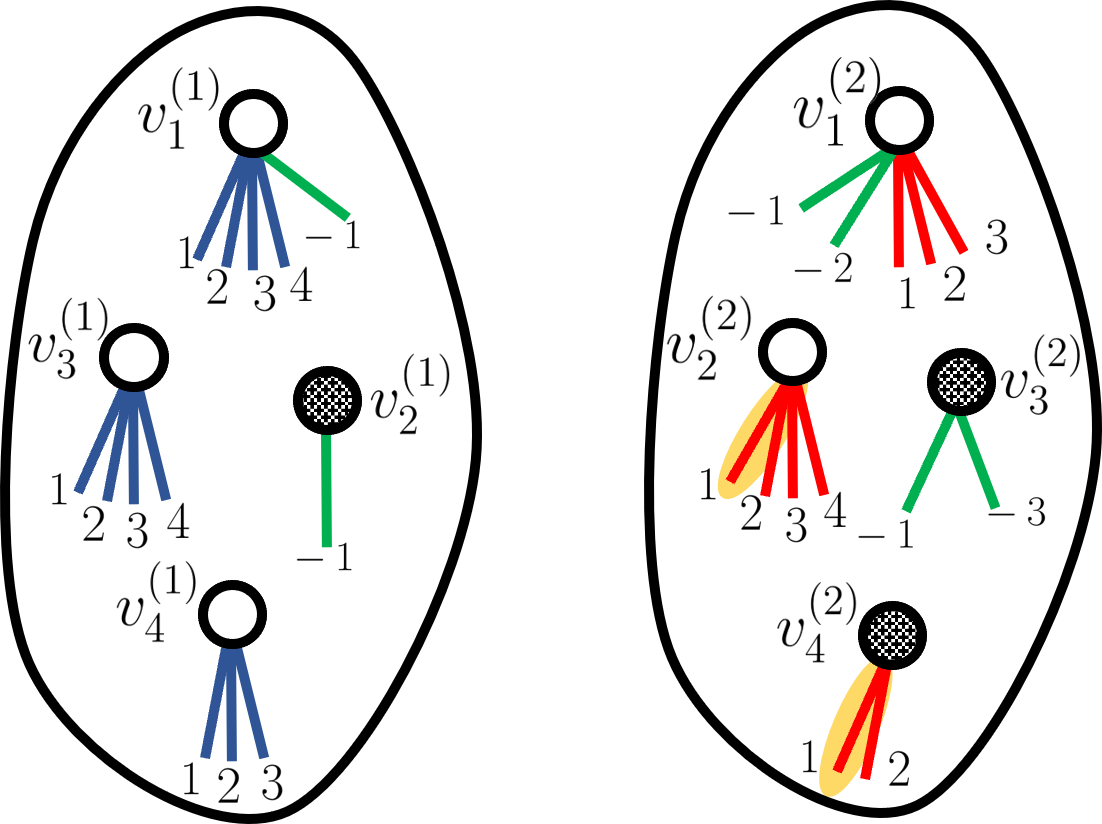}
		\caption{pair $(v_4^{(2)},1)$ and $(v_2^{(2)},1)$.}
	\end{subfigure}\\
	\medskip
	\begin{subfigure}[t]{0.45\linewidth}
		\centering
		\includegraphics*[width = 0.75\textwidth]{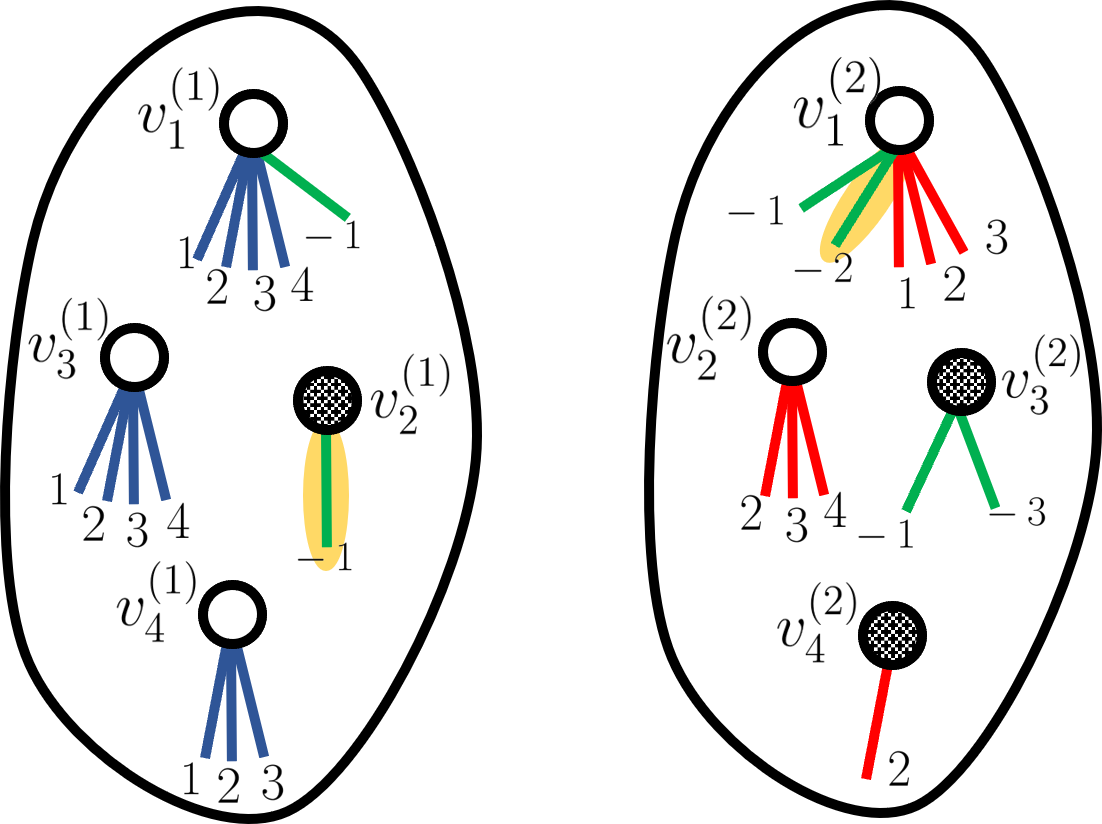}
		\caption{pair $(v_2^{(1)},-1)$ and $({v}_3^{(2)},-1)$.}
	\end{subfigure}\hfill
	\begin{subfigure}[t]{0.45\linewidth}
		\centering
		\includegraphics*[width = 0.75\textwidth]{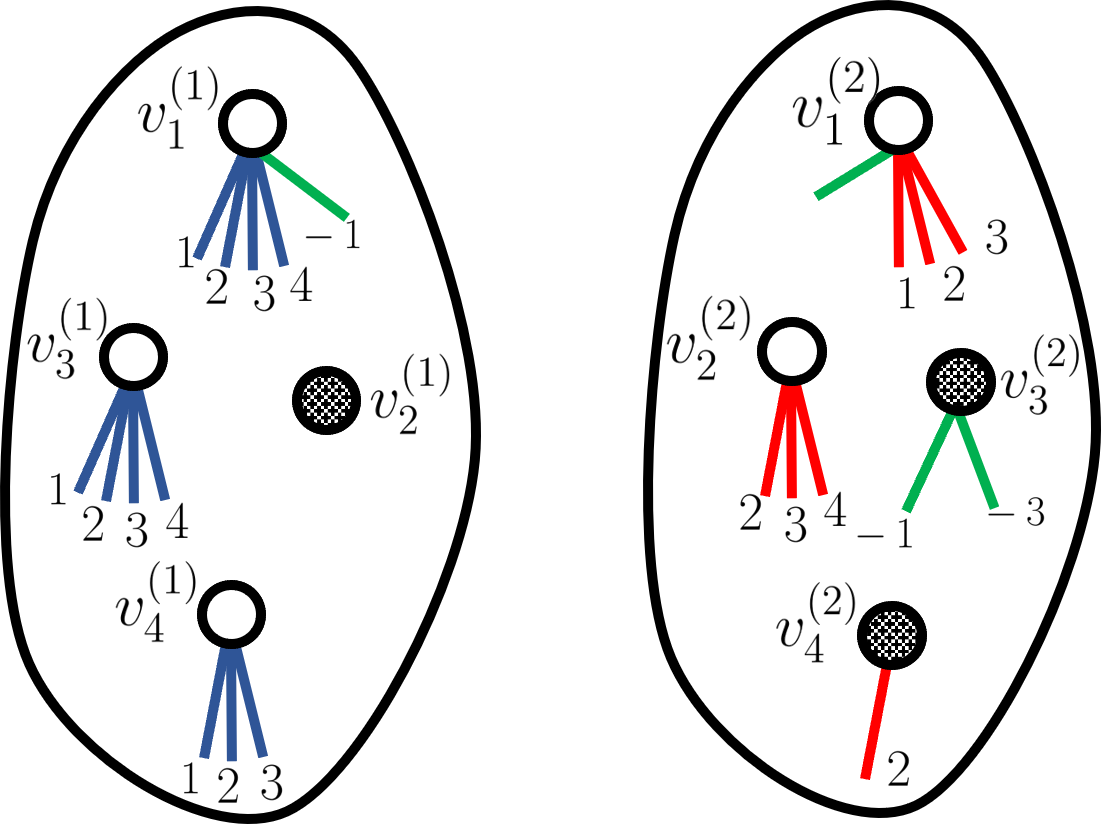}
		\caption{continue coupling.}
	\end{subfigure}
	\caption{Coupling of the truncated process and twisted process given by Figure \ref{fig:coupling_multi_twist}.}
	\label{fig:coupling_multi_trunc}
\end{figure}

Consider a typical time $k$ in the twisted process and suppose that we are given an active half-edge $\widehat{e}(k)$. Note that $\widehat{e}(k)$ can be either a regular or an augmented active half-edge. Suppose that Conditions \ref{def:couplingcond_1}-\ref{def:couplingcond_3} are satisfied. Based on the type of $\widehat{e}(k)$, we have one of the following cases:
\begin{enumerate}[leftmargin=!,itemindent=1em,align=left,label = Case C.$\arabic*$:,ref= C.$\arabic*$]
	\item({\it within community connection}) Suppose that $\widehat{e}(k)$ belongs to community $j\in\{1,2\}$ and it can be paired with another half-edge in the same community. Let $\widehat{T}_{(j\gets j)}(k)$ denote the number of times, after adding the new vertices, that we have removed pair of half-edges from community $j$. Note that by Condition \ref{def:couplingcond_1}, we have $\widehat{T}_{(j\gets j)}(k) < 2\floor{n\chi_\kappa^{(j\gets j)} \lambda_j(n)/2}$. Let $\widehat{E}_{(j\gets j)}(k)$ denote the set of {\em all half-edges}, both augmented and truncated, available in community $j$ at time $k$, to be paired with $\widehat{e}(k)$. Note that
	\begin{align*}
		\left|\widehat{E}_{(j\gets j)}(k)\right| &= \left|{E}_{(j\gets j)}\right| + 2\floor{n\chi_\kappa^{(j\gets j)} \lambda_j(n)/2}  - 2\widehat{T}_{(j\gets j)}(k) - 1 \\
		&\geq \left|{E}_{(j\gets j)}\right| - \widehat{T}_{(j\gets j)}(k) = \left|{E}_{(j\gets j)}\setminus \widetilde{S}_{(j\gets j)}(\widehat{T}(k))\right|
	\end{align*}
	Consider a realization of $\Upsilon_k\left({E}_{(j\gets j)}\setminus \widetilde{S}_{(j\gets j)}(\widehat{T}_{(j\gets j)}(k)),\widehat{E}_{(j\gets j)}(k)\right)$. We pair $\widehat{e}(k)$ with a random half-edge $e$ by first tossing a biased coin, where the ratio of head and tail probabilities are
		\begin{align*}
			\frac{\left|{E}_{(j\gets j)}\setminus \widetilde{S}_{(j\gets j)}(\widehat{T}_{(j\gets j)}(k))\right|}{\left|\widehat{E}_{(j\gets j)}(k)\right| - \left|{E}_{(j\gets j)}\setminus \widetilde{S}_{(j\gets j)}(\widehat{T}_{(j\gets j)}(k))\right|}.
		\end{align*}
	If the outcome is heads, we set $e$ to be
	\begin{align*}
		\Upsilon_k\left({E}_{(j\gets j)}\setminus \widetilde{S}_{(j\gets j)}(\widehat{T}_{(j\gets j)}(k)),\widehat{E}_{(j\gets j)}(k)\right) (\widetilde{e}_{(j\gets j)}(\widehat{T}_{(j\gets j)}(k)+1));
	\end{align*}
	otherwise, we pick $e$ uniformly at random from
	\begin{align*}
		\widehat{E}_{(j\gets j)}(k)\setminus \Upsilon_k\left({E}_{(j\gets j)}\setminus \widetilde{S}_{(j\gets j)}(\widehat{T}_{(j\gets j)}(k)),\widehat{E}_{(j\gets j)}(k)\right).
	\end{align*}
	If $e$ was a regular half-edge, then pick $\widehat{e}(k+1)$ uniformly at random from the set of all available active half-edges. Otherwise, set $\widehat{e}(k+1)$ to be an augmented half-edge in community $j$ that can be paired with half-edges in community $j$.
	\item({\it between community connection}) Suppose that $\widehat{e}(k)$ belongs to community $j\in\{1,2\}$ and it can be paired with another half-edge in community $-j$. Let $\widehat{T}_{(j\gets   -j)}(k)$ denote the number of times, after adding the new vertices, that we have removed pair of half-edges from between the communities, such that the initial active half-edge was in community $j$. Note that by Condition \ref{def:couplingcond_2}, \label{def:coupling_case2}
	\begin{align*}
		\widehat{T}_{(1\gets   2)}(k) < \floor{n\chi_\kappa^{(2\gets 1)} \lambda_m(n)}, \qquad \widehat{T}_{(2\gets   1)}(k) < \floor{n\chi_\kappa^{(1\gets 2)} \lambda_m(n)}.
	\end{align*}
	Let $\widehat{E}_{(j\gets   -j)}(k)$ denote the set of {\em all half-edges}, both augmented and truncated, available in community $-j$ at time $k$, to be paired with $\widehat{e}(k)$.  Note that
	\begin{align*}
		&\left|\widehat{E}_{(j\gets   -j)}(k)\right| \\
		&\myquad[2]= \left|{E}_{(j\gets   -j)}\right| + \floor{n\chi_\kappa^{(1\gets 2)} \lambda_m(n)} +  \floor{n\chi_\kappa^{(2\gets 1)} \lambda_m(n)}  - \widehat{T}_{(1\gets   2)}(k) - \widehat{T}_{(2\gets   1)}(k) \\
		&\myquad[2]> \left|{E}_{(j\gets   -j)}\right| - \widehat{T}_{(j\gets   -j)}(k) =  \left|{E}_{(j\gets   -j)}\setminus \widetilde{S}_{(j\gets -j)}(\widehat{T}_{(j\gets   -j)}(k))\right|.
	\end{align*}
	Consider a realization of $\Upsilon_k\left({E}_{(j\gets   -j)}\setminus \widetilde{S}_{(j\gets -j)}(\widehat{T}_{(j\gets   -j)}(k)),\widehat{E}_{(j\gets   -j)}(k)\right)$. We pair $\widehat{e}(k)$ with a random half-edge $e$ by first tossing a biased coin, where the ratio of head and tail probabilities are
	\begin{align*}
		\frac{\left|{E}_{(j\gets   -j)}\setminus \widetilde{S}_{(j\gets -j)}(\widehat{T}_{(j\gets   -j)}(k))\right|}{\left|\widehat{E}_{(j\gets   -j)}(k)\right| - \left|{E}_{(j\gets   -j)}\setminus \widetilde{S}_{(j\gets -j)}(\widehat{T}_{(j\gets   -j)}(k))\right|}.
	\end{align*}
	If the outcome is heads, we set $e$ to be
	\begin{align*}
		\Upsilon_k\left({E}_{(j\gets   -j)}\setminus \widetilde{S}_{(j\gets -j)}(\widehat{T}_{(j\gets   -j)}(k)),\widehat{E}_{(j\gets   -j)}(k)\right) (\widetilde{e}_{(j\gets   -j)}(\widehat{T}_{(j\gets   -j)}(k)+1));
	\end{align*}
	otherwise, we pick $e$ uniformly at random from
	\begin{align*}
		\widehat{E}_{(j\gets   -j)}(k)\setminus \Upsilon_k\left({E}_{(j\gets   -j)}\setminus \widetilde{S}_{(j\gets -j)}(\widehat{T}_{(j\gets   -j)}(k)),\widehat{E}_{(j\gets   -j)}(k)\right).
	\end{align*}
	If $e$ was a regular half-edge, then pick $\widehat{e}(k+1)$ uniformly at random from the set of all available active half-edges. Otherwise, set $\widehat{e}(k+1)$ to be an augmented half-edge in community $j$, that can be paired with half-edges in community $-j$.
\end{enumerate}
In Figure \ref{fig:coupling_multi_aug}, we illustrate a realization of the sequences $\widetilde{S}_{(1\gets 1)}$, $\widetilde{S}_{(2\gets 2)}$, $\widetilde{S}_{(2\gets 1)}$ and $\widetilde{S}_{(1\gets 2)}$ during the augmented process. We then use these sequences to illustrate a coupling between the augmented and the twisted processes in Figure \ref{fig:coupling_multi_twist}, using the realization of the relabeling functions given in Figure \ref{fig:Upsilonrel}. We also present the natural coupling between the resulted twisted and the truncated processes in Figure \ref{fig:coupling_multi_trunc}.

Next, we study which of Conditions \ref{def:couplingcond_1}-\ref{def:couplingcond_3} are violated first; in particular, we show that with high probability Conditions \ref{def:couplingcond_3} will be violated first. Let $\widehat{k}_{\text{dec}}$ denote the time at which we decouple the twisted and augmented processes. Let $\widehat{T}_{(-j\gets j)}(\widehat{k}_{\text{dec}})$ denote the number of times, after the augmentation, that we have paired an active half-edge in community $-j$ with a random half-edge in community $j$. Based on which condition is violated at time $\widehat{k}_{\text{dec}}$, we have one of the followings:
\begin{enumerate}[leftmargin=!,itemindent=1em,align=left,label = Case V.$\arabic*$:]
	\item({\it Condition \ref{def:couplingcond_1} is violated first}) Suppose that by time $\widehat{k}_{\text{dec}}$, we have fully explored $\widetilde{S}_{(j\gets j)}$ for some $j\in\{1,2\}$. In this case, we have removed $2\floor{n\chi_\kappa^{(j\gets j)} \lambda_j(n)/2}$ pair of half-edges from community $j$ and paired them with each other. In particular, the total number of remaining half-edges in community $j$ of the twisted process for in-community connections at time $\widehat{k}_{\text{dec}}$ is the same as the total number of remaining half-edges in community $j$ of the augmented process for in-community connections at time $\floor{{t}^{(1\gets 2)}_\kappa n}$. Note that for any half-edge removed from community $j$ in the augmented process during this time, we have either removed the same half-edge or wasted an active half-edge in the twisted process. Also, note that removing half-edges from community $-j$ or from between the communities can only increase the total number of active half-edges in community $j$ by activating some of the vertices in community $j$. Hence, we have the following inequality:
	\begin{align*}
		\widetilde{A}_j(\floor{{t}^{(1\gets 2)}_\kappa n}) \geq \widehat{A}_j(\widehat{k}_{\text{dec}}) \geq \widehat{A}_{\text{aug},j}(\widehat{k}_{\text{dec}}),
	\end{align*}
	where $\widehat{A}_{\text{aug},j}(k)$ is the number of remaining augmented half-edges in community $j$ of the twisted process at time $k$ that can be paired with half-edges in the same community.
	\item({\it Condition \ref{def:couplingcond_2} is violated first}) Suppose that by time $\widehat{k}_{\text{dec}}$, we have fully explored $\widetilde{S}_{(j\gets -j)}$ for some $j\in\{1,2\}$. In this case, we have picked $\floor{n \chi^{(-j\gets j)} \lambda_m(n)}$ active half-edges from community $j$ and then paired them with random half-edges that belong to community $-j$. We have further removed $\widehat{T}_{(-j\gets j)}(\widehat{k}_{\text{dec}})$ random half-edges from community $j$, which are paired with active half-edges that belongs to community $-j$. Note that $\widehat{T}_{(-j\gets j)}(\widehat{k}_{\text{dec}})\leq \floor{n \chi^{(j\gets -j)} \lambda_m(n)}$. Following the same coupling as in Case \ref{def:coupling_case2} of the above coupling (tossing a biased coin and removing half-edges accordingly), we remove $\floor{n \chi^{(j\gets -j)} \lambda_m(n)} - \widehat{T}_{(-j\gets j)}(\widehat{k}_{\text{dec}})$ random half-edges from community $j$. We refer to this removal, as the final round of coupling. Let $N_m^{(j)}$ denote the number of active half-edges present at time $\widehat{k}_{\text{dec}}$ that are removed during the final round of coupling, and $N_m^{(\text{aug},j)}$ denote the number of augmented half-edges available at time $\widehat{k}_{\text{dec}}$ that are removed during the final round of coupling. Clearly, the following inequalities hold:
	\begin{align*}
		&\widehat{A}_m^{(j)}(\widehat{k}_{\text{dec}}) > N_m^{(j)}, ~\widehat{A}_m^{(\text{aug},j)}(\widehat{k}_{\text{dec}}) > N_m^{(\text{aug},j)}, \text{ and }\widehat{A}_m^{(j)}(\widehat{k}_{\text{dec}}) - \widehat{A}_m^{(\text{aug},j)}(\widehat{k}_{\text{dec}}) > N_m^{(j)} - N_m^{(\text{aug},j)},
	\end{align*}
	where $\widehat{A}_m^{(\text{aug},j)}(k)$ is the number of augmented half-edges in community $j$ of the twisted process at time $k$ that can be paired with half-edges in community $-j$.	Following the same argument as in the previous case, we have
	\begin{align*}
		\widetilde{A}_m^{(j)}(\floor{{t}^{(1\gets 2)}_\kappa n}) \geq \widehat{A}_m^{(j)}(\widehat{k}_{\text{dec}}) - N_m^{(j)} \geq \widehat{A}_m^{(\text{aug},j)}(\widehat{k}_{\text{dec}}) - N_m^{(\text{aug},j)}.
	\end{align*}
	 \item{(\it Condition \ref{def:couplingcond_3} is violated first)} This is the desired outcome.
\end{enumerate}

Note that during Phase \ref{def:twist_phase1} of the twisted process (and the final round of coupling), augmented half-edges are removed only if they are chosen uniformly at random. In particular, if the random half-edge that is paired with the active half-edge is augmented and both half-edges are in community $j$, then we remove two augmented half-edges from community $j$ (e.g. uniformly at random); similarly, if the random half-edge is in community $j$ and the active half-edge is in community $-j$ and the random half-edge is augmented, then we remove the random half-edge as well as another augmented half-edge from community $-j$ (e.g. uniformly at random). Hence, the number of remaining augmented half-edges in community $j$ after removing $2\floor{n\chi^{(j\gets j)}\lambda_j(n)/2}$ pair of half-edges since time $\floor{t_\kappa n}$, concentrates around
\begin{align*}
	\frac{\left(\mu_{\kappa}^{(j\gets j)}\right)^2 \lambda_j(n)}{2\chi_\kappa^{(j\gets j)}\lambda_j(n) + \left(\mu_{\kappa}^{(j\gets j)}\right)^2 \lambda_j(n)} \times 2\floor{n\chi^{(j\gets j)}\lambda_j(n)/2}
\end{align*}
with high probability. Note that the first term is (almost) the ratio of regular half-edges over all half-edges that are in community $j$ and can be paired with half-edges in the same community at the time of augmentation. Recall that by the analysis of Appendix \ref{app:augment}, the total number of half-edges in community $j$ of the augmented process at time $\floor{t^{(1\gets 2)}_\kappa n}$ concentrates around $n\widetilde{a}_j(t_\kappa^{(1\gets 2)})$, where $\widetilde{a}_j(t_\kappa^{(1\gets 2)}) \leq \kappa^2/2 +  \chi^{(j\gets j)} \lambda_j(n) \zeta_\kappa + O(\kappa^2)$, with high probability. By Conditions \ref{def:mathcalE_i} and \ref{def:mathcalE_ii}, and \eqref{eq:propineq_incomm}, for all small enough $\kappa$ (independent of $n$, so that constants in $O(\kappa^2)$ are small), we have
\begin{align*}
	\frac{\left(\mu_{\kappa}^{(j\gets j)}\right)^2 \lambda_j(n)}{2\chi_\kappa^{(j\gets j)}\lambda_j(n) + \left(\mu_{\kappa}^{(j\gets j)}\right)^2 \lambda_j(n)} \times \chi^{(j\gets j)}\lambda_j(n) &> \kappa^2/2 +  \chi^{(j\gets j)} \lambda_j(n) \zeta_\kappa + O(\kappa^2) \\
	&\geq \widetilde{a}_j(t_\kappa^{(1\gets 2)}),
\end{align*}
which in turn implies that the event that ``Condition \ref{def:couplingcond_1} is violated first'' is a low probability event. Similarly, the number of remaining augmented half-edges between the communities after removing $\floor{n\chi_\kappa^{(1\gets 2)} \lambda_m(n)} +  \floor{n\chi_\kappa^{(2\gets 1)} \lambda_m(n)}$ half-edges uniformly at random since time $\floor{t_\kappa n}$, concentrates around
\begin{align*}
	&\frac{\mu_{\kappa}^{(1\gets 2)} \mu_{\kappa}^{(2\gets 1)}\lambda_m(n) }{ \lambda_m(n)\left(\chi_\kappa^{(1\gets 2)}+ \chi_\kappa^{(2\gets 1)}\right)\! + \! \mu_{\kappa}^{(1\gets 2)} \mu_{\kappa}^{(2\gets 1)}\lambda_m(n) }\! \times \! \left( \floor{n\chi_\kappa^{(1\gets 2)} \lambda_m(n)}\! + \! \floor{n\chi_\kappa^{(2\gets 1)} \lambda_m(n)}\right)\\
	&\myquad[2]\geq n \left(\kappa^2/2 +  \chi^{(-j\gets j)}\lambda_m(n) +\chi^{(j\gets -j)}\lambda_m(n) \zeta_\kappa+ O(\kappa^2)\right) \geq n\times \widetilde{a}^{(-j)}_m(t_\kappa^{(1\gets 2)})
\end{align*}
where the first inequality follows by Conditions \ref{def:mathcalE_i} and \ref{def:mathcalE_ii}, and \eqref{eq:propineq_incomm} for all small enough $\kappa$. Hence, the event that ``Condition \ref{def:couplingcond_2} is violated first'' is a low probability event. Gathering all the above results, the event that ``Condition \ref{def:couplingcond_3} is violated first'' in the above coupling is a high probability event.

To summarize, similar to the case of one community, for any small enough constant $\kappa > 0$ we have defined an event $\Omega_n(\kappa)$ with $\lim_{n\to\infty} P(\Omega_n(\kappa)) = 0$ such that outside $\Omega_n(\kappa)$ a scaled-version of the truncated process hits its natural stopping time at time $\left(t_*\pm O(\epsilon)\right)n$. Letting $n\to\infty$, then $\kappa\to 0$, we have the asymptotic characterization of the truncated process. Letting $\delta\to 0$, and using \eqref{eq:naturalstop}, we get the following theorem for the asymptotic behavior of the Markov process of adoption.

\begin{theorem}\label{thm:odesmaininf}
	Consider the Markov process of adoption. Assume there is a constant $0<d_{\max} \leq \infty$, independent of $n$, such that $K_j(d_j,d_{-j}) = d_j+d_{-j}$ for all $d_j+d_{-j} > d_{\max}$ and $j\in\{1,2\}$. Suppose the degree regularity conditions given in Definition \ref{def:regcon_deg} hold. Let $\bs{\mu}_\infty(t)$ denote the solution of ODEs \eqref{eq:alterode} where the function $\Ffuncbold_\infty(\cdot)$ is given by the right-hand side of \eqref{eq:meanfield_mujj}-\eqref{eq:meanfield_muj-j} (Mean-Field equations). Let $\bs{\mu}_{*,\infty} = \lim_{s\to\infty}\Ffuncbold_\infty^s(\bs{1})$ to be the closest fixed point of $\Ffuncbold_\infty(\cdot)$ to $\bs{1}$ in sup-norm. Define  the function $t_\infty:[0,1]^4\to \mathbb{R}_+$ using equation \eqref{eq:sol_t} as follows:
	\begin{align}
	t_\infty(\bs{\mu}) \coloneqq \frac{\lambda_1}{2}\left(1- \left({\mu^{(1\gets 1)}}\right)^2\right) + \frac{\lambda_2}{2}\left(1- {\left(\mu^{(2\gets 2)}\right)}^2\right) +\lambda_m\left(1-\mu^{(1\gets 2)}\mu^{(2\gets 1)}\right)\label{eq:t(mu)},
	\end{align}
	and let $t_{*,\infty} \coloneqq t_\infty(\bs{\mu}_{*,\infty})$. For $t\leq t_{*,\infty}$, define
	\begin{align}
	&i^{(j)}_{d_j,d_{-j},u_j,u_{-j},\infty}(t) = \prob_{j,m}(d_j,d_{-j})\, Bi(u_j;d_j,1-\mu^{(j\gets j)}_{t,\infty})\,Bi(u_{-j};d_{-j},1-\mu^{(j\gets -j)}_{t,\infty}),\label{eq:soli_final}\\
	&\tau_{j,\infty}(t) = \frac{\lambda_j}{2} \left(1 - \left(\mu^{(j\gets j)}_{t,\infty}\right)^2\right)\\
	&\tau_{m,\infty}(t) = \lambda_m\left(1 - \mu^{(2\gets 1)}_{t,\infty}\mu^{(1\gets 2)}_{t,\infty}\right)
	\end{align}
	where $\bs{\mu}_{t,\infty} \coloneqq t_\infty^{-1}(t)$, and $t_\infty^{-1}:[0,t_{*,\infty}]\to\{\bs{\mu}_\infty(x):\bs{\mu}_\infty\text{ is the solution of the ODEs} \}$ is the inverse of the function $t_\infty(\cdot)$ restricted to the trajectory of the ODEs.
	\begin{enumerate}
		\item Assume $d_{\max} < \infty$. Then, for any $t < t_{*,\infty}$, we have\label{thm:odesmaininf_1}
		\begin{align*}
		& \frac{I^{(j)}_{d_j,d_{-j},u_j,u_{-j}}(tn)}{n} \xrightarrow{P} i^{(j)}_{d_j,d_{-j},u_j,u_{-j},\infty}(t),\\
		& \frac{T_j(tn)}{n} \xrightarrow{P} \tau_{j,\infty}(t), \text{ and }\frac{T_m(tn)}{n} \xrightarrow{P} \tau_{m,\infty}(t).
		\end{align*}
		\item Assume $d_{\max} \leq \infty$, where $d_{\max} = \infty$ is interpreted as no constraint on the threshold function $K_j(d_j,d_{-j})$. Suppose that the Perron-Frobenius eigenvalue of the non-negative matrix $\bs{J}_{\Ffuncbold_\infty(\cdot)}(\bs{\mu}_{*,\infty})$ is smaller than $1$, where $\bs{J}_{\Ffuncbold_\infty(\cdot)}(\bs{\mu}_{*,\infty})$ is the Jacobian matrix of $\Ffuncbold_\infty(\cdot)$ at $\bs{\mu}_{*,\infty}$. Then, at the natural stopping time of the Markov process of adoption $X^n$, we have:\label{thm:odesmaininf_2}
		\begin{align*}
		&\frac{I^{(j)}_{d_j,d_{-j},u_j,u_{-j}}(\stoptime^n)}{n} \xrightarrow{P} i^{(j)}_{d_j,d_{-j},u_j,u_{-j},\infty}(t_{*,\infty}),\\
		&\frac{T_j(\stoptime^n)}{n} \xrightarrow{P} \tau_{j,\infty}(t_{*,\infty}) \text{ and }\frac{T_m(\stoptime^n)}{n} \xrightarrow{P} \tau_{m,\infty}(t_{*,\infty}).
		\end{align*}
	\end{enumerate}
\end{theorem}
%\begin{proof}
%See Appendix~\ref{proof:odesmaininf}.
%\end{proof}
%\begin{remark}
%	We can replace the convergence in probability with asymptotically almost surely, including an $o(1)$ error term, if $(\mathit{i})$we consider the ODE \eqref{eq:alterode} with the function $\bs{F}$ given by \eqref{eq:trunc_Fjj}-\eqref{eq:trunc_Fj-j} instead of \eqref{eq:meanfield_mujj}-\eqref{eq:meanfield_muj-j}, $(\mathit{ii})$with initial values given by \eqref{eq:ode_ic} in Appendix \eqref{app:odederive}, and $(\mathit{iii})$we assume for some constant $d_{\max} > 0$ we have for all $d_j+d_{-j} > d_{\max}$ either $K_j(d_j,d_{-j}) = d_j+d_{-j}$ or $\alpha_j(d_j,d_{-j}) = 1$.
%\end{remark}
\begin{remark}
	As we commented in Point \ref{point:2} of Section \ref{sec:odeapprox}, the natural stopping time of the original Markov process of adoption $X^n$ is bounded between the same quantities for the truncated versions $X_{L,\delta}^n$ and $X_{U,\delta}^n$. However, this bound does not apply to the whole trajectory and the proof of the first part of Theorem \ref{thm:odesmaininf} is restricted to the case $d_{\max} < \infty$. Nonetheless, we conjecture that the same result holds for $d_{\max} = \infty$.
\end{remark}
\begin{remark}\label{rem:approximatedelta}
	Note that if the Perron-Frobenius eigenvalue of $\bs{J}_{\Ffuncbold_\infty(\cdot)}(\bs{\mu}_{*,\infty})$ is greater than or equal to $1$, Theorem \ref{thm:odesmaininf}, provides an asymptotic lower bound for the final proportion of adopters in the Markov process of adoption.
\end{remark}
Tallying all the (scaled) inactive vertices we can determine the total (scaled) number of inactive vertices in community $j$. This is an immediate corollary of Theorem \ref{thm:odesmaininf}.
\begin{corollary}\label{cor:odesmaininf}
Let $I_n(k)$ denote the total number of inactive vertices at time $k$ of the Markov process of adoption $X^n$. Consider the function $\bs{\Phi} = (\Phi_1,\Phi_2)$ given by the right-hand side of \eqref{eq:meanfield_phij}.
\begin{enumerate}
	\item  Assume $d_{\max} < \infty$. Then for all $t<t_{*,\infty}$, we have
	\begin{align*}
	&\frac{I_n(tn)}{n}\xrightarrow{P}\sum_{u_j +u_{-j} \leq K_j(d_j, d_{-j})} i_{d_j,d_{-j},u_j,u_{-j},\infty}^{(j)}(t)\\
	&\myquad[15]=\beta_1\Phi_1(\mu^{(1\gets 1)}_{t,\infty},\mu^{(1\gets 2)}_{t,\infty}) + \beta_2\Phi_2(\mu^{(2\gets 2)}_{t,\infty},\mu^{(2\gets 1)}_{t,\infty}),\allowdisplaybreaks\\
	&\frac{A_j(tn)}{n}\xrightarrow{P}\lambda_j - 2\tau_{j,\infty}(t)-\sum_{u_j+u_{-j} \leq K_j(d_j,d_{-j})} (d_j-u_j)i^{(j)}_{d_j,d_{-j},u_j,u_{-j},\infty}(t),\allowdisplaybreaks\\
	&\frac{A^{(j)}_m(tn)}{n}\xrightarrow{P} \lambda_m -\tau_{m,\infty}(t)-\sum_{u_j+u_{-j} \leq K_j(d_j,d_{-j})} (d_{-j}-u_{-j})i^{(j)}_{d_j,d_{-j},u_j,u_{-j},\infty}(t).
	\end{align*}
	\item  Assume $d_{\max} \leq \infty$ and suppose that $\bs{\mu}_{*,\infty}$ is a stable equilibrium of ODEs \eqref{eq:alterode}. Then, we have
	\begin{align*}
	&\frac{I_n(\stoptime^n)}{n}\xrightarrow{P}\sum_{u_j +u_{-j} \leq K_j(d_j, d_{-j})} i_{d_j,d_{-j},u_j,u_{-j}}^{(j)}(t_{*,\infty})\\
	&\myquad[15]=\beta_1\Phi_1(\mu^{(1\gets 1)}_{*,\infty},\mu^{(1\gets 2)}_{*,\infty}) + \beta_2\Phi_2(\mu^{(2\gets 2)}_{*,\infty},\mu^{(2\gets 1)}_{*,\infty}).
	\end{align*}
\end{enumerate}
\end{corollary}
\begin{proof}
	The proof follows from Theorem \ref{thm:odesmaininf}, and  the fact that for large enough $d_{\max}$, the proportion of vertices with degree higher than $d_{\max}$ is small.
\end{proof}
This machinery can easily be generalized to any finite number of communities. We conclude this section by presenting the generalization to $k$ communities. The degree regularity conditions need not to be revised for this setting. Note that the only assumptions that we used in the proof of Theorem \ref{thm:odesmaininf} are the degree regularity conditions (given in Definition \ref{def:regcon_deg}). As we mentioned earlier, the graph regularity conditions (given in Definition \ref{def:regcon_graph}) are necessary to get a uniform simple random graph using the configuration model with positive probability. We also comment that the result generalizes trivially to the case when the threshold of vertices is random, assuming the distribution of the threshold depends on the community and degrees of vertices. In particular, assuming $K_j(d_j,d_{-j})$ is random, the same formula works after taking expectation with respect to it.

\begin{theorem}\label{thm:odesmaininf_general}
Assume there are $k$ communities, and size of communities are given by $n_1$, $n_2$, $\cdots$, $n_k$ such that $\sum_i n_i = n$. Assume $\lim_{n\to\infty} n_i/n = \beta_i$ for all $i\in\{1,2,\cdots,k\}$. Suppose the degree regularity conditions hold, and define $\lambda_r \coloneqq \lambda_{r,r}\beta_r$ and $\lambda_m^{(r,s)} \coloneqq \lambda_{r,s} \beta_r$ for all $r,s\in\{1,2,\cdots,k\}$ (note that $\lambda_{r,s} \beta_r = \lambda_{s,r}\beta_s$). Assume there is a constant $0<d_{\max} \leq \infty$, independent of $n$, such that the inactive vertices with degree higher than $d_{\max}$ cannot be activated.
Let $X^n$ denote the Markov process of adoption. Let $\bs{\mu}_\infty(t) = (\mu^{(i,j)}_\infty(t))_{i,j\in\{1,2,\cdots,k\}}$ denote the solution of $k^2$-dimensional ODEs.
\begin{align}
\frac{d\bs{\mu}}{dt} = \Ffuncbold_\infty(\bs{\mu}) - \bs{\mu} \myquad[4]\bs{\mu}(0) = \bs{1}, \myquad[4]\bs{\mu} \in [0,1]^{k^2}.
\end{align}
where the function $\Ffuncbold_\infty(\cdot)$ is given by Mean-Field equations. Let $\bs{\mu}_{*,\infty} = \lim_{s\to\infty}\Ffuncbold_\infty^s(\bs{1})$ to be the closest fixed point of $\Ffuncbold_\infty(\cdot)$ to $\bs{1}$ in sup-norm. Define the function $t_\infty(\cdot)$ as follows:
\begin{align}
t_\infty(\bs{\mu}) \coloneqq \sum_{r=1}^{k} \frac{\lambda_r}{2}\left(1- \left(\mu^{(r,r)}\right)^2\right) + \sum_{\substack{r,s=1\\s\neq r}}^{k} \frac{\lambda^{(r,s)}_m}{2}\left(1- \mu^{(r,s)}\mu^{(s,r)}\right)
\end{align}
and let $t_{*,\infty} \coloneqq t_\infty(\bs{\mu}_{*,\infty})$. Now, the result of Theorem \ref{thm:odesmaininf} and Corollary \ref{cor:odesmaininf} holds by using the following functions:
\begin{align*}
&i^{(j)}_{d_1,d_2,\cdots,d_r,u_1,u_2,\cdots,u_r,\infty}(t) = \prob_{j,m}(d_1,d_2,\cdots,d_k) \prod_{r=1}^{k}Bi(u_r;d_r,1-\mu^{(j,r)}_{t,\infty})\\
&\tau_{r,\infty}(t) =  \frac{\lambda_r}{2}\left(1- \left(\mu^{(r,r)}_{t,\infty}\right)^2\right)\\
&\tau^{(r,s)}_{m,\infty}(t) =  \lambda_m^{(r,s)}\left(1- \mu^{(r,s)}_{t,\infty}\mu^{(s,r)}_{t,\infty}\right)
\end{align*}
where $\bs{\mu}_{t,\infty} \coloneqq t_\infty^{-1}(t)$. Note that our notation is slightly different from the case of two communities, as we use $I^{(j)}_{d_1,d_2,\cdots,d_r,u_1,u_2,\cdots,u_r}$ to denote the number of vertices in community $j$ with $d_r$ half-edges in community $r$, such that $u_r$ of them have been already removed.
\end{theorem}
\begin{proof}
	The proof of the generalized $k$ follows by recycling the proof of $k=2$.
\end{proof}

	\section{Contagion Threshold}\label{sec:contagion}
	Recall that by definition, vertex $i$ in community $j$ with $d_{j,i}^{n}$ neighbors in community $j$ and $d_{m,i}^{n}$ neighbors in the other community is an early adopter with probability $\alpha_j(d_{j,i}^{n},d_{m,i}^{n})$. If we assume $K_j(d_j,d_{-j})\equiv\theta(d_j+d_{-j})$ for some $\theta \in (0,1)$, then the largest value of $\theta$ that results in a cascade (i.e., $O(n)$ vertices becoming active) when a small number of vertices ($o(n)$, often taken to be a constant number) are initially seeded is called the contagion threshold; denote it by $\theta_*$. Morris~\cite{Morris2000} showed that $\theta_*\leq 0.5$ and the upper-bound is loose for many graphs. It's argued that the contagion threshold of the graph family can be calculated by choosing $\alpha_j(d_j,d_{-j})\equiv\alpha$, letting $\alpha \to 0$, and varying $\theta$. In this section, we formalize this intuition and characterize the contagion condition for general threshold functions. % (not just linear threshold model.

Let $\bs{\alpha} = \{\alpha_j(d_j,d_{-j})\}_{j,d_j,d_{-j}}$ represent the seeding strategy. Let's rewrite the function $\Ffuncbold_\infty(\bs{\mu})$ as $\Ffuncbold_\infty(\bs{\alpha},\bs{\mu})$ to emphasis the dependency of function $\Ffuncbold_\infty$ over the seeding strategy. Similarly, we write $\mathcal{U}_\infty(\bs{\alpha})$ to denote the largest connected set containing $\bs{1}$ such that $\forall \bs{\mu}\in \mathcal{U}_\infty(\bs{\alpha})$, $\bs{\mu} \geq \Ffuncbold_\infty(\bs{\alpha},\bs{\mu})$. The question of interest is the final proportion of adopters when the seeding affects only a finite population, i.e., the proportion of early adopters goes to 0 as $n\to\infty$. We provide an answer to this question for general threshold functions in the following theorem.
\begin{theorem}\label{thm:contagion}
Consider an arbitrary sequence $\{\bs{\alpha}_s\}_{s=1}^\infty$ that represent a sequence of non-zero seeding strategies that converges to zero in sup-norm, i.e., $\norm{\bs{\alpha}_s}_\infty \rightarrow0$. Let $\zeta_\infty(\bs{0})$ denote the Perron-Frobenius eigenvalue of the non-negative matrix $\bs{J}_{\Ffuncbold_\infty(\bs{0},\cdot)}(\bs{1})$. If $\zeta_\infty(\bs{0}) < 1$, then $\mathcal{U}_\infty(\bs{0}) = \{\bs{1}\}$, and the final proportion of adopters converges to 0 as $\norm{\bs{\alpha}_s}_\infty \rightarrow 0$. If $\zeta_{*,\infty}(\bs{0}) > 1$,
$\{\bs{1}\}$ is in the interior of $\mathcal{U}_\infty(\bs{0})$, and the final proportion of adopters as $\norm{\bs{\alpha}_s}_\infty \rightarrow 0$ is strictly positive, and we have
\begin{align*}
\bs{\mu}_{*,\infty}(\bs{0}) = \lim_{r\to\infty} \Ffuncbold_\infty^r(\bs{0},\bs{u}) = \lim_{s\to\infty} \bs{\mu}_{*,\infty}(\bs{\alpha}_s) \myquad[1] \forall\bs{u}\in \mathcal{U}_\infty(\bs{0}) \cap \{\bs{x} :  \bs{\mu}_{*,\infty}(\bs{0}) \leq \bs{x} \leq \bs{1} \} \setminus \{\bs{1}\},
\end{align*}
where $\bs{\mu}_{*,\infty}(\bs{\alpha}_s)\coloneqq \lim_{r\to\infty}\Ffuncbold_\infty^r(\bs{\alpha}_s,\bs{1})$, and $\bs{\mu}_{*,\infty}(\bs{0})$ is the closest fixed point of $\Ffuncbold_\infty(\bs{0},\bs{1})$ to $\bs{1}$ other than $\bs{1}$ itself.
\end{theorem}
\begin{proof}
See Appendix \ref{proof:contagion}.
\end{proof}
\begin{remark}
	Note that if Perron-Frobenius eigenvalue of $\bs{J}_{\Ffuncbold_\infty(\bs{0},\cdot)}(\bs{\mu}_{*,\infty}(\bs{0}))$ is smaller than $1$, then we can use the same formulas as in Theorem~\ref{thm:odesmaininf} to characterize the asymptotic proportion of adopters. Otherwise, by Remark \ref{rem:approximatedelta}, we get a lower bound for this.
\end{remark}
Note that the elements of the Jacobian matrix $\bs{J}_{\Ffuncbold_\infty(\bs{0},\cdot)}(\bs{1})$ have a simple form (see Appendix \ref{proof:jacobianmatrix}):
\begin{align*}
	&\left. \frac{\partial \Ffunc_{(j\gets j),\infty}(\bs{0},\bs{\mu})}{\partial \mu^{(j\gets j)}} \right|_{\bs{\mu} = \bs{1}} = \sum_{(d_j,d_{-j}):K_j(d_j,d_{-j}) = 0}  \,(d_j-1) \,\prob_{j*,m}(d_j,d_{-j}),
	\allowdisplaybreaks\\
	&\left. \frac{\partial \Ffunc_{(j\gets j),\infty}(\bs{0},\bs{\mu})}{\partial \mu^{(j\gets -j)}}\right|_{\bs{\mu} = \bs{1}} =
	\sum_{(d_j,d_{-j}):K_j(d_j,d_{-j}) = 0}\,d_{-j} \,\prob_{j*,m}(d_j,d_{-j}),
	\allowdisplaybreaks\\
	&\left. \frac{\partial \Ffunc_{(j\gets -j),\infty}(\bs{0},\bs{\mu})}{\partial \mu^{(-j\gets -j)}} \right|_{\bs{\mu} = \bs{1}} =
	\sum_{(d_j,d_{-j}):K_j(d_j,d_{-j}) = 0}\,d_{-j} \,\prob_{-j,m*}(d_{-j},d_j),
	\allowdisplaybreaks\\
	&\left.\frac{\partial \Ffunc_{(j\gets -j),\infty}(\bs{0},\bs{\mu})}{\partial \mu^{(-j\gets j)}}\right|_{\bs{\mu} = \bs{1}} =
	\sum_{(d_j,d_{-j}):K_j(d_j,d_{-j}) = 0}\,(d_j-1) \,\prob_{-j,m*}(d_{-j},d_j),
\end{align*}
There is an interesting intuition behind the Perron-Frobenius eigenvalue $\zeta_\infty(\bs{0})$ of $\bs{J}_{\Ffuncbold_\infty(\bs{0},\cdot)}(\bs{1})$ and the contagion threshold. Let $\mathcal{P}_n$ denote the set of vertices that need only one active neighbor to adopt the new technology. If $\zeta_\infty(\bs{0}) < 1$, then after random pairing of half-edges, $\mathcal{P}_n$ consist of many small components with high probability. On the other hand, $\zeta_\infty(\bs{0}) > 1$ implies that after random pairing of half-edges, $\mathcal{P}_n$ has one giant component with high probability. Hence, if $\zeta_\infty(\bs{0}) > 1$, then activating one of the vertices in the giant component of $\mathcal{P}_n$ will active a large proportion of the population. This has also been reported in \cite{Lelarge2012} for the case of one community.

The discussion on contagion can also be generalized to $k$ communities with the same statement as in Theorem \ref{thm:contagion}. As we mentioned before, this also generalizes to the case of random threshold.
\begin{remark}
	In case of one community, the above analysis yields the same characterization in terms of the derivative of $\Ffunc_\infty(\bs{0},\mu)$ at $\mu = 1$. In particular, contagion happens if $\sum_{d:K(d) = 0} \,(d-1)\prob_{*}(d) > 1$ and does not happen if $\sum_{d:K(d) = 0} \,(d-1)\prob_{*}(d) < 1$, where $\prob_{*}(\cdot)$ is the size-biased distribution of the asymptotic degree distribution $\prob(\cdot)$. This is the same criteria as in \cite[Cascade condition (7)]{Lelarge2012} for the case of one community.
\end{remark}

	\section{Poisson Degree Distributions}\label{sec:poissdeg}
	
We will now specialize our results to Poisson degree distributions. An Erd\H{o}s-R\'enyi random graph is an example of a graph family that asymptotically yields a Poisson degree distribution. The two community stochastic block model is then the appropriate generalization of the Erd\H{o}s-R\'enyi random graph that will asymptotically produce Poisson degree distributions within the community and across the communities. We will show in the following results that under some symmetry assumptions for the threshold and the advertising strategy, the solution of ODEs \eqref{eq:alterode} simplifies considerably. In the case of Poisson degree distribution, we assume
\begin{align*}
\prob_{j,m}(d_j,d_{-j}) &= \mathrm{e}^{-\lambda_{j,j}}\frac{(\lambda_{j,j})^{d_j}}{d_j !}\times\mathrm{e}^{-\lambda_{j,m}}\frac{(\lambda_{j,m})^{d_{-j}}}{d_{-j} !}
\end{align*}
Note that $\prob_{j*,m}(d_j,d_{-j}) = \prob_{j,m}(d_j-1,d_{-j})$ and $\prob_{j,m*}(d_j,d_{-j}) = \prob_{j,m}(d_j,d_{-j}-1)$.

\begin{theorem} \label{thm:poissreduc}
Assume that the threshold of each vertex depends on its community and the total number of its neighbors, i.e. $K_j(d_j,d_{-j}) = K_j(d_j + d_{-j})$. Moreover, assume the advertisement strategy is based on the community affiliation and the total number of neighbors, i.e. $\alpha_j(d_j,d_{-j}) = \alpha_j(d_j + d_{-j})$. Now, if the asymptotic degree distributions are Poisson with parameters $\lambda_{1,1}$, $\lambda_{1,m}$, $\lambda_{2,m}$, and $\lambda_{2,2}$, then the solution of ODEs \eqref{eq:alterode} with the function $\Ffuncbold_\infty(\cdot)$ given by the right-hand side of \eqref{eq:meanfield_mujj}-\eqref{eq:meanfield_muj-j} simplifies as follows:
$\mu_\infty^{(1\gets 1)}(t) = \mu_\infty^{(2\gets 2)}(t)$ and $\mu_\infty^{(2\gets 1)}(t) = \mu_\infty^{(1\gets 2)}(t)$ for all $t\geq0$; that is to say, the dimension of the differential equations reduces to 2.
\end{theorem}
\begin{proof}
See Appendix~\ref{proof:poissreduc}.
\end{proof}

The next theorem concerns general distributions.
\begin{theorem} \label{thm:symmreduc}
Assume both the advertisement strategy and the threshold function are symmetric in the following sense: $\alpha_j(d_j,d_{-j}) = \alpha_{-j}(d_{-j},d_{j})$  and $K_{j}(d_{j},d_{-j}) = K_{-j}(d_{-j},d_{j})$ for all $d_j, d_{-j} \geq 0$ and $j\in\{1,2\}$. Also assume that the asymptotic degree distribution in both communities are the same, i.e. $\prob_{1,m}=\prob_{2,m}$. Then, $\mu^{(1\gets 1)}(x) = \mu^{(2\gets 1)}(x)$ and $\mu^{(2\gets 2)}(x) = \mu^{(1\gets 2)}(x)$ for all $x\geq0$; that is to say, the dimension of the differential equations reduces to 2.
\end{theorem}
\begin{proof}
See Appendix~\ref{proof:symmreduc}.
\end{proof}
Similar generalization holds for the case of $k$ communities: given similar assumptions to Theorem \ref{thm:poissreduc} or Theorem \ref{thm:symmreduc}, in the case of $k$ communities, the dimension of ODEs reduces to $k$. An immediate corollary is the following which asserts that if both the assumptions hold, then the dimension reduces to $1$. This is also true for general $k$.
\begin{corollary}\label{cor:poissym}
Assume the assumptions of Theorems \ref{thm:poissreduc} and Theorems \ref{thm:symmreduc} hold, then the dimension of ODEs \eqref{eq:alterode} reduces to $1$. In particular, the resulted system of ODE is the same as if there was only one community with asymptotic degree distribution given by Poisson$(\lambda_{1,1} + \lambda_{1,m})$.
\end{corollary}

It is interesting to note that given assumptions of Theorems \ref{thm:poissreduc} and Theorems \ref{thm:symmreduc}, the contagion threshold is the same as if there was only one community. The derivation of contagion threshold then matches the ones presented in~\cite{Amini2010,Lelarge2012} for the case of one community.

	\section{Numerical Investigation}\label{sec:numres}
	We present some numerical results using the analysis presented above. The main point is to show how the community structure impacts seeding strategies. A natural question to ask is the following: what is the best seeding strategy given a budget constraint? In this section, we also formalize this question and provide a partial answer to it using a gradient-based heuristic algorithm.

By Theorem \ref{thm:odesmaininf}, we can approximate the state of the Markov process of adoption at its natural stopping time if $\bs{\mu}_{*,\infty}$ is a stable equilibrium point of ODEs \eqref{eq:alterode} where the function $\Ffuncbold(\cdot)$ is replaced with $\Ffuncbold_\infty(\cdot)$. Then, by Corollary \ref{cor:odesmaininf}, our strategy is to pick $\bs{\alpha}_*$ that minimizes $\beta_1 \Phi_1(\bs{\mu}_{*,\infty})+\beta_2\Phi_2(\bs{\mu}_{*,\infty})$.

Specifically, consider the following budget constraint which constraints the expected proportion of early adopters:
\begin{align*}
\sum_{j,d_j,d_{-j}} \beta_j \, \prob_{j,m}(d_j,d_{-j}) \, \alpha_{j}(d_j,d_{-j}) = \mathscr{B},
\end{align*}
where $\mathscr{B} > 0$ is the total available budget. We now formulate the ``optimum seeding strategy'' as follows:
\begin{align*}
\inf_{\bs{\alpha}} &&&\beta_1\Phi_1(\bs{\alpha},\bs{\mu}_{*,\infty})+\beta_2\Phi_2(\bs{\alpha},\bs{\mu}_{*,\infty})\\
\text{subject to} &&& \text{$(i)$ Budget constraint:}\sum_{j,d_j,d_{-j}} \beta_j \, \prob_{j,m}(d_j,d_{-j}) \, \alpha_{j}(d_j,d_{-j}) = \mathscr{B},\\
&&& \text{$(ii)$ Definition of $\bs{\mu}_{*,\infty}$:}~ \bs{\mu}_{*,\infty} = \argmin_{\bs{u}: \Ffuncbold_\infty(\bs{\alpha},\bs{u}) = \bs{u}} \norm{\bs{u} - \bs{1}}_\infty,\\
&&& \text{$(ii)'$ Definition of $\bs{\mu}_{*,\infty}$:}~ \bs{\mu}_{*,\infty} = \lim_{s\to\infty}\Ffuncbold_\infty^s(\bs{\alpha},\bs{1}).
\end{align*}
Note that by Corollary \ref{cor:Fprop_sol}, the conditions $(ii)$ and $(ii)'$ are equivalent. As the constraints are highly nonlinear, we take a heuristic approach to provide a partial answer to this optimization problem. We comment that if $\bs{\mu}_{*,\infty}$ is not a stable equilibrium point, then the process may not be in the vicinity of scaled-time close to $t_{*,\infty}$ as we don't know what happens beyond this point (see also Remark~\ref{rem:approximatedelta}). However, the above formulation is heuristically our best recourse to maximize the contagion.

For the sake of simplicity, let us assume the degrees are uniformly bounded. This assumption is justified by noting that the total number of half-edges associated with high degree vertices is small. Note that the fixed point $\bs{\mu}_{*,\infty}$ depends on the seeding strategy $\bs{\alpha}$; abusing notation, we denote this by $\bs{\mu}_{*,\infty}(\bs{\alpha})$. Recall that $\bs{\mu}_{*,\infty}(\bs{\alpha}) = \bs{F}(\bs{\alpha},\bs{\mu}_{*,\infty}(\bs{\alpha}))$. Using the chain rule, we have
\begin{align*}
\bs{J}_{\bs{\mu}_{*,\infty}(\cdot)}(\bs{\alpha}) = \bs{J}_{\Ffuncbold_\infty(\cdot,\bs{\mu}_{*,\infty}(\bs{\alpha}))}(\bs{\alpha}) +\bs{J}_{\Ffuncbold_\infty(\bs{\alpha},\cdot)}(\bs{\mu}_{*,\infty}(\bs{\alpha})) \,\bs{J}_{\bs{\mu}_{*,\infty}(\cdot)}(\bs{\alpha})
\end{align*}
where $\bs{J}_{\bs{\mu}_{*,\infty}(\cdot)}(\bs{\alpha})$ is the Jacobian matrix of $\bs{\mu}_{*,\infty}(\cdot)$ at $\bs{\alpha}$, $\bs{J}_{\Ffuncbold_\infty(\cdot,\bs{\mu})}(\bs{\alpha})$ is the Jacobian matrix of $\Ffuncbold_\infty(\cdot,\bs{\mu})$ at $\bs{\alpha}$, and $\bs{J}_{\Ffuncbold_\infty(\bs{\alpha},\cdot)}(\bs{\mu})$ is the Jacobian matrix of $\Ffuncbold_\infty(\bs{\alpha},\cdot)$ at $\bs{\mu}$. Hence, assuming the inverse of $\bs{I} - \bs{J}_{\Ffuncbold_\infty(\bs{\alpha},\cdot)}(\bs{\mu}_{*,\infty}(\bs{\alpha}))$ exists (otherwise, we can use the pseudo-inverse), we have
\begin{align}\label{eq:jacobmu*}
\bs{J}_{\bs{\mu}_{*,\infty}(\cdot)}(\bs{\alpha}) = \left(\bs{I} - \bs{J}_{\Ffuncbold_\infty(\bs{\alpha},\cdot)}(\bs{\mu}_{*,\infty}(\bs{\alpha}))\right)^{-1}
\bs{J}_{\Ffuncbold_\infty(\cdot,\bs{\mu}_{*,\infty}(\bs{\alpha}))}(\bs{\alpha}) ,
\end{align}
where $\bs{I}$ is the identity matrix.

The heuristic seeding algorithm is an iterative algorithm that has two stages. The logic behind the algorithm is simple: at $\mathit{Stage~1}$, the algorithm tries to find the best direction for updating the seeding strategy while keeping the budget constraint, and at $\mathit{Stage~2}$, the algorithm validates the choice of the direction.

Fix some $\xi_0 > 0$, and pick $\xi_0 < \xi < 1$ arbitrary. Also pick $\bs{\alpha}$ arbitrary such that it satisfies the budget constraint. The heuristic algorithm is given as follows:

\medskip

\noindent$\mathit{Stage~1}.$ Numerically solve the following linear optimization problem:
\begin{align*}
\inf_{\Delta\bs{\alpha} = [\Delta\bs{\alpha}_{j}(d_j,d_{-j})]_{j,d_j,d_{-j}}} &&&\big\langle\Delta\bs{\alpha}\, ,\,  \sum_j \bs{J}_{\Phi_j(\cdot,\bs{\mu}_{*,\infty}(\bs{\alpha}))}(\bs{\alpha}) +\bs{J}_{\Phi_j(\bs{\alpha},\cdot)}(\bs{\mu}_{*,\infty}(\bs{\alpha}))\,\bs{J}_{\bs{\mu}_{*,\infty}(\cdot)}(\bs{\alpha})\big\rangle_F
\\
\text{subject to} &&&\sum_{j,d_j,d_{-j}} \beta_j \, \prob_{j,m}(d_j,d_{-j}) \, (\alpha_{j}(d_j,d_{-j}) +  \xi \, \Delta\bs{\alpha}_{j}(d_j,d_{-j}))= \mathscr{B},\\
 &&&\alpha_{j}(d_j,d_{-j}) +  \xi \,\Delta\bs{\alpha}_{j}(d_j,d_{-j}) \in [0,1]\text{ for all $j,d_j,d_{-j}$}\\
 &&&\norm{\Delta\bs{\alpha}}_F = 1
\end{align*}
where $\langle A,B \rangle_F\coloneqq \sum{a_{i,j}b_{i,j}}$ is the Frobenius inner product, $\bs{J}_{\Phi_j(\cdot,\bs{\mu})}(\bs{\alpha})$ is the Jacobian matrix of $\Phi_j(\cdot,\bs{\mu})$ at $\bs{\alpha}$, $\bs{J}_{\Phi_j(\bs{\alpha},\cdot)}(\bs{\mu})$ is the Jacobian matrix of $\Phi_j(\bs{\alpha},\cdot)$ at $\bs{\mu}$, and $\norm{\cdot}_F$ is the Frobenius norm.

\smallskip

\noindent$\mathit{Stage~2}.$ Let $\bs{\alpha}_{\text{new}} =[ \alpha_{j}(d_j,d_{-j}) +  \xi \,\Delta\bs{\alpha}_{j}(d_j,d_{-j})]_{j,d_j,d_{-j}}$. If
\begin{align*}
&\beta_1\Phi_1(\bs{\alpha},\bs{\mu}_{*,\infty}(\bs{\alpha}))+\beta_2\Phi_2(\bs{\alpha},\bs{\mu}_{*,\infty}(\bs{\alpha})) >\\ &\myquad[8]\beta_1\Phi_1(\bs{\alpha}_{\text{new}},\bs{\mu}_{*,\infty}(\bs{\alpha}_{\text{new}}))+\beta_2\Phi_2(\bs{\alpha}_{\text{new}},\bs{\mu}_{*,\infty}(\bs{\alpha}_{\text{new}})),
\end{align*}
then update $\bs{\alpha}\leftarrow\bs{\alpha}_{\text{new}}$ and go to $\mathit{Stage~1}.$ Otherwise, update $\xi\leftarrow\xi/2$. If $\xi<\xi_0$ terminate the algorithm, otherwise go to $\mathit{Stage~1}.$

\medskip

Next, we compare different seeding strategies using Theorem \ref{thm:odesmaininf}, Corollary \ref{cor:odesmaininf}, and Corollary \ref{cor:Fprop_sol}. We focus on Poisson degree distributions, owing to analytical simplifications and the fact there are only three parameters to tune. Moreover, for simplicity we assume that the threshold functions are given by $K_j(d_j,d_{-j}) = \theta \times (d_j+d_{-j})-1$ where $\theta = 1/4$. The choice of $1/4$ is motivated by \cite[Figure 2]{Lelarge2012} that shows the contagion threshold of sparse Erd\H{o}s-R\'enyi random graph (with parameter $\lambda$) is below $1/4$. Henceforth, we assume $\beta_1 = \beta_2$.

The vertices that are seeded by the advertisers are early adopters. A few strategies that we consider are:
(1) Random seeding: first, we assume the advertiser does not even know about the existence of two communities. This scenario is named as {\it global seeding}. Second, we assume the advertiser knows the community structure and decides to seed just asymmetrically in the two communities. This advertisement strategy is denoted by {\it local seeding}. (2) Degree-targeted seeding: the advertiser knows the degree distribution of the network and the identity of the vertices that possess a certain degree but does not know the underlying connectivity structure.

\begin{figure}[h]
\centering
\includegraphics[width=0.85\textwidth]{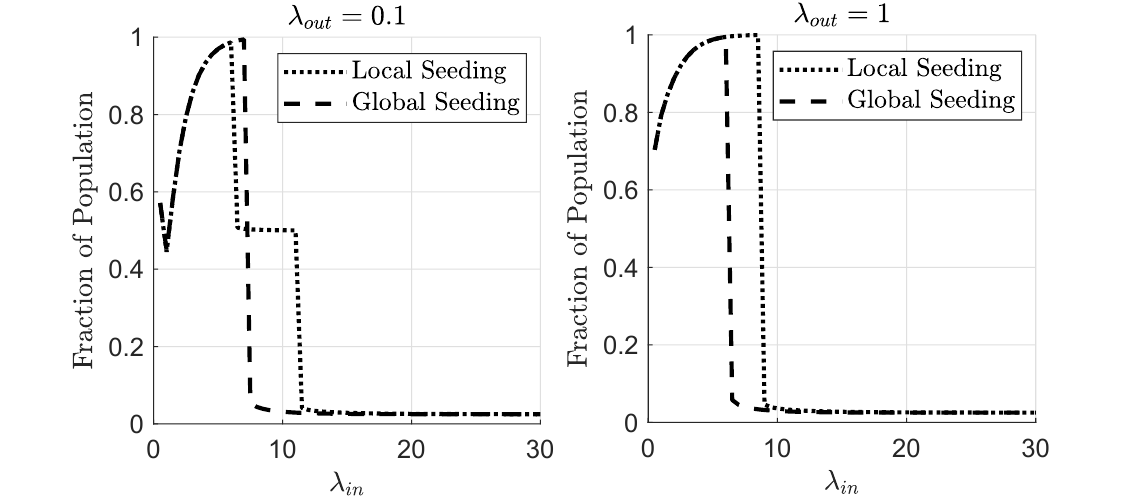}
\caption{Random seeding strategy on symmetric communities. The ratio of of early adopters is 2.5$\%$ of the population. $\lambda_{\text{in}}$ denotes the in-community connectivity, and $\lambda_{\text{out}}$ denotes the out-community connectivity.} \label{fig:randseedpoisssymm}
\end{figure}

\begin{figure}[h]
\centering
\includegraphics[width=0.8\textwidth]{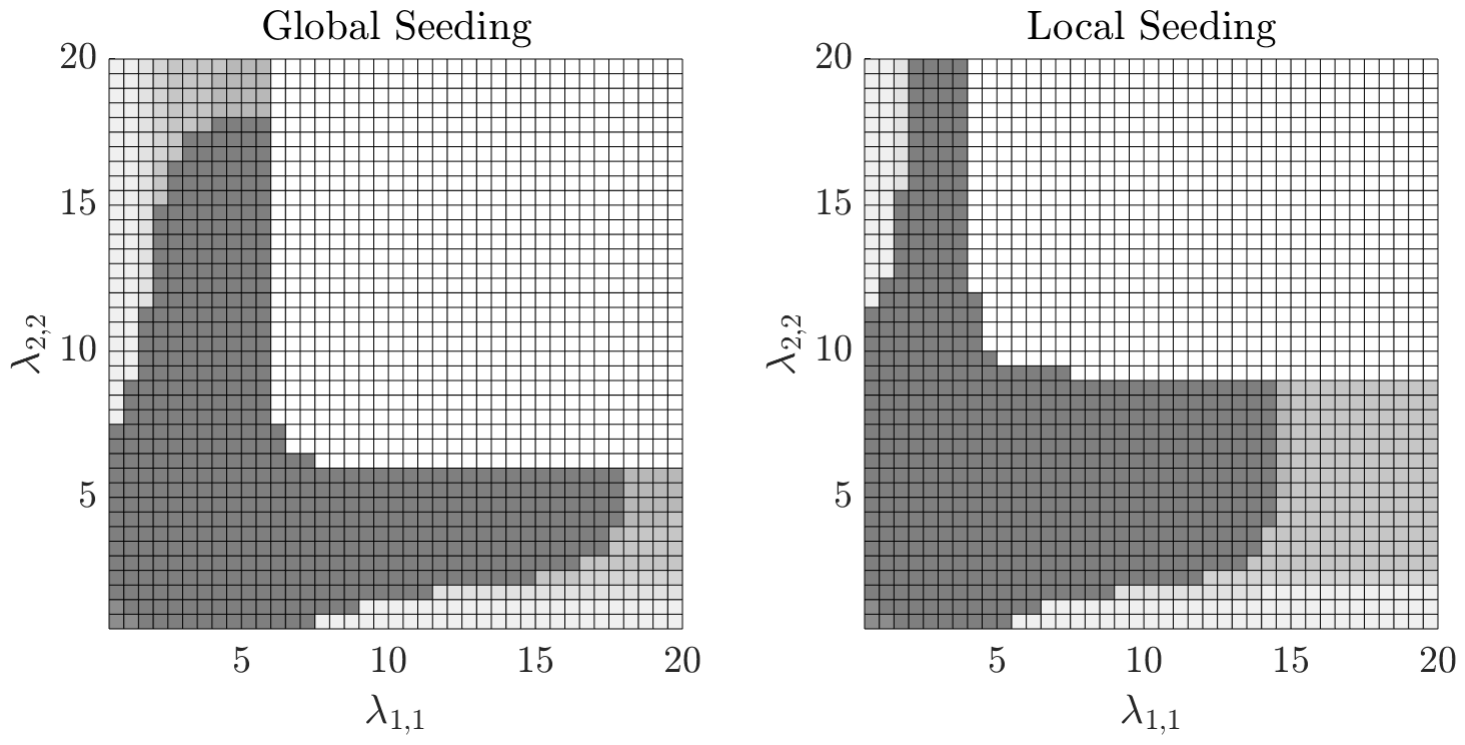}
\caption{Random seeding strategy on asymmetric communities. The ratio of of early adopters is 2.5$\%$ of the population. $\lambda_{1,1}$ and $\lambda_{2,2}$ are parameters of the in-community distributions. $\lambda_{1,2} = \lambda_{2,1} = \lambda_{\text{out}} = 1$ in both cases. Intensity of grayscale indicates the final proportion of adopters: darkest color corresponds to $1$ and lightest color corresponds to $0$.}\label{fig:randseedpoissasymm}
\end{figure}

In Figure \ref{fig:randseedpoisssymm}, we assume that both the in-community and the out-community degree distributions are Poisson distributions with parameters $\lambda_{1,1} = \lambda_{2,2} = \lambda_{\text{in}}$ and $\lambda_{1,2} = \lambda_{2,1} = \lambda_{\text{out}}$. The figure suggests that if the communities are symmetric, and if they are well-connected ($\lambda_{\text{out}} = 1$), then the best strategy is to use the whole budget in one community.
In Figure \ref{fig:randseedpoissasymm} we consider the general case where distributions can have different parameters in the two communities, i.e., $\lambda_{1,1}$ and $\lambda_{2,2}$ need not be equal. We also assume $\lambda_{1,2} = \lambda_{2,1} = \lambda_{\text{out}} = 1$. In this case, the community structure dramatically changes the cascade potential: there are scenarios where global seeding can cause a cascade while local seeding won't, and {\it vice-versa}.

%Next, we consider degree-targeted seeding in Figures \ref{fig:degseedpoiss} and \ref{fig:degseedpoisscomp}.  In general, high-degree vertices can potentially stop a cascade if they are not early adopters; hence it might make sense to seed these vertices in each community. We will consider the following cases: the budget is spread equally in both communities, denoted by $(0.5,0.5)$; the budget is concentrated in community 1, denoted by $(1,0)$; and $25\%$ of budget is in community 1, denoted by $(0.25,0.75)$. The outer-community connectivity is given by $\lambda_{1,2} = \lambda_{2,1} = \lambda_{\text{out}} = 1$. The main observation is the dramatic difference in the proportion of final adopters based on how asymmetric the targeting is. Additionally, a higher inter-community connectivity leads to a bigger cascade.
%Also note that seeding vertices with the highest degree gives better result than random seeding. In Figure \ref{fig:degseedpoisscomp}, we compare the highest degree strategy with the heuristic algorithm proposed in the beginning of the section, using the same parameters as in Figure \ref{fig:degseedpoiss}. The results are dramatically different. We also illustrate the resulting seeding distribution, $\alpha_1(d_1,d_2)$ and $\alpha_2(d_2,d_1)$, for $\lambda_{1,1} = 18$ and  $\lambda_{2,2}=10.5$ in Figure \ref{fig:heurisseed}. These two figures highlights the importance of community structure on the optimum seeding strategy.

\begin{figure}[h]
	\centering
	\includegraphics[width=0.8\textwidth]{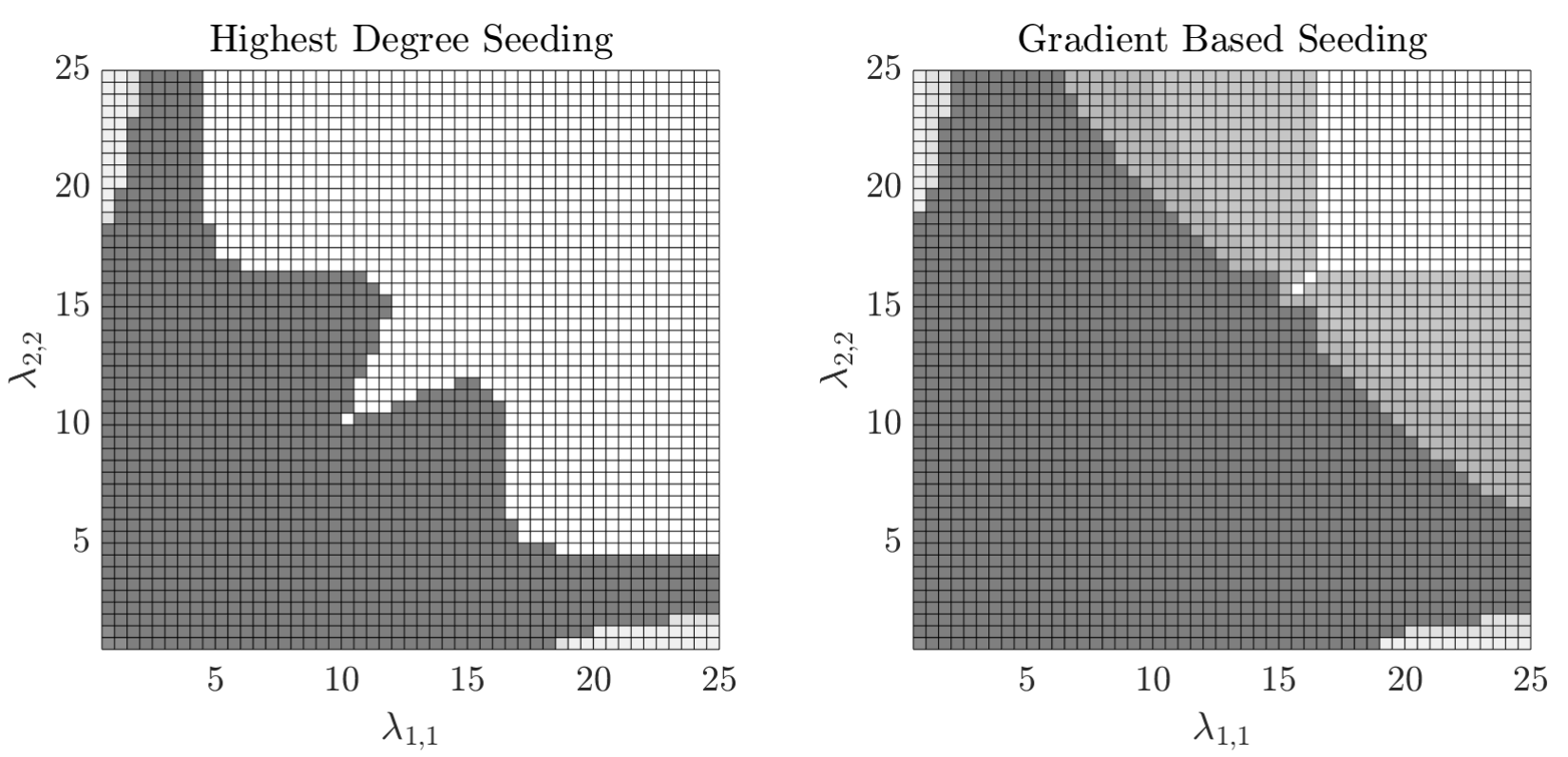}
	\caption{Highest degree seeding strategy versus proposed heuristic seeding strategy on asymmetric communities. The ratio of of early adopters is 2.5$\%$ of the population. $\lambda_{1,1}$ and $\lambda_{2,2}$ are parameters of the in-community distributions. $\lambda_{1,2} = \lambda_{2,1} = \lambda_{\text{out}} = 1$ in all three cases. Intensity of grayscale indicates the final proportion of adopters: darkest color corresponds to $1$ and lightest color corresponds to $0$.}\label{fig:degseedpoisscomp}
\end{figure}

\begin{figure}[h]
	\centering
	\includegraphics[width=0.85\textwidth]{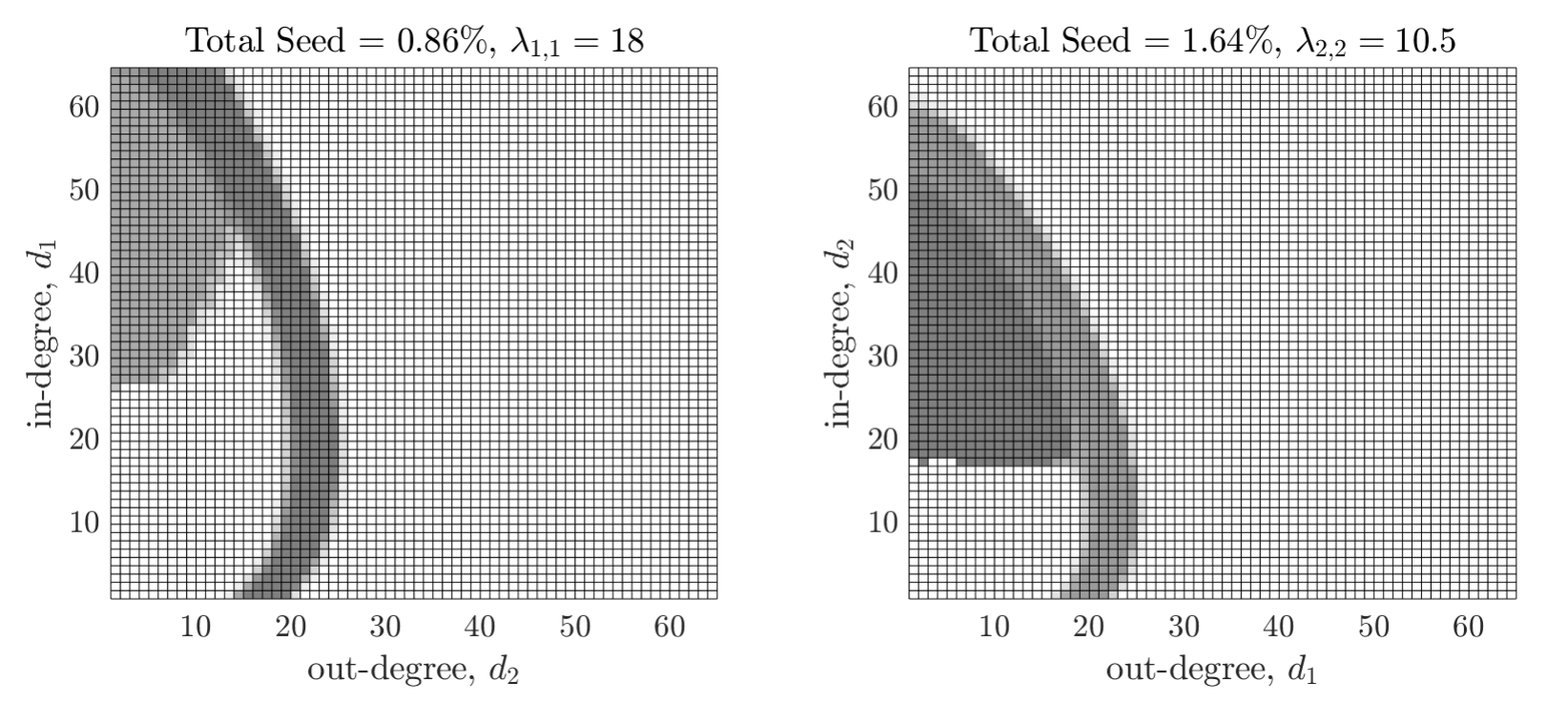}
	\caption{Distribution of $\alpha_1(d_1,d_2)$ (left) and $\alpha_2(d_2,d_1)$ (right) in the proposed heuristic seeding algorithm for $\lambda_{1,1} = 18$,  $\lambda_{2,2}=10.5$ and $\lambda_{\text{out}} = 1$. The ratio of of early adopters is 2.5$\%$ of the population. Intensity of grayscale indicates the value of $\alpha_j(d_j,d_{-j})$: darkest color corresponds to $1$ and lightest color corresponds to $0$.}\label{fig:heurisseed}
\end{figure}

Next, we consider degree-targeted seeding in Figure \ref{fig:degseedpoisscomp}.  In general, high-degree vertices can potentially stop a cascade if they are not early adopters; hence it might make sense to seed these vertices in each community. In Figure \ref{fig:degseedpoisscomp}, we compare the highest degree seeding strategy with the heuristic algorithm proposed at the beginning of the section. The outer-community connectivity is given by $\lambda_{1,2} = \lambda_{2,1} = \lambda_{\text{out}} = 1$, and the heuristic algorithm is initialized with global seeding strategy. The results are dramatically different. We also illustrate the resulting seeding distribution, $\alpha_1(d_1,d_2)$ and $\alpha_2(d_2,d_1)$, for $\lambda_{1,1} = 18$ and  $\lambda_{2,2}=10.5$ in Figure \ref{fig:heurisseed}. These two figures highlight the importance of community structure in the optimum seeding strategy.

%\begin{figure}[H]
%\centering
%\includegraphics[width=1\textwidth]{Figures/AttackHigh-Lou=1.pdf}
%\caption{Highest degree seeding strategy on asymmetric communities.The ratio of of early adopters is 2.5$\%$ of the population. $\lambda_{1,1}$ and $\lambda_{2,2}$ are parameters of the in-community distributions. $\lambda_{1,2} = \lambda_{2,1} = \lambda_{\text{out}} = 1$ in all three cases. Intensity of grayscale indicates the final proportion of adopters.}\label{fig:degseedpoiss}
%\end{figure}

\begin{figure}
	\centering
	\includegraphics[width=0.7\textwidth]{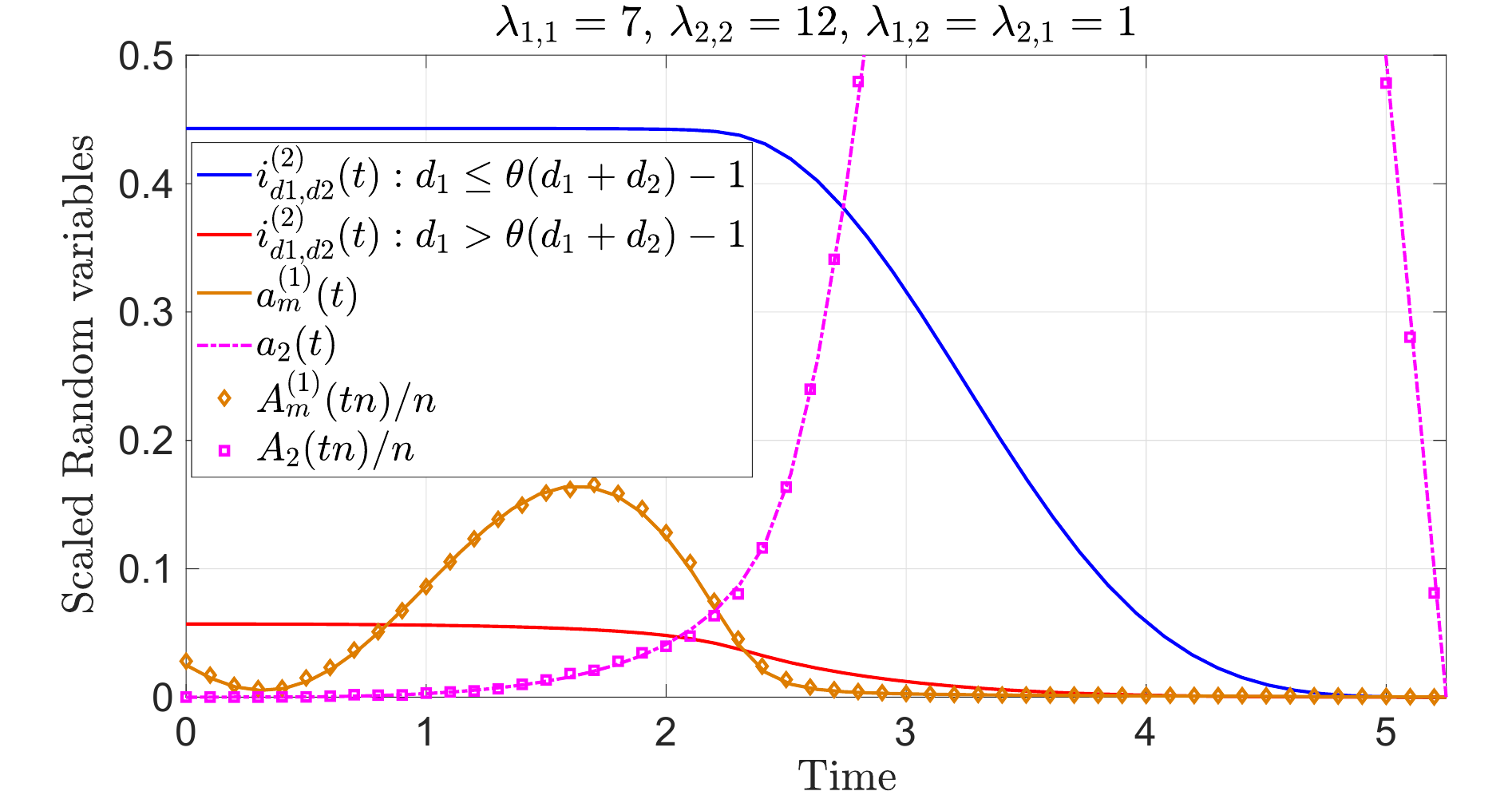}
	\caption{Evolution of cascade for $\lambda_{1,1} = 7$,  $\lambda_{2,2} = 12$, and $\lambda_{1,2} = \lambda_{2,1} = 1$ using local seeding strategy. The total proportion of early adopters is 2.5$\%$.\label{fig:odeloca}}
\end{figure}

Next, we discuss the evolution of cascade using Theorem \ref{thm:odesmaininf} and Corollary \ref{cor:odesmaininf}. Figure \ref{fig:odeloca} illustrates the evolution of active half-edges and inactive vertices in the second community for $\lambda_{1,1} = 7$,  $\lambda_{2,2} = 12$, and $\lambda_{1,2} = \lambda_{2,1} = 1$ when the seeding strategy is to put the whole budget in the first community. Figure \ref{fig:randseedpoissasymm} suggests that global seeding strategy will not result in any cascade. On the other hand, a global cascade emerges following local seeding strategy: it develops in the first community and then moves to the next community; this happens when the inactive vertices in community $2$ with $d_1 > \theta\times(d_1+d_2)-1$ become active, causing a cascade in the second community. In this figure, we also present the total number of active half-edges in the second community for the Markov process of adoption, for $n=20000$, to numerically validate the result of Theorem \ref{thm:odesmaininf}.

\begin{figure}
\centering
\includegraphics[width=0.75\textwidth]{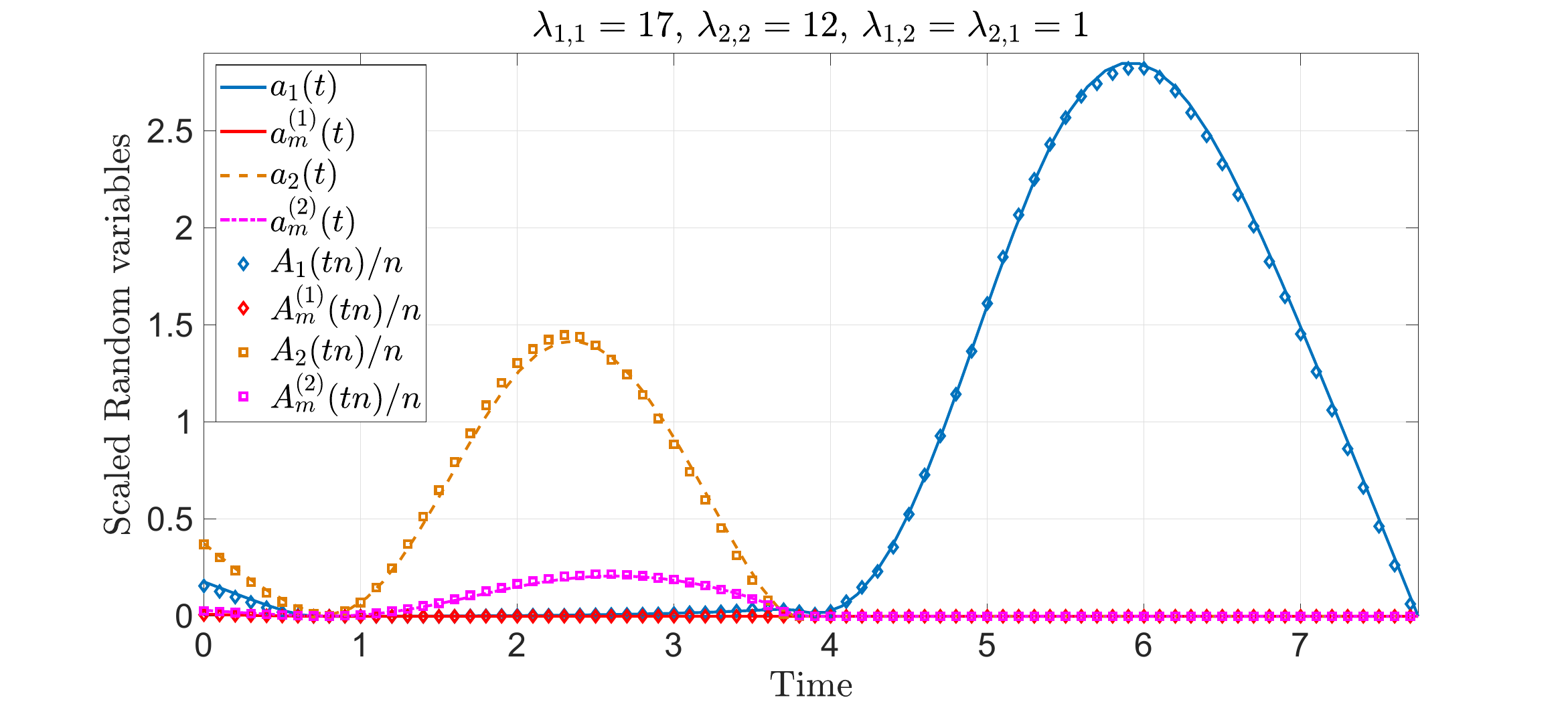}
\caption{Evolution of cascade, $\lambda_{\text{out}} = 1$, $\lambda_{in,1} = 17$ and  $\lambda_{in,2} = 12$;  proportion of early adopters is 2.5$\%$, x axis is time and y axis is the quantity of corresponding scaled variables. The seeding strategy is given by $(0.25,0.75)$. \label{fig:asymseed}}
\end{figure}

Finally, there are scenarios where neither global nor local seeding strategy can cause a cascade. Figure \ref{fig:asymseed} illustrate the evolution of active half-edges for $\lambda_{1,1} = 17$, $\lambda_{2,2} = 12$ and $\lambda_{1,2} = \lambda_{2,1} = 1$, when $25\%$ of budget is used in community $1$ and $75\%$ of the budget is used in community $2$, seeding vertices with the highest degree. Active half-edges in both communities get close to zero; nevertheless, a cascade happens in the second community. This cascade then moves to the first community, and almost all vertices adopt the new technology. This example illustrates the importance of active half-edges $a_m^{(2)}(t)$ in triggering a cascade in the first community. We also include the total number of active half-edges for the Markov process of adoption, for $n=20000$, to numerically validate the result of Theorem \ref{thm:odesmaininf}.

	\section{Open Problems}\label{sec:openprob}
	We close our paper by presenting some open problems:
\begin{enumerate}[wide]
\item What if $\zeta_{*,\infty} \geq 1$? Following the same intuition we presented in Section \ref{sec:contagion}, let $\mathcal{P}_n(k)$ denote the set of pivotal players at time $k$ of the process, i.e., set of remaining inactive vertices at time $k$ that only need one further active neighbor to adopt the new technology. A closer look at the condition $\zeta_{*,\infty} < 1$, suggest that after random pairing of half-edges, $\mathcal{P}_n(\floor{t_* n})$ consist of many small components, with high probability. This is why the cascade cannot grow much further. On the other-hand, if $\zeta_{*,\infty} > 1$, then after random pairing of half-edges, $\mathcal{P}_n(\floor{t_* n})$ has one giant component, with high probability. Hence, it is natural to expect that the cascade won't stop here. However, making this argument rigorous, needs much more work. Intuitively, the cascade process grows up to the closest fixed point of $\mathcal{U}_\infty$ to $\bs{1}$ for which the eigenvalue of the Jacobian matrix is bounded by $1$. We leave this as a conjecture for future work.

\item How many fixed points does $\Ffuncbold_\infty(\cdot)$ have? In a related work, Balogh and Pittel~\cite{Balogh2007} shown that the answer is two for regular random graphs (given some additional conditions). \label{openq:2}

\item What is the trajectory of the contagion process? The evolution of the Markov process of adoption is not the same as the evolution of cascades on the network. Since the random graph model converges locally weakly to $\text{GWMT}_*$ defined in Section \ref{sec:meanfield}, the trajectory of the contagion process on the random graph is related to the trajectory of the contagion process on $\text{GWMT}_*$. Moreover, it is easy to see that the evolution of cascades on $\text{GWMT}_*$ is given by the iterations of function $\Ffuncbold_\infty(\cdot)$ starting from $\bs{1}$ (See~\cite{Como2019} for a related discussion). This intuitive argument also justifies the connection between $\lim_{s\to\infty} \Ffuncbold_\infty^s(\bs{1})$ and the equilibrium of the ODEs \eqref{eq:alterode}. However, a rigorous proof is still missing in the literature.

\item Finally, what is the optimum seeding strategy? We proposed a heuristic seeding strategy in Section \ref{sec:numres} that performs well in different scenarios (Figure \ref{fig:degseedpoisscomp}). However, we have no proof that this seeding strategy is optimum, nor do we have a characterization of its sub-optimality if it is not. This question is also related to the question \ref{openq:2} above.
\end{enumerate}

	\appendix
	\section{When Greedy Maximization Is Bad}\label{app:greed}
	
We now give an example showing how the results in \cite{Kempe2003,Mossel2010} break if $\theta_v$ is assumed to be fixed. We build a network as follows: start from a $3n\times 3n$ torus, i.e. vertex $(i,j)$ with $1\leq i,j\leq 3n$ has four neighbors: $(i+1,j),(i-1,j),(i,j+1),(i,j-1)$ where operations are performed modulo $3n$. Now for each $1\leq j\leq 3n$, and $0\leq k\leq n-1$, we add a vertex $v(j,k)$ connected to the vertices of the torus $(3k+1,j)$, $(3k+2,j)$ and $(3k+3,j)$. Finally each of these vertices $v(j,k)$ are part of a cycle of size $K\geq 3$ with no other common point with the rest of the graph except through $v(j,k)$. In summary, we have $9n^2$ vertices on the torus, and $3n^2$ disjoint cycles of size $K$ which are connected to the torus only through the vertices $v(j,k)$. There is a total of $9n^2+3n^2 K$ vertices. Note that the degree of a vertex on the torus is $5$ (4 neighbors on the torus and 1 on a cycle) as well as for the vertices $v(j,k)$. We take $\theta=2/5$ so that a vertex of degree $d$ becomes active as soon as $\theta d$ of its neighbors are active. In particular a vertex on the torus or a $v(j,k)$ needs only $2$ active neighbors to become active. Moreover, activating a vertex $v(j,k)$ will activate all the $K$ vertices on the cycle. Because of this, it is easy to see that any greedy algorithm with budget $b\leq 3n^2$ will only activate the vertices $v(j,k)$. Note however that by activating the set of vertices on the torus: $(1,1), (1,2),\dots ,(1,3n)$ and $(2,1)$ will result in the global activation of the network. Hence for any $3n+1\leq b \leq 3n^2$, we can find a set activating the $9n^2+3n^2K$ vertices of the networks, whereas the greedy algorithm only activates $Kb$ vertices which is far removed from the optimum solution.

	\section{One-Step Drift}\label{app:onestepdrift}
	Recall that $m_j(n)$ denotes the total number of edges on side $j\in\{1,2\}$, and $m_m(n)$ denotes the total number of edges between the two communities.
The one-step drifts of the random variables associated with the Markov process of adoption are given as follows:
\begin{itemize}
\item One-step drift of $A_j(\cdot)$ for $j\in\{1,2\}$:
\begin{align*}
&\expect[A_j(k+1) - A_j(k)|X^n(k)] =\allowdisplaybreaks\\
&\myquad[4] -\frac{A_j(k)}{A_1(k) + A_2(k) + A^{(1)}_m(k) + A^{(2)}_m(k)}\allowdisplaybreaks\\
&\myquad[4] - \frac{A_j(k)}{A_1(k) + A_2(k) + A^{(1)}_m(k) + A^{(2)}_m(k) }\times\frac{A_j(k)-1}{2m_j(n) - 2T_j(k) -1} \allowdisplaybreaks\\
&\myquad[4] + \frac{A_j(k)}{A_1(k) + A_2(k) + A^{(1)}_m(k) + A^{(2)}_m(k) }\times\allowdisplaybreaks\\
&\myquad[4] \myquad[4] \sum_{u_j+u_{-j} = K_j(d_j,d_{-j})}(d_j - u_j-1)\times\frac{(d_j-u_j)I^{(j)}_{d_j,d_{-j},u_j,u_{-j}}(k)}{2m_j(n) - 2T_j(k) -1}\allowdisplaybreaks\\
&\myquad[4] + \frac{A^{(-j)}_m(k)}{A_1(k) + A_2(k) + A^{(1)}_m(k) + A^{(2)}_m(k) }\times\allowdisplaybreaks\\
&\myquad[4] \myquad[4]  \sum_{u_j+u_{-j} = K_j(d_j,d_{-j})}(d_j - u_j)\times\frac{(d_{-j}-u_{-j})I^{(j)}_{d_j,d_{-j},u_j,u_{-j}}(k)}{m_m(n) - (k-T_1(k)-T_2(k))}
\end{align*}
\item One-step drift of $T_j(\cdot)$ for $j\in\{1,2\}$:
\begin{align*}
&\expect[T_j(k+1) - T_j(k)|X^n(k)] = + \frac{A_j(k)}{A_1(k) + A_2(k) + A^{(1)}_m(k) + A^{(2)}_m(k) }
\end{align*}
\item One-step drift of $A^{(j)}_m(\cdot)$ for $j\in\{1,2\}$:\\
\begin{align*}
&\expect[A^{(j)}_m(k+1) - A^{(j)}_m(k)|X^n(k)] =\allowdisplaybreaks\\
&\myquad[4] -\frac{A^{(j)}_m(k)}{A_1(k) + A_2(k) + A^{(1)}_m(k) + A^{(2)}_m(k)}\allowdisplaybreaks\\
&\myquad[4] -\frac{A^{(-j)}_m(k)}{A_1(k) + A_2(k) + A^{(1)}_m(k) + A^{(2)}_m(k)}\times\frac{A^{(j)}_m(k)}{m_m(n) - (k-T_1(k)-T_2(k))}\allowdisplaybreaks\\
&\myquad[4] + \frac{A^{(-j)}_m(k)}{A_1(k) + A_2(k) + A^{(1)}_m(k) + A^{(2)}_m(k) }\times\allowdisplaybreaks\\
&\myquad[4] \myquad[4] \sum_{u_j+u_{-j} = K_j(d_j,d_{-j})}(d_{-j}- u_{-j}-1)\times\frac{(d_{-j}-u_{-j})I^{(j)}_{d_j,d_{-j},u_j,u_{-j}}(k)}{m_m(n) - (k-T_1(k)-T_2(k))}\allowdisplaybreaks\\
&\myquad[4] + \frac{A_j(k)}{A_1(k) + A_2(k) + A^{(1)}_m(k) + A^{(2)}_m(k) }\times\allowdisplaybreaks\\
&\myquad[4] \myquad[4] \sum_{u_j+u_{-j} = K_j(d_j,d_{-j})}(d_{-j} - u_{-j})\times\frac{(d_j-u_j)I^{(j)}_{d_j,d_{-j},u_j,u_{-j}}(k)}{2m_j(n) - 2T_j(k) - 1}
\end{align*}
\item One-step drift of $I^{(j)}_{d_j,d_{-j},u_j,u_{-j}}(\cdot)$ for $j\in\{1,2\}$:
\begin{align*}
&\expect[I^{(j)}_{d_j,d_{-j},u_j,u_{-j}}(k+1) -I^{(j)}_{d_j,d_{-j},u_j,u_{-j}}(k)|X^n(k)] =\allowdisplaybreaks\\
&\myquad[4] - \frac{A_j(k)}{A_1(k) + A_2(k) + A^{(1)}_m(k) + A^{(2)}_m(k) }\times\frac{(d_j-u_j)I^{(j)}_{d_j,d_{-j},u_j,u_{-j}}(k)}{2m_j(n) - 2T_j(k) -1} \allowdisplaybreaks\\
&\myquad[4]- \frac{A^{(-j)}_m(k)}{A_1(k) + A_2(k) + A^{(1)}_m(k) + A^{(2)}_m(k) }\times\frac{(d_{-j}-u_{-j})I^{(j)}_{d_j,d_{-j},u_j,u_{-j}}(k)}{m_m(n) - (k-T_1(k)-T_2(k))}\allowdisplaybreaks\\
&\myquad[4]+ \frac{A_j(k)}{A_1(k) + A_2(k) + A^{(1)}_m(k) + A^{(2)}_m(k) } \times\frac{(d_j-u_j+1)I^{(j)}_{d_j,d_{-j},u_j-1,u_{-j}}(k)}{2m_j(n) - 2T_j(k) - 1}\allowdisplaybreaks\\
&\myquad[4]+ \frac{A^{(-j)}_m(k)}{A_1(k) + A_2(k) + A^{(1)}_m(k) + A^{(-j)}_m(k) } \times\frac{(d_{-j}-u_{-j}+1)I^{(j)}_{d_j,d_{-j},u_j,u_{-j}-1}(k)}{m_m(n) - (k-T_1(k)-T_2(k))}
\end{align*}
\end{itemize}
Although we've presented the one-step drift of $A_j$ and $A_m^{(j)}$ for $j\in\{1,2\}$, we are not going to use them. As is pointed out in Remark \ref{rem:trackvar}, we only need to keep track of $T_j$ and $I_{d_j,d_{-j},u_j,u_{-j}}^{(j)}$. In particular, using the balance equations \eqref{eq:balanceeq}, we replace all the terms $A_j$ and $A_m^{(j)}$ by an affine function of $T_j$ and $I_{d_j,d_{-j},u_j,u_{-j}}^{(j)}$.

	\section{Derivation of ODEs}\label{app:odederive}
	The scaled variables are supposed to model the behavior of their discrete counterpart, as we mentioned in \eqref{eq:scalevar}. Using Remark \ref{rem:newlambda}, Point \ref{point:2} in Section \ref{sec:odeapprox}, Remark \ref{rem:newrand}, and the one-step drifts in Appendix \ref{app:onestepdrift}, the ODEs are given as follows:\\
\begin{align}
&
\begin{aligned}
& \frac{di^{(j)}_{d_j,d_{-j},u_j,u_{-j}}}{dx} = \\
&\myquad[1]f_{j,d_j,d_{-j},u_j,u_{-j}}(x,\lambda_1(n),\lambda_2(n),\lambda_m(n),\tau_1,\tau_2,w_1,w_2,w^{(1)}_m,w^{(2)}_m,i^{(1)}_{d_1,d_2,u_1,u_2},i^{(2)}_{d_2,d_1,u_2,u_1})\coloneqq\allowdisplaybreaks\\
&\myquad[5]- \frac{a_j(x)}{a_1(x) + a_2(x) + a^{(1)}_m(x) + a^{(2)}_m(x)} \times\frac{(d_j-u_j)i^{(j)}_{d_j,d_{-j},u_j,u_{-j}}(x)}{\lambda_j(n) - 2\tau_j(x) }\allowdisplaybreaks\\
&\myquad[5]- \frac{a^{(-j)}_m(x)}{a_1(x) + a_2(x) + a^{(1)}_m(x) + a^{(2)}_m(x)}\times \frac{(d_{-j}-u_{-j})i^{(j)}_{d_j,d_{-j},u_j,u_{-j}}(x)}{\lambda_m(n) - \tau_m(x)}\allowdisplaybreaks\\
&\myquad[5]+ \frac{a_j(x)}{a_1(x) + a_2(x) + a^{(1)}_m(x) + a^{(2)}_m(x)}
\times \frac{(d_j-u_j+1)i^{(j)}_{d_j,d_{-j},u_j-1,u_{-j}}(x)}{\lambda_j(n) - 2\tau_j(x)  }\allowdisplaybreaks\\
&\myquad[5]+ \frac{a^{(-j)}_m(x)}{a_1(x) + a_2(x) + a^{(1)}_m(x) + a^{(2)}_m(x)}
\times \frac{(d_{-j}-u_{-j}+1)i^{(j)}_{d_j,d_{-j},u_j,u_{-j}-1}(x)}{\lambda_m(n) - \tau_m(x)},
\end{aligned}\label{eq:diffeq_i}\allowdisplaybreaks\\
&
\begin{aligned}
&\frac{d\tau_j}{dx}  = f_{j}(x,\lambda_1(n),\lambda_2(n),\lambda_m(n),\tau_1,\tau_2,w_1,w_2,w^{(1)}_m,w^{(2)}_m,i^{(1)}_{d_1,d_2,u_1,u_2},i^{(2)}_{d_2,d_1,u_2,u_1}) \coloneqq \\
& \myquad[5] \frac{a_j(x)}{a_1(x) + a_2(x) + a^{(1)}_m(x) + a^{(2)}_m(x)},
\end{aligned} \label{eq:diffeq_tau}\allowdisplaybreaks\\
&
\begin{aligned}
&\frac{dw_j}{dx}  = f_{j+2}(x,\lambda_1(n),\lambda_2(n,)\lambda_m(n),\tau_1,\tau_2,w_1,w_2,w^{(1)}_m,w^{(2)}_m,i^{(1)}_{d_1,d_2,u_1,u_2},i^{(2)}_{d_2,d_1,u_2,u_1})
\coloneqq\allowdisplaybreaks\\
&\myquad[5] \frac{a_j(x)}{a_1(x) + a_2(x) + a^{(1)}_m(x) + a^{(2)}_m(x)} \times \frac{-w_j(x)}{{\lambda_j(n) - 2\tau_j(x)}},
\end{aligned}\label{eq:diffeq_wj}\allowdisplaybreaks\\
&
\begin{aligned}
&\frac{dw^{(j)}_m}{dx}  = f_{j+4}(x,\lambda_1(n),\lambda_2(n),\lambda_m(n),\tau_1,\tau_2,w_1,w_2,w^{(1)}_m,w^{(2)}_m,i^{(1)}_{d_1,d_2,u_1,u_2},i^{(2)}_{d_2,d_1,u_2,u_1})
\coloneqq\allowdisplaybreaks\\
&\myquad[5] \frac{a^{(-j)}_m(x)}{a_1(x) + a_2(x) + a^{(1)}_m(x) + a^{(2)}_m(x)} \times \frac{ -w^{(j)}_{m}(x)}{\lambda_m(n) - \tau_m(x)},
\end{aligned}\label{eq:diffeq_wmj}
\end{align}
where $\tau_m(x)\coloneqq x-\tau_1(x)-\tau_2(x)$,
\begin{align}
&a_j(x) \coloneqq -\sum_{\substack{u_j+u_{-j} \leq K_j(d_j,d_{-j})\\d_j + d_{-j} \leq d_{\max}}} (d_j-u_j)i^{(j)}_{d_j,d_{-j},u_j,u_{-j}}(x) + \lambda_j(n) - 2\tau_j(x) - w_j(x), \text{ and }\label{eq:diffeq_aj}\\
&a^{(j)}_m(x) \coloneqq -\sum_{\substack{u_j+u_{-j} \leq K_j(d_j,d_{-j})\\d_j + d_{-j} \leq d_{\max}}} (d_{-j}-u_{-j})i^{(j)}_{d_j,d_{-j},u_j,u_{-j}}(x) + \lambda_m(n) -\tau_m(x)- w^{(j)}_m(x).\label{eq:diffeq_amj}
\end{align}
The initial condition is given as follows: for $j\in\{1,2\}$, we have
\begin{align}
\begin{aligned}
& \tau_j(0)=0, \\
&w_j(0) =\sum_{d_j + d_{-j} > d_{\max}} d_j \frac{I^{(j)}_{d_j,d_{-j},0,0}(0)}{n},\\
&w_m^{(j)}(0) =\sum_{d_j + d_{-j} > d_{\max}} d_{-j}\frac{I^{(j)}_{d_j,d_{-j},0,0}(0)}{n},\allowdisplaybreaks\\
& i_{d_j,d_{-j},u_j,u_{-j}}^{(j)}(0)=
\begin{cases}
\frac{1}{n} I^{(j)}_{d_j,d_{-j},0,0}(0)& \text{if } u_j=u_{-j}=0 \text{ and } d_j+d_{-j} \leq d_{\max} \allowdisplaybreaks\\
0 & \text{otherwise}
\end{cases}.
\end{aligned}
\label{eq:ode_ic}
\end{align}
Note that if $\alpha_j(d_j,d_{-j}) = 1$ for all $d_j+d_{-j} > d_{\max}$, then $w_j(x) = w_m^{(j)}(x) = 0$, $\forall x$. Also, note that the above ODEs depend on the value of $n$ via $\lambda_1(n)$, $\lambda_2(n)$, and $\lambda_m(n)$, as well as the initial condition.

We comment that the above initial condition is random. However, as $n\to\infty$, the corresponding random variables converge to a constant in probability; this is a consequence of the third assumption in the graph regularity conditions.
\begin{lemma}\label{lem:initcond_conv}
	Assume the degree regularity conditions given in Definition \ref{def:regcon_deg} hold. Then, as $n\to\infty$:
	\begin{align*}
	&i^{(j)}_{d_j,d_{-j},0,0}(0)\,d_j/\lambda_j(n) \xrightarrow{P} \prob_{j*,m}(d_j,d_{-j})(1-\alpha_j(d_j,d_{-j})), \allowdisplaybreaks\\ &i^{(-j)}_{d_{-j},d_{j},0,0}(0)\,d_j/\lambda_m(n)\xrightarrow{P}\prob_{-j,m*}(d_{-j},d_{j})(1-\alpha_{-j}(d_{-j},d_{j})),\allowdisplaybreaks\\
	&w_j(0)/\lambda_j(n) \xrightarrow{P} \sum_{d_j+d_{-j} > d_{\max}}\prob_{j*,m}(d_j,d_{-j})(1-\alpha_j(d_j,d_{-j})),\allowdisplaybreaks\\
	&w^{(-j)}_m(0)/\lambda_m(n)\xrightarrow{P}\sum_{d_j+d_{-j} > d_{\max}}\prob_{-j,m*}(d_{-j},d_{j})(1-\alpha_{-j}(d_{-j},d_{j})),\allowdisplaybreaks\\
	& a_j(0) \xrightarrow{P} \beta_j\sum_{d_j,d_{-j}}d_j\,\prob_{j,m}(d_{j},d_{-j})\alpha_{j}(d_{j},d_{-j}),\\
	& a^{(j)}_m(0) \xrightarrow{P} \beta_j\sum_{d_j,d_{-j}}d_{-j}\,\prob_{j,m}(d_{j},d_{-j})\alpha_{j}(d_{j},d_{-j}).
	\end{align*}
\end{lemma}
\begin{proof}
	Recall that $I^{(j)}_{d_j,d_{-j},0,0}(0)=\sum_{l\in N(d_j,d_{-j})} (1-\bs{\alpha}_l^{(j)})$, where $N(d_j,d_{-j})\coloneqq\{l:d^{n}_{j,l} = d_j, d^{n}_{m,l} = d_{-j}\,\text{and } l\in\text{Community }j\}$, and $\{\bs{\alpha}_l^{(j)}\}_{l\in N(d_j,d_{-j})}$ are {\it i.i.d.} Bernoulli random variables with success probability $\alpha_j(d_j,d_{-j}) \in (0,1)$ (note that if $\alpha_j(d_j,d_{-j}) \in \{0,1\}$, there is nothing to prove.). Now, using the Chebyshev's inequality, for any fixed $\epsilon > 0$, we have:
	\begin{align*}
	&\prob\left(\left|\frac{d_j}{n\lambda_{j}(n)}I^{(j)}_{d_j,d_{-j},0,0}(0) - \frac{d_j}{n\lambda_{j}(n)} N(d_j,d_{-j})(1-\alpha_j(d_j,d_{-j}))\right| > \epsilon \right)\allowdisplaybreaks\\
	 &\myquad[5]\leq \left(\frac{d_j}{n\lambda_{j}(n)\epsilon}\right)^2 N(d_j,d_{-j}) \times \alpha_j(d_j,d_{-j})(1-\alpha_j(d_j,d_{-j}))\allowdisplaybreaks\\
	 &\myquad[5]\leq \left(\frac{d_j}{\lambda_{j,j}(n)\epsilon}\right)^2 \times
	 \frac{N(d_j,d_{-j})}{n_j} \times \frac{1}{n_j} =  \frac{o(n_j)}{n_j}\xrightarrow{n\to\infty} 0,
	\end{align*}
	Note that by degree regularity conditions and Remark \ref{rem:newlambda}, we have
	\begin{align*}
		\frac{N(d_j,d_{-j})}{n_j} \xrightarrow{n\to\infty} \prob_{j,m}(d_j,d_{-j}),\myquad[1]\frac{d_j\,N(d_j,d_{-j})}{n\lambda_j(n)} = \frac{d_j\,N(d_j,d_{-j})}{n_j\lambda_{j,j}(n)} \xrightarrow{n\to\infty} \prob_{j*,m}(d_j,d_{-j}).
	\end{align*}

	The proof of the other cases are similar. Note that to prove the convergence in probability of $w_j(0)/\lambda_j(n)$ (and similarly $w^{(-j)}_m(0)/\lambda_m(n)$), we need to invoke the third assumption of degree regularity conditions:
	\begin{align*}
	&\prob\Bigg(\Bigg|\sum_{d_j + d_{-j} > d_{\max}} \frac{d_j}{n \lambda_{m}(n)} I^{(j)}_{d_j,d_{-j},0,0}(0)  \\
	&\myquad[10] - \sum_{d_j + d_{-j} > d_{\max}} \frac{d_j}{n \lambda_{m}(n)} N(d_j,d_{-j})(1-\alpha_j(d_j,d_{-j}))\Bigg| > \epsilon \Bigg)\allowdisplaybreaks\\
	&\myquad[5]\leq \sum_{d_j + d_{-j} > d_{\max}} \left(\frac{d_j}{n\lambda_{m}(n)\epsilon}\right)^2 N(d_j,d_{-j}) \times \alpha_j(d_j,d_{-j})(1-\alpha_j(d_j,d_{-j}))\allowdisplaybreaks\\
	&\myquad[5]\leq \frac{\sum_{i}(d^{n}_{j,i})^2}{n_j}\times \frac{1}{n_j}\times \frac{1}{  (\epsilon\lambda_{m,j}(n))^2} = \frac{o(n_j)}{n_j} \xrightarrow{n\to\infty} 0.
	\end{align*}
	 Finally, the convergence of $a_j(0)$ and $a^{(j)}_m(0)$ in probability follow by the balance equations \eqref{eq:balanceeq_W}.
\end{proof}

	\section{Augmented Process}\label{app:augment}
	Let us fix the initial condition given by \eqref{eq:ode_ic}. As we mentioned for the truncated process, we only need to keep track of the followings to study the evolution of the augmented process: the number of times the algorithm visits each community, and number of inactive vertices in each community with different attributes. Let us denote the corresponding random variables for the augmented process at time $k$ with $\widetilde{I}^{(j)}_{d_j,d_{-j},u_j,u_{-j}}(k)$, $\widetilde{T}_j(k)$, $\widetilde{W}_j(k)$, and $\widetilde{W}_m^{j}(k)$.

	Note that the one-step drifts of these random variables at time $k$ not only depend on $\widetilde{X}^n_\delta$, but also on the value of $k$ itself. In particular, the one-step drifts of these random variables before adding $\widetilde{v}_1$ and $\widetilde{v}_2$, and after removing all the augmented half-edges are same as the ones given by Appendix \ref{app:onestepdrift} and Remark \ref{rem:newrand}; however, the one-step drifts during time $\floor{t_\kappa n},\cdots,\floor{t_\kappa n} + 2\floor{n\chi^{(1\gets 1)}\lambda_1(n)/2}+2\floor{n\chi^{(2\gets 2)} \lambda_2(n)/2}+\floor{n\chi^{(1\gets 2)} \lambda_m(n)} + \floor{n\chi^{(2\gets 1)} \lambda_m(n)}-1$ follow a different formulation.

	Let $k_\kappa \coloneqq \floor{t_\kappa n}$, $k_{\kappa}^{(1\gets 1)} \coloneqq k_\kappa + 2\floor{n\chi^{(1\gets 1)}\lambda_1(n)/2}$, $k_{\kappa}^{(2\gets 1)} \coloneqq k_{\kappa}^{(1\gets 1)} + \floor{n\chi^{(2\gets 1)} \lambda_m(n)}$, $k_{\kappa}^{(2\gets 2)} \coloneqq k_{\kappa}^{(2\gets 1)} + 2 \floor{n \chi^{(2\gets 2)} \lambda_2(n)/2}$, and $k_{\kappa}^{(1\gets 2)} \coloneqq k_{\kappa}^{(2\gets 2)} + \floor{n \chi^{(1\gets 2)} \lambda_m(n)}$. These values correspond to the time in which the augmented process changes its behavior. Next, we describe the one-step drifts of the associated random variables for different times:
	\begin{enumerate}[label = (\roman*)]
		\item $0 \leq k < k_\kappa$: during this interval, the augmented process follows the same path as the truncated process and the one-step drifts of $\widetilde{I}^{(j)}_{d_j,d_{-j},u_j,u_{-j}}$, $\widetilde{T}_j(k)$, $\widetilde{W}_j(k)$, and $\widetilde{W}^{(j)}_m(k)$ are same as the ones given by Appendix \ref{app:onestepdrift} and Remark~\ref{rem:newrand}.
		\item $k_\kappa \leq k < k_{\kappa}^{(1\gets 1)}$: during this interval, we remove random regular half-edges from within community $1$. The one-step drifts are given as follows: \label{case:aug(ii)}
		\begin{align*}
			&\expect[\widetilde{T}_1(k+1) - \widetilde{T}_1(k)|\widetilde{X}^n_{\delta}(k)] = 1/2,\allowdisplaybreaks\\
			&\expect[\widetilde{I}^{(1)}_{d_1,d_{2},u_1,u_{2}}(k+1) -\widetilde{I}^{(1)}_{d_1,d_{2},u_1,u_{2}}(k)|\widetilde{X}^n_{\delta}(k)] = \allowdisplaybreaks\\
			&\myquad[6]- \frac{(d_1-u_1)\widetilde{I}^{(1)}_{d_1,d_{2},u_1,u_{2}}(k)}{2m_1(n) - 2{T}_1(k_\kappa) - (k-k_\kappa) -1} \allowdisplaybreaks\\
			&\myquad[6]+ \frac{(d_1-u_1+1)\widetilde{I}^{(1)}_{d_1,d_{2},u_1-1,u_{2}}(k)}{2m_1(n) - 2{T}_1(k_\kappa) - (k-k_\kappa) -1},\allowdisplaybreaks\\
			&\expect[\widetilde{W}_1(k+1) - \widetilde{W}_1(k)|\widetilde{X}^n_{\delta}(k)] =\frac{ -\widetilde{W}_1(k)}{2m_1(n) - 2{T}_1(k_\kappa) - (k-k_\kappa) -1}.
		\end{align*}
		Note that the one-step drifts of all the other random variables are zero. Also, note that at time $k$, we remove only one half-edge which is not augmented; this is why we increment the value of $\widetilde{T}_j(k)$ by $1/2$ instead of $1$.
		\item $k_{\kappa}^{(1\gets 1)} \leq k < k_{\kappa}^{(2\gets 1)}$: during this interval, we remove random regular half-edges from community $2$ that are between the two communities. The one-step drifts are given as follows:\label{case:aug(iv)}
		\begin{align*}
			&\expect[\widetilde{I}^{(2)}_{d_2,d_{1},u_2,u_{1}}(k+1) -\widetilde{I}^{(2)}_{d_2,d_{1},u_2,u_{1}}(k)|\widetilde{X}^n_{\delta}(k)] = \allowdisplaybreaks\\
			&\myquad[6]- \frac{(d_1-u_1)\widetilde{I}^{(2)}_{d_2,d_{1},u_2,u_{1}}(k)}{m_m(n) - (k_{\kappa}-{T}_1(k_{\kappa})-{T}_2(k_{\kappa})) - (k-k_{\kappa}^{(1\gets 1)}) }\allowdisplaybreaks\\
			&\myquad[6]+ \frac{(d_1-u_1+1)\widetilde{I}^{(2)}_{d_2,d_{1},u_2,u_{1}-1}(k)}{m_m(n) - (k_{\kappa}-{T}_1(k_{\kappa})-{T}_2(k_{\kappa})) - (k-k_{\kappa}^{(1\gets 1)}) }\allowdisplaybreaks\\
			&\expect[\widetilde{W}^{(2)}_{m}(k+1) - \widetilde{W}^{(2)}_{m}(k)|\widetilde{X}^n_{\delta}(k)] = \frac{- \widetilde{W}^{(2)}_{m}(k)}{m_m(n) - (k_{\kappa}-{T}_1(k_{\kappa})-{T}_2(k_{\kappa})) - (k-k_{\kappa}^{(1\gets 1)}) }.
		\end{align*}
		\item $k_{\kappa}^{(2\gets 1)} \leq k < k_{\kappa}^{(2\gets 2)}$: during this interval, we remove random regular half-edges from within community $2$. Similar to \ref{case:aug(ii)} above, the one-step drifts are given as follows:
		\begin{align*}
			&\expect[\widetilde{T}_2(k+1) - \widetilde{T}_2(k)|\widetilde{X}^n_{\delta}(k)] = 1/2,\allowdisplaybreaks\\
			&\expect[\widetilde{I}^{(2)}_{d_2,d_{1},u_2,u_{1}}(k+1) -\widetilde{I}^{(2)}_{d_2,d_{1},u_2,u_{1}}(k)|\widetilde{X}^n_{\delta}(k)] = \allowdisplaybreaks\\
			&\myquad[6]- \frac{(d_2-u_2)\widetilde{I}^{(2)}_{d_2,d_{1},u_2,u_{1}}(k)}{2m_2(n) - 2{T}_2(k_\kappa) - (k-k_{\kappa}^{(2\gets 1)}) -1} \allowdisplaybreaks\\
			&\myquad[6]+ \frac{(d_2-u_2+1)I^{(2)}_{d_2,d_{1},u_2-1,u_{1}}(k)}{2m_2(n) - 2{T}_2(k_\kappa) - (k-k_{\kappa}^{(2\gets 1)}) -1},\allowdisplaybreaks\\
			&\expect[\widetilde{W}_2(k+1) - \widetilde{W}_2(k)|\widetilde{X}^n_{\delta}(k)] =\frac{ -\widetilde{W}_2(k)}{2m_2(n) - 2{T}_2(k_\kappa) - (k-k_{\kappa}^{(2\gets 1)}) -1}.
		\end{align*}
		\item $k_{\kappa}^{(2\gets 2)} \leq k < k_{\kappa}^{(1\gets 2)}$: during this interval, we remove random regular half-edges from community $1$ that are between the two communities. Similar to \ref{case:aug(iv)} above, the one-step drifts are given as follows:
		\begin{align*}
			&\expect[\widetilde{I}^{(1)}_{d_1,d_{2},u_1,u_{2}}(k+1) -\widetilde{I}^{(1)}_{d_1,d_{2},u_1,u_{2}}(k)|\widetilde{X}^n_{\delta}(k)] = \allowdisplaybreaks\\
			&\myquad[6]- \frac{(d_2-u_2)\widetilde{I}^{(1)}_{d_1,d_{2},u_1,u_{2}}(k)}{m_m(n) - (k_{\kappa}-{T}_1(k_{\kappa})-{T}_2(k_{\kappa})) - (k-k_{\kappa}^{(2\gets 2)}) }\allowdisplaybreaks\\
			&\myquad[6]+ \frac{(d_2-u_2+1)\widetilde{I}^{(1)}_{d_1,d_{2},u_1,u_{2}-1}(k)}{m_m(n) - (k_{\kappa}-{T}_1(k_{\kappa})-{T}_2(k_{\kappa})) - (k-k_{\kappa}^{(2\gets 2)}) },\allowdisplaybreaks\\
			&\expect[\widetilde{W}^{(1)}_{m}(k+1) - \widetilde{W}^{(1)}_{m}(k)|\widetilde{X}^n_{\delta}(k)] = \frac{- \widetilde{W}^{(1)}_{m}(k)}{m_m(n) - (k_{\kappa}-{T}_1(k_{\kappa})-{T}_2(k_{\kappa})) - (k-k_{\kappa}^{(2\gets 2)}) }.
		\end{align*}
		\item $k_{\kappa}^{(1\gets 2)} \leq k$: after time $k_{\kappa}^{(1\gets 2)}$, the augmented process proceed normally, and the one-step drift of $\widetilde{I}^{(j)}_{d_j,d_{-j},u_j,u_{-j}}$ and $\widetilde{T}_j(k)$ follows the same formulation as in Appendix \ref{app:onestepdrift}, with one exception: the term $k - \widetilde{T}_1(k) - \widetilde{T}_2(k)$ needs to be replaced with the term $k - \widetilde{T}_1(k) - \widetilde{T}_2(k) - (\floor{n\chi^{(1\gets 1)}\lambda_1(n) /2}+\floor{n\chi^{(2\gets 2)} \lambda_2(n)/2})$.
	\end{enumerate}

	Following similar approach as in Section \ref{sec:odeapprox}, we derive a set of ODEs associated with the one-step drifts given above. Invoking the Wormald's theorem once again \cite[Section 5]{Wormald1999}, we get similar results as in Theorem~\ref{thm:odesfinite} and Corollary \ref{cor:odesfinite} for the augmented process with minimal changes.

	Next, following the same logic as in Sections \ref{sec:odeprobsol}-\ref{sec:odeanalysisfinite}, we relate the above ODEs with a system of four-dimensional ODEs. In particular, one can either follow the same intuitive discussion as in Section \ref{sec:odeprobsol}, or use a plug-in approach as in Lemma \ref{lem:diffeq_sol} to obtained the ODEs that are discussed below.

	Consider the initial condition \eqref{eq:ode_ic}. Let $t_{\kappa}^{(j,j')}$ denote the scaled version of $k_{\kappa}^{(j,j')}$, i.e., for $j,j'\in\{1,2\}$:
	\begin{align*}
		&t_\kappa^{(1\gets 1)} = t_\kappa + \chi^{(1\gets 1)} \lambda_1(n),\myquad[2] t_\kappa^{(2\gets 1)} = t_\kappa^{(1\gets 1)} + \chi^{(2\gets 1)} \lambda_m(n),\\
		&t_\kappa^{(2\gets 2)} = t_\kappa^{(2\gets 1)} + \chi^{(2\gets 2)} \lambda_2(n),\myquad[2] t_\kappa^{(1\gets 2)} = t_\kappa^{(2\gets 2)} + \chi^{(1\gets 2)} \lambda_m(n)
	\end{align*}
	Note that we are only interested in the state of the process at time $k_{\kappa}^{(1\gets 2)}$ (which corresponds to the scaled-time $t_\kappa^{(1\gets 2)}$).	Let $\widetilde{\bs{\mu}}(t) = \left(\widetilde{\mu}^{(1\gets 1)}(t),\widetilde{\mu}^{(1\gets 2)}(t),\widetilde{\mu}^{(2\gets 1)}(t),\widetilde{\mu}^{(2\gets 2)}(t)\right)$ denote the solution of the following system of differential equations up to $t\leq t_\kappa^{(1\gets 2)}$:
	\begin{enumerate}[label = (\roman*)]
		\item For $t < t_\kappa$: for $j\in\{1,2\}$, we have
		\begin{align*}
			&\begin{aligned}
				&\frac{-\widetilde{a}_j(t)}{\widetilde{a}_1(t) + \widetilde{a}_2(t) + \widetilde{a}^{(1)}_m(t) + \widetilde{a}^{(2)}_m(t)}
				=\lambda_j(n) \frac{d\widetilde{\mu}^{(j\gets j)}}{dt}\left(\widetilde{\mu}^{(j\gets j)}(t)\right),
			\end{aligned}
			\\
			&\begin{aligned}
				&\frac{-\widetilde{a}^{(-j)}_m(t)}{\widetilde{a}_1(t) + \widetilde{a}_2(t) + \widetilde{a}^{(1)}_m(t) + \widetilde{a}^{(2)}_m(t)}
				=\lambda_m(n)\frac{d\widetilde{\mu}^{(j\gets -j)}}{dt}\left(\widetilde{\mu}^{(-j\gets j)}(t)\right),
			\end{aligned}
		\end{align*}
		with initial condition $\widetilde{\bs{\mu}}(0) = \bs{1}$ and $\widetilde{\bs{\mu}} \in {\mathcal{D}}_{\varepsilon(\kappa),n}$. Recall that $\bs{\mu}(t_\kappa) = \bs{\mu}_\kappa$. Note that the above ODEs are same as the ODEs associated with the truncated process.

		\item For $t_\kappa \leq t < t_\kappa^{(1\gets 1)}$: \label{part:newode2}
		\begin{align*}
			&\begin{aligned}
				&-1=\lambda_1(n) \mu^{(1\gets 1)}_{\kappa} \times \derivative{\widetilde{\mu}^{(1\gets 1)}}{t}
			\end{aligned}
			\\
			&\begin{aligned}
				0 = \derivative{\widetilde{\mu}^{(1\gets 2)}}{t},\qquad 0 = \derivative{\widetilde{\mu}^{(2\gets 1)}}{t},\qquad 0 = \derivative{\widetilde{\mu}^{(2\gets 2)}}{t}
			\end{aligned}
		\end{align*}
		with initial condition $\widetilde{\bs{\mu}}(t_\kappa) = \bs{\mu}(t_\kappa)$ and $\widetilde{\bs{\mu}} \in [0,1]^4$. Note that the solution of the ODEs associated with the one-step drifts of the augmented process is given by \eqref{eq:sol_i}, \eqref{eq:sol_w}, $\widetilde{\tau}_{2}(t_\kappa) = {\tau}_{2}(t_\kappa)$ and
		\begin{align*}
			\widetilde{\tau}_{1}(t) &= \frac{\lambda_1(n)}{2}\left(1 - \mu^{(1\gets 1)}_{\kappa}\widetilde{\mu}^{(1\gets 1)}(t)\right).
		\end{align*}
		Solving the above ODE, we have $\widetilde{\mu}^{(1\gets 1)}(t_\kappa^{(1\gets 1)}) = \mu^{(1\gets 1)}_{\kappa} - \chi^{(1\gets 1)}/\mu^{(1\gets 1)}_{\kappa} =  \mu^{(1\gets 1)}_{\kappa} - \rho^{(1\gets 1)}(\kappa)$. This can be obtained using the following equalities:
		\begin{gather*}
			\frac{\lambda_1(n)}{2}\left(1 - \mu^{(1\gets 1)}_{\kappa}\widetilde{\mu}^{(1\gets 1)}(t_\kappa^{(1\gets 1)})\right) = \widetilde{\tau}_{1}(t_\kappa^{(1\gets 1)}) = \tau_1(t_\kappa) + \frac{1}{2}(t_\kappa^{(1\gets 1)} - t_\kappa),\\
			t_\kappa^{(1\gets 1)} - t_\kappa = \chi^{(1\gets 1)} \lambda_1(n), \text{ and }\tau_1(t_\kappa) = \frac{\lambda_1(n)}{2}\left(1-\left(\mu^{(1\gets 1)}_{\kappa}\right)^2\right).
		\end{gather*}
		\item For $t_\kappa^{(1\gets 1)} \leq t < t_\kappa^{(2\gets 1)}$: \label{part:newode3}
		\begin{align*}
			&\begin{aligned}
				&-1=\lambda_m(n) \mu^{(1\gets 2)}_{\kappa}  \times \derivative{\widetilde{\mu}^{(2\gets 1)}}{t}
			\end{aligned}
			\\
			&\begin{aligned}
				0 = \derivative{\widetilde{\mu}^{(1\gets 1)}}{t},\qquad 0 = \derivative{\widetilde{\mu}^{(1\gets 2)}}{t},\qquad 0 = \derivative{\widetilde{\mu}^{(2\gets 2)}}{t}
			\end{aligned}
		\end{align*}
		with initial condition $\widetilde{\bs{\mu}}(t_\kappa^{(1\gets 1)})$ given by the part~\ref{part:newode2} and $\widetilde{\bs{\mu}} \in [0,1]^4$. For this time interval, the solution of the ODEs associated with the one-step drifts of the augmented process is given by \eqref{eq:sol_i}, \eqref{eq:sol_w}, $\widetilde{\tau}_{2}(t) = \widetilde{\tau}_{2}(t_\kappa^{(1\gets 1)})$ and $\widetilde{\tau}_{1}(t) = \widetilde{\tau}_{1}(t_\kappa^{(1\gets 1)})$.
		Solving the above ODE, we have $\widetilde{\mu}^{(2\gets 1)}(t_\kappa^{(2\gets 1)}) = \mu^{(2\gets 1)}_{\kappa} - \chi^{(2\gets 1)}/\mu^{(1\gets 2)}_{\kappa} = \mu^{(2\gets 1)}_{\kappa} - \rho^{(2\gets 1)}(\kappa)$.
%		Note that
%		\begin{align*}
%			\mu^{(2\gets 1)}_{\kappa} - \rho(\kappa) \leq \widetilde{\mu}^{(2\gets 1)}(t_\kappa^{(2\gets 1)}) \leq \mu^{(2\gets 1)}_{\kappa} - \rho(\kappa)\min\left(\frac{\mu^{(1\gets 2)}_{\kappa}}{\mu^{(2\gets 1)}_{\kappa}},\frac{\mu^{(2\gets 1)}_{\kappa}}{\mu^{(1\gets 2)}_{\kappa}}\right)
%		\end{align*}
		\item For $t_\kappa^{(2\gets 1)} \leq t < t_\kappa^{(2\gets 2)}$: \label{part:newode4}
		\begin{align*}
			&\begin{aligned}
				&-1=\lambda_2(n) \mu^{(2\gets 2)}_{\kappa} \times \derivative{\widetilde{\mu}^{(2\gets 2)}}{t}
			\end{aligned}
			\\
			&\begin{aligned}
				0 = \derivative{\widetilde{\mu}^{(1\gets 2)}}{t},\qquad 0 = \derivative{\widetilde{\mu}^{(2\gets 1)}}{t},\qquad 0 = \derivative{\widetilde{\mu}^{(1\gets 1)}}{t}
			\end{aligned}
		\end{align*}
		with initial condition $\widetilde{\bs{\mu}}(t_\kappa^{(2\gets 1)})$ given by the part~\ref{part:newode3} and $\widetilde{\bs{\mu}} \in [0,1]^4$. Note that the solution of the ODEs associated with the one-step drifts of the augmented process is given by \eqref{eq:sol_i}, \eqref{eq:sol_w}, $\widetilde{\tau}_{1}(t) = \widetilde{\tau}_{1}(t_\kappa^{(2\gets 1)})$ and
		\begin{align*}
			\widetilde{\tau}_{2}(t) &= \frac{\lambda_2(n)}{2}\left(1 - \mu^{(2\gets 2)}_{\kappa}\widetilde{\mu}^{(2\gets 2)}(t)\right).
		\end{align*}
		Solving the above ODE, we have $\widetilde{\mu}^{(2\gets 2)}(t_\kappa^{(2\gets 2)}) = \mu^{(2\gets 2)}_{\kappa} - \chi^{(2\gets 2)}/\mu^{(2\gets 2)}_{\kappa} =  \mu^{(2\gets 2)}_{\kappa} - \rho^{(2\gets 2)}(\kappa)$.
		\item For $t_\kappa^{(2\gets 2)} \leq t < t_\kappa^{(1\gets 2)}$: \label{part:newode5}
		\begin{align*}
			&\begin{aligned}
				&-1=\lambda_m(n) \mu^{(2\gets 1)}_{\kappa}  \times \derivative{\widetilde{\mu}^{(1\gets 2)}}{t}
			\end{aligned}
			\\
			&\begin{aligned}
				0 = \derivative{\widetilde{\mu}^{(1\gets 1)}}{t},\qquad 0 = \derivative{\widetilde{\mu}^{(2\gets 1)}}{t},\qquad 0 = \derivative{\widetilde{\mu}^{(2\gets 2)}}{t}
			\end{aligned}
		\end{align*}
		with initial condition $\widetilde{\bs{\mu}}(t_\kappa^{(2\gets 2)})$ given by the part~\ref{part:newode4} and $\widetilde{\bs{\mu}} \in [0,1]^4$. For this time interval, the solution of the ODEs associated with the one-step drifts of the augmented process is given by \eqref{eq:sol_i}, \eqref{eq:sol_w}, $\widetilde{\tau}_{2}(t) = \widetilde{\tau}_{2}(t_\kappa^{(2\gets 2)})$ and $\widetilde{\tau}_{1}(t) = \widetilde{\tau}_{1}(t_\kappa^{(2\gets 2)})$.
		%	Note that similar to \eqref{eq:sol_taum}, we have
		%	\begin{align*}
			%		\widetilde{\tau}_m(t) = \lambda_m\left(1 - \mu^{(2\gets 1)}_{\kappa}\widetilde{\mu}^{(1\gets 2)}(t)\right).
			%	\end{align*}
		Solving the above ODE, we have $\widetilde{\mu}^{(1\gets 2)}(t_\kappa^{(1\gets 2)}) = \mu^{(1\gets 2)}_{\kappa} - \chi^{(1\gets 2)}/\mu^{(2\gets 1)}_{\kappa} = \mu^{(1\gets 2)}_{\kappa} - \rho^{(1\gets 2)}(\kappa)$.
	\end{enumerate}
	Following the same logic as in Section \ref{sec:odeanalysisfinite}, it is easy to see that the variables $\widetilde{a}_j(t_\kappa^{(1\gets 2)})$ and $\widetilde{a}^{(-j)}_m(t_\kappa^{(1\gets 2)})$ are defined similar to Lemma~\ref{lem:diffeq_sol}. We can similarly rewrite these variables as in Remark~\ref{rem:atoFcon}:
	\begin{align*}
		&\widetilde{a}_j(t_\kappa^{(1\gets 2)})=\lambda_j(n)\widetilde{\mu}^{(j\gets j)}(t_\kappa^{(1\gets 2)}) \left(\mu^{(j\gets j)}_{\kappa}- {\Ffunc}_{(j\gets j)}\left(\widetilde{\mu}^{(j\gets j)}(t_\kappa^{(1\gets 2)}),\widetilde{\mu}^{(j\gets -j)}(t_\kappa^{(1\gets 2)})\right) \right)\allowdisplaybreaks\\
		&\widetilde{a}^{(-j)}_m(t_\kappa^{(1\gets 2)}) = \\
		&\myquad[3]\lambda_m(n)\widetilde{\mu}^{(-j\gets j)}(t_\kappa^{(1\gets 2)})\left({\mu}^{(-j\gets j)}_\kappa  - {\Ffunc}_{(j\gets -j)}\left(\widetilde{\mu}^{(-j\gets -j)}(t_\kappa^{(1\gets 2)}),\widetilde{\mu}^{(-j\gets j)}(t_\kappa^{(1\gets 2)})\right) \right)\\
		&\myquad[3] + \chi^{(-j\gets j)}\lambda_m(n)
	\end{align*}
	where $\widetilde{\mu}^{(j\gets j)}(t_\kappa^{(1\gets 2)})$ and $\widetilde{\mu}^{(-j\gets j)}(t_\kappa^{(1\gets 2)})$ for $j\in\{1,2\}$ are given as follows:
	\begin{align*}
		\widetilde{\mu}^{(j\gets j)}(t_\kappa^{(1\gets 2)}) \coloneqq {\mu}^{(j\gets j)}_\kappa - \rho^{(j\gets j)}(\kappa)\quad \text{and}\quad \widetilde{\mu}^{(-j\gets j)}(t_\kappa^{(1\gets 2)}) \coloneqq {\mu}^{(-j\gets j)}_\kappa - \rho^{(-j\gets j)}(\kappa).
	\end{align*}

	Next, we simplify the given values for $\widetilde{a}_j(t_\kappa^{(1\gets 2)})$  and $\widetilde{a}^{(-j)}_m(t_\kappa^{(1\gets 2)})$, for $j\in\{1,2\}$. Recall that  $\rho^{(j\gets j)}(\kappa) = \kappa \nu_\kappa^{(j\gets j)}$ and $\rho^{(-j\gets j)}(\kappa) = \kappa \nu_\kappa^{(-j\gets j)}$, for $j\in\{1,2\}$, where the vector $\bs{\nu_\kappa}$ is the non-negative Perron-Frobenius eigenvector corresponding to the Perron-Frobenius eigenvalue $\zeta_\kappa$ of the matrix $\bs{J}_{\Ffuncbold(\cdot)}(\bs{\mu}_{\kappa})$, i.e.,
	\begin{center}
		\begin{tikzpicture}
			\matrix[matrix of math nodes, left delimiter={[},right delimiter={]}, nodes={scale = 1.3, minimum height=5ex, inner sep=3pt}, row sep=1ex] (A)
			{
				\frac{\partial \Ffunc_{(1\gets 1)}(\bs{\mu}_{\kappa})}{\partial \mu^{(1\gets 1)}}  &
				\frac{\partial \Ffunc_{(1\gets 1)}(\bs{\mu}_{\kappa})}{\partial \mu^{(1\gets 2)}}  & 0 & 0 \\

				0 & 0 &  \frac{\partial \Ffunc_{(1\gets 2)}(\bs{\mu}_{\kappa})}{\partial \mu^{(2\gets 1)}}  &
				\frac{\partial \Ffunc_{(1\gets 2)}(\bs{\mu}_{\kappa})}{\partial \mu^{(2\gets 2)}} \\

				\frac{\partial \Ffunc_{(2\gets 1)}(\bs{\mu}_{\kappa})}{\partial \mu^{(1\gets 1)}}  &
				\frac{\partial \Ffunc_{(2\gets 1)}(\bs{\mu}_{\kappa})}{\partial \mu^{(1\gets 2)}} & 0 & 0\\

				0 & 0 &  \frac{\partial \Ffunc_{(2\gets 2)}(\bs{\mu}_{\kappa})}{\partial \mu^{(2\gets 1)}}  &
				\frac{\partial \Ffunc_{(2\gets 2)}(\bs{\mu}_{\kappa})}{\partial \mu^{(2\gets 2)}}\\
			};
			\matrix[matrix of math nodes, left delimiter={[},right delimiter={]}, nodes={scale = 1, minimum height=5ex, inner sep=3pt}, row sep=1ex, right=1.5em of A] (B)
			{
				\nu^{(1\gets 1)}_\kappa\\
				\nu^{(1\gets 2)}_\kappa\\
				\nu^{(2\gets 1)}_\kappa\\
				\nu^{(2\gets 2)}_\kappa\\
			};
			\node[right=0.5em of B] (equal) {$= \zeta_\kappa \bs{\nu}_\kappa$.} ;
%			\matrix[matrix of math nodes, left delimiter={[},right delimiter={]}, nodes={scale = 1.3, minimum height=5ex, inner sep=3pt}, row sep=1ex, right=0.5em of equal] (AB)
%			{
%				\nu^{(1\gets 1)}\\
%				\nu^{(1\gets 2)}\\
%				\nu^{(2\gets 1)}\\
%				\nu^{(2\gets 2)}\\
%			};
			\lineaftercolumn{A}{1}\lineaftercolumn{A}{2}\lineaftercolumn{A}{3}
			\linebelowrow{A}{1}\linebelowrow{A}{2}\linebelowrow{A}{3}
			\linebelowrow{B}{1}\linebelowrow{B}{2}\linebelowrow{B}{3}
%			\linebelowrow{AB}{1}\linebelowrow{AB}{2}\linebelowrow{AB}{3}
		\end{tikzpicture}
	\end{center}

	Using the above equality and first order Taylor approximation of ${\Ffunc}(\cdot)$ at $\bs{\mu}_\kappa$, we can simplify the value of $\widetilde{a}_j(t_\kappa^{(1\gets 2)})$ for $j\in\{1,2\}$ as follows:
	\begin{align*}
		&\widetilde{a}_j(t_\kappa^{(1\gets 2)}) \\
			&\myquad[1] =O(\kappa^2) + \lambda_j(n)\left(\mu^{(j\gets j)}_{\kappa} - \rho^{(j\gets j)}(\kappa)\right)\left(\mu^{(j\gets j)}_{\kappa}- {\Ffunc}_{(j\gets j)}\left(\mu^{(j\gets j)}_{\kappa},\mu^{(j\gets -j)}_{\kappa}\right) \right) \\
			&\myquad[3] +\lambda_j(n)\left(\mu^{(j\gets j)}_{\kappa} - \rho^{(j\gets j)}(\kappa)\right)\\
			&\myquad[6]\left(\frac{\partial \Ffunc_{(j\gets j)}(\bs{\mu}_{\kappa})}{\partial \mu^{(j\gets j)}} \times  \rho^{(j\gets j)}(\kappa) + \frac{\partial \Ffunc_{(j\gets j)}(\bs{\mu}_{\kappa})}{\partial \mu^{(j\gets -j)}} \times \rho^{(j\gets -j)}(\kappa)\right)\\
			&\myquad[1] \leq \lambda_j(n)\mu^{(j\gets j)}_{\kappa} \left(\mu^{(j\gets j)}_{\kappa}- {\Ffunc}_{(j\gets j)}\left(\mu^{(j\gets j)}_{\kappa},\mu^{(j\gets -j)}_{\kappa}\right) \right)\\
			&\myquad[3] + \lambda_j(n)\mu^{(j\gets j)}_{\kappa}\zeta_\kappa \rho^{(j\gets j)}(\kappa) + O(\kappa^2)\\
			&\myquad[1] = a_j(t_\kappa) +  \chi^{(j\gets j)} \lambda_j(n) \zeta_\kappa + O(\kappa^2).
	\end{align*}
	Similarly, we can simplify the value of $\widetilde{a}^{(-j)}_m(t_\kappa^{(1\gets 2)})$ for $j\in\{1,2\}$ as follows:
	\begin{align*}
		&\widetilde{a}^{(-j)}_m(t_\kappa^{(1\gets 2)}) \\
		&\myquad[1]=  \lambda_m(n)\left(\mu^{(-j\gets j)}_{\kappa} - \rho^{(-j\gets j)}(\kappa)\right)\left( \mu^{(j\gets -j)}_{\kappa} - {\Ffunc}_{(j\gets -j)}\left(\mu^{(-j\gets -j)}_{\kappa},\mu^{(-j\gets j)}_{\kappa}\right)\right)\\
		&\myquad[4] +\lambda_m(n)\left(\mu^{(-j\gets j)}_{\kappa} - \rho^{(-j\gets j)}(\kappa)\right)\\
		&\myquad[8]\left(\frac{\partial \Ffunc_{(j\gets -j)}(\bs{\mu}_{\kappa})}{\partial \mu^{(-j\gets -j)}} \times  \rho^{(-j\gets -j)}(\kappa) + \frac{\partial \Ffunc_{(j\gets -j)}(\bs{\mu}_{\kappa})}{\partial \mu^{(-j\gets j)}} \times \rho^{(-j\gets j)}(\kappa)\right)\\
		&\myquad[4]+ \chi^{(-j\gets j)}\lambda_m(n) + O(\kappa^2)\\
		&\myquad[1]\leq\lambda_m(n)\mu^{(-j\gets j)}_{\kappa} \left( \mu^{(j\gets -j)}_{\kappa} - {\Ffunc}_{(j\gets -j)}\left(\mu^{(-j\gets -j)}_{\kappa},\mu^{(-j\gets j)}_{\kappa}\right)\right)\\
		&\myquad[4] +  \lambda_m(n)\mu^{(-j\gets j)}_{\kappa}\zeta_\kappa \rho^{(j\gets -j)}(\kappa) + \chi^{(-j\gets j)}\lambda_m(n)+ O(\kappa^2)\\
		&\myquad[1] = {a}^{(-j)}_m(t_\kappa) +  \chi^{(-j\gets j)}\lambda_m(n) +\chi^{(j\gets -j)}\lambda_m(n) \zeta_\kappa+ O(\kappa^2)
%		&\myquad[1]=\lambda_m\left(\mu^{(-j\gets j)}_{\kappa} - \rho^{(-j\gets j)}(\kappa)\right)\Biggg(\frac{{\mu}^{(-j\gets j)}_\kappa \left(\mu^{(j\gets -j)}_{\kappa} - \rho^{(j\gets -j)}(\kappa)\right)}{\left(\mu^{(-j\gets j)}_{\kappa} - \rho^{(-j\gets j)}(\kappa)\right)}\\
%		&\myquad[10]  - {\Ffunc}_{(j\gets -j)}\left(\left(\mu^{(-j\gets -j)}_{\kappa} - \rho^{(-j\gets -j)}(\kappa)\right),\left(\mu^{(-j\gets j)}_{\kappa} - \rho^{(-j\gets j)}(\kappa)\right)\right) \Biggg) \\
%		&\myquad[3] + \chi^{(-j\gets j)}\lambda_m,\allowdisplaybreaks\\
	\end{align*}

	\section{Proofs of Theorems and Lemmas}\label{app:Proofs}
	\subsection{Proof of Theorem \ref{thm:odesfinite}}\label{proof:odesfinite}
	
By the assumptions of Theorem \ref{thm:odesfinite}, we have
\begin{align*}
{\footnotesize \Big(\frac{0}{n},\frac{T_1(0)}{n},\frac{T_2(0)}{n},\frac{W_1(0)}{n},\frac{W_2(0)}{n},\frac{W^{(1)}_m(0)}{n},\frac{W^{(2)}_m(0)}{n},
\frac{I^{(1)}_{d_1,d_2,u_1,u_2}(0)}{n},\frac{I^{(2)}_{d_2,d_1,u_2,u_1}(0)}{n}\Big) \!\!\in\!\! \widehat{\mathcal{D}}_{\varepsilon,n}}.
\end{align*}
Moreover, it is easy to see that the functions given in Appendix \ref{app:odederive} satisfy a Lipschitz condition on
\begin{align*}
\mathcal{D}_\epsilon \cap \{(t,\tau_1,\tau_2,w_1,w_2,w^{(1)}_m,w^{(2)}_m,i^{(1)}_{d_1,d_2,u_1,u_2},i^{(2)}_{d_2,d_1,u_2,u_1}): t\geq 0\}
\end{align*}
with the same Lipschitz constant (``\textit{Lipschitz hypothesis}''). Also, for $\theta_1 =O(n^{-\eta})$
\begin{align*}
&
\Big| \mathbb{E}(T_j(k+1) - T_j(k)\mid X^n(k)) - \\
&\myquad[1] f_j\Big(\frac{k}{n},\lambda_1(n),\lambda_2(n),\lambda_m(n),\frac{T_1(k)}{n},\frac{T_2(k)}{n},\cdots, \frac{I^{(1)}_{d_1,d_2,u_1,u_2}(k)}{n},\frac{I^{(2)}_{d_2,d_1,u_2,u_1}(k)}{n}\Big)\Big| \leq \theta_1,\allowdisplaybreaks\\
&
\Big| \mathbb{E}(W_j(k+1) - W_j(k)\mid X^n(k)) - \\
&\myquad[1] f_{j+2}\Big(\frac{k}{n},\lambda_1(n),\lambda_2(n),\lambda_m(n),\frac{T_1(k)}{n},\frac{T_2(k)}{n},\cdots, \frac{I^{(1)}_{d_1,d_2,u_1,u_2}(k)}{n},\frac{I^{(2)}_{d_2,d_1,u_2,u_1}(k)}{n}\Big)\Big| \leq \theta_1,\allowdisplaybreaks\\
&
\Big| \mathbb{E}(W^{(j)}_m(k+1) - W^{(j)}_m(k)\mid X^n(k)) - \\
&\myquad[1] f_{j+4}\Big(\frac{k}{n},\lambda_1(n),\lambda_2(n),\lambda_m(n),\frac{T_1(k)}{n},\frac{T_2(k)}{n},\cdots, \frac{I^{(1)}_{d_1,d_2,u_1,u_2}(k)}{n},\frac{I^{(2)}_{d_2,d_1,u_2,u_1}(k)}{n}\Big)\Big| \leq \theta_1,
\end{align*}
and,
\begin{align*}
\Big| &\mathbb{E}(I^{(j)}_{d_j,d_{-j},u_j,u_{-j}}(k+1) - I^{(j)}_{d_j,d_{-j},u_j,u_{-j}}(k)\mid X^n(k)) - \\
&\myquad[1] f_{j,d_j,d_{-j},u_j,u_{-j}}\Big(\frac{k}{n},\lambda_1(n),\lambda_2(n),\lambda_m(n),\frac{T_1(k)}{n},\frac{T_2(k)}{n},\cdots,\frac{I^{(2)}_{d_2,d_1,u_2,u_1}(k)}{n}\Big)\Big| \leq \theta_1,
\end{align*}
for all $k< T_{\mathcal{D}_\epsilon}$, where $T_{\mathcal{D}_\epsilon}$ is the minimum $k>0$ such that
\begin{align*}
{\footnotesize\Big(\frac{k}{n},\frac{T_1(k)}{n},\frac{T_2(k)}{n},\frac{W_1(k)}{n},\frac{W_2(k)}{n},\frac{W^{(1)}_m(k)}{n},\frac{W^{(2)}_m(k)}{n},
\frac{I^{(1)}_{d_1,d_2,u_1,u_2}(k)}{n},\frac{I^{(2)}_{d_2,d_1,u_2,u_1}(k)}{n}\Big)  \!\!\notin\!\! \widehat{\mathcal{D}}_{\varepsilon,n}}
\end{align*}
(``\textit{Trend hypothesis}''). Finally, the changes for each random variable in successive steps of the Markov process of adoption is bounded by $1$ (``\textit{Bounded hypothesis}''). Now, Theorem \ref{thm:odesfinite} follows by the direct application of Wormald's Theorem~\cite[Theorem 5.1]{Wormald1999}.
	\subsection{Proof of Lemma \ref{lem:diffeq_sol}}\label{proof:diffeq_sol}
	
The proof follows by substituting the form of the solution and checking the validity of the corresponding differential equations. Before substituting the form of the solution, let us prove the equality in \eqref{eq:sol_t}. Summing up \eqref{eq:sol_mujj} and \eqref{eq:sol_muj-j} for $j\in\{1,2\}$, we have
\begin{align*}
&\lambda_1(n)\frac{d\mu^{(1\gets 1)}}{dt}\left(\mu^{(1\gets 1)}(t)\right) +\lambda_2(n)\frac{d\mu^{(2\gets 2)}}{dt}\left(\mu^{(2\gets 2)}(t)\right) + \\
&\myquad[3] \lambda_m(n) \frac{d\mu^{(1\gets 2)}}{dt}\left(\mu^{(2\gets 1)}(t)\right) + \lambda_m(n) \frac{d\mu^{(2\gets 1)}}{dt}\left(\mu^{(1\gets 2)}(t)\right) =-1.
\end{align*}
Equality \eqref{eq:sol_t} follows by integrating the both sides of the above equation from $0$ to $t$, where the constant of integration is determined by the initial condition \eqref{eq:sol_mu0}. Note that \eqref{eq:sol_t} suggests $\tau_m(t) = t - \tau_1(t)-\tau_2(t)$, where $\tau_m(t)$ is defined by \eqref{eq:sol_taum}. Also, by \eqref{eq:sol_tau}, \eqref{eq:sol_mujj}, and \eqref{eq:sol_muj-j} we have
\begin{align}
\frac{d\mu^{(j\gets j)}}{dt}\left(\mu^{(j\gets j)}(t)\right)^{-1}
&= \lambda_j(n)\frac{d\mu^{(j\gets j)}}{dt}\left(\mu^{(j\gets j)}(t)\right)\times\left(\lambda_j(n)\mu^{(j\gets j)}(t)^2\right)^{-1} \nonumber\allowdisplaybreaks\\
&=\frac{-a_j(t)}{a_1(t) + a_2(t) + a^{(1)}_m(t) + a^{(2)}_m(t) } \times \frac{1}{\lambda_j(n) - 2\tau_j(t) }\label{eq:diffeq_mujj_tmp},
\end{align}
and
\begin{align}
\frac{d\mu^{(j\gets -j)}}{dt}\left(\mu^{(j\gets -j)}(t)\right)^{-1}&\!=\! \lambda_m(n)\frac{d\mu^{(j\gets -j)}}{dt}\mu^{(-j\gets j)}(t)\left(\lambda_m(n)\mu^{(-j\gets j)}(t)\mu^{(j\gets -j)}(t)\right)^{-1}  \nonumber\allowdisplaybreaks\\
&\!=\! \frac{-a^{(-j)}_m(t)}{a_1(t) + a_2(t) + a^{(1)}_m(t) + a^{(2)}_m(t)} \times \frac{1}{\lambda_m(n)-(t-\tau_1-\tau_2)}.\label{eq:diffeq_muj-j_tmp}
\end{align}

Consider the form of $i^{(j)}_{d_j,d_{-j},u_j,u_{-j}}(t)$ given by \eqref{eq:sol_i}. We have
\begin{align*}
&\frac{di^{(j)}_{d_j,d_{-j},u_j,u_{-j}}}{dt} \\
&\myquad[2]= i^{(j)}_{d_j,d_{-j},0,0}(0) \times \allowdisplaybreaks\\
&\myquad[2]\Bigg\{+(d_j-u_j)\frac{d\mu^{(j\gets j)}}{dt}\dbinom{d_j}{u_j}\left(1-\mu^{(j\gets j)}(t)\right)^{u_j}\left(\mu^{(j\gets j)}(t)\right)^{d_j-u_j-1}\\
&\myquad[2]\myquad[10]\times Bi\left(u_{-j};d_{-j},1-\mu^{(j\gets -j)}(t)\right) \allowdisplaybreaks\\
&\myquad[2]~ + (d_{-j}-u_{-j})\frac{d\mu^{(j\gets -j)}}{dt}\dbinom{d_{-j}}{u_{-j}}\left(1 - \mu^{(j\gets -j)}(t)\right)^{u_{-j}} \times \left(\mu^{(j\gets -j)}(t)\right)^{d_{-j}-u_{-j}-1}\\
&\myquad[2]\myquad[10]\times Bi\left(u_j;d_j,1-\mu^{(j\gets j)}(t)\right) \allowdisplaybreaks\\
&\myquad[2]~ - u_j\frac{d\mu^{(j\gets j)}}{dt}\dbinom{d_j}{u_j}\left(1-\mu^{(j\gets j)}(t)\right)^{u_j-1}\left(\mu^{(j\gets j)}(t)\right)^{d_j-u_j}\\
&\myquad[2]\myquad[10]\times Bi\left(u_{-j};d_{-j},1- \mu^{(j\gets -j)}(t)\right) \allowdisplaybreaks\\
&\myquad[2]~ - u_{-j}\frac{d\mu^{(j\gets -j)}}{dt}\dbinom{d_{-j}}{u_{-j}}\left(1 - \mu^{(j\gets -j)}(t)\right)^{u_{-j}-1} \left(\mu^{(j\gets -j)}(t)\right)^{d_{-j}-u_{-j}}\\
&\myquad[2]\myquad[10]\times Bi\left(u_j;d_j,1-\mu^{(j\gets j)}(t)\right)\Bigg\} .
\end{align*}
Using \eqref{eq:sol_i}, we have
\begin{align*}
\frac{di^{(j)}_{d_j,d_{-j},u_j,u_{-j}}}{dt} &= ~(d_j-u_j) \times \frac{d\mu^{(j\gets j)}}{dt}\left( \mu^{(j\gets j)}(t)\right)^{-1} \times i^{(j)}_{d_j,d_{-j},u_j,u_{-j}}(t)\allowdisplaybreaks\\
&~+(d_{-j}-u_{-j}) \times \frac{d\mu^{(j\gets -j)}}{dt}\left(\mu^{(j\gets -j)}(t)\right)^{-1} \times i^{(j)}_{d_j,d_{-j},u_j,u_{-j}}(t)\allowdisplaybreaks\\
&~-(d_j-u_j+1) \times \frac{d\mu^{(j\gets j)}}{dt}\left( \mu^{(j\gets j)}(t)\right)^{-1} \times i^{(j)}_{d_j,d_{-j},u_j-1,u_{-j}}(t) \allowdisplaybreaks\\
&~-(d_{-j}-u_{-j}+1) \times \frac{d\mu^{(j\gets -j)}}{dt}\left(\mu^{(j\gets -j)}(t)\right)^{-1} \times i^{(j)}_{d_j,d_{-j},u_j,u_{-j}-1}(t).
\end{align*}
Now \eqref{eq:diffeq_i} follows by substituting \eqref{eq:diffeq_mujj_tmp} and \eqref{eq:diffeq_muj-j_tmp} into the above equality. Next, consider the function $\tau_j(t)$ given by \eqref{eq:sol_tau}. It is easy to see that,
\begin{align*}
\frac{d\tau_j}{dt} &= -\lambda_j(n)\mu^{(j\gets j)}\frac{d\mu^{(j\gets j)}}{dt}=\frac{a_j(t)}{a_1(t) + a_2(t) + a^{(1)}_m(t) + a^{(2)}_m(t)}.
\end{align*}
Finally, for the functions $w_j(t)$ and $w^{(j)}_m(t)$ given by \eqref{eq:sol_w}, we have
\begin{align*}
&\frac{dw_j}{dt} = w_j(0) \frac{d\mu^{(j\gets j)}}{dt} = w_j(t)\frac{d\mu^{(j\gets j)}}{dt}\left(\mu^{(j\gets j)}(t)\right)^{-1},\allowdisplaybreaks\\
&\frac{dw^{(j)}_m}{dt} = w^{(j)}_m(0) \frac{d\mu^{(j\gets -j)}}{dt} = w^{(j)}_m(t)\frac{d\mu^{(j\gets -j)}}{dt}\left(\mu^{(j\gets -j)}(t)\right)^{-1}.
\end{align*}
Now, using \eqref{eq:diffeq_mujj_tmp} and \eqref{eq:diffeq_muj-j_tmp}, the equations \eqref{eq:diffeq_wj} and \eqref{eq:diffeq_wmj} follows.
	\subsection{Proof of Lemma \ref{lem:Fprop_increasing}}\label{proof:Fprop_increasing}
	Let $Y(x) \coloneqq Bi(u,d,1-x)$ for $x \in(0,1)$. We have:
\begin{align*}
\frac{\mathrm{d}Y}{\mathrm{d}x} &= \dbinom{d}{u}\left((d-u)\times x^{d-u-1}(1-x)^u-u\times x^{d-u}(1-x)^{u-1}\right)\allowdisplaybreaks\\
&= d\times\left(Bi(u,d-1,1-x) - Bi(u-1,d-1,1-x)\right)
\end{align*}
The proof of the lemma follows by straightforward algebraic calculation.
\begin{align*}
&\frac{\partial \Ffunc_{(j\gets j)}}{\partial \mu^{(j\gets j)}} \\
&\begin{aligned}
&\myquad[2]=\sum_{\substack{u_j+u_{-j} \leq K_j(d_j,d_{-j})\allowdisplaybreaks\\d_j+d_{-j} \leq d_{\max}}}  \,\frac{d_j}{\lambda_j(n)} \,i^{(j)}_{d_j,d_{-j},0,0}(0)\, Bi(u_{-j};d_{-j},1 - \mu^{(j\gets -j)})\times (d_j-1)\allowdisplaybreaks\\
&\myquad[6]\left(Bi(u_j;d_j-2,1 -\mu^{(j\gets j)}) - Bi(u_j-1;d_j-2,1 -\mu^{(j\gets j)})\right)  \allowdisplaybreaks\\
&\myquad[2]= \sum_{\substack{u_{-j} \leq K_j(d_j,d_{-j})\allowdisplaybreaks\\d_j+d_{-j} \leq d_{\max}}}\,\frac{d_j}{\lambda_j(n)} \,i^{(j)}_{d_j,d_{-j},0,0}(0)\, Bi(u_{-j};d_{-j},1 - \mu^{(j\gets -j)})\times (d_j-1)\allowdisplaybreaks\\
&\myquad[6] Bi(K_j(d_j,d_{-j})-u_{-j};d_j-2,1 -\mu^{(j\gets j)}) > 0
\end{aligned}\allowdisplaybreaks\\
&\frac{\partial \Ffunc_{(j\gets j)}}{\partial \mu^{(j\gets -j)}} \\
&\begin{aligned}
&\myquad[2]=\sum_{\substack{u_j+u_{-j} \leq K_j(d_j,d_{-j})\allowdisplaybreaks\\d_j+d_{-j} \leq d_{\max}}}  \,\frac{d_j}{\lambda_j(n)} \,i^{(j)}_{d_j,d_{-j},0,0}(0)\, Bi(u_j;d_j-1,1 - \mu^{(j\gets j)})\times d_{-j}\allowdisplaybreaks\\
&\myquad[6]\left(Bi(u_{-j};d_{-j}-1,1 -\mu^{(j\gets -j)}) - Bi(u_{-j}-1;d_{-j}-1,1 -\mu^{(j\gets -j)})\right)   \allowdisplaybreaks\\
&\myquad[2]= \sum_{\substack{u_j \leq K_j(d_j,d_{-j})\allowdisplaybreaks\\d_j+d_{-j} \leq d_{\max}}}  \,\frac{d_j}{\lambda_j(n)} \,i^{(j)}_{d_j,d_{-j},0,0}(0)\, Bi(u_j;d_j-1,1 - \mu^{(j\gets j)})\times d_{-j}\allowdisplaybreaks\\
&\myquad[6] Bi(K_j(d_j,d_{-j})-u_{j};d_{-j}-1,1 -\mu^{(j\gets -j)}) > 0
\end{aligned}\allowdisplaybreaks\\
&\frac{\partial \Ffunc_{(j\gets -j)}}{\partial \mu^{(-j\gets -j)}}\\
&\begin{aligned}
&\myquad[2]= \sum_{\substack{u_j+u_{-j} \leq K_{-j}(d_{-j},d_{j})\allowdisplaybreaks\\d_j+d_{-j} \leq d_{\max}}}  \,\frac{d_j}{\lambda_m(n)} \,i^{(-j)}_{d_{-j},d_{j},0,0}(0)\, Bi(u_{j};d_{j}-1,1 - \mu^{(-j\gets j)})\times d_{-j}\allowdisplaybreaks\\
&\myquad[6]\left(Bi(u_{-j};d_{-j}-1,1 -\mu^{(-j\gets -j)}) - Bi(u_{-j}-1;d_{-j}-1,1 -\mu^{(-j\gets -j)})\right)  \allowdisplaybreaks\\
&\myquad[2]= \sum_{\substack{u_{j} \leq K_{-j}(d_{-j},d_{j})\allowdisplaybreaks\\d_j+d_{-j} \leq d_{\max}}}\,\frac{d_j}{\lambda_m(n)} \,i^{(-j)}_{d_{-j},d_{j},0,0}(0)\, Bi(u_{j};d_{j}-1,1 - \mu^{(-j\gets j)})\times d_{-j}\allowdisplaybreaks\\
&\myquad[6] Bi(K_{-j}(d_{-j},d_{j})-u_{j};d_{-j}-1,1 -\mu^{(-j\gets -j)}) > 0
\end{aligned}\allowdisplaybreaks\\
&\frac{\partial \Ffunc_{(j\gets -j)}}{\partial \mu^{(-j\gets j)}} \\
&\begin{aligned}
&\myquad[2]=\sum_{\substack{u_j+u_{-j} \leq K_{-j}(d_{-j},d_{j})\allowdisplaybreaks\\d_j+d_{-j} \leq d_{\max}}}  \,\frac{d_j}{\lambda_m(n)} \,i^{(-j)}_{d_{-j},d_{j},0,0}(0)\, Bi(u_{-j};d_{-j},1 - \mu^{(-j\gets -j)})\times (d_{j}-1)\allowdisplaybreaks\\
&\myquad[6]\left(Bi(u_j;d_j-2,1 -\mu^{(-j\gets j)}) - Bi(u_j-1;d_j-2,1 -\mu^{(-j\gets j)})\right)  \allowdisplaybreaks\\
&\myquad[2]= \sum_{\substack{u_{-j} \leq K_{-j}(d_{-j},d_{j})\allowdisplaybreaks\\d_j+d_{-j} \leq d_{\max}}}\,\frac{d_j}{\lambda_j(n)} \,i^{(-j)}_{d_{-j},d_{j},0,0}(0)\, Bi(u_{-j};d_{-j},1 - \mu^{(-j\gets -j)})\times (d_j-1)\allowdisplaybreaks\\
&\myquad[6] Bi(K_{-j}(d_{-j},d_{j})-u_{-j};d_j-2,1 -\mu^{(-j\gets j)}) > 0
\end{aligned}
\end{align*}
%Similar inequalities holds for other cases.
	\subsection{Proof of Lemma \ref{lem:Fprop_feasreg}}\label{proof:Fprop_feasreg}
	Note that $\Ffuncbold(\bs{1})\leq \bs{1}$.
\begin{enumerate}[label=(\roman*)]
\item Fix $\bs{\mu} \in \mathcal{U}$. Consider the closed, convex, and compact set $$\mathcal{S} = \left\{ \bs{x}\in[0,1]^4:\Ffuncbold(\bs{\mu}) \leq \bs{x} \leq \bs{\mu} \right\},$$ where the inequalities are interpreted component-wise. Note that $\mathcal{S}$ is a hyperrectangle. By Lemma \ref{lem:Fprop_increasing}, $\Ffuncbold(\bs{s}) \leq \Ffuncbold(\bs{\mu})\leq\bs{s}$ for all $\bs{s}\in\mathcal{S}$ since $\bs{s}\leq\bs{\mu}$. Hence, $\mathcal{S}$ is a subset of $\mathcal{U}$.
\item Define $\bs{G}(\bs{x}) \coloneqq \bs{x} - \Ffuncbold(\bs{x})$ for $\bs{x}\in[0,1]$. Clearly, the set $\mathcal{A} = \{\bs{x}:\bs{G}(\bs{x})\geq \bs{0}\}$ is a closed set as $\bs{G}(\cdot)$ is a continuous function. Since $\mathcal{U}$ is the largest connected component of $\mathcal{A}$ that contains $\bs{1}=(1,1,1,1)$, $\mathcal{U}$ is closed as well. Now, compactness follows from the fact that $\mathcal{U}$ is bounded.
\item Consider the sequence $\big\{\Ffuncbold^k(\bs{u})\big\}_{k=1}^{\infty}$ for some $\bs{u}\in \mathcal{U}$. Since $\Ffuncbold(\mathcal{U})\subset \mathcal{U}$, we have $\Ffuncbold^k(\bs{u})\in \mathcal{U}$ for all $k$. By compactness of $\mathcal{U}$, this sequence has a subsequence that converges to a point $\bs{u}_*\in \mathcal{U}$. Now the result follows by the fact that $\Ffuncbold^{k+1}(\bs{u}) \leq \Ffuncbold^{k}(\bs{u})$ for all $k\geq 0$, where $\Ffuncbold^{0}(\bs{u})\coloneqq\bs{u}$.
\item Recall that $\Ffuncbold \coloneqq (\Ffunc_{(1\gets 1)},\Ffunc_{(1\gets 2)},\Ffunc_{(2\gets 1)},\Ffunc_{(2\gets 2)})$, and for any $\bs{\mu}\in[0,1]^4$ we use the notation $\bs{\mu}=\left(\mu^{(1\gets 1)},\mu^{(1\gets 2)},\mu^{(2\gets 1)},\mu^{(2\gets 2)}\right)$. Also, recall that $\Ffunc_{(j\gets j)}$ is a function of $\mu^{(j\gets j)}$ and $\mu^{(j\gets -j)}$, and $\Ffunc_{(j\gets -j)}$ is a function of $\mu^{(-j\gets -j)}$ and $\mu^{(-j\gets j)}$, for $j\in\{1,2\}$. Now, by symmetry there are two cases that we need to consider:
 \begin{enumerate}[label=(\alph*)]
\item $u^{(j\gets j)} =u^{(j\gets j)}_*$: If $u^{(j\gets -j)}  > u^{(j\gets -j)}_*$, then by equality $\Ffunc_{(j\gets j)}(u^{(j\gets j)}_*,u^{(j\gets -j)}_*) = u^{(j\gets j)}_*$ and Lemma \ref{lem:Fprop_increasing}, we have $u^{(j\gets j)} < \Ffunc_{(j\gets j)}(u^{(j\gets j)},u^{(j\gets -j)})$. If $u^{(j\gets -j)}  = u^{(j\gets -j)}_*$, then either $u^{(-j\gets j)}  > u^{(-j\gets j)}_*$ or $u^{(-j\gets -j)} > u^{(-j\gets -j)}_*$.
Now, if $u^{(-j\gets -j)}  > u^{(-j\gets -j)}_*$ then by the same argument we have $u^{(j\gets -j)} < \Ffunc_{(j\gets -j)}(u^{(-j\gets -j)},u^{(j\gets -j)}).$
Otherwise, we have $u^{(-j\gets j)}  > u^{(-j\gets j)}_*$ and by the same argument $u^{(-j\gets -j)} < \Ffunc_{(-j\gets -j)}(u^{(-j\gets -j)},u^{(-j\gets j)}).$
\item $u^{(j\gets -j)} =u^{(j\gets -j)}_*$: The argument is exactly the same as the previous case, and we avoid repetition.
\end{enumerate}
\end{enumerate}
	\subsection{Proof of Theorem \ref{thm:alterode_sol}}\label{proof:alterode_sol}
	
By Corollary \ref{cor:Fprop_sol}, $\bs{\mu}_* = \lim_{n\rightarrow \infty} \Ffuncbold^n(\bs{1})$ is a fixed point of $\Ffuncbold$. Let $$\mathcal{N}\coloneqq \mathcal{U} \cap \left\{ \bs{x}\in[0,1]^4:\bs{\mu}_*\leq \bs{x} \leq \bs{1}\right\}.$$ For any arbitrary point $\bs{u}\in \mathcal{N}$, define $\mathcal{S}_{\bs{u}} \coloneqq \left\{ x:\Ffuncbold(\bs{u}) \leq x \leq \bs{u} \right\}$. By the proof of Lemma \ref{lem:Fprop_feasreg} part (i), $\mathcal{S}_{\bs{u}}$ is a subset of $\mathcal{U}$. Moreover, if $\bs{u}\neq \bs{\mu}_*$, then $\Ffuncbold(\bs{u}) \geq \Ffuncbold(\bs{\mu}_*)$ and hence $\mathcal{S}_{\bs{u}}$ is a subset of $\mathcal{N}$. So we have $\Ffuncbold(\mathcal{N})\subset \mathcal{N}$. It is also easy to see that $\mathcal{N}$ is closed and compact, and $\bs{\mu}_*$ is the unique fixed point of $\Ffuncbold$ in $\mathcal{N}$.

Now, consider the ODEs \eqref{eq:alterode}. Note that the initial condition is in $\mathcal{N}$. Moreover, if $\bs{\mu}(t) \in \mathcal{N}$, then $\dot{\bs{\mu}}(t) = \Ffuncbold(\bs{\mu}(t)) - \bs{\mu}(t)$ is directing toward $\mathcal{N}$ as $\bs{\mu}(t)+\delta(\Ffuncbold(\bs{\mu}(t)) - \bs{\mu}(t)) \in \mathcal{S}_{\bs{\mu}(t)} \subset \mathcal{N}$ for all $\delta \in [0,1]$. Hence, $\mathcal{N}$ is a positive invariant set.

Consider the function $V(\bs{\mu})\coloneqq(\bs{\mu} - \bs{\mu}_*) (\bs{\mu} - \bs{\mu}_*)^T$. Note that $\forall\bs{\mu} \in \mathcal{N}\setminus \{\bs{\mu}_*\}$, we have
\begin{align*}
\frac{\dot{V}(\bs{\mu})}{2} &= \frac{1}{2}\nabla V \cdot \frac{d\bs{\mu}}{dt}\\
&= \big(\bs{\mu} -\bs{\mu}_*\big) (\Ffuncbold(\bs{\mu}) - \bs{\mu})^T\\
&=\big(\bs{\mu} - \Ffuncbold(\bs{\mu})+\Ffuncbold(\bs{\mu})-\bs{\mu}_*\big)(\Ffuncbold(\bs{\mu}) - \bs{\mu})^T \\
& = -(\bs{\mu} - \Ffuncbold(\bs{\mu}))(\bs{\mu} - \Ffuncbold(\bs{\mu}))^T + \big(\Ffuncbold(\bs{\mu})-\bs{\mu}_*\big)(\Ffuncbold(\bs{\mu}) - \bs{\mu})^T< 0,
\end{align*}
where the last inequality follows by the fact that $\bs{\mu} \geq \bs{\mu}_*$ implies $\Ffuncbold(\bs{\mu}) \geq \Ffuncbold(\bs{\mu}_*) = \bs{\mu}_*$. Also, note that $\dot{V}(\bs{\mu}_*) = V(\bs{\mu}_*) = 0$. Now, the proof of Theorem \ref{thm:alterode_sol} follows by the LaSalle Invariance Principle~\cite{Lasalle1960}. Specifically, all trajectories with initial value in $\mathcal{N}$ converge to $\bs{\mu}_*$.

Next, we prove that the trajectory of the solution meets the set $$\partial\mathcal{U}\setminus\{\bs{\mu}:\mu^{(j\gets j')}=1 \text{ for some }j,j'\in\{1,2\}\}$$ only at  $\bs{\mu}_*$, where $\partial\mathcal{U}$ is the boundary of $\mathcal{U}$. For sake of contradiction, let us assume that at time $t_0>0$ we have $\bs{\mu}(t_0)\in \partial \mathcal{U}$, $\bs{\mu}(t_0)<\bs{1}$ and $\bs{\mu}(t_0)\neq \bs{\mu}_*$. By the above argument and Lemma \ref{lem:Fprop_feasreg}, we have $\bs{\mu}_* < \bs{\mu}(t_0) < \bs{1}$. Moreover, for all small enough $\delta > 0$ we have $\bs{\mu}(t_0 - \delta) = \bs{\mu}(t_0) - \delta \times (\Ffuncbold(\bs{\mu}(t_0)) - \bs{\mu}(t_0))+ O(\delta^2)$.

	Since $\Ffuncbold$ is continuous, at least one of the components of $\Ffuncbold(\bs{\mu}(t_0)) - \bs{\mu}(t_0)$ is equal to zero (otherwise $\bs{\mu}(t_0)$ was an interior point of $\mathcal{U}$).
	Hence, we have the following cases:
	\begin{enumerate}[label=(\alph*)]
	\item $\mu^{(j\gets j)}(t_0) = \Ffunc_{(j\gets j)}(\mu^{(j\gets j)}(t_0),\mu^{(j\gets -j)}(t_0))$: If $\Ffunc_{(j\gets -j)}(\bs{\mu}(t_0)) < \mu^{(j\gets -j)}(t_0)$, then for all small enough $\delta>0$ we have
	\begin{align*}
		\mu^{(j\gets -j)}(t_0 - \delta) &= \mu^{(j\gets -j)}(t_0) - \delta \mathrm{Const}_{(j\gets -j)} + O(\delta^2),
	\end{align*}
	where $\mathrm{Const}_{(j\gets -j)} \coloneqq \Ffunc_{(j\gets -j)}(\bs{\mu}(t_0)) - \mu^{(j\gets -j)}(t_0)< 0$. Moreover, for all small enough $\delta>0$ we have
	\begin{align*}
		\mu^{(j\gets j)}(t_0 - \delta) &= \mu^{(j\gets j)}(t_0) + O(\delta^2)
		=\Ffunc_{(j\gets j)}(\mu^{(j\gets j)}(t_0),\mu^{(j\gets -j)}(t_0)) + O(\delta^2).
	\end{align*}
	Now, combining the above equalities we have
	\begin{align*}
		&\mu^{(j\gets j)}(t_0 - \delta) = \Ffunc_{(j\gets j)}(\mu^{(j\gets j)}(t_0 - \delta),\mu^{(j\gets -j)}(t_0-\delta) + \delta \mathrm{Const}_{(j\gets -j)}) + O(\delta^2),
	\end{align*}
	which is smaller than $\Ffunc_{(j\gets j)}(\bs{\mu}(t_0-\delta))$ for small enough $\delta> 0$ and contradicts the assumption that $\bs{\mu}(t_0-\delta)\in\mathcal{U}$; hence, $\mu^{(j\gets -j)}(t_0) = \Ffunc_{(j\gets -j)}(\bs{\mu}(t_0))$. Next, following the same argument we show that $\Ffuncbold(\bs{\mu}(t_0)) = \bs{\mu}(t_0)$.

	Assume that $\Ffuncbold(\bs{\mu}(t_0)) \neq \bs{\mu}(t_0)$. Using the same logic as above, we have
	\begin{align*}
		&\mu^{(j\gets -j)}(t_0 - \delta)\\
		&\myquad[1]=\Ffunc_{(j\gets -j)}(\mu^{(-j\gets -j)}(t_0),\mu^{(-j\gets j)}(t_0)) + O(\delta^2)\\
		&\myquad[1]= \Ffunc_{(j\gets -j)}(\mu^{(-j\gets -j)}(t_0 - \delta) + \delta \mathrm{Const}_{(-j\gets -j)},\mu^{(-j\gets j)}(t_0-\delta) + \delta \mathrm{Const}_{(-j\gets j)}) + O(\delta^2),
	\end{align*}
	where $\mathrm{Const}_{(-j\gets -j)} \coloneqq \Ffunc_{(-j\gets -j)}(\bs{\mu}(t_0)) - \mu^{(-j\gets -j)}(t_0)\leq 0$. If either $\mathrm{Const}_{(-j\gets -j)} < 0$ or $\mathrm{Const}_{(-j\gets j)}<0$, then $\mu^{(j\gets -j)}(t_0 - \delta)$ is smaller than $\Ffunc_{(j\gets -j)}(\bs{\mu}(t_0-\delta))$ for small enough $\delta > 0$, which contradicts the assumption that $\bs{\mu}(t_0-\delta)\in\mathcal{U}$. Hence,
	$\mathrm{Const}_{(-j\gets -j)} = \mathrm{Const}_{(-j\gets j)} = 0$, and we have $\Ffuncbold(\bs{\mu}(t_0)) = \bs{\mu}(t_0)$.
	\item $\mu^{(j\gets -j)}(t_0) = \Ffunc_{(j\gets -j)}(\mu^{(-j\gets -j)}(t_0),\mu^{(-j\gets j)}(t_0))$: The argument is exactly the same as the previous case, and we avoid repetition.
\end{enumerate}

	\subsection{Jacobian Matrix in Section \ref{sec:odeanalysisinf}}\label{proof:jacobianmatrix}
	Following the derivation in Appendix \ref{proof:Fprop_increasing} and using the asymptotic values for the initial condition given in Lemma \ref{lem:initcond_conv}, we have:
\begin{align*}
	&\frac{\partial \Ffunc_{(j\gets j),\infty}(\bs{\mu}_{*,\infty})}{\partial \mu^{(j\gets j)}}=\\*
	&\myquad[2]\begin{aligned}
		 &\sum_{\substack{u_{-j} \leq K_j(d_j,d_{-j})\allowdisplaybreaks\\d_j+d_{-j} \leq d_{\max}}}\,\prob_{j*,m}(d_j,d_{-j})(1-\alpha_j(d_j,d_{-j}))\, Bi(u_{-j};d_{-j},1 - {\mu}_{*,\infty}^{(j\gets -j)}) \\
		&\myquad[8] \times (d_j-1)\,Bi(K_j(d_j,d_{-j})-u_{-j};d_j-2,1 -{\mu}_{*,\infty}^{(j\gets j)})
	\end{aligned},\allowdisplaybreaks\\
	&\frac{\partial \Ffunc_{(j\gets j),\infty}(\bs{\mu}_{*,\infty})}{\partial \mu^{(j\gets -j)}} =\\*
	&\myquad[2]\begin{aligned}
		&\sum_{\substack{u_j \leq K_j(d_j,d_{-j})\allowdisplaybreaks\\d_j+d_{-j} \leq d_{\max}}}\,\prob_{j*,m}(d_j,d_{-j})(1-\alpha_j(d_j,d_{-j}))\, Bi(u_j;d_j-1,1 - {\mu}_{*,\infty}^{(j\gets j)})\\
		&\myquad[8] \times d_{-j}\,Bi(K_j(d_j,d_{-j})-u_{j};d_{-j}-1,1 -{\mu}_{*,\infty}^{(j\gets -j)})
	\end{aligned},\allowdisplaybreaks\\
	&\frac{\partial \Ffunc_{(j\gets -j),\infty}(\bs{\mu}_{*,\infty})}{\partial \mu^{(-j\gets -j)}} =\\*
	&\myquad[2]\begin{aligned}
		&\sum_{\substack{u_{j} \leq K_{-j}(d_{-j},d_{j})\allowdisplaybreaks\\d_j+d_{-j} \leq d_{\max}}}\prob_{-j,m*}(d_{-j},d_{j})(1-\alpha_{-j}(d_{-j},d_{j}))\, Bi(u_{j};d_{j}-1,1 - {\mu}_{*,\infty}^{(-j\gets j)})\\
		&\myquad[8] \times d_{-j}\,Bi(K_{-j}(d_{-j},d_{j})-u_{j};d_{-j}-1,1 -{\mu}_{*,\infty}^{(-j\gets -j)})
	\end{aligned},\allowdisplaybreaks\\
	&\frac{\partial \Ffunc_{(j\gets -j),\infty}(\bs{\mu}_{*,\infty})}{\partial \mu^{(-j\gets j)}}=\\
	&\myquad[2]\begin{aligned}
		&\sum_{\substack{u_{-j} \leq K_{-j}(d_{-j},d_{j})\allowdisplaybreaks\\d_j+d_{-j} \leq d_{\max}}}\,\prob_{-j,m*}(d_{-j},d_{j})(1-\alpha_{-j}(d_{-j},d_{j}))\, Bi(u_{-j};d_{-j},1 - {\mu}_{*,\infty}^{(-j\gets -j)})\allowdisplaybreaks\\
		&\myquad[8] \times (d_j-1)\,Bi(K_{-j}(d_{-j},d_{j})-u_{-j};d_j-2,1 -{\mu}_{*,\infty}^{(-j\gets j)})
	\end{aligned}.
\end{align*}
	\subsection{Proof of Theorem \ref{thm:contagion}}\label{proof:contagion}
	By Taylor expansion of $\Ffuncbold_\infty(\bs{0},\bs{u})$ at $\bs{u} = \bs{1}$, it is easy to see that if $\zeta_\infty(\bs{0}) < 1$, then $\mathcal{U}_\infty(\bs{0}) = \{\bs{1}\}$, and if $\zeta_\infty(\bs{0}) > 1$, then $\bs{1} - \kappa \bs{\nu}(\bs{0}) \in \mathcal{U}_\infty(\bs{0})$ for all small enough $\kappa > 0$, where $\bs{\nu}(\bs{0})$ is the Perron-Frobenius eigenvector of $\bs{J}_{\Ffuncbold_\infty(\bs{0},\cdot)}(\bs{1})$ corresponding to eigenvalue $\zeta_\infty(\bs{0})$. It is also easy to see that for any $\bs{\alpha} \neq \bs{0}$, we have $\Ffuncbold_\infty(\bs{\alpha},\bs{u}) < \Ffuncbold_\infty(\bs{0},\bs{u})$, and hence, $\mathcal{U}_\infty(\bs{0})\subset \mathcal{U}_\infty(\bs{\alpha})$.

Next, we show a simple observation that if $\norm{\bs{\alpha}_s}_\infty \rightarrow0$, then we have $\Ffuncbold_\infty(\bs{\alpha}_s,\bs{u}) \rightarrow\Ffuncbold_\infty(\bs{0},\bs{u})$ for all $\bs{u}\in[0,1]^4$.
\begin{lemma}\label{lem:temp}
Assume $\{\bs{\alpha}_s\}_{s=1}^\infty$ converges to zero in sup-norm. Then, we have $$\norm{\Ffuncbold_\infty(\bs{\alpha}_s,\bs{u}) - \Ffuncbold_\infty(\bs{0},\bs{u})}_\infty \to 0$$ uniformly over $\bs{u}\in[0,1]^4$.
\end{lemma}
\begin{proof}
Fix $\delta > 0$. Pick $s_0(\delta)\in\mathbb{N}$ large enough such that $\norm{\bs{\alpha}_s}_\infty < \delta$ for all $s>s_0(\delta)$. It is easy to see that $\Ffuncbold_\infty(\bs{0},\bs{u}) \geq \Ffuncbold_\infty(\bs{\alpha}_s,\bs{u})  > \Ffuncbold_\infty(\delta\bs{1},\bs{u})$ for all $\bs{u}\in[0,1]^4$. Now, for any $\bs{u}\in[0,1]^4$ we have
\begin{align*}
\begin{aligned}
\norm{ \Ffuncbold_\infty(\bs{0},\bs{u}) - \Ffuncbold_\infty(\bs{\alpha}_s,\bs{u})}_\infty&\leq  \norm{\Ffuncbold_\infty(\bs{0},\bs{u}) - \Ffuncbold_\infty(\delta\bs{1},\bs{u}) }_\infty \\
&= \delta \norm{\Ffuncbold_\infty(\bs{0},\bs{u})}_\infty
\end{aligned}
 \label{eq_pr:ContConvFLemm}
\end{align*}
\end{proof}

Following the proof of Lemma \ref{lem:temp}, for any $\delta > 0$ we have $\Ffuncbold_\infty(\bs{0},\bs{u}) > \Ffuncbold_\infty(\bs{\alpha}_s,\bs{u})  > \Ffuncbold_\infty(\delta\bs{1},\bs{u})$ for all $s>s_0(\delta)$. Hence, $\mathcal{U}_\infty(\bs{0}) \subset\mathcal{U}_\infty(\bs{\alpha}_s) \subset \mathcal{U}_\infty(\delta\bs{1})$ for all $s>s_0(\delta)$. Now, using the fact that
$\cap_{k=1}^\infty \mathcal{U}_\infty(1/k\times\bs{1}) = \mathcal{U}_\infty(\bs{0})$, we have $\lim_{s\to\infty}\mathcal{U}_\infty(\bs{0}) \cap \mathcal{U}_\infty(\bs{\alpha}_s) =\mathcal{U}_\infty(\bs{0})$.

By Corollary \ref{cor:Fprop_sol}, $\bs{\mu}_{*,\infty}(\bs{\alpha}_s) = \lim_{r\to\infty}\Ffuncbold_\infty^r(\bs{\alpha}_s,\bs{1})$ is the closest fixed point of $\Ffuncbold_\infty(\bs{\alpha}_s,\cdot)$ to $\bs{1}$ in sup-norm. Define $\bs{\mu}_{*,\infty}(\bs{0}) \in \mathcal{U}(\bs{0})$ as follows: if $\mathcal{U}(\bs{0})$ is a singleton, set $\bs{\mu}_{*,\infty}(\bs{0}) \coloneqq \{\bs{1}\}$, otherwise, set $\bs{\mu}_{*,\infty}(\bs{0})$ to be the closest fixed point of $\Ffuncbold_\infty(\bs{0},\cdot)$ to $\bs{1}$ other than $\bs{1}$ itself. Note that for all $\bs{u}\in  \mathcal{U}_\infty(\bs{0}) \cap \{\bs{x} :  \bs{\mu}_{*,\infty}(\bs{0}) \leq \bs{x} \leq \bs{1} \} \setminus \{\bs{1}\}$, using the same argument as in the proof of Lemma \ref{lem:Fprop_feasreg}, we have $\lim_{r\to\infty} \Ffuncbold_\infty^r(\bs{0},\bs{u}) =  \bs{\mu}_{*,\infty}(\bs{0}) $. Now, the sequence $\{\bs{\mu}_{*,\infty}(\bs{\alpha}_s)\}_{s=0}^\infty$ is sandwiched between the sequence $\{\bs{\mu}_{*,\infty}(1/k\times \bs{1})\}_{s=0}^\infty$ and $\bs{\mu}_{*,\infty}(\bs{0})$. Hence, if $\mathcal{U}_\infty(\bs{0})$ is a singleton, then the final proportion of adopters converges to 0 as $\bs{\alpha_s}\to \bs{0}$. Otherwise, the final proportion of adopters is strictly positive, and we have
\begin{align*}
 \lim_{s\to\infty} \bs{\mu}_{*,\infty}(\bs{\alpha}_s) = \lim_{r\to\infty} \Ffuncbold_\infty(\bs{0},\bs{u}) \qquad \forall\mu\in \mathcal{U}_\infty(\bs{0}) \cap \{\bs{x} :  \bs{\mu}_{*,\infty}(\bs{0}) \leq \bs{x} \leq \bs{1} \} \setminus \{\bs{1}\}.
\end{align*}
	\subsection{Proof of Theorem \ref{thm:poissreduc}}\label{proof:poissreduc}
	
By definition of $\Ffuncbold(\cdot)$ given by the right-hand side of \eqref{eq:meanfield_mujj}-\eqref{eq:meanfield_muj-j}, and the fact that $\prob_{j*,m}(d_j,d_{-j}) = \prob_{j,m}(d_j-1,d_{-j})$ and $\prob_{j,m*}(d_j,d_{-j}) = \prob_{j,m}(d_j,d_{-j}-1)$, we have
\begin{align*}
&\Ffunc_{(j\gets j),\infty}(\mu^{(j\gets j)},\mu^{(j\gets -j)}) \\
&\myquad[2]= \sum_{u_j+u_{-j} \leq K_j(d_j + d_{-j})} \prob_{j,m}(d_j-1,d_{-j})(1-\alpha_j(d_j + d_{-j}))\times \allowdisplaybreaks\\
&\myquad[2]~\myquad[6] Bi(u_j;d_j-1,1-\mu^{(j\gets j)})Bi(u_{-j};d_{-j},1-\mu^{(j\gets -j)}) \allowdisplaybreaks\\
&\myquad[2]= \sum_{u_j+u_{-j} \leq K_j(d_j + d_{-j}+1)} \prob_{j,m}(d_j,d_{-j})(1-\alpha_j(d_j + d_{-j}+1)) \times\allowdisplaybreaks\\
&\myquad[2]~ \myquad[6] Bi(u_j;d_j,1-\mu^{(j\gets j)})Bi(u_{-j};d_{-j},1-\mu^{(j\gets -j)}) \allowdisplaybreaks\\
\
&\Ffunc_{(-j\gets j),\infty}(\mu^{(j\gets j)},\mu^{(j\gets -j)}) \\
&\myquad[2]= \sum_{u_j+u_{-j} \leq K_j(d_j + d_{-j})}  \prob_{j,m}(d_j,d_{-j}-1)(1-\alpha_j(d_j + d_{-j}))\times\allowdisplaybreaks\\
&\myquad[2]~ \myquad[6] Bi(u_{-j};d_{-j}-1,1-\mu^{(j\gets -j)})Bi(u_j;d_j,1-\mu^{(j\gets j)}) \allowdisplaybreaks\\
&\myquad[2]= \sum_{u_j+u_{-j} \leq K_j(d_j + d_{-j}+1)} \prob_{j,m}(d_j,d_{-j})(1-\alpha_j(d_j + d_{-j}+1))\times \allowdisplaybreaks\\
&\myquad[2]~ \myquad[6] Bi(u_{-j};d_{-j},1-\mu^{(j\gets -j)})Bi(u_j;d_j,1-\mu^{(j\gets j)})\allowdisplaybreaks\\
\
\end{align*}
Hence, for all $\bs{\mu}\in[0,1]^4$, we have $\Ffunc_{(1\gets 1),\infty}(\mu^{(1\gets 1)},\mu^{(1\gets 2)}) = \Ffunc_{(2\gets 1),\infty}(\mu^{(1\gets 1)},\mu^{(1\gets 2)})$ and $\Ffunc_{(2\gets 2),\infty}\allowbreak(\mu^{(2\gets 2)},\mu^{(2\gets 1)}) = \Ffunc_{(1\gets 2),\infty}(\mu^{(2\gets 2)},\mu^{(2\gets 1)})$. Now, if $\mu^{(1\gets 1)} = \mu^{(2\gets 1)}$ and $\mu^{(2\gets 2)} = \mu^{(1\gets 2)}$, then we have
\begin{align*}
&\Ffunc_{(1\gets 1),\infty}(\mu^{(1\gets 1)},\mu^{(1\gets 2)})  - \mu^{(1\gets 1)}= \Ffunc_{(2\gets 1),\infty}(\mu^{(1\gets 1)},\mu^{(1\gets 2)}) - \mu^{(2\gets 1)}\allowdisplaybreaks\\
&\Ffunc_{(2\gets 2),\infty}(\mu^{(2\gets 2)},\mu^{(2\gets 1)}) - \mu^{(2\gets 2)}= \Ffunc_{(1\gets 2),\infty}(\mu^{(2\gets 2)},\mu^{(2\gets 1)}) - \mu^{(1\gets 2)}
\end{align*}
Since these equalities hold at time $0$ of the ODEs \eqref{eq:alterode} with the function $\Ffuncbold_\infty(\cdot)$, they are satisfied on the whole trajectory; that is to say, $\mu_\infty^{(1\gets 1)}(t) = \mu_\infty^{(2\gets 1)}(t)$ and $\mu_\infty^{(2\gets 2)}(t) = \mu_\infty^{(1\gets 2)}(t)$ for all $t\geq0$, where $\bs{\mu}_\infty(t)$ is the solution of the ODEs.
	\subsection{Proof of Theorem \ref{thm:symmreduc}}\label{proof:symmreduc}
	
If $\mu^{(1\gets 1)} = \mu^{(2\gets 2)}$ and $\mu^{(2\gets 1)} = \mu^{(1\gets 2)}$ then we have
\begin{align*}
	&\Ffunc_{(1\gets 1),\infty}(\mu^{(1\gets 1)},\mu^{(1\gets 2)}) = \Ffunc_{(2\gets 2),\infty}(\mu^{(2\gets 2)},\mu^{(2\gets 1)}) \text{ and }\\
	&\Ffunc_{(2\gets 1),\infty}(\mu^{(1\gets 1)},\mu^{(1\gets 2)}) = \Ffunc_{(1\gets 2),\infty}(\mu^{(2\gets 2)},\mu^{(2\gets 1)})
\end{align*}
Using the same argument as in the proof of Theorem \ref{thm:poissreduc}, we have $\mu_\infty^{(1\gets 1)}(t) = \mu_\infty^{(2\gets 1)}(t)$ and $\mu_\infty^{(2\gets 2)}(t) = \mu_\infty^{(1\gets 2)}(t)$ for all $t\geq0$, where $\bs{\mu}_\infty(t)$ is the solution of the ODEs.

	\section*{Acknowledgements}
	M. Moharrami was partially suppoerted by the NSF under grants CNS 1422211, ECCS 1446521, IIS 1538827, AST 1516075 and Rackham Graduate Predoctoral Fellowship. The majority of the work was done while the first author was at the University of Michigan. He also acknowledges support by NSF CCF 1934986 during his stay as a Postdoc at the University of Illinois at Urbana-Champaign.
	V. Subramanian was partially supported by NSF via grants IIS 1538827, AST 1516075, CCF 2008130, ECCS 2038416 and CNS 1955777. He also acknowledges support by INRIA for a visit to Paris in May 2015 and discussions with Venkat Anantharam and Sanjay Shakkottai on the limiting behavior of the cascade process.
	M. Liu was partially supported by the NSF under grants CNS 1422211 and ECCS 1446521.
	M. Lelarge was partially supported by ANR via grant ANR-11-JS02-005-01.

	\bibliographystyle{plain}
	\bibliography{references}
\end{document}